\documentclass{article}[12pt]
\usepackage{a4}
\usepackage[T1]{fontenc}
\usepackage[latin1]{inputenc}
\usepackage[francais]{babel}
\usepackage{hyperref}
\usepackage[dvips]{graphics}
\usepackage{amsmath}
\usepackage{amssymb}
\usepackage{amsopn}
\usepackage{amsbsy}
\usepackage{amstext}
\usepackage{latexsym}
\usepackage{amsfonts}
\usepackage{array}
\usepackage{epsfig}
\usepackage[mathscr]{eucal}

 \newtheorem{theorem}{Th\'eor\`eme}[section]
 \newtheorem{corollary}[theorem]{Corollaire}
 \newtheorem{lemma}[theorem]{Lemme}
 \newtheorem{proposition}[theorem]{Proposition}

 \newtheorem{remarque}[theorem]{Remarque}
 \newtheorem{thm*}{Th\'eor\`eme}
 \numberwithin{equation}{section}
\def\EMdash{\leavevmode\hbox to 7.5mm{\vrule height .63ex depth -.59ex
    width 5.4mm\hfill}}

\def\Kbar{{\, \overline{K}}}
\def\hbar{{\, \overline{h}}}

\def\card{{\rm{card}}}
\def\cl{{\rm{cl}}}
\def\plongement{\hookrightarrow}
\def\Xsoul{{\underline{X}}}
\def\Ysoul{{\underline{Y}}}

\def\Asoul{{\underline{A}}}

\def\alphasoul{{\underline{\alpha}}}

\def\l{{\ell}}

\def\rest#1{\vert_{#1}}              
\def\EMts{\mspace{.3mu}}  
\def\nb#1{{\left\vert{\EMts\EMts #1 \EMts\EMts}\right\vert}}        
\def\:{\!\!:}

\def\0{{\mathbf{0}}}

\let\epsilon=\varepsilon

\begin{document}

\title{APPROXIMATION DIOPHANTIENNE SUR UNE COURBE ELLIPTIQUE:}
\author{Bakir FARHI \\
\footnotesize{D\'epartement de Math\'ematiques, Universit\'e du
Maine,} \\ \footnotesize{Avenue Olivier Messiaen, 72085 Le Mans
Cedex 9, France.} \\ \footnotesize{Bakir.Farhi@univ-lemans.fr}}
\date{}

\maketitle \tableofcontents \newpage
\section{Introduction}
Notre but ici est de donner un analogue du
th\'eor\`eme 2 de [Fa2] dans le cas particulier o\`u la variété
abelienne $A$ est une courbe elliptique et la sous-variété $E$ de
$A$ est réduite au point à l'infini de cette courbe elliptique
(voir ci-dessous). L'intér\^et est de fournir une preuve plus
élémentaire (utilisant les polyn\^omes au lieu des sections) et en
plus donnant une estimation explicite pour le nombre de points
exceptionnels en question. \\ Soit $K$ un corps de nombres, $A$
une variété abelienne définie sur $K$ et $E$ une $K$-sous-variété
de $A$. Soit aussi $w$ une place de $K$, pour $x \in A(K)$, on
peut définir sa hauteur multiplicative $H(x)$ (apr\`es avoir
choisi un diviseur ample sur $A$) comme on peut définir aussi la
distance $w$-adique $d_w(x , E)$ de $x$ à $E$ (voir §$2$).
Rappelons d'abord le
théor\`eme 2 de [Fa2]: \\
{\bf{Théor\`eme}} ({\rm{G. Faltings}}) {\it{Pour tout $\epsilon >
0$ et pour presque\footnote{Le mot presque veut dire ici ``à
l'exception d'un nombre fini de points''.}tout point $K$-rationnel
$x \in A - E$ on a: $d_w(x , E) \geq {H(x)}^{- \epsilon}$.}}
\section{Notations et résultats}
Soit $E \stackrel{i}{\hookrightarrow}{\mathbb P}_2$ une courbe
elliptique définie sur un corps de nombres $K$, plongée (à la
Weierstrass) dans l'espace projectif ${\mathbb P}_2$ et d'équation
projective: $$Y^2 Z = 4 X^3 - g_2 X Z^2 - g_3 Z^3 ~~~~(g_2 , g_3
\in K )~, $$ nous prenons le point à l'infini $\mathbf 0$,
représenté dans ${\mathbb P}_2$ par les coordonnées $(0 , 1 , 0)$,
comme élément neutre de $E$. Soit aussi $\widetilde{E}$ la courbe
affine $\widetilde{E} := E \cap \{Y \neq 0\}$ qu'on peut plonger
dans ${\mathbb A}_2$ gr\^ace au morphisme:
$$
\begin{array}{ccc}
\pi\!: ~~~~\widetilde E & \longrightarrow & {\mathbb A}_2 \\
~~~~~~(x : y : z) & \longmapsto & \left(\frac{x}{y} ,
\frac{z}{y}\right)
\end{array}
$$
$\pi(\widetilde E)$ est alors d'équation: $\widetilde{G}(X , Z) =
0$ avec $\widetilde{G}(X , Z) := Z + g_2 X Z^2 + g_3 Z^3 - 4 X^3$.
Posons aussi $\widetilde{\Delta}(X , Z ) := \frac{\partial
\widetilde{G}}{\partial Z} (X , Z) = 3g_3 Z^2 + 2 g_2 X Z + 1$,
$\widetilde{G}$ et $\widetilde{\Delta}$ sont donc des polyn\^omes
de $K[X , Z]$. Par ailleurs notons respectivement par $A(E)$ et
$A(\widetilde E)$ les anneaux de coordonnées de $E$ et $\widetilde
E ,$ par $K(E)$ et $K(\widetilde E)$ les corps de fractions de
$A(E)$ et $A(\widetilde E)$ et par ${K(E)}_0$ les éléments
homog\`enes de degré $0$ de $K(E) .$
\\ Maintenant, pour toute place $v$ de $K$, on note par $M_v$ et $m_v$ les deux r\'eels positifs:
\begin{eqnarray*}
M_v & := & \max\left\{1 , {\mid g_2 \mid}_v , {\mid g_3 \mid}_v\right\} \\
m_v & := & \log M_v
\end{eqnarray*}
et par $\eta$ le r\'eel positif:
$$\eta ~:=~ h(1 : g_2 : g_3) ~=~ \sum_{v \in M_K} \frac{[K_v : {\mathbb Q}_v]}{[K : \mathbb Q]} m_v .$$
On note aussi par $c_v$ la constante absolue:
$$c_v ~:=~ \begin{cases}
0 & \text{si $v$ est finie} \\
16 & \text{si $v$ est infinie}
\end{cases}.$$
Si $P$ est un polyn\^ome \`a une ou plusieurs ind\'etermin\'ees
\`a coefficients dans $K$, désignons par $H_v(P)$ (resp $L_v(P)$)
le maximum (resp la somme) des valeurs absolues $v$-adiques de
tous les coefficients de $P$. Si de plus $P$ est non identiquement
nul, on désigne par $h_v(P)$ et $\ell_v(P)$ les deux nombres
r\'eels:
\begin{eqnarray*}
h_v(P) & := & \log H_v(P) \\
\ell_v(P) & := & \log L_v(P).
\end{eqnarray*}
Plus g\'en\'eralement, si $F = {(P_i)}_{i \in I}$ est une famille
finie de polyn\^omes \`a une ou plusieurs ind\'eterminées, \`a
coefficients dans $K$ et non tous identiquement nuls, on désigne
par $H_v(F) , L_v(F) , h_v(F)$ et $\ell_v(F)$ les nombres r\'eels:
\begin{eqnarray*}
H_v(F) & := & \max_{i \in I} H_v(P_i) \\
L_v(F) & := & \max_{i \in I} L_v(P_i) \\
h_v(F) & := & \max_{i \in I} h_v(P_i) ~=~ \log H_v(F) \\
\ell_v(F) & := & \max_{i \in I} \ell_v(P_i) ~=~ \log L_v(F)
\end{eqnarray*}
et on entend par {\sl{hauteur de Gauss-Weil}} de $F$ le nombre
r\'eel:
$$\widetilde{h}(F) ~:=~ \sum_{v \in M_K} \frac{[K_v : {\mathbb Q}_v]}{[K : \mathbb Q]} h_v(F) .$$
En voici dans ce qui suit quelques propriétés qu'on utilisera
souvent dans ce qui va suivre. \\ \\
$\underline{\text{Quelques propriétés des hauteurs et longueurs locales}}$: \\

Soient $P_1 , \dots , P_n$ $(n \in {\mathbb N}^*)$ des polyn\^omes de $K[X_1 , \dots , X_d]$ $(d \in \mathbb N)$ et $Q$ un polyn\^ome de $K[Y_1 , \dots Y_n]$. On a: \\
\begin{itemize}
\item[1)] Pour toute place finie $v$ de $K$:
\begin{align}
H_v(P_1 + \dots + P_n) &\leq \max_{1 \leq i \leq n} H_v(P_i) , \notag \\
H_v(P_1 \dots P_n) &\leq \prod_{i = 1}^{n} H_v(P_i) . \notag
\end{align}
\item[2)] Pour toute place infinie $v$ de $K$:
\begin{align}
H_v(P_1 + \dots + P_n) &\leq n H_v\!\left(\{P_1 , \dots , P_n\}\right) \leq n \prod_{i = 1}^{n} \max\{1 , H_v(P_i)\} , \notag \\
H_v(P_1 \dots P_n) &\leq \prod_{i = 1}^{n - 1} \mathcal{N}(P_i) . \prod_{i = 1}^{n} H_v(P_i) \notag \\
\intertext{(o\`u on a noté $\mathcal{N}$ l'application associant à
tout polyn\^ome -à coefficients complexes et en un certain nombre
d'indétérminées- le nombre de mon\^omes intervenant dans son
écriture canonique). De plus:}
L_v(P_1 + \dots + P_n) &\leq \sum_{i = 1}^{n} L_v(P_i) , \notag \\
L_v(P_1 \dots P_n) &\leq \prod_{i = 1}^{n} L_v(P_i) . \notag
\end{align}
\item[3)] Si $v$ est une place finie de $K$:
\begin{align}
H_v\left(Q(P_1 , \dots , P_n)\right) &\leq H_v(Q) \prod_{i = 1}^{n} {H_v(P_i)}^{{d°}_{Y_i}Q} \notag \\
\intertext{et si $v$ est une place infinie de $K$:} L_v\left(Q(P_1
, \dots , P_n)\right) &\leq L_v(Q) \prod_{i = 1}^{n}
{L_v(P_i)}^{{d°}_{Y_i}Q} \notag .
\end{align}
\item[4)] Pour tout $\alphasoul = (\alpha_1 , \dots , \alpha_n)
\in {\mathbb N}^n$ et toute place $v$ de $K$:
\begin{align}
H_v\!\!\left(\frac{1}{\alpha_1! \dots \alpha_n!} \frac{\partial^{\nb{\alphasoul}} Q}{\partial Y_{1}^{\alpha_1} \dots \partial Y_{n}^{\alpha_n}}\right) &\leq H_v(Q) \notag \\
\intertext{si $v$ est finie et:} H_v\!\!\left(\frac{1}{\alpha_1!
\dots \alpha_n!} \frac{\partial^{\nb{\alphasoul}} Q}{\partial
Y_{1}^{\alpha_1} \dots \partial Y_{n}^{\alpha_n}}\right) &\leq
\prod_{i = 1}^{n} \binom{{d°}_{Y_i} Q}{\alpha_i} . H_v(Q) \notag
\end{align}
si $v$ est infinie.
\end{itemize}
$~~$ \\

On d\'esigne, par ailleurs, par ${\mbox{dist}}_v$ la distance projective $v$-adique sur ${\mathbb P}_2(K)$ d\'efinie de la mani\`ere suivante: \\
Pour tous points $\mathbf{p}$ et $\mathbf{q}$ de ${\mathbb
P}_2(K)$ repr\'esent\'es respectivement par les deux syst\`emes
projectifs de $K^3$: $\underline p = (p_0 , p_1 , p_2)$ et
$\underline q = (q_0 , q_1 , q_2)$, on définit:
\begin{equation*}
{\mbox{dist}}_v(\mathbf p , \mathbf q) ~:=~ \frac{\max\left({\mid
p_0 q_1 - q_0 p_1\mid}_v , {\mid p_0 q_2 - q_0 p_2\mid}_v , {\mid
p_1 q_2 - q_1 p_2 \mid}_v\right)}{\max\left({\mid p_0 \mid}_v ,
{\mid p_1 \mid}_v , {\mid p_2 \mid}_v\right) . \max\left({\mid q_0
\mid}_v , {\mid q_1 \mid}_v , {\mid q_2 \mid}_v\right)} .
\end{equation*}
Cette distance ${\mbox{dist}}_v$ a les deux particularit\'es
int\'eressantes suivantes:
\begin{description}
\item[i)] $\forall \mathbf x \in {\mathbb P}_2(K)$ on a:
${\mbox{dist}}_v(\mathbf x , \mathbf 0) \leq 1 ,$ \item[ii)]
$\forall \mathbf x \in {\mathbb P}_2(K)$ on a: \item[]
${\mbox{dist}}_v(\mathbf x , \mathbf 0) < 1 \Rightarrow \mathbf x
\in {\mathbb P}_2(K) \setminus \{Y = 0\}$ et si $\underline x = (x
, 1 , z) \in K^3$ est un \item[] repr\'esentant de $\mathbf x$
alors ${\mbox{dist}}_v(\mathbf x , \mathbf 0) = \max\left({\mid x
\mid}_v , {\mid z \mid}_v\right) .$
\end{description}
En effet, lorsque $\underline x = (x , y , z) \in K^3$ est un
repr\'esentant d'un point $\mathbf x$ de ${\mathbb P}_2(K)$, on a
bien:
$${\mbox{dist}}_v(\mathbf x , \mathbf 0) = \frac{\max\left({\mid x \mid}_v , {\mid z \mid}_v\right)}{\max\left({\mid x \mid}_v , {\mid y \mid}_v , {\mid z \mid}_v\right)} .$$
Cette derni\`ere identité entra{\sf\^\i}ne imm\'ediatement les
deux propri\'et\'es i) et ii) pr\'ec\'edentes pour la distance
${\mbox{dist}}_v$.
\\ Nos résultats principaux sont les suivants:
\begin{theorem}[Premier théor\`eme principal]\label{c.1}
Soit $E$ une courbe elliptique définie sur un corps de nombres $K$
de degré $D$, plongée dans ${\mathbb P}_2$ à la Weierstrass,
d'équation projective $Y^2 Z = 4 X^3 - g_2 X Z^2 - g_3 Z^3 ~ (g_2
, g_3 \in K)$ et d'élément neutre (en tant que groupe) le point à
l'infini $\mathbf 0$ représenté dans ${\mathbb P}_2$ par les
coordonnées projectives $(0 : 1 : 0)$. On désigne par $r$ le rang
de Mordell-Weil de $E(K)$ que l'on suppose non nul. Soient aussi
$S$ un ensemble fini de places de $K$, $m_v , c_v ~ (v \in M_K)$
et $\eta$ les réels positifs définis précédemment et
${(\lambda_v)}_{v \in S}$ une famille de réels positifs
satisfaisant: $$\sum_{v \in S} \frac{[K_v : {\mathbb Q}_v]}{[K :
\mathbb Q]} \lambda_v = 1 .$$ Soit enfin $\epsilon$ un réel
strictement positif $< \frac{1}{15788}$. Alors l'ensemble des
points $\mathbf x$ de $E(K)$ satisfaisant le syst\`eme
d'inégalités simultanées:
$${\mbox{dist}}_v(\mathbf x , \mathbf 0) <
\exp\!\left\{\!\!- \lambda_v \!\left(\epsilon h(\mathbf x) +
56(\eta + 5) \epsilon^{- 1 - \frac{183}{\log \mid \log \epsilon
\mid}}\right) \!\!- 2 m_v - c_v \!\right\} ~ (v \in S)$$ est de
cardinal majoré par:
$$
34 \epsilon^{- 1/2} {\mid \!\log \epsilon \!\mid}^{3/2}
\!\!\left(\log \!\mid \!\log \epsilon \!\mid\right)^{- 1/2}
\left[499 \epsilon^{- 1/2} \exp\!\left(\sqrt{\mid \!\log \epsilon
\!\mid . \log \!\mid \!\log \epsilon \!\mid}\right)\right]^{r} .
$$
\end{theorem}
\begin{theorem}[Deuxi\`eme théor\`eme principal]\label{c.2}
Dans la situation du théor\`eme \ref{c.1}, en remplaçant
l'hypoth\`ese $\epsilon < \frac{1}{15788}$ par $\epsilon \leq e^{-
4/r}$; l'ensemble des points $\mathbf x$ de $E(K)$ satisfaisant le
syst\`eme d'inégalités simultanées:
\begin{equation*}
{\mbox{dist}}_v(\mathbf x , \mathbf 0) < \exp\!\left\{\!\!-
\lambda_v \!\left(\epsilon h(\mathbf x) + (\eta + 5)
e^{\left(\frac{r}{4} \mid \log \epsilon \mid + 2\right) \left(\log
\mid \log \epsilon \mid + \log r + 16\right)}\!\right) \!- 2 m_v -
c_v \right\} (v \in S)
\end{equation*}
est de cardinal majoré par:
$$2 r^2 \epsilon^{- 1/2} {\mid\! \log \epsilon \!\mid}^{2} \!\left(\log r + \log \!\mid\! \log \epsilon \!\mid\! + 82\right)\!
\left(499 \epsilon^{- 1/2}\right)^{\!\!r} .$$
\end{theorem}
\begin{theorem}[Troisi\`eme théor\`eme principal]\label{c.3}
Sous les m\^emes hypoth\`eses que le théor\`eme \ref{c.1} et en
désignant de plus par ${E(K)}_{\mbox{tor}}$ le sous-groupe des
points de torsion de $E(K)$, par $\widehat{h}$ la hauteur de
Néron-Tate sur $E$ définie au paragraphe §$13.3$ et par
${\widehat{h}}_{\min}$ la plus petite valeur non nulle des
hauteurs de Néron-Tate des points de $E(K)$; l'ensemble des points
$\mathbf x$ de $E(K)$ satisfaisant le syst\`eme d'inégalités
simultanées:
$${\mbox{dist}}_v(\mathbf x , \mathbf 0) < \exp\left\{- \lambda_v \epsilon h(\mathbf x) - 2 m_v - c_v \right\} ~~~~~~~~~~ (v \in S)$$
est de cardinal majoré par:
\begin{equation*}
\begin{split}
\sharp {E(K)}_{\mbox{tor}} \!\left(1 + \frac{15 (\eta + 4)^{\frac{1}{2}} \epsilon^{- \frac{1}{2} - \frac{92}{\log \mid \log \epsilon \mid}}}{{{\widehat{h}}^{\frac{1}{2}}}_{\min}}\right)^{\!\!\!r} + 34 \epsilon^{- 1/2} {\mid \!\log \epsilon \!\mid}^{3/2}\!\!\left(\log \!\mid \!\log \epsilon \!\mid\right)^{- 1/2} & \\
&\!\!\!\!\!\!\!\!\!\!\!\!\!\!\!\!\!\!\!\!\!\!\!\!\!\!\!\!\!\!\!\!\!\!\!\!\!\!\!\!\!\!\!\!\!\!\!\!\!\!\!\!\!\!\!\!\!\!\!\!\!\!\!\!\!\!\!\!\!\!\!\!\!\!\!\!\!\!\!\times
\left[499 \epsilon^{- 1/2} \exp\!\left(\sqrt{\mid \!\log \epsilon
\!\mid . \log \!\mid \!\log \epsilon \!\mid}\right)\right]^{r} .
\end{split}
\end{equation*}
\end{theorem}
\begin{corollary}[du théor\`eme \ref{c.1}]\label{c.18}
Soit $E$ une courbe elliptique définie sur un corps de nombres $K$
de degré $D$, plongée dans ${\mathbb P}_2$ à la Weierstrass,
d'équation projective $Y^2 Z = 4 X^3 - g_2 X Z^2 - g_3 Z^3 ~ (g_2
, g_3 \in K)$ et d'élément neutre (en tant que groupe) le point à
l'infini $\mathbf 0$ représenté dans ${\mathbb P}_2$ par les
coordonnées projectives $(0 : 1 : 0)$. Soient aussi $r$ le rang de
Mordell-Weil de $E(K)$ que l'on suppose non nul, $S$ un ensemble
fini de places de $K$ et $\epsilon$ un réel strictement positif $<
\frac{1}{15788}$. Alors, l'ensemble des points $\mathbf x$ de
$E(K)$ satisfaisant l'inégalité:
$$\prod_{v \in S} {{\mbox{dist}}_v(\mathbf x , \mathbf 0)}^{\frac{[K_v : {\mathbb Q}_v]}{[K : \mathbb Q]}} ~\leq~ e^{- \epsilon h(\mathbf x) - 57 (\eta + 5) \\\epsilon^{- 1 - \frac{183}{\log \mid \log \epsilon \mid}}}$$
est de cardinal majoré par:
\begin{equation*}
\begin{split}
34 . 5^{card (S)}\!\left(\frac{2}{\epsilon}\right)^{\!\!1/2}
\!\!\left(\log\left(\frac{2}{\epsilon}\right)\right)^{\!\!3/2}
\!\!\left(\log \log\left(\frac{2}{\epsilon}\right)\right)^{\!\!- 1/2} & \\
&\!\!\!\!\!\!\!\!\!\!\!\!\!\!\!\!\!\!\!\!\!\!\!\!\!\!\!\!\!\!\!\!\!\!\!\!\!\!\!\!\!\!\!\!\!\!\!\!\!\!\!\!\times
\left[499\!\left(\frac{2}{\epsilon}\right)^{\!\!1/2}
\!\!\exp\!\left(\sqrt{\log\!\left(\frac{2}{\epsilon}\right)
\log\log\!\left(\frac{2}{\epsilon}\right)}\right)\right]^{r} .
\end{split}
\end{equation*}
\end{corollary}
\begin{corollary}[du théor\`eme \ref{c.2}]\label{c.19}
Dans la situation du corollaire \ref{c.18}, en remplaçant
l'hypoth\`ese $\epsilon < \frac{1}{15788}$ par $\epsilon \leq e^{-
4/r}$; l'ensemble des points $\mathbf x$ de $E(K)$ satisfaisant
l'inégalité:
$$\prod_{v \in S} {{\mbox{dist}}_v(\mathbf x , \mathbf 0)}^{\frac{[K_v : {\mathbb Q}_v]}{[K : \mathbb Q]}} ~\leq~ \exp\!\left\{\!\!- \epsilon h(\mathbf x) - (\eta + 5) e^{\left(\frac{r}{4} \mid \log \epsilon \mid + 2\right)
\left(\log \mid \log \epsilon \mid + \log r +
17\right)}\!\right\}$$ est de cardinal majoré par:
$$2 . 5^{card (S)} r^2 \sqrt{\frac{2}{\epsilon}} \left(\log\left(\frac{2}{\epsilon}\right)\right)^{\!\!2} \!\left(\log r + \log\log\left(\frac{2}{\epsilon}\right) + 82\right)\!
\left(499\sqrt{\frac{2}{\epsilon}}\right)^{\!\!r} .$$
\end{corollary}
\begin{corollary}[du théor\`eme \ref{c.3}]\label{c.20}
Sous les hypoth\`eses du corollaire \ref{c.18}, l'ensemble des
points $\mathbf x$ de $E(K)$ satisfaisant l'inégalité:
$$\prod_{v \in S} {{\mbox{dist}}_v(\mathbf x , \mathbf 0)}^{\frac{[K_v : {\mathbb Q}_v]}{[K : \mathbb Q]}} ~\leq~ e^{- \epsilon h(\mathbf x) - 2 \eta - 16}$$
est de cardinal majoré par:
$$
5^{card (S)} \left\{\sharp {E(K)}_{\mbox{tor}}\!\left(1 + \frac{15
(\eta +
4)^{\frac{1}{2}}\left(\frac{2}{\epsilon}\right)^{\frac{1}{2} +
\frac{92}{\log\log\left(\frac{2}{\epsilon}\right)}}}{{{\widehat{h}}_{\min}}^{\!\frac{1}{2}}}\right)^{\!\!\!r}
\right.
$$
$$+ 34 \left(\frac{2}{\epsilon}\right)^{\!\!\frac{1}{2}} \!\left(\log\!\left(\frac{2}{\epsilon}\right)\right)^{\!\!\frac{3}{2}} \!\left(\log\log\!\left(\frac{2}{\epsilon}\right)\right)^{\!\!\!- \frac{1}{2}}\left.\phantom{\left(1 + \frac{15 (\eta + 4)^{\frac{1}{2}}\left(\frac{2}{\epsilon}\right)^{\frac{1}{2} + \frac{92}{\log\log\left(\frac{2}{\epsilon}\right)}}}{{{\widehat{h}}_{\min}}^{\!\frac{1}{2}}}\right)^{\!\!\!r}}
\!\!\!\!\!\!\!\!\!\!\!\!\!\!\!\!\!\!\!\!\!\!\!\!\!\!\!\!\!\!\!\!\!\!\!\!\!\!\!\!\!\!\!\!\!\!\!\!\!\!\!\!\!\!\!\!\!\!\!\!\!\!\!\!\!\!\!\!\!\!\!\!\!\!\!\!\!\!\!\!\!\!\!\!\!\!\!\!\!\!\left[499
\left(\frac{2}{\epsilon}\right)^{\!\!\frac{1}{2}}
\exp\!\left(\sqrt{\log\!\left(\frac{2}{\epsilon}\right)
\log\log\!\left(\frac{2}{\epsilon}\right)}\right)\right]^{\!\!r}\right\}
.
$$
\end{corollary}

Avant de se lancer dans les détails, décrivons grosso-modo les différentes étapes nous permettant d'aboutir aux résultats: \\
Nous commençons (§$3$ qui suit) par paramétriser $E$ au voisinage
d'un point quelconque $(x : y : z)$ d'une certaine carte de $E$,
en prenant $t = \frac{x}{y}$ comme param\`etre et en exprimant
$\frac{z}{y}$ comme fonction enti\`ere en $t$. Au §$4$ nous
introduisons un entier $m \geq 2$ et des entiers strictement
positifs $a_1 , \dots , a_{m - 1}$, à l'aide desquels nous
plongeons $E^m$ dans $E^m \times E^{m - 1}$ comme suit:
$$
\begin{array}{rcl}
E^m & \stackrel{\psi_{\underline a}}{\hookrightarrow} & E^m \times
E^{m - 1} \\
 (x_1 , \ldots , x_m) & \longmapsto & (x_1 , \ldots , x_m , a_1 x_1 - x_m , \ldots , a_{m - 1} x_{m - 1} - x_m)
\end{array}~,$$
puis nous plongeons $E^m \times E^{2 m - 1}$ dans ${\mathbb
P}_{2}^{2 m - 1}$ (plongement de Weierstrass) et nous appelons
$\varphi_{\underline a}$ le plongement composé. Nous calculons au
lemme $2.4.1$ -en utilisant le théor\`eme de Wirtinger- les
différents multidegrés de $\varphi_{\underline a}(E^m)$. Nous
déduisons naturellement de la paramétrisation locale de $E$, une
paramétrisation locale $\Omega_{\underline a}$ de
$\varphi_{\underline a}(E^m) \plongement {\mathbb P}_{2}^{2 m -
1}$ sur une certaine carte $\mathcal{\alpha}$ contenant
${\{\mathbf 0\}}^{2 m - 1}$. Nous avons besoin pour celà d'un
syst\`eme complet de familles de formes représentant l'addition
sur $E$, lequel est donné dans [La-Ru], et d'une famille de forme
représentant la multiplication d'un point de $E$ par un entier
positif donné. Pour cette derni\`ere, nous n'avons pas trouvé de
référence donnant des estimations totalement explicites des degrés
et hauteurs de ces formules, nous avons donc repris les calculs en
suivant [La3], ce qui nous a amené au théor\`eme $2.13.2$
(formulaire). Nous avons défini des opérateurs de dérivations
$\partial^{(i_1 , \dots , i_m)}$$(i_1 , \dots , i_m \in \mathbb
N)$ sur l'anneau des coordonnées de $\varphi_{\underline a}(E^m)$
tels que pour toute forme $P_1$ sur $\varphi_{\underline a}(E^m)$,
l'annulation du coefficient $u_{1}^{i_1} \dots u_{m}^{i_m}$ (o\`u
$u_1 , \dots , u_m$ sont les param\`etres) dans la série
$\Omega_{\underline a}(P_1)$ en un certain point équivaut à
l'annulation de la forme $\partial^{(i_1 , \dots , i_m)}P_1$ au
m\^eme point. Nous estimons ensuite dans le corollaire $2.5.2$ les
degrés et hauteurs des dérivées d'une forme donnée sur
$\varphi_{\underline a}(E^m)$, en fonction du degré et de la
hauteur de cette forme et de $E$. Au §$6$ nous introduisons des
param\`etres positifs $0 < \epsilon_0 < 1 / 2$, $\epsilon_1 > 0$
et $\delta \in {\mathbb N}^*$ (destiné à tendre vers l'infini)
avec $\epsilon_0 \delta \in {\mathbb N}^*$ assez grand et:
\begin{gather}
\frac{m - 1}{m!} \left(\frac{7}{3}\right)^m
\frac{\epsilon_{1}^{m}}{\epsilon_0 (m + \epsilon_0) (1 +
\epsilon_0)^{m - 2}} \leq \frac{1}{2} . \tag{$2'$}
\end{gather}
Nous construisons par le lemme de Siegel usuel, une forme non
identiquement nulle $P$ sur $\varphi_{\underline a}(E^m)$ de
multidegré $(\epsilon_0 \delta a_{1}^{2} , \dots , \epsilon_0
\delta a_{m}^{2} , \delta , \dots , \delta)$, s'annulant en
${\{\mathbf 0\}}^{2 m - 1}$ avec une multiplicité définie par le
dessous d'escalier de ${\mathbb N}^m$:
$$T_{\delta} := \left\{(\tau_1 , \dots , \tau_m) \in {\mathbb N}^m / \frac{\tau_1}{a_{1}^{2}} + \dots + \frac{\tau_{m - 1}}{a_{m - 1}^{2}} +
\frac{\tau_m}{m - 1} \leq 7 \epsilon_1 \delta\right\} $$
et qui soit de hauteur $\ll \delta a_{1}^{2}$ (o\`u dans tout ce qui suit $\ll$ veut dire inférieur ou égal à une constante multiplicative pr\`es
qui ne dépend que de $E$). \\
En effet, dans le syst\`eme linéaire de Siegel, les inconnues sont les coefficients de la forme $P$ à construire et le nombre d'équations est égal
au cardinal de l'ensemble $T_{\delta}$. Le nombre d'inconnues est alors la valeur de la fonction de Hilbert de l'idéal
$\mathfrak{I}(\varphi_{\underline a}(E^m))$ en $(\epsilon_0 \delta a_{1}^{2} , \dots , \epsilon_0 \delta a_{m}^{2} , \delta , \dots , \delta)$;
comme $\delta$ et $\epsilon_0 \delta$ sont supposés assez grands, cette valeur co{\sf\"\i}ncide avec la valeur d'un polyn\^ome de
$Q[\Xsoul_1 , \dots , \Xsoul_m , \Ysoul_1 , \dots ,$ \\ $\Ysoul_{m - 1}]$ en $(\epsilon_0 \delta a_{1}^{2} , \dots , \epsilon_0 \delta a_{m}^{2} , \delta , \dots , \delta)$ dont la partie homog\`ene dominante est connue explicitement en fonction des multidegrés de $\varphi_{\underline a}(E^m)$ lesquels sont calculés par le lemme $2.4.1$. On estime le nombre d'inconnues, puis $\sharp T_{\delta}$ est grossi\`erement estimé par ${\rm{vol}} ~\!\!T_{\delta}$, qu'on calcule facilement, et afin d'appliquer le lemme de Siegel, on vérifie gr\^ace à $(2')$ que le nombre d'équations est strictement inférieur au nombre d'inconnues. \\
Au §$7$, nous introduisons un nouveau param\`etre $0 < \alpha < 1$
et des points ${\mathbf x}_1 , \dots , {\mathbf x}_m$ de $E(K)$,
ordonnés par ordre croissant de leurs hauteurs, que nous supposons
contenus dans un petit c\^one d'angle $\leq \arccos(1 - \alpha /
4)$ de l'espace euclidien $E(K) \otimes_{\mathbb Z} \mathbb R$ et
de hauteurs assez espacées. Nous posons $a_i := [\nb{{\mathbf
x}_m} / \nb{{\mathbf x}_i}]$$(i = 1 , \dots , m)$ de sorte que les
points ${\mathbf y}_i = a_i {\mathbf x}_i - {\mathbf x}_m$$(i = 1
, \dots , m - 1)$ soient de hauteurs assez petites en comparaison
avec les points ${\mathbf x}_i$ (la géométrie euclidienne nous
donne plus précisement $\widehat{h}({\mathbf y}_i) \leq \alpha
(a_{i}^{2} \widehat{h}({\mathbf x}_i) + \widehat{h}({\mathbf
x}_m))$). Nous posons $\mathbf x = ({\mathbf x}_1 , \dots ,
{\mathbf x}_m)$ et $\mathbf y = ({\mathbf x}_1 , \dots , {\mathbf
x}_m , {\mathbf y}_1 , \dots , {\mathbf y}_{m - 1})$ et nous
considérons la forme translatée $\tau_{- \mathbf y}^* P$ de notre
forme $P$ construite au §$6$, par le point $- \mathbf y$. Celle-ci
s'annule en $\mathbf y$ avec la multiplicité définie par le
dessous d'escalier $T_{\delta}$ et nous pouvons estimer les degrés
et les hauteurs des formes dérivées $\partial^{(i_1 , \dots ,
i_m)} \tau_{- \mathbf y}^{*} P$. Au §$8$, nous supposons que les
param\`etres $\epsilon_0 , \epsilon_1$ et $\alpha$ sont liés par
la relation:
\begin{gather}
4 m (\epsilon_0 + 2 \alpha) \leq \epsilon \epsilon_1 , \tag{$1'$}
\end{gather}
Nous introduisons un ensemble fini $S$ de places de $K$ et
${(\lambda_v)}_{v \in S}$ une famille de réels positifs
satisfaisants:
$$\sum_{v \in S} \frac{[K_v : {\mathbb Q}_v]}{[K : \mathbb Q]} \lambda_v = 1 .$$
Nous supposons que les points ${\mathbf x}_i$ satisfont le
syst\`eme d'inégalités simultanées:
\begin{gather}
\forall v \in S ,~ d_v({\mathbf x}_i , \mathbf 0) \ll e^{-
\lambda_v \epsilon h({\mathbf x}_i)} , \tag{$\rm{S.S}$}
\end{gather}
qu'ils sont contenu dans un petit c\^one d'angle $\leq \arccos(1 -
\alpha / 4)$ de $E(K) \otimes_{\mathbb Z} \mathbf R$ et que
$\widehat{h}({\mathbf x}_1) \gg \frac{1}{\epsilon \epsilon_1}$.
Dans l'extrapolation, l'hypoth\`ese principale $(\rm{S.S})$,
l'hypoth\`ese $(1')$ et l'hypoth\`ese $\widehat{h}({\mathbf x}_1)
\gg \frac{1}{\epsilon \epsilon_1}$ entra{\sf\^\i}nent que la forme
$\tau_{- \mathbf y}^{*}P$ s'annule en $(\mathbf 0 , \dots ,
\mathbf 0)$ avec la multiplicité définie par le dessous d'escalier
$T_{\delta / 2}$ (ceci se démontre en comparant les coefficients
des deux séries $\tau_{- \mathbf y}^{*}P({\underline x}_1 , \dots
, {\underline x}_m)$ et $\tau_{- \mathbf y}^{*}P(\underline 0 ,
\dots , \underline 0)$ et en utilisant la formule du produit pour
chaque coefficient de la série $\tau_{- \mathbf y}^{*}P(\underline
0 , \dots , \underline 0)$ correspondant à un exposant $\underline
i \in T_{\delta / 2}$). Or, ceci est équivalent à dire que notre
forme $P$ s'annule en $(- \mathbf y)$ avec la multiplicité définie
par le dessous d'escalier $T_{\delta / 2}$. Et, en retirant $P$ en
une forme $Q$ sur $E^m$, $Q$ s'annule en $(- \mathbf x)$ avec la
m\^eme multiplicité et de plus $h(Q) \ll \delta a_{1}^{2}$. Au
§$9$ nous arrivons à ce que nous appelons ``inégalité à la
Vojta''; en supposant de plus que les points ${\mathbf x}_1 ,
\dots , {\mathbf x}_m$ sont de hauteurs un peu plus grandes
$\widehat{h}({\mathbf x}_1) \geq (6 m^2 / \epsilon_1)^m$ et sont
un peu plus espacés $\widehat{h}({\mathbf x}_i) \gg m (6 m^2 /
\epsilon_1)^m \widehat{h}({\mathbf x}_{i - 1})$ (cette condition
d'espacement de hauteurs donne l'hypoth\`ese principale du
théor\`eme du produit, puisque les $a_i$ sont inversement
proportionels aux hauteurs des points ${\mathbf x}_i$) et en
appliquant le théorème du produit de \cite{300}, on obtient une
contradiction avec le fait que $\widehat{h}({\mathbf x}_1)$ soit
assez grand. Nous déduisons alors le théor\`eme $2.9.1$ qui anonce
que pour $\epsilon_0 , \epsilon_1$ et $\alpha$ des réels positifs
satisfaisant les contraintes $(1')$ et $(2')$ et pour des points
${\mathbf x}_1 , \dots , {\mathbf x}_m$ contenus dans un petit
c\^one d'angle $\leq \arccos(1 - \alpha / 4)$, de hauteurs:
\\ $\widehat{h}({\mathbf x}_m) \geq \dots \geq \widehat{h}({\mathbf
x}_1) \geq (6 m^2 / \epsilon_1)^m$ et satisfaisant le syst\`eme
simultané $(\rm{S.S})$, on a l'une au moins des inégalités:
$\widehat{h}({\mathbf x}_i) < m (6 m^2 / \epsilon_1)^m
\widehat{h}({\mathbf x}_{i - 1})$. Apr\`es cela, nous choisissons
les param\`etres $\epsilon_0 , \epsilon_1$ et $\alpha$ en fonction
de $\epsilon$ et $m$ de façon à satisfaire les contraintes $(1')$
et $(2')$. \`A une constante absolue ($< 1$) multiplicative pr\`es
on prend: $\epsilon_0 \simeq \epsilon^{\frac{m}{m - 1}}$,
$\epsilon_1 \simeq m \epsilon^{\frac{1}{m - 1}}$ et $\alpha \simeq
\epsilon^{\frac{m}{m - 1}}$ et nous obtenons l'inégalité de la
hauteur à la Vojta qui s'énonce: Pour ${\mathbf x}_1 , \dots ,
{\mathbf x}_m$ des points de $E(K)$ contenus dans un petit c\^one
de $E(K) \otimes_{\mathbb Z} \mathbb R$ d'angle $\leq \arccos(1 -
\frac{1}{30976} \epsilon^{\frac{m}{m - 1}})$, qui sont de hauteurs
$\geq (cte . m)^m \epsilon^{- \frac{m}{m - 1}}$ (avec $cte$
désigne une constante absolue) et satisfaisant $(\rm{S.S})$, on a
l'une au moins des inégalités $(1< i \leq m)$:
\begin{gather}
\widehat{h}({\mathbf x}_i) < (cte . m)^m \epsilon^{- \frac{m}{m -
1}} \widehat{h}({\mathbf x}_{i - 1}) . \tag{$\rm{I.V}$}
\end{gather}
L'inégalité de la hauteur à la Mumford est grosso-modo l'inégalité
dans l'autre sens pour $m = 2$. Elle s'énonce: Pour ${\mathbf
x}_1$ et ${\mathbf x}_2$ deux points de $E(K)$ contenus dans un
petit c\^one de $E(K) \otimes_{\mathbb Z} \mathbb R$ d'angle $\leq
\arccos(1 - \epsilon / 8)$, qui sont de hauteurs
$\widehat{h}({\mathbf x}_2) \geq \widehat{h}({\mathbf x}_1) \gg
\frac{1}{\epsilon}$ et qui satisfont $(\rm{S.S})$, on a:
\begin{gather}
\widehat{h}({\mathbf x}_2) \geq \left(1 + \frac{1}{3}
\sqrt{\epsilon}\right) \widehat{h}({\mathbf x}_1) .
\tag{$\rm{I.M}$}
\end{gather}
On obtient ce théor\`eme (théor\`eme $2.10.1$) en suivant les m\^emes étapes que pour le théor\`eme $2.9.1$, cependant la preuve est ici beaucoup plus simple du fait qu'on n'utilise ni le lemme de Siegel pour construire les fonctions auxiliaires, ni un lemme de zéros à la fin pour conclure. Voilà bri\`evement comment on fait: \\
Nous procédons par l'absurde, nous supposons que ${\mathbf x}_1$
et ${\mathbf x}_2$$({\mathbf x}_1 \neq {\mathbf x}_2)$ satisfont
toutes les hypoth\`eses du théor\`eme mais ne satisfont pas
$(\rm{I.M})$, ainsi le point $\mathbf y = {\mathbf x}_1 - {\mathbf
x}_2$ sera de hauteur tr\`es petite en comparaison avec ${\mathbf
x}_1$ et ${\mathbf x}_2$ (plus précisément la géométrie
euclidienne nous donne $\widehat{h}(\mathbf y) <
\frac{\epsilon}{2} \min\{\widehat{h}({\mathbf x}_1) ,
\widehat{h}({\mathbf x}_2)\}$). Soit $\underline D = (D_0 , D_1 ,
D_2)$ une famille de formes représentant la différence sur $E$
dans une certaine carte de $E^2$ contenant $\{\mathbf 0\} \times
\{\mathbf 0\}$ et $\{{\mathbf x}_1\} \times \{{\mathbf x}_2\}$.
Nous introduisons les deux fonctions auxiliaires sur $E^2$:
\begin{align}
Q_1(\Xsoul_1 , \Xsoul_2) &:= D_0\left(\underline{D}(\Xsoul_1 , \Xsoul_2) , \underline{y}\right) \notag \\
Q_2(\Xsoul_1 , \Xsoul_2) &:= D_2\left(\underline{D}(\Xsoul_1 ,
\Xsoul_2) , \underline{y}\right) \notag
\end{align}
(o\`u $\underline{y}$ désigne un représentant dans ${\mathbb{P}}_2$ du point $\mathbf{y}$). \\
$Q_1$ et $Q_2$ s'annulent clairement en $({\mathbf x}_1 , {\mathbf x}_2)$. Dans l'extrapolation,
si on suppose que l'une des formes $Q_j$$(j = 1 , 2)$ ne s'annule pas en $(\mathbf 0 , \mathbf 0)$,
en appliquant la formule du produit au nombre non nul $Q_j(\underline 0 , \underline 0)$ de $K$ et en
tenant compte de $(\rm{S.S})$, on aboutit à une contradiction avec le fait que $\widehat{h}(\mathbf y)$
est tr\`es petit relativement à $\widehat{h}({\mathbf x}_1)$ et $\widehat{h}({\mathbf x}_2)$. On en déduit
que $Q_1$ et $Q_2$ doivent s'annuler toutes les deux en $(\mathbf 0 , \mathbf 0)$, ce qui entra{\sf\^\i}ne
$\mathbf y = \mathbf 0$, puis ${\mathbf x}_1 = {\mathbf x}_2$ qui est une contradiction. D'o\`u le théor\`eme
$2.10.1$. \\
En mettant ensemble les deux inégalités $(\rm{I.V})$ et $(\rm{I.M})$, on a un décompte de l'ensemble des points de
$E(K)$ satisfaisant $(\rm{S.S})$ se situant dans un petit c\^one de $E(K) \otimes_{\mathbb Z} \mathbb R$ et qui
sont de hauteurs assez grandes. En effet, en supposant qu'on a $\l$ tels points ${\mathbf x}_1 , \dots ,
{\mathbf x}_{\l}$, ordonnés selon l'ordre croissant de leurs hauteurs, on partage ces $\l$ points en $m$ paquets de $k$
points ($k := [\frac{\l - 1}{m - 1}]$, en oubliant éventuellement quelques un des derniers points) et on consid\`ere
dans chacun de ces paquets le point de plus petite hauteur. On désigne ces derniers par ${\mathbf y}_1 , \dots ,
{\mathbf y}_m$. D'apr\`es le théor\`eme $2.9.1$, $(\rm{I.V})$ est satisfaite pour un certain ${\mathbf y}_j$
$(j \in \{2 , \dots , m\})$ et d'apr\`es le théor\`eme $2.10.1$, $(\rm{I.M})$ est satisfaite pour chaque point
${\mathbf x}_i$ tels que les points ${\mathbf x}_i$ et ${\mathbf x}_{i - 1}$ soient dans l'intervalle
$[{\mathbf y}_{j - 1} , {\mathbf y}_j]$. Ainsi $\widehat{h}({\mathbf y}_j)$ est majoré et minoré en fonction de
$\widehat{h}({\mathbf y}_{j - 1})$ et cette comparaison donne une majoration pour $k$ puis pour $\l$. \\
Pour conclure à nos théor\`emes principaux (théor\`emes $2.2.1$,
$2.2.2$ et $2.2.3$), nous recouvrons l'espace euclidien $E(K)
\otimes_{\mathbb Z} \mathbb R \simeq {\mathbb R}^r$ (o\`u $r$ est
le rang du groupe de Mordell-Weil de $E(K)$) par un nombre fini de
petits c\^ones, et pour faire le décompte des points de petites
hauteurs de $E(K)$ satisfaisant $(\rm{S.S})$, soit on affaiblie
$(\rm{S.S})$ de façon à ce que le point $\mathbf 0$ soit l'unique
point de $E(K)$ de petite hauteur sous $(\rm{S.S})$ ou bien on
compte tous les points à petite hauteur sans tenir compte de
$(\rm{S.S})$.
\section{Paramétrisation locale de $\widetilde E$}
On paramétrise localement la courbe $\widetilde E$ en utilisant
des séries enti\`eres. Le th\'eor\`eme de Bézout montre que
l'ensemble des points $(t_1 , t_2 )$ de $\widetilde{E}(K)$ pour
lesquels $\widetilde{\Delta}(t_1 , t_2) = 0$ est fini et comporte
aux maximum 6 points. Ainsi la proposition qui suit fournit pour
tout point g\'en\'eral donn\'e $(t_1 , t_2)$ de $\widetilde{E}(K)$
(c'est-\`a-dire pour tout point $(t_1 , t_2)$ de
$\widetilde{E}(K)$ satisfaisant $\widetilde{\Delta}(t_1 , t_2)
\neq 0$) une paramétrisation de $\widetilde E$ au voisinage de ce
point (en spécialisant $(X , Z)$ en $(t_1 , t_2)$ dans le
morphisme $\widetilde{\tau}$ ci-dessous).
\begin{proposition}\label{c.4}
 Il existe un monomorphisme d'anneaux: $$ \widetilde{\tau} : A(\widetilde E)
\longrightarrow K(\widetilde E)[[t]]$$ tel que: $$\widetilde{\tau}
(X)
 = X + t ~~\text{et}~~ \widetilde{\tau} (Z) = Z + \sum_{\ell = 1}^{\infty} \frac {\widetilde{{\partial}^{\ell} Z}}
 {{\widetilde{\Delta}}^{2 \ell - 1}} t^{\ell}$$ o\`u
 les $\widetilde{{\partial}^
 {\ell} Z} , \ell \in {\mathbb N}^*$, sont des polyn\^omes de
 $ K[X , Z]$ satisfaisant pour $T \geq 1$ et pour toute place $v$ de $K$:
  $$\max (d°\widetilde{{\partial}^{\ell} Z}~ ; ~\ell = 0 , \cdots , T)~\leq~3T - 1$$
 $$h_v(\widetilde{{\partial}^{\ell} Z}~;~\ell = 0 , \cdots , T) ~\leq~ \begin{cases}
(2 T - 1)m_v & \text{si $v$ est finie} \\
(2 T - 1)(m_v + 4) & \text{si $v$ est infinie}
\end{cases}.$$
Par conséquent pour tout $T \geq 1$ on a:
$$\widetilde h (\widetilde{{\partial}^{\ell} Z}~;~\ell = 0 , \cdots , T)~\leq~(2T - 1)(\eta + 4) .$$
\end{proposition}
{\bf Démonstration.---}
On consid\`ere -durant toute cette démonstration- $\widetilde{G}$ et $\widetilde{\Delta}$ comme des éléments de l'anneau de polyn\^omes $K[\Psi , \Upsilon]$. \\
La relation liant $X$ et $Z$ (considérés comme éléments de
$A(\widetilde{E})$) est évidemment: $\widetilde{G}(X , Z) = 0$,
donc pour que $\widetilde{\tau}$ soit un monomorphisme, il faut et
il suffit qu'on ait: $\widetilde{G}(\widetilde{\tau}(X) ,
\widetilde{\tau}(Z)) = \widetilde{\tau}(\widetilde{G}(X , Z)) =
0$, c'est-à-dire: $\widetilde{G}(X + t , Z + \sum_{\ell =
1}^{\infty} \frac{\widetilde{\partial^{\ell}
Z}}{{\widetilde{\Delta}}^{2 \ell - 1}}) = 0$. Ce qui donne par le
developpement de Taylor:
$$
\sum_{(h , j) \in {\mathbb N}^2 \setminus \{(0 ,
0)\}}^{}\frac{1}{h!j!} \frac{{\partial}^{h + j}
\widetilde{G}}{{\partial \Psi}^h {\partial \Upsilon}^j} (X , Z)
t^h \left(\sum_{\ell = 1}^{\infty}
\frac{\widetilde{{\partial}^{\ell} Z}}{{\widetilde{\Delta}}^{2
\ell - 1}} t^{\ell}\right)^j = 0 .
$$
En développant la série du membre de gauche de cette relation et
en annulant chacun de ses coefficients, on obtient la relation:
\begin{equation}
{\partial}^{\ell} Z ~=~ - \sum_{\begin{array}{c}
   \scriptstyle (h , j) \in {\mathbb N}^2 \setminus \{(0 , 0) , (0 , 1)\} \\
   \scriptstyle k_1 , \ldots , k_j \in {\mathbb N}^* \\
   \scriptstyle h + k_1 + \cdots + k_j = \ell
   \end{array}} {\widetilde{\Delta}}^{2h + j - 2} \frac{1}{h! j!} \frac{{\partial}^{h + j} \widetilde{G}}
   {{\partial \Psi}^h {\partial \Upsilon}^j}(X , Z) \widetilde{{\partial}^{k_1} Z} \dots \widetilde{{\partial}^{k_j} Z} . \label{10.1}
\end{equation}
Cette derni\`ere relation nous permet de calculer les polyn\^omes
$\widetilde{\partial^{\ell} Z}$ de proche en proche. Par exemple,
pour $\ell = 1$, elle donne:
$$
\widetilde{\partial^1 Z} = - \frac{\widetilde{\partial
G}}{\partial \Psi} = 12 \Psi^2 - g_2 \Upsilon^2 .
$$
Gr\^ace à la relation (\ref{10.1}), on démontre aisément (par récurrence sur $\ell$) les estimations de la proposition \ref{c.4} pour les degrés et les hauteurs logarithmiques $v$-adiques des polyn\^omes $\widetilde{\partial^{\ell} Z}$, dans le cas o\`u $v$ est une place finie de $K$ (remarquer que lorsque $v$ est une place finie de $K$, on a: $h_v(\widetilde{\Delta}) \leq m_v$ et $\forall (h , j) \in {\mathbb N}^2$: $h_v(\frac{1}{h! j!} \frac{\partial^{h + j} \widetilde{G}}{{\partial \Psi}^h {\partial \Upsilon}^j}) \leq m_v$). \\
Pour obtenir l'estimation de la proposition \ref{c.4} pour les
hauteurs logarithmiques $v$-adiques des polyn\^omes
$\widetilde{\partial^{\ell} Z}$, quand $v$ est une place infinie
de $K$; nous sommes amenés à utiliser un autre type de relation,
qui est:
\begin{equation}
\begin{split}
\widetilde{{\partial}^{\ell + 1} Z} = \frac{1}{\ell + 1}
\left[\frac{\partial (\widetilde{{\partial}^{\ell} Z})}{\partial
\Psi} {\widetilde{\Delta}}^2 - \frac{\partial
(\widetilde{{\partial}^{\ell} Z})}
{\partial \Upsilon} \frac{\partial{\widetilde{G}}}{\partial \Psi} \widetilde{\Delta} \right.& \\
&\!\!\!\!\!\!\!\!\!\!\!\!\!\!\!\!\!- \left. (2 \ell - 1)
\widetilde{{\partial}^{\ell} Z} \left(\widetilde{\Delta} .
\frac{\partial{\widetilde{\Delta}}} {\partial \Psi} -
\frac{\partial{\widetilde{G}}}{\partial \Psi} .
\frac{\partial{\widetilde{\Delta}}}{\partial
\Upsilon}\right)\right]
\end{split} \label{10.2}
\end{equation}
pour tout $\ell \geq 1$. Cette relation s'établit en montrant
d'abord par récurrence que si deux fonctions enti\`eres (sur un
ouvert de $\mathbb C$) $F_1$ et $F_2$ sont liées par une équation
du type $P(F , G) = 0$, pour un certain polyn\^ome $P$ de
$\mathbb{C}[\Psi , \Upsilon]$ vérifiant $\frac{\partial
P}{\partial \Upsilon} \not\equiv 0$, alors on a $\frac{1}{\ell!}
\frac{d^{\ell} F_2}{{d F_1}^{\ell}} = f_{\ell}(F_1 , F_2)$ pour
tout $\ell \geq 1$, o\`u ${(f_{\ell})}_{\ell \geq 1}$ est la suite
de fonctions de deux variables définie par:
$$
\left\{\!\!
\begin{array}{rcl}
f_1 & = & - \frac {\partial{P}}{\partial \Psi} /
\frac {\partial{P}}{\partial \Upsilon} \\
f_{\ell + 1}(\Psi , \Upsilon) & = & \frac{1}{\ell + 1}
\left[\frac{\partial{f_{\ell}}}
{\partial \Psi}(\Psi , \Upsilon) + f_1(\Psi , \Upsilon) \frac{\partial{f_{\ell}}}{\partial \Upsilon}(\Psi , \Upsilon)\right] ,~~\forall \ell \geq 1 \\
\end{array}
\right. .$$
Ensuite, on applique ce fait aux deux fonctions $F_1 = \widetilde{\tau}(X) = X + t$ et $F_2 = \widetilde{\tau}(Z)$ en remarquant qu'on a: $\frac{1}{{\ell}!} \frac{d^{\ell} \widetilde{\tau}(Z)}{{d \widetilde{\tau}(X)}^{\ell}} = \frac{\widetilde{{\partial}^{\ell} Z}}{{\widetilde{\Delta}}^{2 \ell - 1}} , ~~\forall \ell \geq 1$. \\
Posons maintenant,
\begin{equation*}
\begin{split}
R &:= \widetilde{\Delta} \frac{\partial \widetilde{\Delta}}{\partial \Psi} - \frac{\partial \widetilde{G}}{\partial \Psi} \frac{\partial \widetilde{\Delta}}{\partial \Upsilon} \\
&= 24 g_2 \Psi^3 + 72 g_3 \Psi^2 \Upsilon + 2 g_{2}^{2} \Psi
\Upsilon^2 + 2 g_2 \Upsilon .
\end{split}
\end{equation*}
Etant donné une place infinie $v$ de $K$ et un entier positif
$\ell \geq 1$, les majorations -citées au §$2$- de la hauteur
$H_v$ d'une somme, d'un produit de polyn\^omes ou encore de
dérivées de ployn\^omes, permettent d'estimer le deuxi\`eme membre
de la relation (\ref{10.2}) en fonction de $\ell , m_v$ et
$H_v(\widetilde{\partial^{\ell} Z})$. On obtient:
\begin{equation}
\begin{split}
H_v(\widetilde{\partial^{\ell + 1} Z}) &\leq \frac{1}{\ell + 1}\left[H_v\!\left(\!\frac{\partial(\widetilde{\partial^{\ell} Z})}{\partial \Psi}\!\right) {H_v(\widetilde{\Delta})}^2 \mathcal{N}(\widetilde{\Delta}) \right. \\
&\quad\left.+ H_v\!\left(\!\frac{\partial(\widetilde{\partial^{\ell} Z})}{\partial \Upsilon}\!\right) \!H_v\!\left(\!\frac{\partial \widetilde{G}}{\partial \Psi}\!\right) \!H_v(\widetilde{\Delta}) \mathcal{N}\!\left(\!\frac{\partial \widetilde{G}}{\partial \Psi}\!\right) \!+ (2 \ell - 1) H_v(\widetilde{\partial^{\ell} Z}) H_v(R) \mathcal{N}(R)\right] \\
&\leq \frac{1}{\ell + 1}\left[(3 \ell - 1) H_v(\widetilde{\partial^{\ell} Z}) (3 M_v)^2 . 3 + (3 \ell - 1) H_v(\widetilde{\partial^{\ell} Z}) . 12 M_v . 3 M_v . 2 \right. \\
&\quad \left.+ (2 \ell - 1) H_v(\widetilde{\partial^{\ell} Z}) . 72 M_{v}^{2} . 4\right] \\
&\leq (30 M_v)^2 H_v(\widetilde{\partial^{\ell} Z}) ,
\end{split} \label{10.3}
\end{equation}
o\`u dans cette série d'inégalités, on a désigné par $\mathcal N$ l'application associant à tout polyn\^ome, le nombre de mon\^ome intervenant dans son écriture canonique. De plus, on a majoré $H_v(\partial(\widetilde{\partial^{\ell} Z}) / \partial \Psi)$ et $H_v(\partial(\widetilde{\partial^{\ell} Z}) / \partial \Upsilon)$ par $d°\widetilde{\partial^{\ell} Z} . H_v(\widetilde{\partial^{\ell} Z}) \leq (3 \ell - 1) H_v(\widetilde{\partial^{\ell} Z})$ (d'apr\`es l'estimation de la proposition \ref{c.4} -déja démontrée- pour les degrés des polyn\^omes $\widetilde{\partial^{\ell} Z}$), $H_v(\widetilde{\Delta})$ par $3 M_v$, $H_v(\partial{\widetilde{G}} / \partial{\Psi})$ par $12 M_v$, $H_v(R)$ par $72 M_{v}^{2}$, $\mathcal{N}(\widetilde{\Delta})$ par $3$, $\mathcal{N}(\partial{\widetilde{G}} / \partial{\Psi})$ par $2$ et $\mathcal{N}(R)$ par $4$. \\
Par suite, comme: $H_v(\widetilde{\partial^1 Z}) = H_v(12 \Psi^2 -
g_2 \Upsilon^2) \leq 30 M_v$, la récurrence sur $\ell$, utilisant
(\ref{10.3}), donne:
$$
H_v(\widetilde{\partial^{\ell} Z}) ~\!\leq~\! (30 M_v)^{2 \ell -
1} ~~~~ (\forall \ell \geq 1) .
$$
En prenant finalement les logarithmes des deux termes de cette
derni\`ere inégalité et en majorant $\log{30}$ par $4$, on aboutit
à l'estimation restante de la proposition \ref{c.4}. Ce qui
ach\`eve cette démonstration.  $~~~~\blacksquare$\vspace{1mm}\\
{\bf Avertissement.---} On posera pour tout $\ell \in \mathbb N$,
$$f(\ell) ~:=~ \max(0 , 2 \ell - 1) ~~\mbox{et}~~ g(\ell) := 2 f(\ell) - \ell = \max(0 , 3 \ell - 2) .$$ On posera aussi $\widetilde{\partial^{\ell}
(X)} = \left\{
\begin{array}{ccc}
X & \mbox {si} & \ell = 0 \\ 1 & \mbox {si} & \ell = 1 \\ 0 &
\mbox {si} & \ell \geq 2 \\
\end{array}
\right.$ et $\widetilde{{\partial}^0 Z} = Z$. La proposition
\ref{c.4} donne alors: $$\widetilde{\tau}(X) = \sum_{\ell = 0
}^{\infty} \widetilde{\partial^{\ell} X}~t^{\ell} ~~\mbox{et}~~
\widetilde{\tau}(Z) = \sum_{\ell = 0}^{\infty}
\frac{\widetilde{\partial^{\ell}
Z}}{{\widetilde{\Delta}}^{f(\ell)}} t^{\ell} .$$ Avec toutes ces
notations on peut énoncer :
\begin{corollary}\label{c.5}
 Pour tout mon\^ome $\mathfrak{m} := {X}^{\alpha_1}{Z}^{\alpha_2} $ on
 a:
 $$\widetilde{\tau}(\mathfrak m ) = \sum_{\ell = 0}^{\infty} \frac {\widetilde{\partial^{\ell}
  \mathfrak m}}{{\widetilde{\Delta}}^{f(\ell)}} t^{\ell}$$
 avec $\widetilde{\partial^{\ell} \mathfrak m} \in k[X , Z]~$ (pour tout $\ell
 \in \mathbb N $). De plus, pour $\delta \in \mathbb N$, $T \in {\mathbb N}^*$ et $v$ une place de $K$ on a :
 $$\max \left(d° \widetilde{\partial^{\ell} \mathfrak m} ~;~ d° \mathfrak m \leq \delta ~,~ \ell = 0 ,
 \ldots , T \right)~\leq ~ 3 T + \delta - 2$$
$$h_v \left(\widetilde{\partial^{\ell} \mathfrak m} ~;~ d° \mathfrak m \leq \delta ~,~ \ell = 0 , \ldots , T \right) ~\leq ~ \begin{cases}
2 m_v T & \text{si $v$ est finie} \\
2 (m_v + 6) T + \delta & \text{si $v$ est infinie}
\end{cases}.$$
Par conséquent pour tout $\delta \in \mathbb N$ et $T \in {\mathbb
N}^*$ on a:
 $$\widetilde h \left(\widetilde{\partial^{\ell} \mathfrak m} ~;~ d° \mathfrak m \leq \delta ~,~
  \ell = 0 ,
 \ldots , T \right) ~\leq ~ 2 (\eta + 6) T + \delta .$$
 \end{corollary}
 {\bf Démonstration.---}
 On a $\mathfrak m := X^{\alpha_1} Z^{\alpha_2}$ donc d'apr\`es la
 proposition \ref{c.4}:
 $$\widetilde{\tau}(\mathfrak m) = {\widetilde{\tau}(X)}^{\alpha_1} {\widetilde{\tau}(Z)}^{\alpha_2} = (X + t)^{\alpha_1}
 \left(\sum_{\ell = 0}^{\infty} \frac{\widetilde{{\partial}^{\ell} Z}}
 {{\widetilde{\Delta}}^{f(\ell)}} (X , Z)\right)^{\!\alpha_2} .$$
 Le développement de cette derni\`ere expression donne une série en
 $t$ qui s'identifie à la série $\sum_{\ell = 0}^{\infty} \frac{\widetilde{{\partial}^{\ell} \mathfrak m}}
 {{\widetilde{\Delta}}^{f(\ell)}}
 t^{\ell}$ pour:
 $$\widetilde{{\partial}^{\ell} \mathfrak m} = \sum
 \widetilde{{\partial}^{{\ell}_{11}}X} \dots \widetilde{{\partial}^{{\ell}_{1
 {\alpha_1}}}X} {\widetilde{\Delta}}^{f(\ell) - \left(f({\ell}_{21}) + \dots + f({\ell}_{2 {\alpha
 2}})\right)} \widetilde{{\partial}^{{\ell}_{21}}Z} \dots
 \widetilde{{\partial}^{{\ell}_{2{\alpha}_2}}Z}$$ ($\forall \ell \geq 0$),
 o\`u la somme $\sum$ porte sur tous les uplets $({\ell}_{11} , \dots , {\ell}_{1 {\alpha}_1} , {\ell}_{21} , \dots
, {\ell}_{2 {\alpha}_2})$ de ${\mathbb N}^{\alpha_1 + \alpha_2}$
satisfaisant ${\ell}_{11} + \dots + {\ell}_{1 {\alpha}_1} +
{\ell}_{21} + \dots + {\ell}_{2 {\alpha}_2} = \ell$ et
${\ell}_{11} ,\dots ,{\ell}_{1{\alpha}_1} \in \{0 , 1\}$. Ces
derniers $\widetilde{{\partial}^{\ell} \mathfrak m}$ sont
effectivement des polyn\^omes de $K[X , Z]$ dont on estime les
degrés et hauteurs gr\^ace à la proposition \ref{c.4}. On obtient
ainsi les estimation du corollaire \ref{c.5}.
$~~~~\blacksquare$\vspace{1mm}
\section{Plongements éclatants}
Dans tout ce qui suit $m \geq 2$ désignera un entier positif. Pour
tout $\underline{a} = (a_1 , \ldots , a_m)$ \\ $\in {\mathbb N}^m$ avec
$a_1 \geq 1 , \ldots , a_{m - 1} \geq 1$ et $a_m = 1$ on
consid\`ere le morphisme (plongement) ${\psi}_{\underline{a}}$ de
$E^m$ dans $E^m \times E^{m - 1}$ défini par: $$
\begin{array}{rcl}
E^m & \stackrel{\psi_{\underline a}}{\hookrightarrow} & E^m \times
E^{m - 1} \\
 (x_1 , \ldots , x_m) & \longmapsto & (x_1 , \ldots , x_m , a_1 x_1 - x_m , \ldots , a_{m - 1} x_{m - 1} - x_m)
\end{array}~.$$
En plongeant de nouveau $E^m \times E^{m - 1}$ dans ${{\mathbb
P}_2}^{2 m - 1}$ par $i^{2 m - 1}$ on obtient le {\sl{plongement
éclatant}}: $$\varphi_{\underline a} := i^{2 m - 1} \circ
\psi_{\underline a} : E^m \hookrightarrow {{\mathbb P}_2}^{2 m -
1} .$$ Nous nous intéressons dans le lemme suivant aux multidegrés
de $\varphi_{\underline a}(E^m)$ dans ${{\mathbb P}_2}^{2 m - 1}
.$
\begin{lemma}\label{c.6}
Soient $I$ un sous-ensemble de $\{1 , \ldots , m\}$ et $J$ un
sous-ensemble de $\{1 , \ldots , m - 1\}$ tels que: $\mbox{card}~I
+ \mbox{card}~J = m $. Notons par $(I , J)$ le $(2 m - 1)$-uplet
$(I , J) := ({\mathbb I}_I (1) , \ldots , {\mathbb I}_I (m) ,
{\mathbb I}_J (1) , \ldots , {\mathbb I}_J (m - 1))$ o\`u
${\mathbb I}_I$ et ${\mathbb I}_J$ désignent les fonctions
caractéristiques de $I$ et $J$ respectivement. Le multidegré
$d_{(I , J)} ({\varphi}_{\underline a}(E^m))$ est nul sauf si : $m
\in I$ et $I \cap J = \emptyset$ ou bien $m \not\in I$ et $I \cap
J$ est un singleton et dans ces deux cas il vaut:
 $$d_{(I , J)} ({\varphi}_{\underline a}(E^m)) = \prod_{i \in J \setminus I}\!\!\!{a_i}^2 . {d(E)}^m
 .$$
\end{lemma}
{\bf Démonstration.---} Dans toute cette démonstration on
identifie $E$ au tore complexe $\mathbb C / \Lambda$ pour un
réseau $\Lambda = \mathbb{Z} + \tau \mathbb{Z} ~ (\tau \in
\mathbb{C} , \mbox{IM} \tau > 0)$ de $\mathbb{C}$. L'invariance du
multidegré par translation permet d'écrire pour tout $u = (u_1 ,
\ldots , u_m , v_1 , \ldots , v_{m - 1}) \in E^{2 m - 1}$:
$$d_{(I , J)} ({\varphi}_{\underline
a}(E^m)) = d_{(I , J)}\left(i^{2 m - 1} \circ \psi_{\underline
a}(E^m)\right) = d_{(I , J)}\left(i^{2 m - 1} \circ \tau_u \circ
\psi_{\underline a }(E^m)\right)$$ o\`u $\tau_u$ désigne la
translation par $u$ dans $E^{2 m - 1} .$ D'autre part, en posant
$i_1 = {\mathbb I}_I(1) , \ldots , i_m = {\mathbb I}_I (m)$ et
$j_1 = {\mathbb I}_J (1) , \ldots , j_{m - 1} = {\mathbb I}_J (m -
1)$ on a, d'apr\`es le théor\`eme de Wirtinger multiprojectif:
$$d_{(I , J)}\left(i^{2 m - 1} \circ \tau_u \circ \psi_{\underline a
}(E^m)\right) = d_{(i_1 , \ldots , i_m , j_1 , \ldots , j_{m -
1})} \left(i^{2 m - 1} \circ \tau_u \circ \psi_{\underline a
}(E^m)\right)$$ $$= \int_{i^{2 m - 1 } \circ \tau_u \circ
\psi_{\underline a} (E^m)} \!\!\!\!\!\!\!\!\!{\Omega_{{\mathbb
P}_2}}^{\!\!\! \wedge i_1} \!\!\wedge \ldots \wedge
{\Omega_{{\mathbb P}_2}}^{\!\!\! \wedge i_m} \!\!\wedge
{\Omega_{{\mathbb P}_2}}^{\!\!\! \wedge j_1} \!\!\wedge \ldots
\wedge {\Omega_{{\mathbb P}_2}}^{\!\!\! \wedge j_{m - 1}}$$ o\`u
$\Omega_{{\mathbb P}_2}$ désigne la forme de Fubini-study sur
${\mathbb P}_2$. En faisant un changement de variable dans cette
derni\`ere intégrale et en moyennant $(2 m - 1)$ fois par la
mesure de Haar normalisée $\nu_E$ sur $E$, puis en inversant les
signes d'intégration on aura les égalités successives suivantes:
\begin{equation*}
\begin{split}
&\int_{i^{2 m - 1 } \circ \tau_u \circ \psi_{\underline a} (E^m)}
\!\!\!\!\!\!\!\!\!{\Omega_{{\mathbb P}_2}}^{\!\!\! \wedge i_1}
\!\!\wedge \ldots \wedge {\Omega_{{\mathbb P}_2}}^{\!\!\! \wedge
i_m} \!\!\wedge {\Omega_{{\mathbb P}_2}}^{\!\!\! \wedge j_1}
\!\!\wedge \ldots \wedge {\Omega_{{\mathbb P}_2}}^{\!\!\!
\wedge j_{m - 1}} \\
&= \!\int_{\psi_{\underline a} (E^m)}\!\!\!\!\!\!\!\!\!(i^{2 m - 1
} \circ \tau_u)^*\!\!\left({\Omega_{{\mathbb P}_2}}^{\!\!\! \wedge
i_1} \!\!\wedge \ldots \wedge {\Omega_{{\mathbb P}_2}}^{\!\!\!
\wedge i_m} \!\!\wedge {\Omega_{{\mathbb P}_2}}^{\!\!\! \wedge
j_1} \!\!\wedge
\ldots \wedge {\Omega_{{\mathbb P}_2}}^{\!\!\! \wedge j_{m - 1}} \right)\!\!\left(Z_1 , \dots , Z_{2 m - 1}\right) \\
&= \int_{E^m}(i^{2 m - 1 } \circ
\tau_u)^*\!\left({\Omega_{{\mathbb P}_2}}^{\!\!\! \wedge i_1}
\!\!\wedge \ldots \wedge {\Omega_{{\mathbb P}_2}}^{\!\!\! \wedge
i_m} \!\!\wedge {\Omega_{{\mathbb P}_2}}^{\!\!\! \wedge j_1}
\!\!\wedge \ldots
\wedge {\Omega_{{\mathbb P}_2}}^{\!\!\! \wedge j_{m - 1}} \right) \!\!\left(\phantom{a_{m - 1} Z_{m - 1} - Z_m} \!\!\!\!\!\!\!\!\!\!\!\!\!\!\!\!\!\!\!\!\!\!\!\!\!\!\!\!\!\!\!\!\!\!\!\!\!\!\!Z_1 , \dots , Z_m ,\right. \\
&~~~~~~~~~~~~~~~~~~~~~~~~~~~~~~~~~~~~~~~~~~~~~~~~~~~~~\phantom{= \int_{E^m}}\left.a_1 Z_1 - Z_m , \dots , a_{m - 1} Z_{m - 1} - Z_m \right) \\
&= \int_{E^{2 m - 1 }}\!\!\left[\int_{E^m}(i^{2 m - 1 } \circ
\tau_u)^*\!\left({\Omega_{{\mathbb P}_2}}^{\!\!\! \wedge i_1}
\!\!\wedge \ldots \wedge {\Omega_{{\mathbb P}_2}}^{\!\!\! \wedge
i_m} \!\!\wedge {\Omega_{{\mathbb P}_2}}^{\!\!\! \wedge j_1}
\!\!\wedge
\ldots \wedge {\Omega_{{\mathbb P}_2}}^{\!\!\! \wedge j_{m - 1}}\!\right) \!\!\left(\phantom{a_{m - 1} Z_{m - 1} - Z_m} \!\!\!\!\!\!\!\!\!\!\!\!\!\!\!\!\!\!\!\!\!\!\!\!\!\!\!\!\!\!\!\!\!\!\!\!\!\!\!Z_1 , \dots ,\right.\right. \\
&~~~~~\!\!\left.\phantom{\int_{E^{2 m - 1
}}\!\!\int_{E^m}}\left.Z_m , a_1 Z_1 - Z_m , \dots , a_{m - 1} Z_{m - 1} - Z_m\right) \right] \wedge d \nu_E (u_1) \wedge \ldots \wedge d \nu_E (u_m) \\
&~~~~~~~~~~~~~~~~~~~~~~~~~~~~~~~~~~~~~~~~~~~~~~~~~~~~~\phantom{\int_{E^{2
m - 1
}}\!\!\left[\int_{E^m}\right.}\wedge d \nu_E (v_1) \wedge \ldots \wedge d \nu_E (v_{m - 1}) \\
&= \int_{E^m }\!\!\left[\int_{E^{2 m - 1}}(i^{2 m - 1 } \circ
\tau_u)^*\!\left({\Omega_{{\mathbb P}_2}}^{\!\!\! \wedge i_1} \!\!\wedge \ldots \wedge {\Omega_{{\mathbb P}_2}}^{\!\!\! \wedge i_m} \!\!\wedge {\Omega_{{\mathbb P}_2}}^{\!\!\! \wedge j_1} \!\!\wedge \ldots \wedge {\Omega_{{\mathbb P}_2}}^{\!\!\! \wedge j_{m - 1}}\right) \!\!\left(\phantom{a_{m - 1} Z_{m - 1} - Z_m} \!\!\!\!\!\!\!\!\!\!\!\!\!\!\!\!\!\!\!\!\!\!\!\!\!\!\!\!\!\!\!\!\!\!\!\!\!\!\!Z_1 , \dots ,\right.\right. \\
&~~~~~~~\!\!\phantom{\int_{E^m }\!\!\int_{E^{2 m - 1}}} \left.Z_m
, a_1 Z_1 - Z_m , \dots , a_{m - 1} Z_{m - 1} - Z_m \right) \wedge
d \nu_E
(u_1) \wedge \ldots \wedge d \nu_E (u_m) \\
&~~~~~~~~~~~~~~~~~~~~~~~~~~~~~~~~~~~~~~~~~~~~~~~~~~~~\left.\phantom{\int_{E^m
}\!\!\int_{E^{2 m - 1}}} \wedge d \nu_E (v_1) \wedge \ldots \wedge
d \nu_E (v_{m - 1})\right]
\end{split}
\end{equation*}
\begin{equation*}
\begin{split}
&= \int_{E^m} \!\!\left[\int_{E} (i \circ \tau_{u_1})^*
 \!\left({\Omega_{{\mathbb P}_2}}^{\!\!\!\wedge i_1}\right)(Z_1) \wedge d \nu_E (u_1) \wedge \ldots \wedge \!\!\int_{E} (i \circ \tau_{u_m})^*
 \!\left({\Omega_{{\mathbb P}_2}}^{\!\!\!\wedge i_m}\right)(Z_m) \right. \\
&~~~~~~~~~~~~\!\phantom{\int_{E^m}\int_{E}}\wedge d \nu_E (u_m) \wedge \int_{E} (i \circ \tau_{v_1})^* \!\left({\Omega_{{\mathbb P}_2}}^{\wedge j_1}\right)(a_1 Z_1 - Z_m) \wedge d \nu_E (v_1) \wedge \dots \\
&~~~~~~~~~~~~~\!\!\!\left.\phantom{\int_{E^m}\int_{E}}\wedge
\int_{E} (i \circ \tau_{v_{m - 1}})^*
 \!\left({\Omega_{{\mathbb P}_2}}^{\wedge {j_{m - 1}}}\right)(a_{m - 1} Z_{m - 1} - Z_m) \wedge d \nu_E (v_{m - 1})\right] .
\end{split}
\end{equation*}
 Or, d'apr\`es la proposition 3.1 de [Da-Ph] (adaptée au plongement de
 Weierstrass), cette derni\`ere intégrale vaut:
\begin{equation*}
\begin{split}
&\int_{{D_{\!E}}^m} \!\!\left(\!\frac{H(d Z_1 , d Z_1)}{- 2
i}\!\right)^{\!\!\!\wedge i_1} \!\!\wedge \ldots \wedge
\left(\!\frac{H(d Z_m , d Z_m)}{- 2
 i}\!\right)^{\!\!\!\wedge i_m} \!\!\!\!\wedge \!\left(\!\frac{H(a_1 d Z_1 - d Z_m , a_1 d Z_1 - d Z_m)}{- 2 i}\!\right)^{\!\!\!\wedge j_1} \\
&~~~~~~~~~~~~~~~~~~~~~~~~~~\phantom{\int_{{D_{\!E}}^m}}\wedge
\ldots \wedge \left(\frac{H(a_{m - 1} d Z_{m - 1} - d Z_m , a_{m -
1} d Z_{m - 1} - d Z_m)}
 {- 2 i}\right)^{\!\!\!\wedge j_{m - 1}} ,
\end{split}
\end{equation*}
o\`u $D_{\!E}$ est un domaine fondamental du réseau $\Lambda
\subset \mathbb C$ et $H$ désigne la forme de Riemann associée à
ce réseau. En remarquant finalement que pour tout $t \in \{1 ,
\ldots , m -
 1\}$: $H(a_t d Z_t - d Z_m , a_t d Z_t - d Z_m) = {a_t}^2 H(d Z_t , d Z_t) - 2 a_t H(d Z_t , d Z_m) +
 H(d Z_m , d Z_m)$ et que $\forall \alpha , \beta , \gamma , \delta \in \{1 , \ldots , m\}$: $H(d Z_{\alpha} ,
 d Z_{\beta}) \wedge H(d Z_{\gamma} , d Z_{\delta}) =$ \\ $({\mbox{IM} \tau})^{-2} d Z_{\alpha}
   \wedge d \overline{Z_{\beta}} \wedge d Z_{\gamma} \wedge d
   \overline{Z_{\delta}}$ est nul d\`es que $\alpha = \gamma$ ou
   $\beta = \delta$ , l'inégalité précédente devient:
   \\ $\bullet$ si $m \in I$ et $I \cap J = \emptyset$ ou bien $m \not\in
   I$ et $I \cap J = \{k\}$ pour un certain $k \in \{1 , \ldots , m -
   1\}$:
$$\prod_{i \in J \setminus I} \!\!\!{a_i}^2 . \int_{{D_{\!E}}^m}
\left(\frac{H(d Z_1 , d Z_1)}{- 2 i}\right)
 \wedge \ldots \wedge \left(\frac{H(d Z_m , d Z_m)}{- 2
 i}\right) = \prod_{i \in J \setminus I}\!\!\! {a_i}^2 .
 {d(E)}^m$$(en utilisant une autre fois le théor\`eme de Wirtinger)\\ $\bullet$ nulle sinon.
$~~~~\blacksquare$\vspace{1mm}
\section{Paramétrisation locale de l'image de $E^m$ dans ${\mathbb{P}}_{2}^{2 m - 1}$}
Commençons d'abord par l'homogénisation du monomorphisme
$\widetilde{\tau}$ de la proposition \ref{c.4} qui consiste à
définir un monomorphisme $\tau$ de $A(E) := \frac{K[X , Y ,
Z]}{I(E)}$ dans ${K(E)}_0 [[t]]$ \`a partir duquel on retrouve
notre monomorphisme $\widetilde{\tau}$ en spécialisant la variable
$Y$ à $1$. On doit alors définir $\tau$ de la mani\`ere suivante:
$\forall P_1 \in A(E):$
$$\tau(P_1) := \widetilde{\tau}({\widetilde P}_1) \left(\frac{X}{Y} , \frac{Z}{Y}\right)$$
avec ${\widetilde P}_1$ désigne l'élément de $A(\widetilde E) :$
$${\widetilde P}_1 (X , Z) := P_1 (X , 1 , Z).$$
Posons aussi par définition:
$$\Delta \left(X , Y , Z \right) := Y^2 \widetilde{\Delta}\left(\frac{X}{Y} , \frac{Z}{Y}\right) = 3 g_3 Z^2 + 2 g_2 X Z
 + Y^2$$
 et pour un mon\^ome $\mathfrak m := X^{\alpha_1} Y^{\alpha_2}
 Z^{\alpha_3}$ de $A(E)$ et un entier $\ell \in \mathbb N :$
 $$\partial^{\ell} \mathfrak m := Y^{g(\ell) + d° \mathfrak m} \widetilde{\partial^{\ell} \widetilde{\mathfrak m}}
  \left(\frac{X}{Y} , \frac{Z}{Y}\right)$$
  avec $\widetilde{\mathfrak m}$ désigne le mon\^ome de $A(\widetilde{E}) :$
  $$\widetilde{\mathfrak m}\left(X , Z\right) := \mathfrak m \left(X , 1 , Z\right) = X^{\alpha_1} Z^{\alpha_3} .$$
  Ces derniers $\partial^{\ell} \mathfrak m , \ell \in \mathbb N$
  sont -d'apr\`es le corollaire \ref{c.5}- des formes de $K[X , Y , Z]$ de
  degrés: $d° \partial^{\ell} \mathfrak m = g(\ell) + d° \mathfrak m ~(\forall \ell \in \mathbb
  N)$ et chaque forme $\partial^{\ell} \mathfrak m ~(\ell \in \mathbb
  N)$ a évidemment les m\^emes coefficients que le polyn\^ome $\widetilde{\partial^{\ell} \widetilde{\mathfrak
  m}}$ de $K[X , Z]$, donc \`a fortiori une famille finie de formes
  ${(\partial^{\ell} \mathfrak m)}_{\ell , \mathfrak m}$ a la
  m\^eme hauteur logarithmique que la famille finie correspondante
  des polyn\^omes ${(\widetilde{\partial^{\ell} \widetilde{\mathfrak m}})}_{\ell , \mathfrak
  m} .$ Ainsi du corollaire \ref{c.5} découle immédiatement le corollaire
  suivant:
 \begin{corollary}\label{c.7}
Il existe un monomorphisme d'anneau $\tau$ de $A(E) := \frac{K[X ,
Y , Z]}{I(E)}$ dans ${K(E)}_0 [[t]]$ associant à tout mon\^ome
$\mathfrak m := X^{\alpha_1} Y^{\alpha_2} Z^{\alpha_3}$ de $A(E)
:$
$$\tau(\mathfrak m) = \frac{1}{Y^{d° \mathfrak m}} \left(\sum_{\ell = 0}^{\infty} \frac{Y^{\ell} .
\partial^{\ell} \mathfrak m}{\Delta^{f(\ell)}} t^{\ell}\right)$$
o\`u les $\partial^{\ell} \mathfrak m , ~\ell \in \mathbb N$ sont
des formes de $K[X , Y , Z]$ et pour $\delta \in \mathbb N$ et $T
\in {\mathbb N}^*$ on a:
$$\max\left(d° \partial^{\ell} \mathfrak m ~;~ d° \mathfrak m \leq \delta ~,~ \ell = 0 , \ldots , T\right)
~\leq~ 3 T + \delta - 2$$
$$h_v \left(\partial^{\ell} \mathfrak m ~;~ d° \mathfrak m \leq \delta ~,~ \ell = 0 , \ldots , T\right)
~\leq~ \begin{cases}
2 m_v T & \text{si $v$ est finie} \\
2 (m_v + 6) T + \delta & \text{si $v$ est infinie}
\end{cases} .$$
Par conséquent, pour tout $\delta \in \mathbb N$ et $T \in
{\mathbb N}^*$ on a:
$$\widetilde{h}\left(\partial^{\ell} \mathfrak m ~;~ d° \mathfrak m \leq \delta ~,~ \ell = 0 , \ldots , T\right)
~\leq~ 2 (\eta + 6) T + \delta .$$
 \end{corollary}
 Maintenant, étant donné un $m$-uplet fixé $\underline a = (a_1 , \dots , a_m) \in {\mathbb
 N}^m$ tel que $a_i \geq 5$ pour $i = 1 , \dots , m - 1$ et $a_m = 1
 ,$ posons:
 $$B := \frac{K\left[{\underline X}_1 , \dots , {\underline X}_m , {\underline Y}_1 , \dots , {\underline Y}_{m - 1}
 \right]}{I\left(\varphi_{\underline a} (E^m)\right)}$$
 l'anneau des coordonnées de $\varphi_{\underline a} (E^m)$ et $K(\varphi_{\underline a}
 (E^m))$ son corps de fractions, avec ${\underline X}_i := (X_{i 0} , X_{i 1} , X_{i
 2})$ et ${\underline Y}_i := (Y_{i 0} , Y_{i 1} , Y_{i 2}) ,~i = 1 , \dots , m - 1 .$
 Soit, par ailleurs, $\underline D$ une famille de formes
 bihomog\`enes représentant la différence dans $E$ au voisinage de
 $\{\mathbf{0}\} \times \{\mathbf{0}\} .$ D'apr\`es le théor\`eme \ref{c.8} du formulaire, $\underline D$ peut \^etre prise constituée de formes de bidegré $(2 , 2)$
 et de hauteur logarithmique locale $h_v$ (resp de longueur logarithmique locale $\ell_v$) -pour une place finie (resp infinie) $v$ sur $K$- majorée par: $h_v(\underline D) \leq 3 m_v$ (resp $\ell_v(\underline D) \leq 3 m_v + 7$) et de hauteur de Gauss-Weil majorée par: $\widetilde{h}(\underline D) \leq 3 \eta + 5$. Soit aussi, pour tout entier $k \geq 1$, ${\underline F}^{(k)}$
 une famille de formes homog\`enes représentant la multiplication
 par $k$ dans $E$. D'apr\`es le théor\`eme \ref{c.9} du formulaire, ${\underline F}^{(k)}$ peut \^etre prise constituée de formes de degré $k^2$
 chacune, de hauteur logarithmique locale $v$-adique (resp de longueur logarithmique locale $v$-adique) -pour une place finie (resp infinie) $v$ sur $K$- majorée par $\frac{3}{2} m_v k^2$ (resp $\frac{3}{2}(m_v + 3) k^2$) et de hauteur de Gauss-Weil majorée par: $\widetilde{h}({\underline F}^{(k)}) \leq \frac{3}{2}(\eta + 3) k^2$.
\`A partir du monomorphisme $\tau ,$ on en déduit un monomorphisme
 de paramétrisation locale pour la sous-variété $\varphi_{\underline
 a}(E^m)$ de ${{\mathbb P}_2}^{2 m - 1}$ qu'on notera $\Omega_{\underline
 a}$ et qu'on définit par:
 $$\Omega_{\underline a} : B \rightarrow K\left(\varphi_{\underline a}(E^m)\right)[[\underline u]]$$
 $$~~~~~~~~~~\Omega_{\underline a}({\underline X}_1) := \tau_1 ({\underline X}_1) , \dots ,
 \Omega_{\underline a}({\underline X}_m) := \tau_m ({\underline X}_m);$$
$$~~~~\Omega_{\underline a}({\underline Y}_1) := \underline D
\left({\underline F}^{(a_1)}(\tau_1 ({\underline X}_1)) , \tau_m
({\underline X}_m)\right) $$
$$\vdots~~~~~~~~~~~~~~~~~~~~~~~~~~$$
$$~~~~~~~~~~~\Omega_{\underline a}({\underline Y}_{m - 1}) := \underline D
\left({\underline F}^{(a_{m - 1})}(\tau_{m - 1} ({\underline X}_{m
- 1})) , \tau_m ({\underline X}_m)\right)$$ o\`u $\tau_i : A(E)
\rightarrow {K(E)}_0 [u_i] ~(i = 1 , \dots , m)$ désigne le
monomorphisme $\tau$ du corollaire \ref{c.7} pour la $i$ \`eme
composante $E$ de la sous-variété $E^m$ de ${{\mathbb P}_2}^m$ et
$\underline u := (u_1 , \dots , u_m) .$ On a le lemme suivant:
\begin{lemma}\label{c.10}
Avec toutes les notations précédentes, pour tout mon\^ome \\
$\mathfrak m = {\underline X}_{1}^{{\underline{\alpha}}_1} \dots
{\underline X}_{m}^{{\underline{\alpha}}_m} {\underline
Y}_{1}^{{\underline{\beta}}_1} \dots {\underline Y}_{m -
1}^{{\underline{\beta}}_{m - 1}}$ on a:
$$\Omega_{\underline a}(\mathfrak m) = \frac{1}{X_{1 1}^{d_1} \dots X_{m 1}^{d_m}}
\sum_{i_1 \geq 0 , \dots , i_m \geq 0} \frac{X_{1 1}^{i_1} \dots
X_{m 1}^{i_m} \partial^{(i_1 , \dots , i_m)} \mathfrak
m}{{\Delta({\underline X}_1)}^{f(i_1)} \dots {\Delta({\underline
X}_m)}^{f(i_m)}} u_{1}^{i_1} \dots u_{m}^{i_m}$$ o\`u $d_1 , \dots
d_m$ désignent les entiers positifs:
\begin{equation*}
\begin{split}
d_i &:=
\mid{\underline{\alpha}}_i\mid + 2 {a^2}_{\!\!\!i} \mid{\underline{\beta}}_i\mid ~~~~\text{pour $i = 1 , \dots , m - 1$} \\
\text{et}~~ d_m &:= \mid{\underline{\alpha}}_m\mid + 2
\left(\mid{\underline{\beta}}_1\mid + \dots +
\mid{\underline{\beta}}_{m - 1}\mid\right) ;
\end{split}
\end{equation*}
les $\partial^{(i_1 , \dots , i_m)} \mathfrak m$ sont des formes
de $K[{\underline X}_1 , \dots , {\underline X}_m]$ de degrés
majorés par:
$${d°}_{\!\!\!\!{\underline X}_1} \partial^{(i_1 , \dots , i_m)} \mathfrak m ~\leq~ d_1 + g(i_1)$$
$$\vdots~~~~~~~~~~~~~~~~~~~~$$
$${d°}_{\!\!\!\!{\underline X}_m} \partial^{(i_1 , \dots , i_m)} \mathfrak m ~\leq~ d_m + g(i_m)$$
et pour tout $\delta_1 , \dots , \delta_m , {{\delta}'}_{\!\!1} ,
\dots , {{\delta}'}_{\!\!m - 1} , T \in \mathbb N ,$ la famille
des formes $\partial^{(i_1 , \dots , i_m)} \mathfrak m ,$ \\
${d°}_{\!\!\!\!{\underline X}_1} \mathfrak m \leq \delta_1 , \dots
, {d°}_{\!\!\!\!{\underline X}_m} \mathfrak m \leq \delta_m ,
{d°}_{\!\!\!\!{\underline Y}_1} \mathfrak m \leq
{{\delta}'}_{\!\!1} , \dots , {d°}_{\!\!\!\!{\underline Y}_{m -
1}} \leq {{\delta}'}_{\!\!m - 1} , i_1 + \dots + i_m \leq T$ est
de hauteur logarithmique locale $h_v$ majorée par:
$$\left[2 T + \sum_{i = 1}^{m - 1} 3 ({a^2}_{\!\!\!i} + 1) {\delta'}_{\!\!i}\right] m_v$$ lorsque $v$ est finie et elle est de longueur logarithmique locale $\ell_v$ majorée par:
$$\left[2 T + \sum_{i = 1}^{m - 1} 3 ({a^2}_{\!\!\!i} + 1) {\delta'}_{\!\!i}\right] m_v + 12 T + (\delta_1 + \dots + \delta_m) + \sum_{i = 1}^{m - 1} (11 {a^2}_{\!\!\!i} + 9) {\delta'}_{\!\!i}$$ lorsque $v$ est infinie. Par conséquent, elle est de hauteur de Gauss-Weil majorée par:
$$\left[2 T + \sum_{i = 1}^{m - 1} 3 ({a^2}_{\!\!\!i} + 1) {\delta'}_{\!\!i}\right] \eta + 12 T + (\delta_1 + \dots + \delta_m) + \sum_{i = 1}^{m - 1} (11 {a^2}_{\!\!\!i} + 9) {\delta'}_{\!\!i} .$$
\end{lemma}
{\bf Démonstration.---} Pour tout mon\^ome $\mathfrak m =
{{\underline X}_1}^{{\underline{\alpha}}_1} \dots {{\underline
X}_m}^{{\underline{\alpha}}_m} {{\underline
Y}_1}^{{\underline{\beta}}_1} \dots {{\underline Y}_{m -
1}}^{{\underline{\beta}}_{m - 1}} ,$ on remarque que la série
$\Omega_{\underline a} (\mathfrak m)$ s'obtient en substituant
respectivement dans la forme:
$$R({\underline X}_1 , \dots , {\underline X}_m) := {{\underline
X}_1}^{{\underline{\alpha}}_1} \dots {{\underline
X}_m}^{{\underline{\alpha}}_m} \prod_{i = 1}^{m - 1} \underline D
\left({\underline F}^{(a_i)} ({\underline X}_i) , {\underline
X}_m\right)^{{\underline{\beta}}_i} $$ ${\underline X}_1 , \dots ,
{\underline X}_m$ par $\tau_1({\underline X}_1) , \dots ,
\tau_m({\underline X}_m) .$ \\ Il est clair que cette forme
multihomog\`ene $R$ de $K[{\underline X}_1 , \dots , {\underline
X}_m]$ est de multidegré:
$$
\left(\mid{\underline{\alpha}}_1\mid + 2 {a^2}_{\!\!\!1}
\mid{\underline{\beta}}_1\mid , \dots ,
\mid{\underline{\alpha}}_{m - 1}\mid + 2 {a^2}_{\!\!\!m - 1}
\mid{\underline{\beta}}_{m - 1}\mid ,
\mid{\underline{\alpha}}_m\mid + 2 \sum_{i = 1}^{m - 1}
\mid{\underline{\beta}}_i\mid\right)
$$
et un simple calcul montre que pour toute place $v$ sur $K$ on a:
\begin{align}
h_v(R) &\leq \left[\sum_{i = 1}^{m - 1} 3({a^2}_{\!\!\!} + 1) \mid {\underline \beta}_i \mid \right] m_v \notag \\
\intertext{lorsque $v$ est finie et:} \ell_v(R) &\leq
\left[\sum_{i = 1}^{m - 1} 3({a^2}_{\!\!\!} + 1) \mid {\underline
\beta}_i \mid \right] m_v + \sum_{i = 1}^{m - 1} (9
{a^2}_{\!\!\!i} + 7) \mid {\underline \beta}_i \mid \notag
\end{align}
lorsque $v$ est infinie. \\
De plus le nombre de mon\^omes que contient $R$ est majoré par: \\ \\
$\displaystyle \binom{2 \mid{\underline{\beta}}_1\mid
{a^2}_{\!\!\!1} + 2}{2} \kern-3pt \dots \kern-3pt \binom{2
\mid{\underline{\beta}}_{m - 1}\mid {a^2}_{\!\!\!m - 1} + 2}{2}
\kern-3pt \binom{2 \mid{\underline{\beta}}_1\mid + \dots + 2
\mid{\underline{\beta}}_{m - 1}\mid + 2}{2}$. \\ \\ Il ne reste
qu'à appliquer le lemme \ref{c.7} à chaque mon\^ome intervenant
dans $R$ pour conclure. En effet, en écrivant:
$$R({\underline X}_1 , \dots , {\underline X}_m) =: \!\!\!\sum_{\underline \gamma := ({\underline \gamma}_1 , \dots , {\underline \gamma}_m) \in U} \!\!\!\!\!\!\!r(\underline \gamma) {{\underline X}_1}^{{\underline \gamma}_1} \dots {{\underline X}_m}^{{\underline \gamma}_m}$$
pour un certain sous-ensemble fini $U$ de $({\mathbb N}^3)^m$ et
certains nombres $r(\underline \gamma) , \underline \gamma \in U$
de $K$, on a:
\begin{eqnarray*}
\Omega_{\underline a}(\mathfrak m) & = & R \left(\tau({\underline X}_1) , \dots , \tau({\underline X}_m)\right) \\
& = & \sum_{\underline \gamma \in U} r(\underline \gamma)
\tau({\mathfrak m}_1) \dots \tau({\mathfrak m}_m)
\end{eqnarray*}
avec ${\mathfrak m}_i := {{\underline X}_i}^{{\underline
\gamma}_i}$ pour $i = 1 , \dots , m$. En substituant maintenant
les $\tau({\mathfrak m}_i) ~ (i = 1 , \dots , m)$ par leurs
expressions données par le corollaire \ref{c.7}, on
aura:\begin{eqnarray*}
\Omega_{\underline a}(\mathfrak m) & \!\!\!\!= & \!\!\!\!\sum_{\underline \gamma \in U} \left[r(\underline \gamma) \frac{1}{X_{11}^{\mid {\underline \gamma}_1 \mid}} \left(\sum_{i_1 = 0}^{\infty} \frac{X_{11}^{i_1} \partial^{i_1} {\mathfrak m}_1}{{\Delta({\underline X}_1)}^{f(i_1)}} u_{1}^{i_1}\right) \dots \frac{1}{X_{m1}^{\mid {\underline \gamma}_m \mid}} \left(\sum_{i_m = 0}^{\infty} \frac{X_{m1}^{i_m} \partial^{i_m} {\mathfrak m}_m}{{\Delta({\underline X}_m)}^{f(i_m)}} u_{m}^{i_m}\right)\right] \\
& \!\!\!\! = & \!\!\!\!\frac{1}{{\underline X}_{11}^{d_1} \dots
{\underline X}_{m1}^{d_m}} \sum_{i_1 \geq 0 , \dots , i_m \geq 0}
\frac{{\underline X}_{11}^{i_1} \dots {\underline X}_{m1}^{i_m}
\partial^{(i_1 , \dots , i_m)}\mathfrak m}{{\Delta({\underline
X}_1)}^{f(i_1)} \dots {\Delta({\underline X}_m)}^{f(i_m)}}
u_{1}^{i_1} \dots u_{m}^{i_m}
\end{eqnarray*}
avec $\displaystyle d_i := \mid {\underline \gamma}_i \mid =
\begin{cases}
\mid {\underline \alpha}_i \mid + 2 {a^2}_{\!\!\!i} \mid {\underline \beta}_i \mid & \!\!\!\text{pour $i = 1 , \dots , m - 1$} \\
\mid {\underline \alpha}_m \mid + 2 \left(\mid {\underline
\beta}_1 \mid + \dots + \mid {\underline \beta}_{m - 1}
\mid\right) & \!\!\!\text{pour $i = m$}
\end{cases}$ \\ \\
et $\displaystyle \partial^{(i_1 , \dots , i_m)} \mathfrak m := \sum_{\underline \gamma \in U} r(\underline \gamma) \partial^{i_1} {\mathfrak m}_1 \dots \partial^{i_m} {\mathfrak m}_m ,~~~~~~~~~~ \forall (i_1 , \dots , i_m) \in {\mathbb N}^m$. \\ \\
Ainsi, on a bien la formule du lemme \ref{c.10} pour
$\Omega_{\underline a}(\mathfrak m)$ et de plus, concernant les
degrés on a:
$${d°}_{\!\!\!\!{\underline X}_j} \partial^{(i_1 , \dots , i_m)} \mathfrak m = \max_{\underline \gamma \in U} {d°}_{\!\!\!\!{\underline X}_j} \partial^{i_j} {\mathfrak m}_j ~~~~~~ \forall j = 1 , \dots , m$$
et concernant les hauteurs et les longueurs logarithmiques locales
on a bien pour tout $(i_1 , \dots , i_m) \in {\mathbb N}^m$:
\begin{align*}
h_v(\partial^{(i_1 , \dots , i_m)} \mathfrak m) &\leq h_v(\partial^{i_1} {\mathfrak m}_1) + \dots + h_v(\partial^{i_m} {\mathfrak m}_m) + h_v(R) \\
\intertext{lorsque $v$ est une place finie sur $K$ et:}
\ell_v(\partial^{(i_1 , \dots , i_m)} \mathfrak m) &\leq
\ell_v(\partial^{i_1} {\mathfrak m}_1) + \dots +
\ell_v(\partial^{i_m} {\mathfrak m}_m) + \ell_v(R)
\end{align*}
lorsque $v$ est une place infinie sur $K$. Le reste suit de
l'application des estimations du lemme \ref{c.7}. La démonstration
est achevée.
$~~~~\blacksquare$\vspace{1mm} \\
$~$ \\
Plus généralement on a le corollaire suivant qui est une
conséquence immédiate du lemme \ref{c.10} précédent:
\begin{corollary}\label{c.11}
Pour toute forme multihomog\`ene non identiquement nulle \\ $P_1
\in K[{\underline X}_1 , \dots , {\underline X}_m , {\underline
Y}_1 , \dots , {\underline Y}_{m - 1}]$ de multidegré $(\delta_1 ,
\dots , \delta_m , {{\delta}'}_{\!\!1} , \dots ,
{{\delta}'}_{\!\!m - 1})$ on a:
$$\Omega_{\underline a}(P_1) = \frac{1}{{X_{1 1}}^{d_1} \dots {X_{m 1}}^{d_m}}
\sum_{i_1 \geq 0 , \dots , i_m \geq 0} \frac{{X_{1 1}}^{i_1} \dots
{X_{m 1}}^{i_m} \partial^{(i_1 , \dots , i_m)}
P_1}{{\Delta({\underline X}_1)}^{f(i_1)} \dots {\Delta({\underline
X}_m)}^{f(i_m)}} {u_1}^{i_1} \dots {u_m}^{i_m}$$ avec $d_1 , \dots
d_m$ désignent les entiers positifs:
\begin{eqnarray*}
d_1 & := & \delta_1 + 2 {a^2}_{\!\!\!1} {{\delta}'}_{\!\!1} , \\
\vdots \\
d_{m - 1} & := & \delta_{m - 1} + 2 {a^2}_{\!\!\!m - 1} {{\delta}'}_{\!\!m - 1} \\
\mbox{et}~~~~ d_m & := & \delta_m + 2 ({{\delta}'}_{\!\!1} + \dots
+ {{\delta}'}_{\!\!m - 1}) ;
\end{eqnarray*}
les $\partial^{(i_1 , \dots , i_m)} P_1$ sont des formes de
$K[{\underline X}_1 , \dots , {\underline X}_m]$ de degrés majorés
par:
\begin{eqnarray*}
{d°}_{\!\!\!\!{\underline X}_1} \partial^{(i_1 , \dots , i_m)} P_1
& \leq & \delta_1 +
2 {a^2}_{\!\!\!1} {{\delta}'}_{\!\!1} + g(i_1) \\
\vdots \\
{d°}_{\!\!\!\!{\underline X}_{m - 1}} \partial^{(i_1 , \dots ,
i_m)} P_1 & \leq & \delta_{m - 1} +
2 {a^2}_{\!\!\!m - 1} {{\delta}'}_{\!\!m - 1} + g(i_{m - 1}) \\
{d°}_{\!\!\!\!{\underline X}_m} \partial^{(i_1 , \dots , i_m)} P_1
& \leq & \delta_m + 2 ({{\delta}'}_{\!\!1} + \dots +
{{\delta}'}_{\!\!m - 1}) + g(i_m)
\end{eqnarray*}
et, pour tout $T \in \mathbb N$ et $v$ une place de $K$, la
famille des formes $\partial^{(i_1 , \dots , i_m)} P_1 , i_1 +
\dots + i_m \leq T$ est de hauteur logarithmique locale $h_v$
majorée par:
$$\left[2 T + \sum_{i = 1}^{m - 1} 3({a^2}_{\!\!\!i} + 1) {\delta'}_{\!\!i} \right] m_v + h_v(P_1)$$
lorsque $v$ est finie et elle est de longueur logarithmique locale
$\ell_v$ majorée par:
$$\left[2 T + \sum_{i = 1}^{m - 1} 3({a^2}_{\!\!\!i} + 1) {\delta'}_{\!\!i} \right] m_v + 12 T + (\delta_1 + \dots + \delta_m) + \sum_{i = 1}^{m - 1} (11 {a^2}_{\!\!\!i} + 9) + \ell_v(P_1)$$
lorsque $v$ est infinie. Par conséquent, elle est de hauteur de Gauss-Weil majorée par: \\ \\
$\displaystyle \left[2 T + \sum_{i = 1}^{m - 1} 3({a^2}_{\!\!\!i} + 1) {\delta'}_{\!\!i} \right] \eta + 12 T + (\delta_1 + \dots + \delta_m) + \sum_{i = 1}^{m - 1} (11 {a^2}_{\!\!\!i} + 9) + \widetilde{h}(P_1)$ \\ \\
$\displaystyle ~~~~~~~~~~~~~~~~~~~~~~~~~~~~~~~~~~~~~~~~~~~~+ 2
\sum_{i = 1}^{m} \log\left(\delta_i + 1\right) + 2 \sum_{i = 1}^{m
- 1} \log\left({\delta'}_{\!\!i} + 1\right)$.
\end{corollary}
{\bf Démonstration.---} Ecrivons $$P_1 =: \sum_{\mathfrak m \in
\Lambda} \rho_{\mathfrak m}.\mathfrak m$$ o\`u $\Lambda$ est
l'ensemble fini des mon\^omes unitaires de $K[{\underline X}_1 ,
\dots , {\underline X}_m , {\underline Y}_1 , \dots , {\underline
Y}_{m - 1}]$ de multidegrés $(\delta_1 , \dots , \delta_m ,
{{\delta}'}_{\!\!1} , \dots , {{\delta}'}_{\!\!m - 1})$ et les
$\rho_{\mathfrak m} ~(\mathfrak m \in \Lambda)$ sont des nombres
de $K$ (ce sont les coefficients de $P_1$). On a d'apr\`es le
lemme \ref{c.10} précédent:
\begin{eqnarray*}
\Omega_{\underline a} (P_1) & = & \sum_{\mathfrak m \in \Lambda} \rho_{\mathfrak m}.\Omega_{\underline a} (\mathfrak m) \\
& = & \frac{1}{{X_{1 1}}^{d_1} \dots {X_{m 1}}^{d_m}} \sum_{i_1
\geq 0 , \dots , i_m \geq 0} \frac{{X_{1 1}}^{i_1} \dots {X_{m
1}}^{i_m} \partial^{(i_1 , \dots , i_m)} P_1}{{\Delta({\underline
X}_1)}^{f(i_1)} \dots {\Delta({\underline X}_m)}^{f(i_m)}}
{u_1}^{i_1} \dots {u_m}^{i_m}
\end{eqnarray*}
$$\mbox{avec} ~~~~~~~~~~~~~~~~~~~~~~~~\partial^{(i_1 , \dots , i_m)} P_1 := \sum_{\mathfrak m \in \Lambda} \rho_{\mathfrak m}.\partial^{(i_1 , \dots , i_m)}
\mathfrak m .~~~~~~~~~~~~~~~~~~~~~~~~~~~~~~~~~~~~~~~~~~$$ La suite
du corollaire \ref{c.11} s'obtient des estimations du lemme
\ref{c.10} pour les degrés, hauteurs et longueurs des
$\partial^{(i_1 , \dots , i_m)} \mathfrak m$, en majorant de plus
-afin d'estimer la hauteur de Gauss-Weil des $\partial^{(i_1 ,
\dots , i_m)}P_1$- les $\ell_v(P_1) ~(v \in M_{K}^{\infty})$ par:
\begin{equation*}
\begin{split}
\ell_v(P_1) &\leq h_v(P_1) + \log \mbox{card}~\!\Lambda \\
&\leq h_v(P_1) + \sum_{i = 1}^{m} \log\binom{\delta_i + 2}{2} + \sum_{i = 1}^{m - 1} \log\binom{{\delta'}_{\!\!i} + 2}{2} \\
&\leq h_v(P_1) + 2 \sum_{i = 1}^{m} \log\left(\delta_i + 1\right)
+ 2 \sum_{i = 1}^{m - 1} \log\left({\delta'}_{\!\!i} + 1\right) .
\end{split}
\end{equation*}
La démonstration est achevée.  $~~~~\blacksquare$\vspace{1mm}
\section{La fonction auxiliaire}
Soient, pour toute la suite de ce texte, $\epsilon_0$ et
$\epsilon_1$ deux réels strictement positifs assez petits et
$\delta \in {\mathbb N}^*$ tel que $\epsilon_0 \delta$ soit un
entier positif assez grand. Le param\`etre $\delta$ est destiné
\`a tendre vers l'infini. Nous supposons que $\epsilon_0$ et
$\epsilon_1$ vérifient:
$$\epsilon_0 \leq \frac{1}{2} ~~~\mbox{et}~~~ \frac{m - 1}{m!} \left(\frac{7}{3}\right)^m \frac{{\epsilon_1}^m}{\epsilon_0
(m + \epsilon_0) (1 + \epsilon_0)^{m - 2}} \leq \frac{1}{2} .$$
Soit aussi $\underline a = (a_1 , \dots , a_m)$ un $m$-uplet de
${\mathbb N}^m$ tel que: \\ \\
$\displaystyle a_1 \geq a_2 \geq \dots \geq a_{m - 1} \geq a_m =
1$ et $\displaystyle {a^2}_{\!\!\!1} + \dots + {a^2}_{\!\!\!m - 1}
+ m - 1 \leq
\left(1 + \frac{1}{24}\right){a^2}_{\!\!\!1}$. \\ \\
Ceci entra{\sf\^\i}ne qu'on a:
$${a^2}_{\!\!\!1} \geq \frac{48}{25} (m - 1) \geq m - 1 .$$
Posons $T_{\delta}$ le simplexe de ${\mathbb N}^m$:
$$T_{\delta} := \left\{(\tau_1 , \dots , \tau_m) \in {\mathbb N}^m ~/~ \frac{\tau_1}{{a^2}_{\!\!\!1}} + \dots +
\frac{\tau_{m - 1}}{{a^2}_{\!\!\!m - 1}} + \frac{\tau_m}{m - 1}
\leq 7 \epsilon_1 \delta \right\}$$ et ${\mathfrak q}_{\delta}$
l'idéal de $K[[u_1 , \dots , u_m]]$:
$${\mathfrak q}_{\delta} := \left\{f \in K[[u_1 , \dots , u_m]] ~/~ f_{i_1 , \dots , i_m} = 0 ~~ \forall (
i_1 , \dots , i_m) \in T_{\delta}\right\}~~~$$ o\`u $f_{i_1 ,
\dots , i_m}$ désigne le coefficient de ${u_1}^{i_1} \dots
{u_m}^{i_m}$ dans la série $f$ de $K[[u_1 , \dots , u_m]]$.
\\ Dans un premier temps, nous allons construire une forme
multihomog\`ene non identiquement nulle $P$ de $K[{\underline X}_1
, \dots , {\underline X}_m , {\underline Y}_1 , \dots ,
{\underline Y}_{m - 1}] / I\left(\varphi_{\underline a}
(E^m)\right)$ de multidegré \\ $(\epsilon_0 \delta {a^2}_{\!\!\!1} ,
\dots , \epsilon_0 \delta {a^2}_{\!\!\!m} , \delta , \dots ,
\delta)$ qui soit de hauteur relativement petite par rapport à son
degré et qui satisfasse la condition d'annulation:
$$\Omega_{\underline a} (P) (\underline 0 , \dots , \underline 0) \in {\mathfrak q}_{\delta}$$
avec $\underline 0 := (0 , 1 , 0) .$ Cette condition
d'appartenance à ${\mathfrak q}_{\delta}$ s'interpr\^ete comme
syst\`eme linéaire sur $K$ en les coefficients de $P$. En effet,
en choisissant une base $\mathcal M$ du $K$-espace vectoriel:
$${\left(K[{\underline X}_1 , \dots , {\underline X}_m ,
{\underline Y}_1 , \dots , {\underline Y}_{m - 1}] /
I(\varphi_{\underline a} (E^m))\right)}_{(\epsilon_0 \delta
{a^2}_{\!\!\!1} , \dots , \epsilon_0 \delta {a^2}_{\!\!\!m} ,
\delta , \dots , \delta)}$$ et en écrivant $P =: \sum_{\mathfrak m
\in \mathcal M} P_{\mathfrak m}.\mathfrak m$, la condition
$\Omega_{\underline a} (P) (\underline 0 , \dots , \underline 0)
\in {\mathfrak q}_{\delta}$ se traduit par le syst\`eme
d'équations:
\begin{equation}
\sum_{\mathfrak m \in \mathcal M} P_{\mathfrak m}.{\partial}^{(i_1
, \dots , i_m)} \mathfrak m (\underline 0 , \dots , \underline 0)
= 0 ~~~~\mbox{pour}~ (i_1 , \dots , i_m) \in T_{\delta} .
\label{3.1}
\end{equation}
En choisissant maintenant un ordre pour chacun des deux ensembles
finis $T_{\delta}$ et $\mathcal M$ et en posant $M$ la matrice: $M
:= {\left({\partial}^{(i_1 , \dots , i_m)} \mathfrak m (\underline
0 , \dots , \underline 0)\right)}_{(i_1 , \dots , i_m) \in
T_{\delta} , \mathfrak m \in \mathcal M}$ et $\overrightarrow{V}$
le vecteur: $\overrightarrow{V} := {\left(P_{\mathfrak
m}\right)}_{\mathfrak m \in \mathcal M}$, le syst\`eme (\ref{3.1})
devient: $M.\overrightarrow{V} = \overrightarrow{0}.$ On est donc
ramené \`a chercher un vecteur non nul $\overrightarrow{V}$ \`a
coordonnées dans $K$, de hauteur assez petite et satisfaisant
$M.\overrightarrow{V} = \overrightarrow{0}.$ On utilise pour cela
la version suivante du lemme de Siegel, d\^ue à E. Bombieri (voir
[Bom1]).
\begin{lemma}[lemme de Siegel]\label{b.10}
 Pour toute matrice $A \in {\mathcal M}_{m \times n}(K)$ avec $m < n$ il existe
 $\underline x \in K^n \backslash \{\underline 0 \}$ satisfaisant $A \underline x = \underline 0 $ et tel que
 $$\widetilde h (\underline x) \leq e \left(\widetilde h (A) + \log n \right)
  + (1 + e) c_{\mbox S}(K)$$ o\`u $e := \frac{m}{n - m}$ est l'exposant de Dirichlet de $A$, $c_{\mbox S}(K)$ est une constante dépendant seulement du corps de nombres $K$ et $\widetilde h (A)$ désigne ici la hauteur de Gauss-Weil du vecteur formé de tous les coefficients de $A$.
 \end{lemma}
L'application du lemme de Siegel ci-dessus nous donne la
proposition suivante:
\begin{proposition}\label{c.12}
Sous toutes les hypoth\`eses précédentes, il existe une forme $P
\in K[{\underline X}_1 , \dots , {\underline X}_m , {\underline
Y}_1 , \dots , {\underline Y}_{m - 1}] /
I\left(\varphi_{\underline a} (E^m)\right)$ non nulle, de
multidegré $(\epsilon_0 \delta {a^2}_{\!\!\!1} , \dots ,$ \\ $
\epsilon_0 \delta {a^2}_{\!\!\!m} , \delta , \dots ,
\delta)$, satisfaisant aux équations (\ref{3.1}) et de hauteur
majorée par:
$$\widetilde{h}(P) \leq [14 (\eta + 6) \epsilon_1 + 4 \eta + 12] \delta
{a^2}_{\!\!\!1} + o(\delta) .$$
\end{proposition}
{\bf Démonstration.---} L'application du lemme de Siegel ci-dessus
à notre matrice $M$ qui est à éléments dans $K$ de format
$\mbox{card}~ T_{\delta} \times \mbox{card}~ \mathcal M$ (on
vérifiera ci-dessous que $\mbox{card}~ T_{\delta} <
\mbox{card}~\mathcal M$ pour $\delta$ assez grand) donne
l'existence d'une forme non identiquement nulle $P$ de
$K[{\underline X}_1 , \dots , {\underline X}_m , {\underline Y}_1
, \dots , {\underline Y}_{m - 1}] / I\left(\varphi_{\underline a}
(E^m)\right)$ de multidegré $(\epsilon_0 \delta {a^2}_{\!\!\!1} ,
\dots , \epsilon_0 \delta {a^2}_{\!\!\!m} , \delta , \dots ,
\delta)$ satisfaisant le syst\`eme d'équations (\ref{3.1}) et de
hauteur $\widetilde{h} (P) \leq e (\widetilde{h} (M) + \log
\mbox{card}~\mathcal M) + (1 + e) c_{\mbox{S}}(K)$ avec $e$ est
l'exposant de Dirichlet du syst\`eme (\ref{3.1}): $e :=
\frac{\mbox{card}~T_{\delta}}{\mbox{card}~\mathcal M -
\mbox{card}~T_{\delta}}$ et $c_{\mbox{S}}(K)$ est une constante ne
dépendant que du corps de nombres $K$. L'estimation de la
proposition \ref{c.12} pour $\widetilde{h}(P)$ suit ainsi des
estimations suivantes pour $\widetilde{h}(M) , e$ et $\log
\mbox{card}~\mathcal M$:
\begin{eqnarray}
\widetilde{h}(M) & \leq & [14 (\eta + 6) \epsilon_1 + 4 \eta + 12]
\delta
{a^2}_{\!\!\!1} \label{3.2} \\
e & \leq & 1 + o(1) \label{3.3} \\
\log \mbox{card}~\mathcal M & = & o(\delta) \label{3.4}
\end{eqnarray}
que nous démontrons ci-apr\`es. Commençons d'abord par démontrer
l'estimation (\ref{3.2}) pour $\widetilde{h}(M)$: la hauteur
$\widetilde{h}(M)$ est par définition la hauteur de Gauss-Weil du
vecteur formé de tous les éléments de $M$. Or, on remarque que
pour tout $(i_1 , \dots , i_m) \in T_{\delta}$ et tout $\mathfrak
m \in \mathcal M$ l'élément de $M$ correspendant ${\partial}^{(i_1
, \dots , i_m)}\mathfrak m (\underline 0 , \dots , \underline 0)$
est un coefficient de la forme ${\partial}^{(i_1 , \dots ,
i_m)}\mathfrak m$ et donc on a:
$$\widetilde{h}(M) ~\leq~ \widetilde{h}\left({\partial}^{(i_1 ,
\dots , i_m)}\mathfrak m ~;~ (i_1 , \dots , i_m) \in T_{\delta}
~,~ \mathfrak m \in \mathcal M\right).$$ La hauteur du membre de
droite de cette derni\`ere inégalité est majorée -gr\^ace au lemme
\ref{c.10}-
par: \\ \\
$\displaystyle \widetilde{h}\left({\partial}^{(i_1 , \dots ,
i_m)}\mathfrak m ~;~ (i_1 , \dots , i_m) \in T_{\delta} ~,~
\mathfrak m \in \mathcal M\right)
~\leq $ \\ \\
$\displaystyle (2 \eta + 12).\max_{(i_1 , \dots , i_m) \in
T_{\delta}} (i_1 + \dots + i_m) + \epsilon_0 \delta
{a^2}_{\!\!\!1} + \dots + \epsilon_0 \delta {a^2}_{\!\!\!m} + (3 \eta + 11) \delta \sum_{i = 1}^{m - 1} ({a^2}_{\!\!\!i} + 1)$. \\ \\
Comme maintenant, d'apr\`es nos hypoth\`eses, on a:
$$\epsilon_0 \leq \frac{1}{2} ~~\mbox{et}~~ \sum_{i = 1}^{m - 1} ({a^2}_{\!\!\!i} + 1) \leq \left(1 + \frac{1}{24}\right) {a^2}_{\!\!\!1}$$
et pour tout $(i_1 , \dots , i_m) \in T_{\delta}$:
$$i_1 + \dots + i_m ~\leq~ \left(\frac{i_1}{{a^2}_{\!\!\!1}} + \dots + \frac{i_{m - 1}}{{a^2}_{\!\!\!m - 1}} +
\frac{i_m}{m - 1}\right) {a^2}_{\!\!\!1} ~\leq~ 7 \epsilon_1
\delta {a^2}_{\!\!\!1} ,$$ alors:
$$\widetilde{h}\left({\partial}^{(i_1 ,
\dots , i_m)}\mathfrak m ~;~ (i_1 , \dots , i_m) \in T_{\delta}
~,~ \mathfrak m \in \mathcal M\right) ~\leq~ [14 (\eta + 6)
\epsilon_1 + 4 \eta + 12] \delta {a^2}_{\!\!\!1} .$$ D'o\`u cette
m\^eme majoration aussi pour $\widetilde{h}(M) ,$ c'est-\`a-dire:
$$\widetilde{h}(M) ~\leq~ [14 (\eta + 6) \epsilon_1 + 4 \eta + 12] \delta {a^2}_{\!\!\!1} ,$$ ce qui est l'estimation
(\ref{3.2}).\\ Afin d'établir les estimations (\ref{3.3}) et
(\ref{3.4}), fournissons des estimations pour les deux cardinaux
$\mbox{card}~\!T_{\delta}$ et $\mbox{card}~\!\mathcal M$:
\\ $\bullet$ pour $\mbox{card}~\!T_{\delta}$ on a (d'apr\`es le lemme $2.14.5$ de l'appendice):
$$\mbox{card}~\!T_{\delta} \leq \left(1 + o(1)\right) \mbox{vol}(T_{\delta}) = \left(1 + o(\delta)\right)
\frac{\left(7 \epsilon_1 \delta {a^2}_{\!\!\!1}\right)\dots
\left(7 \epsilon_1 \delta {a^2}_{\!\!\!m - 1}\right)7 \epsilon_1
\delta (m - 1)}{m!}$$ c'est-\`a-dire:
\begin{equation}
\mbox{card}~\!T_{\delta} ~\leq~ \frac{(m - 1) 7^m
{{\epsilon}^m}_{\!\!\!\!\!1}}{m!} {\delta}^m {a^2}_{\!\!\!1} \dots
{a^2}_{\!\!\!m}
 ~+~ o({\delta}^m) \label{3.5}
\end{equation}
$\bullet$ et pour $\mbox{card}~\!\mathcal M$ on a:
 \begin{eqnarray*}
\mbox{card}~\!\mathcal M & = & \dim {\left(\frac{K[{\underline
X}_1 , \dots , {\underline X}_m , {\underline Y}_1 , \dots ,
{\underline Y}_{m - 1}]}{I\left(\varphi_{\underline
a}(E^m)\right)}\right)}_{(\epsilon_0 \delta {a^2}_{\!\!\!1} ,
\dots , \epsilon_0 \delta {a^2}_{\!\!\!m} , \delta , \dots ,
\delta)} \\
 & = & \mathcal H \left(I\left(\varphi_{\underline a}(E^m)\right) ~;~ \epsilon_0 \delta {a^2}_{\!\!\!1} ,
\dots , \epsilon_0 \delta {a^2}_{\!\!\!m} , \delta , \dots ,
\delta\right)
 \end{eqnarray*}
 o\`u $\mathcal H$ désigne la fonction de Hilbert multihomog\`ene.
 Mais comme on sait que lorsque $\delta_1 , \dots , \delta_m , {{\delta}'}_{\!\!1} , \dots , {{\delta}'}_{\!\!m - 1}$
  sont des entiers positifs assez grands, $\mathcal H \!\left(I\left(\varphi_{\underline a}(E^m)\right) ~\!\!;~\!
  \delta_1 , \dots ,\right.$ \\ $\left.\delta_m , {{\delta}'}_{\!\!1} , \dots , {{\delta}'}_{\!\!m -
  1}\right)$ co\"{\i}ncide avec un polyn\^ome en $\delta_1 , \dots , \delta_m , {{\delta}'}_{\!\!1}
  , \dots , {{\delta}'}_{\!\!m - 1}$ dont la partie homog\`ene
  dominante vaut:
  $$\sum_{
  \begin{array}{c}
  \scriptstyle \{i_1 , \dots , i_m , j_1 , \dots , j_{m - 1}\}
  \subset \{0 , 1\} \\
  \scriptstyle i_1 + \dots + i_m + j_1 + \dots + j_{m - 1} = m
  \end{array}
  }\!\!\!\!\!\!\!\!\!\!\!\!\!\!\!\!\!\!\!\!d_{(i_1 , \dots , i_m , j_1 , \dots , j_{m - 1})} \left(\varphi_{\underline a} (E^m)\right) \frac{
  {{\delta}^{i_1}}_{\!\!\!\!\!1} \dots {{\delta}^{i_m}}_{\!\!\!\!\!m} {{{\delta}'}^{j_1}}_{\!\!\!\!\!\!\!1}
   \dots {{{\delta}'}^{j_{m - 1}}}_{\!\!\!\!\!\!\!\!\!\!\!\!\!m - 1}}{i_1! \dots i_m! j_1! \dots j_{m - 1}!}$$
   et que c'est bien notre cas ici puisque l'entier $\epsilon_0
   \delta$ est supposé assez grand, alors on a:
   $$\mbox{card}~\mathcal M
   =\!\!\!\!\!\!\!\!\!\!\!\!\!\!\!\!\!\!\!\!
\sum_{
  \begin{array}{c}
  \scriptstyle \{i_1 , \dots , i_m , j_1 , \dots , j_{m - 1}\}
  \subset \{0 , 1\} \\
  \scriptstyle i_1 + \dots + i_m + j_1 + \dots + j_{m - 1} = m
  \end{array}
  }\!\!\!\!\!\!\!\!\!\!\!\!\!\!\!\!\!\!\!\!d_{(i_1 , \dots , i_m , j_1
  , \dots , j_{m - 1})} \left(\varphi_{\underline a} (E^m)\right) \frac{{(\epsilon_0 \delta {a^2}_{\!\!\!1})}^{i_1}
  \dots {(\epsilon_0 \delta {a^2}_{\!\!\!m})} {\delta}^{j_1 + \dots + j_{m - 1}}}{i_1! \dots i_m! j_1! \dots j_{m - 1}!}$$
  $$= \sideset{}{^*}\sum d_{(I , J)} \left(\varphi_{\underline a}(E^m)\right) \prod_{i \in I} \!{a^2}_{\!\!\!i} . {\delta}^m
   {\epsilon_0}^{card ~\!\!I} + o({\delta}^m)~~~~~~~~~~~~~~~~$$
   o\`u la somme $\sum^{*}$ porte sur tous les couples d'ensembles
   $(I , J)$ tels que $I \subset \{1 , \dots , m\} ,$ \\ $J\subset\{1 , \dots , m -
   1\}$ et $\mbox{card}~\!I + \mbox{card}~\!J = m$, et la notation $d_{(I , J)} \left(\varphi_{\underline a}
   (E^m)\right)$ est celle du lemme \ref{c.6}.\\ En utilisant maintenant
   le lemme \ref{c.6}, on en déduit que:
   $$\mbox{card}~\!\mathcal M = {\sum}' \left(\prod_{i \in J \setminus I} {a^2}_{\!\!\!i}\right) {d(E)}^m. \prod_{i \in I}
   \!{a^2}_{\!\!\!i}. {\delta}^m. {\epsilon_0}^{card ~\!\!I} + o({\delta}^m)$$
   o\`u la somme ${\sum}'$ porte seulement sur les couples $(I ,
   J)$ de la somme $\sum^{*}$ vérifiant en plus: soit $m \in I$ et
   $I \cap J = \emptyset$ ou bien $m \not\in I$ et $I \cap J$ est
   un singleton. D'o\`u:
\begin{equation*}
\begin{split}
card ~\!\!\mathcal{M} &= {\sum}' {a^2}_{\!\!\!1} \dots {a^2}_{\!\!\!m} {d(E)}^m {\delta}^m \epsilon_{0}^{card ~\!\!I} + o({\delta}^m) \\
&= 3^m {\delta}^m {a^2}_{\!\!\!1} \dots {a^2}_{\!\!\!m}
\left({\sum}' \epsilon_{0}^{card ~\!\!I}\right) + o({\delta}^m) .
\end{split}
\end{equation*}
   Comme on a:
   \begin{eqnarray*}
   {\sum}'{\epsilon_0}^{card ~\!\!I} & = &  \sum_{k = 1}^{m}\left(\binom{m - 1}{k - 1} + k \binom{m - 1}{k}\right) {{\epsilon}^k}_{\!\!\!0} \\
    & = & \epsilon_0 {\left(1 + \epsilon_0\right)}^{m - 1} +
    \epsilon_0 (m - 1) {\left(1 + \epsilon_0\right)}^{m - 2} \\
     & = & \epsilon_0 (m + \epsilon_0) {(1 + \epsilon_0)}^{m - 2}
    \end{eqnarray*}
    on déduit enfin que:
\begin{equation}
\mbox{card}~\!\mathcal M = 3^m \epsilon_0 (m + \epsilon_0){(1 +
\epsilon_0)}^{m - 2} {\delta}^m {a^2}_{\!\!\!1} \dots
{a^2}_{\!\!\!m} + o({\delta}^m) \label{3.6}
\end{equation}
     Etablissons maintenant les deux estimations (\ref{3.3}) et (\ref{3.4}): nous remarquons que l'estimation (\ref{3.4}) découle immédiatement
     de (\ref{3.6}), par ailleurs pour établir l'estimation (\ref{3.3}), nous majorons gr\^ace à (\ref{3.5}) et (\ref{3.6}) l'exposant de Dirichlet $e$ par:
     $$e ~\leq~ \frac{\frac{(m - 1) 7^m {{\epsilon}^m}_{\!\!\!1}}{m!} {\delta}^m {a^2}_{\!\!\!1} \dots {a^2}_{\!\!\!m}
      + o({\delta}^m)}{3^m \epsilon_0 (m + \epsilon_0) {(1 + \epsilon_0)}^{m - 2} {\delta}^m {a^2}_{\!\!\!1} \dots
      {a^2}_{\!\!\!m} - \frac{(m - 1) 7^m {{\epsilon}^m}_{\!\!\!1}}{m!} {\delta}^m {a^2}_{\!\!\!1} \dots {a^2}_{\!\!\!m}
      + o({\delta}^m)}.~~~~$$
      En posant provisoirement:
\begin{eqnarray*}
\epsilon_2 & := & \frac{\frac{(m - 1) 7^m
{{\epsilon}^m}_{\!\!\!1}}{m!} {\delta}^m {a^2}_{\!\!\!1}
      \dots {a^2}_{\!\!\!m}}{3^m \epsilon_0 (m + \epsilon_0) {(1 + \epsilon_0)}^{m - 2}
{\delta}^m {a^2}_{\!\!\!1} \dots
      {a^2}_{\!\!\!m}} \\
& = & \frac{m - 1}{m!} {\left(\frac{7}{3}\right)}^{\!\!m}
\!\!\!\!\frac{{{\epsilon}^m}_{\!\!\!1}}{
      \epsilon_0 (m + \epsilon_0) {(1 + \epsilon_0)}^{m - 2}} ~\leq~ \frac{1}{2}~~~~ \mbox{(d'apr\`es nos hypoth\`eses),}
\end{eqnarray*}
on a enfin:
$$e ~\leq~ \frac{\epsilon_2 + o(1)}{1 - \epsilon_2 + o(1)} ~=~ \frac{\epsilon_2}{1 - \epsilon_2} + o(1) ~\leq~ 1 + o(1)$$
ce qui démontre l'estimation (\ref{3.3}) et ach\`eve cette
démonstration.
$~~~~\blacksquare$\vspace{1mm} \\
$~$ \\
En posant maintenant:
\begin{equation*}
Q({\underline X}_1 , \dots , {\underline X}_m) :=
P\!\left({\underline X}_1 , \dots , {\underline X}_m , \underline
D \!\left({\underline F}^{(a_1)}({\underline X}_1) , {\underline
X}_m\right) , \dots , \underline D \!\left({\underline F}^{(a_{m -
1})}({\underline X}_{m - 1}) , {\underline X}_m\right)\right)
\end{equation*}
on obtient le corollaire suivant:
\begin{corollary}\label{c.13}
Sous toutes les hypoth\`eses de cette section, notre forme $Q$ de
\\ $ K[{\underline X}_1 , \dots , {\underline X}_m] / I(E^m)$ est
non nulle, de multidegré $((2 + \epsilon_0) \delta {a^2}_{\!\!\!1}
, \dots , (2 + \epsilon_0) \delta {a^2}_{\!\!\!m - 1} ,$ \\ $(2 m
- 2 + \epsilon_0) \delta)$, satisfait ${\tau}^m (Q) (\underline 0
, \dots , \underline 0) \in {\mathfrak q}_{\delta}$ et pour toute
place $v$ sur $K$ on a:
\begin{description}
\item[$\bullet$] si $v$ est finie: $~~~~~~~~~\! h_v(Q) ~\leq~
h_v(P) + 4 m_v \delta {a^2}_{\!\!\!1}$ \item[$\bullet$] et si $v$
est infinie: $~~~~ \ell_v(Q) ~\leq~ \ell_v(P) + (4 m_v + 10)
\delta {a^2}_{\!\!\!1}$.
\end{description}
Par conséquent $Q$ est de hauteur de Gauss-Weil majorée par:
$$\left[14 (\eta + 6) \epsilon_1 + 8 \eta + 22\right] \delta {a^2}_{\!\!\!1} + \circ(\delta) .$$
\end{corollary}
{\bf Démonstration.---} Il est immédiat -du fait que $P$ est une
forme non identiquement nulle de $K[{\underline X}_1 , \dots ,
{\underline X}_m , {\underline Y}_1 , \dots , {\underline Y}_{m -
1}]$ et de multidegré $(\epsilon_0 \delta {a^2}_{\!\!\!1} , \dots
, \epsilon_0 \delta {a^2}_{\!\!\!m} , \delta , \dots , \delta)$-
que $Q$ est une forme non identiquement nulle de $K[{\underline
X}_1 , \dots , {\underline X}_m]$ de multidegré $((2 + \epsilon_0)
\delta {a^2}_{\!\!\!1} , \dots , (2 + \epsilon_0) \delta
{a^2}_{\!\!\!m - 1} , (2 m - 2 + \epsilon_0) \delta)$. De plus, on
a clairement $\tau^m(Q)(\underline 0 , \dots , \underline 0) =
\Omega_{\underline a}(P)(\underline 0 , \dots , \underline 0) \in
{\mathfrak q}_{\delta}$ (par construction m\^eme de notre forme
$P$). Par ailleurs, si $v$ est une place finie de $K$, on a:
\begin{eqnarray*}
h_v(Q) & \leq & \delta \left(2 h_v({\underline F}^{(a_1)}) + h_v(\underline D)\right) + \dots + \delta \left(2 h_v({\underline F}^{(a_{m - 1})}) + h_v(\underline D)\right) + h_v(P) \\
& \leq & \delta (3 m_v {a^2}_{\!\!\!1} + 3 m_v) + \dots + \delta (3 m_v {a^2}_{\!\!\!m - 1} + 3 m_v) + h_v(P) \\
& \leq & 3 m_v \delta \left({a^2}_{\!\!\!1} + \dots + {a^2}_{\!\!\!m - 1} + m - 1\right) + h_v(P) \\
& \leq & h_v(P) + 3 m_v (1 + 1/24) \delta {a^2}_{\!\!\!1} \\
& \leq & h_v(P) + 4 m_v \delta {a^2}_{\!\!\!1} ,
\end{eqnarray*}
ce qui est l'estimation du corollaire \ref{c.13} pour $h_v(Q)$
dans le cas $v$ finie.
$~$ \\

Si maintenant $v$ est une place infinie de $K$, on a le m\^eme
genre d'inégalité pour les longueurs logarithmiques locales,
c'est-à-dire:
\begin{equation*}
\begin{split}
\ell_v(Q) &\leq~ \delta \left(2 \ell_v({\underline F}^{(a_1)}) + \ell_v(\underline D)\right) + \dots + \delta \left(2 \ell_v({\underline F}^{(a_{m - 1})}) + \ell_v(\underline D)\right) + \ell_v(P) \\
& \!\!\!\!\!\!\!\!\!\!\!\!\!\!\!\leq~ \delta \left(3 (m_v + 3){a^2}_{\!\!\!1} + 3 m_v + 7 \right) + \dots + \delta \left(3 (m_v + 3){a^2}_{\!\!\!m - 1} + 3 m_v + 7 \right) + \ell_v(P) \\
& \!\!\!\!\!\!\!\!\!\!\!\!\!\!\!\leq~ 3 (m_v + 3) \delta \left({a^2}_{\!\!\!1} + \dots + {a^2}_{\!\!\!m - 1}\right) + (3 m_v + 7) (m - 1) \delta + \ell_v(P) \\
& \!\!\!\!\!\!\!\!\!\!\!\!\!\!\!\leq~ \ell_v(P) + (3 m_v + 9) \delta \left({a^2}_{\!\!\!1} + \dots + {a^2}_{\!\!\!m - 1} + m - 1\right) \\
& \!\!\!\!\!\!\!\!\!\!\!\!\!\!\!\leq~ \ell_v(P) + (3 m_v + 9) \delta (1 + 1/24){a^2}_{\!\!\!1} \\
& \!\!\!\!\!\!\!\!\!\!\!\!\!\!\!\leq~ \ell_v(P) + (4 m_v + 10)
\delta {a^2}_{\!\!\!1} ,
\end{split}
\end{equation*}
ce qui est aussi l'estimation du corollaire \ref{c.13} pour
$\ell_v(Q)$ (dans le cas $v$ infinie).

Finalement, pour aboutir à l'estimation du corollaire \ref{c.13},
on majore dans le cas $v$ infinie $\ell_v(P)$ par $h_v(P) +
o(\delta)$, ce qui nous permet d'écrire:
\begin{equation*}
\begin{split}
\forall v \in M_{K}^{\infty} : ~ h_v(Q) \leq \ell_v(Q) &\leq \ell_v(P) + (4 m_v + 10) \delta {a^2}_{\!\!\!1} \\
&\leq h_v(P) + (4 m_v + 10) \delta {a^2}_{\!\!\!1} + o(\delta)
.~~~~~~~~~~~~~~~~~
\end{split}
\end{equation*}
Puis, en reportant dans la définition de $\widetilde{h}(Q)$ ces
majorations des $h_v(Q)$ (selon les cas $v$ finie ou $v$ infinie),
on obtient:
$$\widetilde{h}(Q) \leq \widetilde{h}(P) + (4 \eta + 10) \delta {a^2}_{\!\!\!1} + o(\delta) ,$$
ce qui nous am\`ene -en utilisant la majoration de
$\widetilde{h}(P)$ donnée par la proposition \ref{c.12}- à
l'estimation du corollaire \ref{c.13} pour $\widetilde{h}(Q)$. La
démonstration est achevée.  $~~~~\blacksquare$\vspace{1mm}
\section{Translatée de la fonction auxiliaire}
Soient pour toute la suite de ce texte $\alpha$ un réel
strictement positif assez petit et ${\mathbf x}_1 , {\mathbf x}_2
, \dots , {\mathbf x}_m$ des points de $E(K)$, différents des
points de $2$-torsion, représentés respectivement dans ${\mathbb
P}_2$ par les syst\`emes de coordonnées projectives: ${\underline
x}_1 = (x_1 : 1 : z_1) , \dots , {\underline x}_m = (x_m : 1 :
z_m).$ \\ Désignons par $\widehat{h}$ la hauteur normalisée sur
$E(\overline{\mathbb Q})$ de la hauteur logarithmique absolue et
par $<. , .>$ le produit scalaire associé à la norme de Néron-Tate
$\nb{.} := \sqrt{\widehat{h}}$ dans l'espace euclidien
$E(\overline{\mathbb Q}) \otimes \mathbb R$. \\  Nous introduisons
à partir de cette section les deux hypoth\`eses ${\mbox H}_{\mbox
1}$ et ${\mbox H}_{\mbox 2}$ suivantes:
$${\mbox H}_{\mbox 1}: \left\{
\begin{array}{lr}
\bullet ~\forall i , j \in \{1 , 2 , \dots , m\} , i \neq j : &
\frac{<{\mathbf x}_i , {\mathbf x}_j>}{\mid{\mathbf
x}_i\mid.\mid{\mathbf
x}_j\mid} \geq 1 - \frac{\alpha}{4} \\
\bullet ~\forall i \in \{1 , 2 , \dots , m - 1\} : &
\frac{\mid{\mathbf
x}_m\mid}{\mid{\mathbf x}_i\mid } \geq \frac{1}{\sqrt{\alpha}} + 1 \\
\bullet ~\forall i \in \{2 , \dots , m\} : & \widehat{h}({\mathbf
x}_i) \geq 49 \widehat{h}({\mathbf x}_{i - 1}) \\
\bullet ~\widehat{h}({\mathbf x}_1) \geq 7 \eta + 37
\end{array}
\right.~~~~~~~~~~~$$
\begin{equation*}
\begin{split}
{\mbox H}_{\mbox 2}: \forall v \in S ~\text{et}~ \forall i \in \{1 , \dots , m\} : {\mid x_i \mid}_v &\leq \begin{cases} 1 & \text{si $v$ est finie} \\ 1 - 2/e & \text{si $v$ est infinie} \end{cases} \\
\text{et}~ {\mid z_i \mid}_v &\leq \begin{cases} 1 & \text{si $v$ est finie} \\
1/2 & \text{si $v$ est infinie} \end{cases} .
\end{split}
\end{equation*}
Posons pour tout $i = 1 , \dots , m$:
$$a_i := \left[\frac{\mid{\mathbf x}_m\mid}{\mid{\mathbf x}_i\mid}\right]$$
o\`u $[.]$ désigne la partie enti\`ere.\\ Le lemme \ref{c.31} de
l'appendice montre que les entiers positifs $a_1 , \dots , a_m$
ainsi définis vérifient bien les hypoth\`eses du §$6$ précédent.
Par conséquent, la proposition \ref{c.12} et le corollaire
\ref{c.13} de ce dernier restent valables pour le présent
paragraphe. En particulier, le lemme \ref{c.31} de l'appendice
montre qu'on a:
\begin{equation}
{a^2}_{\!\!\!1} + \dots + {a^2}_{\!\!\!m} \leq \left(1 +
\frac{1}{48}\right) {a^2}_{\!\!\!1} , \label{3.7}
\end{equation}
il montre aussi qu'on a: $a_1 \geq 7 (m - 1) \geq 7$ ce qui
entra\^ine: $m - 1 \leq \frac{a_1}{7} \leq {(\frac{a_1}{7})}^2$
c'est-\`a-dire:
\begin{equation}
m - 1 \leq \frac{{a^2}_{\!\!\!1}}{49} . \label{3.8}
\end{equation}
Le lemme \ref{c.31} montre aussi qu'on a:
\begin{equation}
m - 1 \leq \left(1 + \frac{1}{48}\right) \alpha {a^2}_{\!\!\!1} .
\label{3.9}
\end{equation}
Par ailleurs le lemme \ref{c.30} de l'appendice montre qu'on a
pour tout $i = 2 , \dots , m$:
\begin{equation}
{a^2}_{\!\!\!i} h({\mathbf x}_i) \leq 2 {a^2}_{\!\!\!1} h({\mathbf
x}_1) . \label{3.10}
\end{equation}
Soient aussi ${\mathbf y}_1 , \dots , {\mathbf y}_{m - 1}$ les $m
- 1$ points de $E(K)$ définis par:
$$
~~~~~~~~~~~~~~~~~~~~~~~~~~~~~~~~~~~~~~~~~{\mathbf y}_i := a_i
{\mathbf x}_i - {\mathbf x}_m~~~~~~~~~~~~~~~~~~~~~~~~~~~~i = 1 ,
\dots , m - 1$$ et $\mathbf y$ le point de $E^{2 m - 1}$:
$$~~~~~~~~~~~~~~\mathbf y := ({\mathbf x}_1 , \dots , {\mathbf x}_m , {\mathbf y}_1 , \dots , {\mathbf y}_{m - 1}).$$
Le lemme \ref{c.32} de l'appendice montre alors qu'on a pour tout
$i = 1 , \dots , m - 1$:
$$~\widehat{h}(a_i {\mathbf x}_i - {\mathbf x}_m) \leq \alpha \left({a^2}_{\!\!\!i} \widehat{h}({\mathbf x}_i) +
\widehat{h}({\mathbf x}_m)\right),~~~~~~~~~~$$ cette derni\`ere
inégalité entra\^ine -d'apr\`es le théor\`eme \ref{c.14} du
formulaire- qu'on a pour tout $i = 1 , \dots , m - 1$:
\begin{equation}
h({\mathbf y}_i) \leq \alpha \left[{a^2}_{\!\!\!i} h({\mathbf
x}_i) + h({\mathbf x}_m) + (3\!/\!4 . \eta + 5)({a^2}_{\!\!\!i} +
1)\right] + 3\!/\!2 . \eta + 8 . \label{3.11}
\end{equation}
Soit enfin pour tout $i = 1 , \dots , m : {\underline y}_i =:
(y_{i 0} : y_{i 1} : y_{i 2})$ un représentant dans ${\mathbb
P}_2$ du point ${\mathbf y}_i$. \\
Désignons maintenant par ${\underline D}^{(1)} , \dots ,
{\underline D}^{(m)} , {\underline D}^{(m + 1)} , \dots ,
{\underline D}^{(2 m - 1)}$ des formes de $K[\underline X ,
\underline Y]$, représentant la différence dans $E$ au voisinage
de $\{{\mathbf x}_1\} \times \{{\mathbf x}_1\} , \dots ,
\{{\mathbf x}_m\} \times \{{\mathbf x}_m\} , \{{\mathbf y}_1\}
\times \{{\mathbf y}_1\} , \dots , \{{\mathbf y}_{m - 1}\} \times
\{{\mathbf y}_{m - 1}\}$ respectivement, qu'on choisit de bidegrés
$(2 , 2)$ et de hauteurs logarithmiques locales $h_v$ (resp de
longueurs logarithmiques locales $\ell_v$) -pour une place finie
(resp infinie) $v$ de $K$- majorées par $3 m_v$ (resp $3 m_v + 7$)
et de hauteurs de Gauss-Weil majorées par $3 \eta + 5$ (ceci étant
possible d'apr\`es le théor\`eme \ref{c.8} du formulaire). En
désignant toujours par $P$ notre forme de $\frac{K[{\underline
X}_1, \dots , {\underline X}_m , {\underline Y}_1 , \dots ,
{\underline Y}_{m - 1}]}{I\left(\varphi_{\underline a}
(E^m)\right)}$ construite \`a la proposition \ref{c.12}, soit
${\tau^*}_{\!\!\!\!\!\!- {\mathbf y}}P$ la forme définie par:
\begin{equation*}
\begin{split}
{\tau^*}_{\!\!\!\!\!\!- {\mathbf y}}P := P \left({\underline
D}^{(1)} ({\underline X}_1 , {\underline x}_1)
 , \dots , {\underline D}^{(m)} ({\underline X}_m , {\underline x}_m) ,
 {\underline D}^{(m + 1)} ({\underline Y}_1 , {\underline y}_1) , \dots , \right. & \\
&\left.\!\!\!\!\!\!\!\!\!\!\!\!\!\!\!\!\!\!\!\!\!\!\!{\underline
D}^{(2 m - 1)} ({\underline Y}_{m - 1} , {\underline y}_{m -
1})\right).
\end{split}
\end{equation*}
\begin{lemma}\label{c.15}
Les coefficients de ${\tau^*}_{\!\!\!\!\!\!- \mathbf y}P \in
K[{\underline X}_1 , \dots , {\underline X}_m , {\underline Y}_1 ,
\dots , {\underline Y}_{m - 1}]$ sont des valeurs de formes de
$K[{\underline X'}_{\!\!1} , \dots , {\underline X'}_{\!\!m} ,
{\underline Y'}_{\!\!1} , \dots , {\underline Y'}_{\!\!m - 1}]$ au
point $({\underline x}_1 , \dots , {\underline x}_m , {\underline
y}_1 , \dots , {\underline y}_{m - 1})$. Ces formes sont de
multidegrés $(2 \epsilon_0 \delta {a^2}_{\!\!\!1} , \dots , 2
\epsilon_0 \delta {a^2}_{\!\!\!m} , 2 \delta , \dots , 2 \delta)$
et la famille ${\mathcal F}_1$ constituée de toutes ces formes vérifie pour toute place $v$ de $K$: \\
$\bullet$ si $v$ est finie:
\begin{align}
h_v({\mathcal F}_1) &\leq 2 m_v \delta {a^2}_{\!\!\!1} + h_v(P) \notag \\
\intertext{$\bullet$ et si $v$ est infinie:}
\ell_v({\mathcal F}_1) &\leq (2 m_v + 4) \delta {a^2}_{\!\!\!1} + \ell_v(P) . \notag \\
\intertext{De plus ${\mathcal F}_1$ est de hauteur de Gauss-Weil
majorée par:} \widetilde{h}({\mathcal F}_1) &\leq [14 (\eta + 6)
\epsilon_1 + 6 \eta + 16] \delta {a^2}_{\!\!\!1} + o(\delta) .
\notag
\end{align}
\end{lemma}
{\bf Démonstration.---} Soit $\overline{P} \in K[{\underline X}_1
, \dots , {\underline X}_m , {\underline Y}_1 , \dots ,
{\underline Y}_{m - 1} , {\underline X'}_{\!\!1} , \dots ,
{\underline X'}_{\!\!m} , {\underline Y'}_{\!\!1} , \dots ,$ \\
${\underline Y'}_{\!\!m - 1}]$ la forme égale à:
$$P\left({\underline D}^{(1)}({\underline X}_1 , {\underline
X'}_{\!\!1}) , \dots , {\underline D}^{(m)}({\underline X}_m ,
{\underline X'}_{\!\!m}) , {\underline D}^{(m + 1)}({\underline
Y}_1 , {\underline Y'}_{\!\!1}) , \dots , {\underline D}^{(2 m -
1)}({\underline Y}_{m - 1} , {\underline Y'}_{\!\!m - 1})\right)
.$$ Il est ainsi clair qu'on a:
$${\tau^*}_{\!\!\!\!\!\!- \mathbf y}P = \overline{P}\left({\underline X}_1 , \dots , {\underline X}_m ,
{\underline Y}_1 , \dots , {\underline Y}_{m - 1} , {\underline
x}_1 , \dots , {\underline x}_m , {\underline y}_1 , \dots ,
{\underline y}_{m - 1}\right)$$ et donc les coefficients de
${\tau^*}_{\!\!\!\!\!\!- \mathbf y}P$ sont les valeurs des
coefficients des mon\^omes en $\underline X , \underline Y$ de
$\overline{P}$ -vus comme formes de $K[{\underline X'}_{\!\!1} ,
\dots , {\underline X'}_{\!\!m} , {\underline Y'}_{\!\!1} , \dots
, {\underline Y'}_{\!\!m - 1}]$- au point $({\underline x}_1 ,
\dots , {\underline x}_m , {\underline y}_1 ,$ \\ $\dots ,
{\underline y}_{m - 1})$. Par conséquent, les formes dont il
s'agit dans le lemme \ref{c.15} sont tout simplement ces
coefficients de $\overline{P}$. Ces derniers sont effectivement
des formes de $K[{\underline X'}_{\!\!1} , \dots , {\underline
X'}_{\!\!m} , {\underline Y'}_{\!\!1} , \dots , {\underline
Y'}_{\!\!m - 1}]$ de multidegrés $({d°}_{\!\!\!\!{\underline
X'}_{\!\!1}} \overline{P} , \dots , {d°}_{\!\!\!\!{\underline
X'}_{\!\!m}} \overline{P} , {d°}_{\!\!\!\!{\underline Y'}_{\!\!1}}
\overline{P} , \dots , {d°}_{\!\!\!\!{\underline Y'}_{\!\!m - 1}}
\overline{P})$ et la famille ${\mathcal F}_1$ constituant toutes
ces formes vérifie clairement pour toute place $v$ de $K$:
$h_v({\mathcal F}_1) = h_v(\overline P)$ et $\ell_v({\mathcal
F}_1) \leq \ell_v(\overline P)$. Elle vérifie -par conséquent-
aussi $\widetilde{h}({\mathcal F}_1) = \widetilde{h}(\overline
P)$. La démonstration du lemme \ref{c.15} s'ach\`eve en démontrant
que les degrés de $\overline{P}$ par rapport à ${\underline
X'}_{\!\!1} , \dots , {\underline X'}_{\!\!m} , {\underline
Y'}_{\!\!1} , \dots , {\underline Y'}_{\!\!m - 1}$ sont
respectivement $2 \epsilon_0 \delta {a^2}_{\!\!\!1} , \dots , 2
\epsilon_0 \delta {a^2}_{\!\!\!m} ,$ \\ $2 \delta , \dots , 2 \delta$
et que les estimations du lemme \ref{c.15} concernant les hauteurs
et les longueurs logarithmiques locales ainsi que la hauteur de
Gauss-Weil sont valables m\^eme quand ${\mathcal F}_1$ est
remplacée par $\overline P$. Pour les degrés, comme la forme $P$
est de multidegré $(\epsilon_0 \delta {a^2}_{\!\!\!1} , \dots ,
\epsilon_0 \delta {a^2}_{\!\!\!m} , \delta , \dots , \delta)$ et
que les formes ${\underline D}^{(i)}$ $(i = 1 , \dots , 2 m - 1)$
sont toutes de bidegré $(2 , 2)$, on a clairement:
$${d°}_{\!\!\!\!{\underline X}_i} \overline{P} = {d°}_{\!\!\!\!{\underline X'}_{\!\!i}} \overline{P} = 2 \epsilon_0
\delta {a^2}_{\!\!\!i}~~~~\mbox{pour} ~i = 1 , \dots , m
~~~~\mbox{et}$$
$$~~~~{d°}_{\!\!\!\!{\underline Y}_i} \overline{P} = {d°}_{\!\!\!\!{\underline Y'}_{\!\!i}} \overline{P} = 2
\delta ~~~~~~~~~\mbox{pour} ~i = 1 , \dots , m - 1 .~~~~$$ Par
ailleurs quand $v$ est une place finie de $K$ on a:
\begin{equation*}
\begin{split}
h_v(\overline{P}) &\leq \epsilon_0 \delta {a^2}_{\!\!\!1} h_v \!\left({\underline D}^{(1)}\right) + \dots + \epsilon_0 \delta {a^2}_{\!\!\!m} h_v \!\left({\underline D}^{(m)}\right) \\
&\quad + \delta h_v \!\left({\underline D}^{(m + 1)}\right) + \dots + \delta h_v \!\left({\underline D}^{(2 m - 1)}\right) + h_v(P) \\
&\leq 3 m_v \epsilon_0 \delta \left({a^2}_{\!\!\!1} + \dots + {a^2}_{\!\!\!m}\right) + 3 m_v \delta (m - 1) + h_v(P) \\
&\leq 3\!/\!2 (1 + 1\!/\!48) m_v \delta {a^2}_{\!\!\!1} + 3\!/\!49
. m_v \delta {a^2}_{\!\!\!1} + h_v(P)
\end{split}
\end{equation*}
o\`u cette derni\`ere inégalité est obtenue en majorant
$\epsilon_0$ par $1/2$ et en utilisant (\ref{3.7}) et (\ref{3.8}).
D'o\`u:
$$h_v(\overline P) \leq 2 m_v \delta {a^2}_{\!\!\!1} + h_v(P) .$$
Ce qui entra{\sf\^\i}ne l'estimation du lemme \ref{c.15} pour
$h_v({\mathcal F}_1)$ (lorsque $v$ est finie). Si maintenant $v$
est une place infinie de $K$, on a:
\begin{equation*}
\begin{split}
\ell_v(\overline{P}) &\leq \epsilon_0 \delta {a^2}_{\!\!\!1} \ell_v \!\left({\underline D}^{(1)}\right) + \dots + \epsilon_0 \delta {a^2}_{\!\!\!m} \ell_v \!\left({\underline D}^{(m)}\right) \\
&\quad + \delta \ell_v \!\left({\underline D}^{(m + 1)}\right) + \dots + \delta \ell_v \!\left({\underline D}^{(2 m - 1)}\right) + \ell_v(P) \\
&\leq (3 m_v + 7) \epsilon_0 \delta \left({a^2}_{\!\!\!1} + \dots + {a^2}_{\!\!\!m}\right) + (3 m_v + 7) \delta (m - 1) + \ell_v(P) \\
&\leq (2 m_v + 4) \delta {a^2}_{\!\!\!1} + \ell_v(P)
\end{split}
\end{equation*}
o\`u cette derni\`ere inégalité est obtenue (comme dans le cas $v$
finie) en majorant $\epsilon_0$ par $1/2$ et en utilisant
(\ref{3.7}) et (\ref{3.8}). Ce qui entra{\sf\^\i}ne aussi
l'estimation du lemme \ref{c.15} pour $\ell_v({\mathcal F}_1)$
(lorsque $v$ est infinie). Il résulte de ces deux estimations pour
$h_v(\overline P)$ (quand $v$ est finie) et $\ell_v(\overline P)$
(quand $v$ est infinie) qu'on a:
$$\widetilde{h}(\overline P) \leq (2 \eta + 4) \delta {a^2}_{\!\!\!1} + \widetilde{h}(P) + o(\delta) .$$
En reportant finalement $\widetilde{h}(P)$ par son estimation
donnée à la proposition \ref{c.12}, on a:
$$\widetilde{h}(\overline P) \leq [14 (\eta + 6) \epsilon_1 + 6 \eta + 16] \delta {a^2}_{\!\!\!1} + o(\delta) .$$
Ce qui entra{\sf\^\i}ne l'estimation du lemme \ref{c.15} pour
$\widetilde{h}({\mathcal F}_1)$ et ach\`eve cette démonstration.
 $~~~~\blacksquare$\vspace{1mm}
\begin{lemma}\label{c.16}
Pour tout $(i_1 , \dots , i_m) \in {\mathbb N}^m$, les
coefficients de $\partial^{(i_1 , \dots , i_m)}
({\tau^*}_{\!\!\!\!\!\!- \mathbf y}P)$ \\ $\in K[{\underline X}_1
, \dots , {\underline X}_m]$ sont des valeurs de formes de
$K[{\underline X'}_{\!\!} , \dots , {\underline X'}_{\!\!m} ,
{\underline Y'}_{\!\!1} , \dots , {\underline Y'}_{\!\!m - 1}]$ au
point $({\underline x}_1 , \dots , {\underline x}_m , {\underline
y}_1 , \dots , {\underline y}_{m - 1})$. Ces formes sont de
multidegrés $(2 \epsilon_0 \delta {a^2}_{\!\!\!1} , \dots , 2
\epsilon_0 \delta {a^2}_{\!\!\!m} ,$ \\ $ 2 \delta , \dots , 2
\delta)$ et pour $T \in \mathbb N$, la famille ${\mathcal F}_2$
constituée de telles formes, correspondant \`a tous les $m$-uplets
$(i_1 , \dots , i_m) \in
{\mathbb N}^m$ tels que $i_1 + \dots + i_m \leq T$ satisfait pour toute place $v$ de $K$:\\
\indent{$\bullet$} si $v$ est finie:
\begin{align}
h_v({\mathcal F}_2) &\leq 2 m_v T + 9 m_v \delta {a^2}_{\!\!\!1} + h_v(P) \notag \\
\intertext{\indent{$\bullet$} et si $v$ est infinie:}
\ell_v({\mathcal F}_2) &\leq (2 m_v + 12)T + (9 m_v + 28) \delta {a^2}_{\!\!\!1} + h_v(P) + o(\delta) . \notag \\
\intertext{De plus ${\mathcal F}_2$ est de hauteur de Gauss-Weil
majorée par:} \widetilde{h}({\mathcal F}_2) &\leq (2 \eta + 12)T +
[14(\eta + 6)\epsilon_1 + 13 \eta + 40] \delta {a^2}_{\!\!\!1} +
o(\delta) . \notag
\end{align}
\end{lemma}
{\bf Démonstration.---}
 En écrivant:
 $${\tau^*}_{\!\!\!\!\!\!- \mathbf y} P =: \sum_{\mathfrak m \in \Gamma} \rho_{\mathfrak m} . \mathfrak m$$
o\`u $\Gamma$ désigne un ensemble fini de mon\^omes unitaires de
 $K[{\underline X}_1 , \dots , {\underline X}_m , {\underline Y}_1 , \dots , {\underline Y}_{m -
 1}]$ de multidegrés $(2 \epsilon_0 \delta {a^2}_{\!\!\!1} , \dots , 2 \epsilon_0 \delta {a^2}_{\!\!\!m} , 2 \delta ,
 \dots , 2 \delta)$, on a pour tout $(i_1 , \dots , i_m) \in {\mathbb
 N}^m$:
 $$\partial^{(i_1 , \dots , i_m)} ({\tau^*}_{\!\!\!\!\!\!- \mathbf y}P) = \sum_{\mathfrak m \in \Gamma}
 \rho_{\mathfrak m} . \partial^{(i_1 , \dots , i_m)} \mathfrak m .$$
D'o\`u, en gardant les notations du lemme \ref{c.15}:
$$h_v({\mathcal F}_2) \leq h_v({\mathcal F}_1) + h_v \!\left(\partial^{(i_1 , \dots , i_m)}\mathfrak m , \mathfrak m \in \Gamma , i_1 + \dots + i_m \leq T \right)$$
lorsque $v$ est une place finie de $K$ et:
$$\ell_v({\mathcal F}_2) \leq \ell_v({\mathcal F}_1) + \ell_v \!\left(\partial^{(i_1 , \dots , i_m)}\mathfrak m , \mathfrak m \in \Gamma , i_1 + \dots + i_m \leq T \right) + \log \mbox{card}~\!\Gamma$$
lorsque $v$ est une place infinie de $K$. En utilisant maintenant les estimations du lemme \ref{c.10} pour les hauteurs et longueurs logarithmiques
locales de la famille \\ $\{\partial^{(i_1 , \dots , i_m)}\mathfrak m , \mathfrak m \in \Gamma , i_1 + \dots + i_m \leq T\}$ et les estimations du lemme \ref{c.15} pour les hauteurs et longueurs logarithmiques locales de la famille ${\mathcal F}_1$, on a: \\
$\bullet$ si $v$ est une place finie de $K$:
\begin{equation*}
\begin{split}
h_v({\mathcal F}_2) &\leq 2 m_v \delta {a^2}_{\!\!\!1} + h_v(P) + \left[2 T + \sum_{i = 1}^{m - 1} 6 ({a^2}_{\!\!\!i} + 1) \delta \right] m_v \\
&\leq 2 m_v T + 9 m_v \delta {a^2}_{\!\!\!1} + h_v(P)
\end{split}
\end{equation*}
o\`u cette derni\`ere inégalité est obtenue gr\^ace à la majoration: ${a^2}_{\!\!\!1} + \dots + {a^2}_{\!\!\!m - 1} + m - 1 \leq (1 + 1/24) {a^2}_{\!\!\!1}$. L'estimation du lemme \ref{c.16} pour $h_v({\mathcal F}_2)$ -lorsque $v$ est finie- est alors démontrée. \\
$\bullet$ Et si $v$ est une place infinie de $K$:
\begin{equation*}
\begin{split}
\ell_v({\mathcal F}_2) &\leq (2 m_v + 4) \delta {a^2}_{\!\!\!1} + \ell_v(P) + \left[2 T + \sum_{i = 1}^{m - 1} 6({a^2}_{\!\!\!i} + 1) \delta \right] m_v + 12 T \\
&\quad + 2 \epsilon_0 \delta\left({a^2}_{\!\!\!1} + \dots + {a^2}_{\!\!\!m}\right) + \sum_{i = 1}^{m - 1} (22 {a^2}_{\!\!\!i} + 18) \delta + \log \mbox{card}~\!\Gamma \\
&\!\!\!\!\!\!\!\!\!\!\leq (2 m_v + 12) T + [2 m_v + 4 + 6 (1 + 1\!/\!24) m_v + 22 (1 + 1\!/\!24) + 1 + 1\!/\!48] \delta {a^2}_{\!\!\!1} \\
&\!\!\!\!\!\!\!\!\!\!\quad + h_v(P) + o(\delta) \\
&\!\!\!\!\!\!\!\!\!\!\leq (2 m_v + 12) T + (9 m_v + 28) \delta
{a^2}_{\!\!\!1} + h_v(P) + o(\delta)
\end{split}
\end{equation*}
o\`u l'avant derni\`ere inégalité est obtenue en utilisant les
majorations: ${a^2}_{\!\!\!1} + \dots + {a^2}_{\!\!\!m - 1} + m -
1 \leq (1 + 1/24) {a^2}_{\!\!\!1}$, $\epsilon_0 \leq 1/2$,
${a^2}_{\!\!\!1} + \dots + {a^2}_{\!\!\!m} \leq (1 + 1/48)
{a^2}_{\!\!\!1}$ et $\ell_v(P) \leq h_v(P) + o(\delta)$ et en
remarquant que $\log \mbox{card}~\!\Gamma = o(\delta)$. Ce qui
montre alors l'estimation du lemme \ref{c.16} pour
$\ell_v({\mathcal F}_2)$ (dans le cas $v$ infinie) aussi. On
déduit ainsi pour la hauteur de Gauss-Weil de la famille
${\mathcal F}_2$ la majoration suivante:
$$\widetilde{h}({\mathcal F}_2) \leq (2 \eta + 12) T + [14 (\eta + 6) \epsilon_1 + 13 \eta + 40] \delta {a^2}_{\!\!\!1} + o(\delta) .$$
Ceci ach\`eve cette démonstration. $~~~~\blacksquare$\vspace{1mm}
\begin{lemma}\label{c.22}
Pour tout $T \in \mathbb N$, la famille des formes $\partial^{(i_1
, \dots , i_m)} ({\tau^*}_{\!\!\!\!\!\!-\mathbf y}P) , i_1 + \dots
+ i_m \leq T$ est de hauteur logarithmique locale $h_v$ majorée
par:
\begin{equation*}
\begin{split}
&2 m_v T + 9 m_v \delta {a^2}_{\!\!\!1} + h_v(P) + 2 \epsilon_0 \delta \left({a^2}_{\!\!\!1} h_v({\underline x}_1) + \dots + {a^2}_{\!\!\!m} h_v({\underline x}_m)\right) \\
&\!\!\!\!\!\quad + 2 \delta \left(h_v({\underline y}_1) + \dots +
h_v({\underline y}_{m - 1})\right)
\end{split}
\end{equation*}
lorsque $v$ est une place finie de $K$ et par:
\begin{equation*}
\begin{split}
&(2 m_v + 12) T + (9 m_v + 28) \delta {a^2}_{\!\!\!1} + h_v(P) + 2 \epsilon_0 \delta \left({a^2}_{\!\!\!1} h_v({\underline x}_1) + \dots + {a^2}_{\!\!\!m} h_v({\underline x}_m)\right) \\
&\!\!\!\!\!\quad + 2 \delta \left(h_v({\underline y}_1) + \dots +
h_v({\underline y}_{m - 1})\right) + o(\delta)
\end{split}
\end{equation*}
lorsque $v$ est une place infinie de $K$. \\
De plus elle est de hauteur de Gauus-Weil majorée par:
\begin{equation*}
\begin{split}
&(2 \eta + 12) T + [4 m (\epsilon_0 + 2 \alpha) h({\mathbf x}_1) + 14 (\eta + 6) \epsilon_1 + (5 \eta + 27) \alpha + 13 \eta + 40] \delta {a^2}_{\!\!\!1} \\
&\!\!\!\!\!\quad + o(\delta) .
\end{split}
\end{equation*}
\end{lemma}
{\bf Démonstration.---} Chaque coefficient $c$ de l'une des formes
$\partial^{(i_1 , \dots , i_m)} ({\tau^*}_{\!\!\!\!\!\!- \mathbf
y}P) ,$ \\ $i_1 + \dots + i_m \leq T$ s'écrit d'apr\`es le lemme
\ref{c.16}:
$$c = {\sum}_1 \gamma_c({\underline{\alpha}}_1 , \dots , {\underline{\alpha}}_m , {\underline{\beta}}_1
, \dots , {\underline{\beta}}_{m - 1}) {{\underline
x}_1}^{{\underline{\alpha}}_1} \dots {{\underline
x}_m}^{{\underline{\alpha}}_m} {{\underline
y}_1}^{{\underline{\beta}}_1} \dots {{\underline y}_{m -
1}}^{{\underline{\beta}}_{m - 1}}$$ o\`u la somme $\sum_1$ porte
sur tous les triplets ${\underline{\alpha}}_1 , \dots ,
{\underline{\alpha}}_m , {\underline{\beta}}_1 , \dots ,
{\underline{\beta}}_{m - 1}$ de ${\mathbb N}^3$ satisfaisant:
$\mid{\underline{\alpha}}_1 \mid = 2 \epsilon_0 \delta
{a^2}_{\!\!\!1} , \dots , \mid{\underline{\alpha}}_m \mid = 2
\epsilon_0 \delta {a^2}_{\!\!\!m} , \mid{\underline{\beta}}_1 \mid
= \mid{\underline{\beta}}_2 \mid = \dots =
\mid{\underline{\beta}}_{m - 1}\mid = 2 \delta$ et pour tout $c$,
la famille $\{\gamma_c({\underline{\alpha}}_1 , \dots ,
{\underline{\alpha}}_m , {\underline{\beta}}_1 , \dots ,
{\underline{\beta}}_{m - 1}) , {\underline{\alpha}}_1 , \dots ,
{\underline{\alpha}}_m , {\underline{\beta}}_1 , \dots ,
{\underline{\beta}}_{m - 1}\}$, formée de nombres de $K$, est de hauteur logarithmique locale $v$-adique $\leq h_v({\mathcal F}_2) \leq 2 m_v T + 9 m_v \delta {a^2}_{\!\!\!1} + h_v(P)$ lorsque $v$ est une place finie de $K$ et de longueur logarithmique locale $v$-adique $\leq \ell_v({\mathcal F}_2) \leq (2 m_v + 12) T + (9 m_v + 28) \delta {a^2}_{\!\!\!1} + h_v(P) + o(\delta)$ lorsque $v$ est une place infinie de $K$. On a alors pour toute place $v$ de $K$: \\
$\bullet$ si $v$ est finie:
\begin{equation}
\begin{split}
h_v(c) &\leq h_v({\mathcal F}_2) + 2 \epsilon_0 \delta {a^2}_{\!\!\!1} h_v({\underline x}_1) + \dots + 2 \epsilon_0 \delta {a^2}_{\!\!\!m} h_v({\underline x}_m) \\
&\quad + 2 \delta \left(h_v({\underline y}_1) + \dots + h_v({\underline y}_{m - 1})\right) \\
&\leq 2 m_v T + 9 m_v \delta {a^2}_{\!\!\!1} + h_v(P) + 2 \epsilon_0 \delta \left({a^2}_{\!\!\!1} h_v({\underline x}_1) + \dots + {a^2}_{\!\!\!m} h_v({\underline x}_m)\right) \\
&\quad + 2 \delta \left(h_v({\underline y}_1) + \dots +
h_v({\underline y}_{m - 1})\right) \label{3.12}
\end{split}
\end{equation}
$\bullet$ et si $v$ est infinie:
\begin{equation}
\begin{split}
h_v(c) &\leq \ell_v({\mathcal F}_2) + 2 \epsilon_0 \delta {a^2}_{\!\!\!1} h_v({\underline x}_1) + \dots + 2 \epsilon_0 \delta {a^2}_{\!\!\!m} h_v({\underline x}_m) \\
&\quad + 2 \delta \left(h_v({\underline y}_1) + \dots + h_v({\underline y}_{m - 1})\right) \\
&\leq (2 m_v + 12) T + (9 m_v + 28) \delta {a^2}_{\!\!\!1} + h_v(P) \\
&\quad + 2 \epsilon_0 \delta \left({a^2}_{\!\!\!1} h_v({\underline x}_1) + \dots + {a^2}_{\!\!\!m} h_v({\underline x}_m)\right) \\
&\quad + 2 \delta \left(h_v({\underline y}_1) + \dots +
h_v({\underline y}_{m - 1})\right) + o(\delta) . \label{3.21}
\end{split}
\end{equation}
Les estimations locales du lemme \ref{c.22} suivent clairement de
(\ref{3.12}) et (\ref{3.21}) puisque ces derni\`eres sont valables
pour tout coefficient $c$ de l'une des formes $\partial^{(i_1 ,
\dots , i_m)}({\tau^*}_{\!\!\!\!\!\!- \mathbf y}P) , i_1 + \dots +
i_m \leq T$. On déduit aussi de (\ref{3.12}) et (\ref{3.21}) que
la hauteur de Gauss-Weil de la famille des formes
$\{\partial^{(i_1 , \dots , i_m)}({\tau^*}_{\!\!\!\!\!\!- \mathbf
y}P) , i_1 + \dots + i_m \leq T\}$ est majorée par:
\begin{equation}
\begin{split}
\widetilde{h}\left(\partial^{(i_1 , \dots , i_m)}({\tau^*}_{\!\!\!\!\!\!- \mathbf y}P) , i_1 + \dots + i_m \leq T\right) &\leq (2 \eta + 12) T \\
&\!\!\!\!\!\!\!\!\!\!\!\!\!\!\!\!\!\!\!\!\!\!\!\!\!\!\!\!\!\!\!\!\!\!\!\!\!+ (9 \eta + 28) \delta {a^2}_{\!\!\!1} + \widetilde{h}(P) + 2 \epsilon_0 \delta \left({a^2}_{\!\!\!1} h({\mathbf x}_1) + \dots + {a^2}_{\!\!\!m} h({\mathbf x}_m)\right) \\
&~~~~~~\!\!\!+ 2 \delta \left(h({\mathbf y}_1) + \dots +
h({\mathbf y}_{m - 1})\right) + o(\delta) . \label{3.22}
\end{split}
\end{equation}
En utilisant maintenant les majorations (\ref{3.10}), on a:
\begin{equation}
{a^2}_{\!\!\!1} h({\mathbf x}_1) + \dots + {a^2}_{\!\!\!m}
h({\mathbf x}_m) \leq 2 m {a^2}_{\!\!\!1} h({\mathbf x}_1)
\label{3.23}
\end{equation}
et en utilisant les majorations (\ref{3.11}), on a:
\begin{equation}
\begin{split}
h({\mathbf y}_1) + \dots &+ h({\mathbf y}_{m - 1}) \leq \alpha \left[\left({a^2}_{\!\!\!1} h({\mathbf x}_1) + \dots + {a^2}_{\!\!\!m - 1} h({\mathbf x}_{m - 1})\right) \right. \\
&\quad\left. + (m - 1) h({\mathbf x}_m) + (3\!/\!4 \eta + 5)({a^2}_{\!\!\!1} + \dots + {a^2}_{\!\!\!m - 1} + m - 1)\right] \\
&\quad + (3\!/\!2 \eta + 8)(m - 1) \\
&\leq \alpha \left[(4 m - 5) {a^2}_{\!\!\!1} h({\mathbf x}_1) + (3\!/\!4 \eta + 5) (1 + 1\!/\!24) {a^2}_{\!\!\!1}\right] \\
&\quad + (3\!/\!2 \eta + 8)(1 + 1\!/\!48) \alpha {a^2}_{\!\!\!1} \\
&\leq \alpha \left[4 m h({\mathbf x}_1) + 5\!/\!2 \eta +
27\!/\!2\right] {a^2}_{\!\!\!1} \label{3.24}
\end{split}
\end{equation}
o\`u l'avant derni\`ere inégalité est obtenue gr\^ace aux
majorations (\ref{3.10}), la majoration ${a^2}_{\!\!\!1} + \dots +
{a^2}_{\!\!\!m - 1} + m - 1 \leq (1 + 1\!/\!24) {a^2}_{\!\!\!1}$
et la majoration (\ref{3.9}). En reportant finalement les deux
majorations (\ref{3.23}) et (\ref{3.24}) ainsi que l'estimation de
$\widetilde{h}(P)$ donnée par la proposition \ref{c.12} dans
(\ref{3.22}) on obtient finalement:
\begin{equation*}
\begin{split}
\widetilde{h}\left(\partial^{(i_1 + \dots + i_m)}({\tau^*}_{\!\!\!\!\!\!- \mathbf y}P) , i_1 + \dots + i_m \leq T\right) \leq (2 \eta + 12) T + [4 m (\epsilon_0 + 2 \alpha) h({\mathbf x}_1) \\
+ 14(\eta + 6)\epsilon_1 + (5 \eta + 27) \alpha + 13 \eta + 40]
\delta {a^2}_{\!\!\!1} + o(\delta)
\end{split}
\end{equation*}
ce qui ach\`eve cette démonstration.
$~~~~\blacksquare$\vspace{1mm}
\section{Extrapolation}
Soit pour toute la suite de ce texte $\epsilon$ un réel
strictement positif assez petit. Les param\`etres $\epsilon_0 ,
\epsilon_1$ et $\alpha$ seront choisis \`a la fin en fonction de
$\epsilon , D , m$ et $\eta$. Par ailleurs, dans le but de
démontrer les théor\`emes principaux \ref{c.1}, \ref{c.2} et
\ref{c.3}, introduisons $S$ un ensemble fini de places sur $K$ et
${(\lambda_v)}_{v \in S}$ une famille de réels positifs
satisfaisants:
$$\sum_{v \in S} \frac{[K_v : {\mathbb Q}_v]}{[K : \mathbb Q]} \lambda_v ~=~ 1 .$$ Outre les hypoth\`eses des §§$6$ et $7$ précédents, on fait, dans le présent paragraphe, les hypoth\`eses supplémentaires suivantes: \\
$\displaystyle {\bf{H}}_3 \!: ~~ 4 m (\epsilon_0 + 2 \alpha) \leq \epsilon \epsilon_1$. \\
$\displaystyle {\bf{H}}_4 \!: ~~ \widehat{h}({\mathbf x}_1) \geq \frac{1}{\epsilon \epsilon_1}\left\{(28 \eta + 176) \epsilon_1 + (7 \eta + 36) \alpha + (18 \eta + 54)\right\} + \frac{3}{4} \eta + 5 .$ \\
$\displaystyle {\bf{H}}_5 \!: ~~ $(l'hypoth\`ese principale) \\
$$\forall i \in \{1 , \dots , m\} ~\!,~\! \forall v \in S : ~{\mbox{dist}}_v({\mathbf x}_i , \mathbf 0) ~<~ e^{- \lambda_v \epsilon h({\mathbf x}_i) - 2 m_v - c_v} ,$$
o\`u $c_v$ est une constante dépendant de $v$ déja introduite au
§$2$, qui vaut $0$ si $v$ est finie et $16$ si $v$ est infinie.
\begin{proposition}\label{c.23}
On a:
$$\Omega_{\underline a}({\tau^*}_{\!\!\!\!\!\!- \mathbf y}P)({\underline x}_1 , \dots , {\underline x}_m) ~= \!\!\!\sum_{i_1 \geq 0 , \dots , i_m \geq 0}\!\!\!f_{i_1 , \dots , i_m} {u^{i_1}}_{\!\!\!\!1} \dots {u^{i_m}}_{\!\!\!\!\!m}$$
o\`u les $f_{i_1 , \dots , i_m} := \frac{\partial^{(i_1 , \dots ,
i_m)} ({\tau^*}_{\!\!\!\!\!\!- \mathbf y}P)({\underline x}_1 ,
\dots , {\underline x}_m)}{{\Delta({\underline x}_1)}^{f(i_1)}
\dots {\Delta({\underline x}_m)}^{f(i_m)}} ,~ (i_1 , \dots , i_m)
\in {\mathbb N}^m$ sont des nombres de $K$, nuls pour $\underline
i = (i_1 , \dots , i_m) \in T_{\delta}$ et vérifiant pour tout
$\underline i \in {\mathbb N}^m$ et pour toute place $v \in S$:
\begin{equation*}
\begin{split}
{\mid f_{\underline i} \mid}_v & \leq \exp\left\{\phantom{\left(h_v({\underline y}_1) + \dots + h_v({\underline y}_{m - 1})\right)} \!\!\!\!\!\!\!\!\!\!\!\!\!\!\!\!\!\!\!\!\!\!\!\!\!\!\!\!\!\!\!\!\!\!\!\!\!\!\!\!\!\!\!\!\!\!\!\!\!\!\!\!\!\!\!\!\!\!\!\!\!\!\!\!\!\!9 m_v \delta {a^2}_{\!\!\!1} + h_v(P) + 2 \epsilon_0 \delta \left({a^2}_{\!\!\!1} h_v({\underline x}_1) + \dots + {a^2}_{\!\!\!m} h_v({\underline x}_m)\right)\right. \\
& \quad \left.+ 2 \delta \left(h_v({\underline y}_1) + \dots +
h_v({\underline y}_{m - 1})\right)\right\} \left(e^{2
m_v}\right)^{i_1 + \dots + i_m}
\end{split}
\end{equation*}
si $v$ est finie et:
\begin{equation*}
\begin{split}
{\mid f_{\underline i} \mid}_v & \leq \exp\left\{\phantom{\left(h_v({\underline y}_1) + \dots + h_v({\underline y}_{m - 1})\right)} \!\!\!\!\!\!\!\!\!\!\!\!\!\!\!\!\!\!\!\!\!\!\!\!\!\!\!\!\!\!\!\!\!\!\!\!\!\!\!\!\!\!\!\!\!\!\!\!\!\!\!\!\!\!\!\!\!\!\!\!\!\!\!\!\!\!(9 m_v + 28)\delta {a^2}_{\!\!\!1} + h_v(P) + 2 \epsilon_0 \delta \left({a^2}_{\!\!\!1} h_v({\underline x}_1) + \dots + {a^2}_{\!\!\!m} h_v({\underline x}_m)\right)\right. \\
& \quad \left.+ 2 \delta \left(h_v({\underline y}_1) + \dots +
h_v({\underline y}_{m - 1})\right) + o(\delta)\right\} \left(e^{2
m_v + 13}\right)^{i_1 + \dots + i_m}
\end{split}
\end{equation*}
si $v$ est infinie.
\end{proposition}
{\bf Démonstration.---} Le fait que les $f_{\underline i} ~
(\underline i \in {\mathbb N}^m)$ sont nuls pour $\underline i \in
T_{\delta}$ vient de ce que notre forme $P$ satisfait -par
construction m\^eme- la condition d'annulation:
$\Omega_{\underline a}(P)(\underline 0 , \dots , \underline 0) \in
{\mathfrak q}_{\delta}$. Par ailleurs, écrivons pour un $(i_1 ,
\dots , i_m) \in {\mathbb N}^m$:
$$\partial^{(i_1 , \dots , i_m)}
({\tau^*}_{\!\!\!\!\!\!- \mathbf y}P)({\underline X}_1 , \dots ,
{\underline X}_m) ~= \!\!\!\sum_{\mid {\underline{\alpha}}_1 \mid
= d_1 , \dots , \mid {\underline{\alpha}}_m \mid = d_m}
\!\!\!\!\!\!\lambda({\underline{\alpha}}_1 , \dots ,
{\underline{\alpha}}_m) {{\underline
X}_1}^{{\underline{\alpha}}_1} \dots {{\underline
X}_m}^{{\underline{\alpha}}_m}$$ o\`u $d_1 , \dots , d_m$ sont les
degrés de la forme $\partial^{(i_1 , \dots ,
i_m)}({\tau^*}_{\!\!\!\!\!\!- \mathbf y}P)$ par rapport \`a
${\underline X}_1 , \dots , {\underline X}_m$ respectivement et la
famille $\!\!\{\lambda({\underline{\alpha}}_1 , \dots ,
{\underline{\alpha}}_m) , {\underline{\alpha}}_1 , \dots ,
{\underline{\alpha}}_m\}$, constituée de nombres de $K$, de
hauteur logarithmique $v$-adique (pour une place quelconque $v$ de
$K$) majorée par le lemme \ref{c.22}. En spécialisant
$({\underline X}_1 , \dots , {\underline X}_m)$ en $({\underline
x}_1 , \dots , {\underline x}_m)$ on a:
$$\partial^{(i_1 , \dots , i_m)}
({\tau^*}_{\!\!\!\!\!\!- \mathbf y}P)({\underline x}_1 , \dots ,
{\underline x}_m) ~= \!\!\!\sum_{\mid {\underline{\alpha}}_1 \mid
= d_1 , \dots , \mid {\underline{\alpha}}_m \mid = d_m}
\!\!\!\!\!\!\lambda({\underline{\alpha}}_1 , \dots ,
{\underline{\alpha}}_m) {{\underline
x}_1}^{{\underline{\alpha}}_1} \dots {{\underline
x}_m}^{{\underline{\alpha}}_m} .$$
D'o\`u, pour toute place $v \in S$: \\
$\displaystyle {\mid \partial^{(i_1 , \dots , i_m)}
({\tau^*}_{\!\!\!\!\!\!- \mathbf y}P)({\underline x}_1 , \dots ,
{\underline x}_m) \mid}_v$
$$ \leq \begin{cases}
\max\left\{{\mid \lambda({\underline{\alpha}}_1 , \dots , {\underline{\alpha}}_m) \mid}_v , {\underline{\alpha}}_1 , \dots , {\underline{\alpha}}_m\right\}.\prod_{i = 1}^{m} \max_{\mid {\underline{\alpha}}_i \mid = d_i}{\mid {{\underline x}_i}^{{\underline{\alpha}}_i} \mid}_v & \!\!\text{si $v \nmid \infty$ } \\
\max\left\{{\mid \lambda({\underline{\alpha}}_1 , \dots ,
{\underline{\alpha}}_m) \mid}_v , {\underline{\alpha}}_1 , \dots ,
{\underline{\alpha}}_m\right\}.\prod_{i = 1}^{m} \sum_{\mid
{\underline{\alpha}}_i \mid = d_i}{\mid {{\underline
x}_i}^{{\underline{\alpha}}_i} \mid}_v & \!\!\text{si $v \mid
\infty$}
\end{cases}\!\!.
$$ En utilisant le lemme \ref{c.22} pour majorer $\max\{{\mid \lambda({\underline{\alpha}}_1 , \dots , {\underline{\alpha}}_m) \mid}_v , {\underline{\alpha}}_1 , \dots , {\underline{\alpha}}_m\}$ et le lemme \ref{c.35} pour majorer les deux quantités $\prod_{i = 1}^{m}\max_{\mid {\underline{\alpha}}_i \mid = d_i}{\mid {{\underline x}_i}^{{\underline{\alpha}}_i} \mid}_v $ (lorsque $v$ est finie) et $\prod_{i = 1}^{m} \sum_{\mid {\underline{\alpha}}_i \mid = d_i}{\mid {{\underline x}_i}^{{\underline{\alpha}}_i} \mid}_v$ (lorsque $v$ est infinie) \footnote{Le lemme \ref{c.35} de l'appendice majore les deux quantités $\prod_{i = 1}^{m}\max_{\mid {\underline{\alpha}}_i \mid = d_i}{\mid {{\underline x}_i}^{{\underline{\alpha}}_i} \mid}_v $ (lorsque $v$ est finie) et $\prod_{i = 1}^{m} \sum_{\mid {\underline{\alpha}}_i \mid = d_i}{\mid {{\underline x}_i}^{{\underline{\alpha}}_i} \mid}_v$ (lorsque $v$ est infinie) par $1$ et $e^m \leq e^{o(\delta)}$ respectivement.} on aura pour toute place $v \in S$:
\begin{equation}
\begin{split}
{\mid \partial^{(i_1 , \dots , i_m)}({\tau^*}_{\!\!\!\!\!\!- \mathbf y}P)({\underline x}_1 , \dots , {\underline x}_m) \mid}_v  & \leq\exp\left\{\phantom{\left(h_v({\underline y}_1) + \dots + h_v({\underline y}_{m - 1})\right)} \!\!\!\!\!\!\!\!\!\!\!\!\!\!\!\!\!\!\!\!\!\!\!\!\!\!\!\!\!\!\!\!\!\!\!\!\!\!\!\!\!\!\!\!\!\!\!\!\!\!\!\!\!\!\!\!\!\!\!\!\!\!\!\!\!\!2 m_v (i_1 + \dots + i_m) + 9 m_v \delta {a^2}_{\!\!\!1}\right. \\
& \quad + h_v(P) + 2 \epsilon_0 \delta \left({a^2}_{\!\!\!1} h_v({\underline x}_1) + \dots + {a^2}_{\!\!\!m} h_v({\underline x}_m)\right) \\
& \quad \left.+ 2 \delta \left(h_v({\underline y}_1) + \dots +
h_v({\underline y}_{m - 1})\right)\right\} \label{3.30}
\end{split}
\end{equation}
si $v$ est finie et:
\begin{equation}
\begin{split}
{\mid \partial^{(i_1 , \dots , i_m)}({\tau^*}_{\!\!\!\!\!\!- \mathbf y}P)({\underline x}_1 , \dots , {\underline x}_m) \mid}_v  & \leq\exp\left\{\phantom{\left(h_v({\underline y}_1) + \dots + h_v({\underline y}_{m - 1})\right)} \!\!\!\!\!\!\!\!\!\!\!\!\!\!\!\!\!\!\!\!\!\!\!\!\!\!\!\!\!\!\!\!\!\!\!\!\!\!\!\!\!\!\!\!\!\!\!\!\!\!\!\!\!\!\!\!\!\!\!\!\!\!\!\!\!\!(2 m_v + 12) (i_1 + \dots + i_m) \right. \\
& \!\!\!\!\!\!\!\!\!\!\!\!\!\!\!\!\!\!\!\!\!\!\!\!\!\!\!\!\!\!\!\!\!\!\!\!\!\!\!\! + (9 m_v + 28) \delta {a^2}_{\!\!\!1} + h_v(P) + 2 \epsilon_0 \delta \left({a^2}_{\!\!\!1} h_v({\underline x}_1) + \dots + {a^2}_{\!\!\!m} h_v({\underline x}_m)\right) \\
&
\!\!\!\!\!\!\!\!\!\!\!\!\!\!\!\!\!\!\!\!\!\!\!\!\!\!\!\!\!\!\!\!\!\!\!\!\!\!\!\!\left.+
2 \delta \left(h_v({\underline y}_1) + \dots + h_v({\underline
y}_{m - 1})\right) + \circ(\delta)\right\} \label{3.31}
\end{split}
\end{equation}
si $v$ est infinie. \\
D'autre part, on a aussi d'apr\`es le lemme \ref{c.35} de
l'appendice pour toute place $v \in S$:
\begin{equation}
\frac{1}{{\mid {\Delta({\underline x}_1)}^{f(i_1)} \dots
{\Delta({\underline x}_m)}^{f(i_m)} \mid}_v} ~\leq~ \begin{cases}
1 & \text{si $v \nmid \infty$} \\
{\sqrt{e}}^{f(i_1) + \dots + f(i_m)} \leq e^{i_1 + \dots + i_m} &
\text{si $v \mid \infty$}
\end{cases}. \label{3.32}
\end{equation}
La proposition \ref{c.23} suit finalement des trois majorations
(\ref{3.30}), (\ref{3.31}) et (\ref{3.32}) puisqu'on a pour tout
$(i_1 , \dots , i_m) \in {\mathbb N}^m$ et toute place $v \in S$:
$${\mid f_{i_1 , \dots , i_m} \mid}_v ~=~ {\mid \partial^{(i_1 , \dots , i_m)}({\tau^*}_{\!\!\!\!\!\!- \mathbf y}P)({\underline x}_1 , \dots , {\underline x}_m) \mid}_v.\frac{1}{{\mid {\Delta({\underline x}_1)}^{f(i_1)} \dots {\Delta({\underline x}_m)}^{f(i_m)} \mid}_v} .$$
La démonstration est achevée.  $~~~~\blacksquare$\vspace{1mm}
\begin{corollary}\label{c.36}
Pour toute place $v \in S$, la série de la proposition \ref{c.23}
est absolument convergente en valeur absolue $v$-adique d\`es que:
$${\mid u_i \mid}_v ~<~ \begin{cases}
e^{-2 m_v} & \text{si $v$ est finie} \\
e^{-2 m_v - 13} & \text{si $v$ est infinie}
\end{cases} ~~~~ i = 1 , \dots , m .$$
En particulier, pour toute place $v \in S$, la série sus-citée
converge absolument en valeur absolue $v$-adique au voisinage du
point $(- x_1 , \dots , - x_m)$.
\end{corollary}
{\bf Démonstration.---} Etant donné une place $v \in S$, la
condition suffisante pour la convergence absolue $v$-adique de la
série de la proposition \ref{c.23} est claire et entraine
-d'apr\`es l'hypoth\`ese principale ${\mbox{H}}_5$ et la propriété
$ii)$ du §$2$ pour la distance ${\mbox{dist}}_v$- la convergence
absolue $v$-adique de cette m\^eme série au voisinage du point $(-
x_1 , \dots , - x_m)$.  $~~~~\blacksquare$\vspace{1mm}
\begin{lemma}\label{c.25}
On a:
$$\Omega_{\underline a}({\tau^*}_{\!\!\!\!\!\!- \mathbf y}P)(\underline 0 , \dots , \underline 0) ~=~ \sum_{i_1 \geq 0 , \dots , i_m \geq 0} g_{i_1 , \dots , i_m} {u^{i_1}}_{\!\!\!\!1} \dots {u^{i_m}}_{\!\!\!\!\!m}$$
o\`u les $g_{i_1 , \dots , i_m} := \partial^{(i_1 , \dots , i_m)}
({\tau^*}_{\!\!\!\!\!\!- \mathbf y}P)(\underline 0 , \dots
\underline 0)$, $(i_1 , \dots , i_m) \in {\mathbb N}^m$ sont des
nombres algébriques de $K$ satisfaisant, pour toute place $v \in
S$, l'identité au sens $v$-adique: $\forall \underline{\ell} =
(\ell_1 , \dots , \ell_m) \in {\mathbb N}^m$:
$$g_{\underline \ell} ~=~ \!\!\!\sum_{
\begin{array}{c}
\scriptstyle{i_1 \geq \ell_1 , \dots , i_m \geq \ell_m} \\
\scriptstyle{\underline i \not\in T_{\delta}}
\end{array}} \!\!\!\!\!\!f_{i_1 , \dots , i_m} \binom{i_1}{\ell_1} \dots \binom{i_m}{\ell_m} (- x_1)^{i_1 - \ell_1} \dots (- x_m)^{i_m - \ell_m}$$
dont la série du membre droit est $v$-adiquement absolument
convergente.
\end{lemma}
{\bf Démonstration.---} Le fait que:
$$\Omega_{\underline a}({\tau^*}_{\!\!\!\!\!\!- \mathbf y}P)(\underline 0 , \dots , \underline 0) ~=~ \sum_{i_1 \geq 0 , \dots , i_m \geq 0} g_{i_1 , \dots , i_m} {u^{i_1}}_{\!\!\!\!1} \dots {u^{i_m}}_{\!\!\!\!\!m}$$
avec les $g_{i_1 , \dots , i_m} := \partial^{(i_1 , \dots , i_m)}
({\tau^*}_{\!\!\!\!\!\!- \mathbf y}P)(\underline 0 , \dots
\underline 0)$, (pour $(i_1 , \dots , i_m) \in {\mathbb N}^m$
vient immédiatement de l'application du corollaire \ref{c.11} pour
$P_1 = {\tau^*}_{\!\!\!\!\!\!- \mathbf y}P$, en spécialisant
$({\underline X}_1 , \dots , {\underline X}_m)$ \`a (\underline 0
, \dots , \underline 0). Par ailleurs soit $v$ une place fixée de
$S$ et posons $\mathcal S$ la série de la proposition \ref{c.23}:
$$\mathcal S ~:=~ \Omega_{\underline a}({\tau^*}_{\!\!\!\!\!\!- \mathbf y}P)({\underline x}_1 , \dots , {\underline x}_m) ~=~ \sum_{i_1 \geq 0 , \dots , i_m \geq 0} f_{i_1 , \dots , i_m} {u^{i_1}}_{\!\!\!\!1} \dots {u^{i_m}}_{\!\!\!\!\!m} .$$
Comme d'apr\`es le corollaire \ref{c.36}, $\mathcal S$ est
$v$-adiquement convergente au voisinage de $(- x_1 , \dots , -
x_m)$, alors quand le point $(u_1 , \dots , u_m)$ de ${\mathbb
C}_{v}^{m}$ est suffisamment proche de $(0 , \dots , 0)$ on a au
sens $v$-adique:
\begin{equation*}
\begin{split}
\Omega_{\underline a}({\tau^*}_{\!\!\!\!\!\!- \mathbf y}P)(\underline 0 , \dots , \underline 0)(u_1 , \dots , u_m) \\
& \!\!\!\!\!\!\!\!\!\!\!\!\!\!\!\!\!\!\!\!\!\!\!\!\!\!\!\!\!\!\!\!\!\!\!\!\!\!\!\!\!\!\!\!\!\!\!\!\!\!=~ \Omega_{\underline a}({\tau^*}_{\!\!\!\!\!\!- \mathbf y}P)({\underline x}_1 , \dots , {\underline x}_m)(- x_1 + u_1 , \dots , - x_m + u_m) \\
& \!\!\!\!\!\!\!\!\!\!\!\!\!\!\!\!\!\!\!\!\!\!\!\!\!\!\!\!\!\!\!\!\!\!\!\!\!\!\!\!\!\!\!\!\!\!\!\!\!\!=~ \mathcal S (- x_1 + u_1 , \dots , - x_m + u_m) \\
&
\!\!\!\!\!\!\!\!\!\!\!\!\!\!\!\!\!\!\!\!\!\!\!\!\!\!\!\!\!\!\!\!\!\!\!\!\!\!\!\!\!\!\!\!\!\!\!\!\!\!=\!\!\!\sum_{i_1
\geq 0 , \dots , i_m \geq 0} \frac{1}{i_1! \dots i_m!}
\frac{\partial^{(i_1 + \dots + i_m)} \mathcal S}{\partial
{u_1}^{\!i_1} \dots \partial {u_m}^{\!i_m}}(- x_1 , \dots , -
x_m).{u^{i_1}}_{\!\!\!\!1} \dots {u^{i_m}}_{\!\!\!\!\!m} .
\end{split}
\end{equation*}
D'o\`u, par identification, pour tout $(i_1 , \dots , i_m) \in
{\mathbb N}^m$:
\begin{equation}
g_{i_1 , \dots , i_m} ~=~ \frac{1}{i_1! \dots i_m!}
\frac{\partial^{(i_1 + \dots + i_m)} \mathcal S}{\partial
{u_1}^{\!i_1} \dots \partial {u_m}^{\!i_m}}(- x_1 , \dots , - x_m)
. \label{3.33}
\end{equation}
Du fait que la série $\mathcal S$ converge absolument -au sens
$v$-adique- au voisinage du point $(- x_1 , \dots , - x_m)$, ses
dérivées partielles de tout ordre le seront aussi et on peut
intervertir, en dérivant $\mathcal S$, la sommation et la
dérivation. On obtient ainsi : pour tout $\underline{\ell} =
(\ell_1 , \dots , \ell_m) \in {\mathbb N}^m$:
\begin{equation*}
\begin{split}
g_{\underline{\ell}} & =~ \frac{1}{\ell_1! \dots \ell_m!} \frac{\partial^{\mid \underline{\ell} \mid} \mathcal S}{\partial {u_1}^{\!\ell_1} \dots \partial {u_m}^{\!\ell_m}}(- x_1 , \dots , - x_m) \\
& =~\!\!\!\sum_{i_1 \geq \ell_1 , \dots , i_m \geq \ell_m}
\!\!\!\!\!\!f_{i_1 , \dots , i_m} \binom{i_1}{\ell_1} \dots \binom{i_m}{\ell_m} (- x_1)^{i_1 - \ell_1} \dots (- x_m)^{i_m - \ell_m} \\
& =~\!\!\!\sum_{
\begin{array}{c}
\scriptstyle{i_1 \geq \ell_1 , \dots , i_m \geq \ell_m} \\
\scriptstyle{\underline i \not\in T_{\delta}}
\end{array}} \!\!\!\!\!\!f_{i_1 , \dots , i_m} \binom{i_1}{\ell_1} \dots \binom{i_m}{\ell_m} (- x_1)^{i_1 - \ell_1} \dots (- x_m)^{i_m - \ell_m}
\end{split}
\end{equation*}
o\`u cette derni\`ere égalité tient du fait que les $f_{i_1 ,
\dots , i_m}$ sont nuls pour $(i_1 , \dots , i_m) \in T_{\delta}$.
La démonstration est achevée.  $~~~~\blacksquare$\vspace{1mm}
\begin{proposition}\label{c.26}
On a: $g_{\underline{\ell}} = 0 ,~~ \forall \underline{\ell} \in
T_{\delta / 2}$.
\end{proposition}
{\bf Démonstration.---} Procédons par l'absurde, c'est-\`a-dire
supposons que pour un certain $\underline{\ell} \in T_{\delta /
2}$ on ait $g_{\underline{\ell}} \neq 0$. Ainsi, comme
$g_{\underline{\ell}} \in K$, $g_{\underline{\ell}}$ doit
satisfaire la formule du produit:
\begin{equation}
\sum_{v \in M_K} \frac{[K_v : {\mathbb Q}_v]}{[K : \mathbb Q]}
\log {\mid g_{\ell} \mid}_v ~=~ 0 . \label{3.34}
\end{equation}
Nous allons montrer que (\ref{3.34}) ne peut pas avoir lieu et ceci en majorant son membre gauche par une quantité strictement négative. Majorons, pour toute place $v$ sur $K$, le nombre $\log {\mid g_{\ell} \mid}_v$. Pour ce faire on distingue les quatres cas suivants: \\
$\underline{{\mbox{1}}^{\mbox{er}} \mbox{cas:}}$ (si $v \in S$ et $v$ est finie) \\
Dans ce premier cas on utilise l'identité du lemme \ref{c.25} pour
majorer $\log {\mid g_{\ell} \mid}_v$. Cette derni\`ere montre
qu'on a:
$${\mid g_{\underline{\ell}} \mid}_v ~\leq~ \max_{\begin{array}{c}
\scriptstyle{\underline i \geq \underline{\ell}} \\
\scriptstyle{\underline i \not\in T_{\delta}}
\end{array}} \left\{{\mid f_{\underline i} \mid}_v.{{\mid x_1 \mid}_v}^{\! i_1 - \ell_1} \dots {{\mid x_m \mid}_v}^{\! i_m - \ell_m}\right\}$$
o\`u $\underline i \geq \underline{\ell}$ veut dire $i_n \geq
\ell_n ~ \forall n \in \{1 , \dots , m\}$. Puis en utilisant la
majoration de la proposition \ref{c.23} pour les ${\mid
f_{\underline i} \mid}_v ~ (\underline i \in {\mathbb N}^m)$ on
obtient:
\begin{equation*}
\begin{split}
{\mid g_{\underline{\ell}} \mid}_v & \leq~ \exp\!\left\{\phantom{\left(h_v({\underline y}_1) + \dots + h_v({\underline y}_{m - 1})\right)} \!\!\!\!\!\!\!\!\!\!\!\!\!\!\!\!\!\!\!\!\!\!\!\!\!\!\!\!\!\!\!\!\!\!\!\!\!\!\!\!\!\!\!\!\!\!\!\!\!\!\!\!\!\!\!\!\!\!\!\!\!\!\!\!\!\!9 m_v \delta {a^2}_{\!\!\!1} + h_v(P) + 2 \epsilon_0 \delta \left({a^2}_{\!\!\!1} h_v({\underline x}_1) + \dots + {a^2}_{\!\!\!m} h_v({\underline x}_m)\right)\right. \\
& \quad \left.+ 2 \delta \!\left(h_v({\underline y}_1) + \dots + h_v({\underline y}_{m - 1})\right)\!\right\} \\
& \quad \times \max_{\begin{array}{c}
\scriptstyle{\underline i \geq \underline{\ell}} \\
\scriptstyle{\underline i \not\in T_{\delta}}
\end{array}} \!\!\!\left\{(e^{2 m_v})^{i_1 + \dots + i_m} {{\mid x_1 \mid}_v}^{\! i_1 - \ell_1} \dots {{\mid x_m \mid}_v}^{\! i_m - \ell_m}\right\} .
\end{split}
\end{equation*}
En faisant maintenant -dans le maximum du membre de droite de
cette derni\`ere inégalité- le changement d'indice $\underline j =
\underline i - \underline{\ell}$ et en remarquant que pour tout
$\underline i \geq \underline{\ell}$: $\underline i \not\in
T_{\delta} \Rightarrow \underline j \not\in T_{\delta / 2}~$ (car
$\underline{\ell} \in T_{\delta / 2}$), on a -a fortiori-
l'inégalité:
\begin{equation*}
\begin{split}
{\mid g_{\underline{\ell}} \mid}_v & \leq~ \exp\!\left\{\phantom{\left(h_v({\underline y}_1) + \dots + h_v({\underline y}_{m - 1})\right)} \!\!\!\!\!\!\!\!\!\!\!\!\!\!\!\!\!\!\!\!\!\!\!\!\!\!\!\!\!\!\!\!\!\!\!\!\!\!\!\!\!\!\!\!\!\!\!\!\!\!\!\!\!\!\!\!\!\!\!\!\!\!\!\!\!\!9 m_v \delta {a^2}_{\!\!\!1} + h_v(P) + 2 \epsilon_0 \delta \left({a^2}_{\!\!\!1} h_v({\underline x}_1) + \dots + {a^2}_{\!\!\!m} h_v({\underline x}_m)\right)\right. \\
& \quad \left.+ 2 \delta \left(h_v({\underline y}_1) + \dots + h_v({\underline y}_{m - 1})\right)\!\right\}.(e^{2 m_v})^{\mid \underline{\ell} \mid} \\
& \quad \times \max_{\underline j \not\in T_{\delta / 2}}
\left\{(e^{2 m_v})^{\mid \underline j \mid} \prod_{k = 1}^{m}
{{\mid x_k \mid}_v}^{\!j_k}\right\} .
\end{split}
\end{equation*}
Or, d'apr\`es notre hypoth\`ese principale ${\mbox{H}}_5$ et la
propriété $ii)$ du §$2$ pour la distance ${\mbox{dist}}_v$, on a
bien pour tout $k \in \{1 , \dots , m\}$:
$${\mid x_k \mid}_v ~<~ e^{- \lambda_v \epsilon h({\mathbf x}_k) - 2 m_v} .$$
On obtient, gr\^ace \`a ces derni\`eres inégalités:
\begin{equation*}
\begin{split}
{\mid g_{\underline{\ell}} \mid}_v & \leq~ \exp\!\left\{\phantom{\left(h_v({\underline y}_1) + \dots + h_v({\underline y}_{m - 1})\right)} \!\!\!\!\!\!\!\!\!\!\!\!\!\!\!\!\!\!\!\!\!\!\!\!\!\!\!\!\!\!\!\!\!\!\!\!\!\!\!\!\!\!\!\!\!\!\!\!\!\!\!\!\!\!\!\!\!\!\!\!\!\!\!\!\!\!9 m_v \delta {a^2}_{\!\!\!1} + h_v(P) + 2 \epsilon_0 \delta \left({a^2}_{\!\!\!1} h_v({\underline x}_1) + \dots + {a^2}_{\!\!\!m} h_v({\underline x}_m)\right)\right. \\
& \quad \left.+ 2 \delta \left(h_v({\underline y}_1) + \dots + h_v({\underline y}_{m - 1})\right)\!\right\}.(e^{2 m_v})^{\mid \underline{\ell} \mid} \\
& \quad \times \max_{\underline j \not\in T_{\delta / 2}} e^{-
\epsilon \lambda_v \left\{\sum_{k = 1}^{m} j_k h({\mathbf
x}_k)\right\}} .
\end{split}
\end{equation*}
D'o\`u, en passant aux logarithmes:
\begin{equation*}
\begin{split}
\log {\mid g_{\underline{\ell}} \mid}_v & \leq~ 9 m_v \delta {a^2}_{\!\!\!1} + h_v(P) + 2 \epsilon_0 \delta \left({a^2}_{\!\!\!1} h_v({\underline x}_1) + \dots + {a^2}_{\!\!\!m} h_v({\underline x}_m)\right) \\
& \quad + 2 \delta \left(h_v({\underline y}_1) + \dots + h_v({\underline y}_{m - 1})\right) + 2 m_v \mid \underline{\ell} \mid \\
& \quad - \epsilon \lambda_v \min_{\underline j \not\in T_{\delta
/ 2}} \left\{\sum_{k = 1}^{m} j_k h({\mathbf x}_k)\right\} .
\end{split}
\end{equation*}
Remarquons maintenant que pour tout $\underline j \not\in
T_{\delta / 2}$ (c'est-\`a-dire: $\frac{j_1}{{a^2}_{\!\!\!1}} +
\dots + \frac{j_{m - 1}}{{a^2}_{\!\!\!m -1}} + \frac{j_m}{m - 1} >
\frac{7}{2} \epsilon_1 \delta$) on a:
\begin{equation}
\begin{split}
& j_1 h({\mathbf x}_1) + \dots + j_m h({\mathbf x}_m) \\
& \geq~ \frac{j_1}{{a^2}_{\!\!\!1}} \left({a^2}_{\!\!\!1} h({\mathbf x}_1)\right) + \dots + \frac{j_{m - 1}}{{a^2}_{\!\!\!m -1}} \left({a^2}_{\!\!\!m -1} h({\mathbf x}_{m - 1})\right) + \frac{j_m}{m - 1} \left({a^2}_{\!\!\!m} h({\mathbf x}_m)\right) \\
& \geq~ \frac{j_1}{{a^2}_{\!\!\!1}} \left(\frac{1}{2} {a^2}_{\!\!\!1} h({\mathbf x}_1)\right) + \dots + \frac{j_{m - 1}}{{a^2}_{\!\!\!m -1}} \left( \frac{1}{2} {a^2}_{\!\!\!1} h({\mathbf x}_1)\right) + \frac{j_m}{m - 1} \left(\frac{1}{2}{a^2}_{\!\!\!1} h({\mathbf x}_1)\right) \\
& \geq~ \frac{1}{2} {a^2}_{\!\!\!1} h({\mathbf x}_1) \left(\frac{j_1}{{a^2}_{\!\!\!1}} + \dots + \frac{j_{m - 1}}{{a^2}_{\!\!\!m -1}} + \frac{j_m}{m - 1}\right) \\
& \geq~ \frac{7}{4} \epsilon_1 \delta {a^2}_{\!\!\!1} h({\mathbf
x}_1) . \label{3.35}
\end{split}
\end{equation}
En utilisant (\ref{3.35}) et en majorant $\mid \underline{\ell}
\mid$ par:
$$\mid \underline{\ell} \mid ~\leq~ {a^2}_{\!\!\!1} \left(\frac{\ell_1}{{a^2}_{\!\!\!1}} + \dots + \frac{\ell_{m - 1}}{{a^2}_{\!\!\!m - 1}} + \frac{\ell_m}{m - 1}\right) ~\leq~ {a^2}_{\!\!\!1}.\frac{7}{2} \epsilon_1 \delta ~=~ \frac{7}{2} \epsilon_1 \delta {a^2}_{\!\!\!1}$$
(car $\underline{\ell} \in T_{\delta / 2}$), on aura finalement:
 \begin{equation*}
\begin{split}
\mathbf{(a)}~~~~
\log {\mid g_{\underline{\ell}} \mid}_v & \leq~ \left[- \frac{7}{4} \lambda_v \epsilon \epsilon_1 h({\mathbf x}_1) + (7 \epsilon_1 + 9) m_v\right] \!\!\delta {a^2}_{\!\!\!1} + h_v(P) \\
& \!\!\!\!\!\!\!\!\!\!\!\!\!\!\!\!\!\!\! + 2 \epsilon_0 \delta
\!\left({a^2}_{\!\!\!1} h_v({\underline x}_1) + \dots +
{a^2}_{\!\!\!m} h_v({\underline x}_m)\right) + 2 \delta
\!\left(h_v({\underline y}_1) + \dots + h_v({\underline y}_{m -
1})\right) .
\end{split}
\end{equation*}
$\underline{{\mbox{2}}^{\mbox{\`eme}} \mbox{cas:}}$ (si $v \in S$ et $v$ est infinie) \\
On suit exactement la m\^eme méthode que celle du premier cas,
sauf qu'ici le fait que $v$ est infinie apporte un petit
changement sur la majoration de $\log {\mid g_{\underline{\ell}}
\mid}_v$ par rapport au premier cas. L'identité du lemme
\ref{c.25} donne la majoration:
$${\mid g_{\underline{\ell}} \mid}_v ~\leq \!\!\sum_{\begin{array}{c}
\scriptstyle{\underline i \geq \underline{\ell}} \\
\scriptstyle{\underline i \not\in T_{\delta}}
\end{array}} \!\!\!{\mid f_{\underline i} \mid}_v .  2^{i_1 + \dots + i_m} . {{\mid x_1 \mid}_v}^{\! i_1 - \ell_1} \dots {{\mid x_m \mid}_v}^{\! i_m - \ell_m} .$$
En utilisant ensuite la majoration de la proposition \ref{c.23}
pour les ${\mid f_{\underline i} \mid}_v ~ (\underline i \in
{\mathbb N}^m)$ on obtient:
\begin{equation*}
\begin{split}
{\mid g_{\underline{\ell}} \mid}_v & \leq~ \exp\!\left\{\phantom{\left(h_v({\underline y}_1) + \dots + h_v({\underline y}_{m - 1})\right)} \!\!\!\!\!\!\!\!\!\!\!\!\!\!\!\!\!\!\!\!\!\!\!\!\!\!\!\!\!\!\!\!\!\!\!\!\!\!\!\!\!\!\!\!\!\!\!\!\!\!\!\!\!\!\!\!\!\!\!\!\!\!\!\!\!\!(9 m_v + 28) \delta {a^2}_{\!\!\!1} + h_v(P) + 2 \epsilon_0 \delta \left({a^2}_{\!\!\!1} h_v({\underline x}_1) + \dots + {a^2}_{\!\!\!m} h_v({\underline x}_m)\right)\right. \\
& \quad \left.+ 2 \delta \left(h_v({\underline y}_1) + \dots + h_v({\underline y}_{m - 1})\right) + o(\delta)\right\} \\
& \quad \times \!\!\!\sum_{\begin{array}{c}
\scriptstyle{\underline i \geq \underline{\ell}} \\
\scriptstyle{\underline i \not\in T_{\delta}}
\end{array}} \!\!\!\left(2 e^{2 m_v + 13}\right)^{i_1 + \dots + i_m} {{\mid x_1 \mid}_v}^{\! i_1 - \ell_1} \dots {{\mid x_m \mid}_v}^{\! i_m - \ell_m} .
\end{split}
\end{equation*}
En faisant le changement d'indice $\underline j = \underline i -
\underline{\ell}$ dans la somme intervenant dans le deuxi\`eme
membre de cette derni\`ere inégalité, gr\^ace \`a la remarque
faite au premier cas, on obtient -a fortiori- l'inégalité:
\begin{equation*}
\begin{split}
{\mid g_{\underline{\ell}} \mid}_v & \leq~ \exp\!\left\{\phantom{\left(h_v({\underline y}_1) + \dots + h_v({\underline y}_{m - 1})\right)} \!\!\!\!\!\!\!\!\!\!\!\!\!\!\!\!\!\!\!\!\!\!\!\!\!\!\!\!\!\!\!\!\!\!\!\!\!\!\!\!\!\!\!\!\!\!\!\!\!\!\!\!\!\!\!\!\!\!\!\!\!\!\!\!\!\!(9 m_v + 28) \delta {a^2}_{\!\!\!1} + h_v(P) + 2 \epsilon_0 \delta \left({a^2}_{\!\!\!1} h_v({\underline x}_1) + \dots + {a^2}_{\!\!\!m} h_v({\underline x}_m)\right)\right. \\
& \quad \left.+ 2 \delta \left(h_v({\underline y}_1) + \dots + h_v({\underline y}_{m - 1})\right) + o(\delta)\right\}.\left(2 e^{2 m_v + 13}\right)^{\mid \underline{\ell} \mid} \\
& \quad \times \!\!\!\sum_{\underline j \not\in T_{\delta / 2}}
\!\!\!\left\{\left(2 e^{2 m_v + 13}\right)^{\mid \underline j
\mid} \prod_{k = 1}^{m} {{\mid x_k \mid}_v}^{\!j_k}\right\} .
\end{split}
\end{equation*}
Or, d'apr\`es notre hypoth\`ese principale ${\bf{H}}_5$ et la
propriété $ii)$ du §$2$ pour la distance ${\mbox{dist}}_v$, on a
bien pour tout $k \in \{1 , \dots , m\}$:
$${\mid x_k \mid}_v ~<~ e^{- \lambda_v \epsilon h({\mathbf x}_k) - 2 m_v - 16} .$$
Gr\^ace \`a ces derni\`eres inégalités on obtient:
\begin{equation*}
\begin{split}
{\mid g_{\underline{\ell}} \mid}_v & \leq~ \exp\!\left\{(9 m_v + 28) \delta {a^2}_{\!\!\!1} + h_v(P) + 2 \epsilon_0 \delta \left({a^2}_{\!\!\!1} h_v({\underline x}_1) + \dots + {a^2}_{\!\!\!m} h_v({\underline x}_m)\right)\right. \\
& \quad \left.+ 2 \delta \left(h_v({\underline y}_1) + \dots + h_v({\underline y}_{m - 1})\right) + o(\delta)\right\}.\left(2 e^{2 m_v + 13}\right)^{\mid \underline{\ell} \mid} \\
& \quad \times \!\!\!\sum_{\underline j \not\in T_{\delta / 2}}
\!\!\! \left(\frac{2}{e^3}\right)^{\mid \underline j \mid} \!e^{-
\epsilon \lambda_v \left\{\sum_{k = 1}^{m} j_k h({\mathbf
x}_k)\right\}} .
\end{split}
\end{equation*}
D'une part, d'apr\`es (\ref{3.35}): pour tout $\underline j
\not\in T_{\delta / 2}$ on a: $\sum_{k = 1}^{m} j_k h({\mathbf
x}_k) \geq \frac{7}{4} \epsilon_1 \delta {a^2}_{\!\!\!1}
h({\mathbf x}_1)$ et d'autre part, comme:
$$\underline j \not\in T_{\delta / 2} \Longrightarrow \mid \underline j \mid = j_1 + \dots + j_m \geq \frac{j_1}{{a^2}_{\!\!\!1}} + \dots + \frac{j_{m - 1}}{{a^2}_{\!\!\!m - 1}} + \frac{j_m}{m - 1} > \frac{7}{2} \epsilon_1 \delta ,$$
on a:
\begin{eqnarray*}
\sum_{\underline j \not\in T_{\delta / 2}} \!\!\! \left(\frac{2}{e^3}\right)^{\mid \underline j \mid} & \leq & \sum_{\underline j / \mid \underline j \mid \geq \frac{7}{2} \epsilon_1 \delta} \!\!\! \left(\frac{2}{e^3}\right)^{\mid \underline j \mid} \\
& \leq & \sum_{j \geq \frac{7}{2} \epsilon_1 \delta} \!\!\! \binom{j + m - 1}{m - 1} \left(\frac{2}{e^3}\right)^j ~~~~ \text{(en posant $j = \mid \underline j \mid $)} \\
& \leq & \sum_{j \geq \frac{7}{2} \epsilon_1 \delta} \!\!\! 2^{j + m - 1}\left(\frac{2}{e^3}\right)^j \\
& \leq & 2^{m - 1} \!\!\!\sum_{j \geq \frac{7}{2} \epsilon_1 \delta} \!\!\!\left(\frac{4}{e^3}\right)^j \\
& \leq & 2^{m - 1} \left(\frac{4}{e^3}\right)^{\frac{7}{2}
\epsilon_1 \delta} \frac{1}{1 - \frac{4}{e^3}}
\end{eqnarray*}
puis:
\begin{eqnarray*}
\sum_{\underline j \not\in T_{\delta / 2}} \!\!\! \left(\frac{2}{e^3}\right)^{\mid \underline j \mid} & \leq & 2^m \left(\frac{4}{e^3}\right)^{\frac{7}{2} \epsilon_1 \delta} ~~~~~~~~~~~~~~~~~ \text{(car $\frac{1}{1 - \frac{4}{e^3}} < 2$)} \\
& \leq & e^{o(\delta)} .
\end{eqnarray*}
La derni\`ere estimation pour ${\mid g_{\underline{\ell}} \mid}_v$
entra{\sf\^\i}ne alors qu'on a:
\begin{equation*}
\begin{split}
{\mid g_{\underline{\ell}} \mid}_v & \leq~ \exp\!\left\{\phantom{\left(h_v({\underline y}_1) + \dots + h_v({\underline y}_{m - 1})\right)} \!\!\!\!\!\!\!\!\!\!\!\!\!\!\!\!\!\!\!\!\!\!\!\!\!\!\!\!\!\!\!\!\!\!\!\!\!\!\!\!\!\!\!\!\!\!\!\!\!\!\!\!\!\!\!\!\!\!\!\!\!\!\!\!\!\!(9 m_v + 28) \delta {a^2}_{\!\!\!1} + h_v(P) + 2 \epsilon_0 \delta \left({a^2}_{\!\!\!1} h_v({\underline x}_1) + \dots + {a^2}_{\!\!\!m} h_v({\underline x}_m)\right)\right. \\
& \!\!\!\!\!\!\!\!\!\!\left.+ 2 \delta \left(h_v({\underline y}_1)
+ \dots + h_v({\underline y}_{m - 1})\right) +
o(\delta)\right\}.\left(2 e^{2 m_v + 13}\right)^{\mid
\underline{\ell} \mid} \!.~\! e^{- \frac{7}{4} \lambda_v \epsilon
\epsilon_1 \delta {a^2}_{\!\!\!1} h({\mathbf x}_1)} .
\end{split}
\end{equation*}
D'o\`u, en passant aux logarithmes et en majorant $\mid
\underline{\ell} \mid$ par $\frac{7}{2} \epsilon_1 \delta
{a^2}_{\!\!\!1}$ (comme dans le premier cas) et $\frac{7}{2} (13 +
\log 2)$ par $48$, on aura finalement:
\begin{equation*}
\begin{split}
\mathbf{(b)}~~~~
\log {\mid g_{\underline{\ell}} \mid}_v & \leq~ \left[- \frac{7}{4} \lambda_v \epsilon \epsilon_1 h({\mathbf x}_1) + (7 \epsilon_1 + 9) m_v + 48 \epsilon_1 + 28 \right] \!\!\delta {a^2}_{\!\!\!1} + h_v(P) \\
&
\!\!\!\!\!\!\!\!\!\!\!\!\!\!\!\!\!\!\!\!\!\!\!\!\!\!\!\!\!\!\!\!\!\!\!
+ 2 \epsilon_0 \delta \!\left({a^2}_{\!\!\!1} h_v({\underline
x}_1) + \dots + {a^2}_{\!\!\!m} h_v({\underline x}_m)\right) + 2
\delta \!\left(h_v({\underline y}_1) + \dots + h_v({\underline
y}_{m - 1})\right) + o(\delta) .
\end{split}
\end{equation*}
$\underline{{\mbox{3}}^{\mbox{\`eme}} \mbox{cas:}}$ (si $v \not\in S$ et $v$ est finie) \\
Dans les deux cas restants, on majore na{\sf\"\i}vement le nombre
$\log {\mid g_{\underline{\ell}} \mid}_v$, c'est-\`a-dire on
utilise l'identité de définition:
$$g_{\underline{\ell}} ~:=~ \partial^{\underline{\ell}} \left({\tau^*}_{\!\!\!\!\!\!- \mathbf y}P\right)(\underline 0 , \dots , \underline 0)$$
qui montre que $g_{\underline{\ell}}$ est un coefficient de la
forme $\partial^{\underline{\ell}}({\tau^*}_{\!\!\!\!\!\!- \mathbf
y}P)$, par suite on a évidemment:
$$\log {\mid g_{\underline{\ell}} \mid}_v ~=~ h_v(g_{\underline{\ell}}) ~\leq~ h_v\left(\partial^{\underline{\ell}}({\tau^*}_{\!\!\!\!\!\!- \mathbf y}P)\right) .$$ En majorant maintenant $h_v\left(\partial^{\underline{\ell}}({\tau^*}_{\!\!\!\!\!\!- \mathbf y}P)\right)$ \`a l'aide du lemme \ref{c.22} on obtient (comme $v$ est finie):
\begin{equation*}
\begin{split}
\log {\mid g_{\underline{\ell}} \mid}_v & \leq~ 2 m_v \!\mid \underline{\ell} \mid + 9 m_v \delta {a^2}_{\!\!\!1} + h_v(P) + 2 \epsilon_0 \delta \left({a^2}_{\!\!\!1} h_v({\underline x}_1) + \dots + {a^2}_{\!\!\!m} h_v({\underline x}_m)\right) \\
& \quad + 2 \delta \left(h_v({\underline y}_1) + \dots +
h_v({\underline y}_{m - 1})\right) .
\end{split}
\end{equation*}
Et en majorant finalement $\mid \underline{\ell} \mid$ par
$\frac{7}{2} \epsilon_1 \delta {a^2}_{\!\!\!1}$ -comme dans le
premier et le second cas- il vient:
\begin{equation*}
\begin{split}
\mathbf{(c)}~~
\log {\mid g_{\underline{\ell}} \mid}_v & \leq~ (7 \epsilon_1 + 9) m_v \delta {a^2}_{\!\!\!1} + h_v(P) + 2 \epsilon_0 \delta \!\left({a^2}_{\!\!\!1} h_v({\underline x}_1) + \dots + {a^2}_{\!\!\!m} h_v({\underline x}_m)\right) \\
& \quad + 2 \delta \!\left(h_v({\underline y}_1) + \dots +
h_v({\underline y}_{m - 1})\right) .
\end{split}
\end{equation*}
$\underline{{\mbox{4}}^{\mbox{\`eme}} \mbox{cas:}}$ (si $v \not\in S$ et $v$ est infinie) \\
Comme dans le troisi\`eme cas, on a:
\begin{equation*}
\begin{split}
\log {\mid g_{\underline{\ell}} \mid}_v & \leq~ h_v\left(\partial^{\underline{\ell}}({\tau^*}_{\!\!\!\!\!\!- \mathbf y}P)\right) \\
& \!\!\!\!\!\!\!\!\!\!\!\!\!\!\!\!\!\!\!\!\!\!\!\!\!\leq~ (2 m_v + 12) \!\mid \underline{\ell} \mid + (9 m_v + 28) \delta {a^2}_{\!\!\!1} + h_v(P) \\
& \!\!\!\!\!\!\!\!\!\!\!\!\!\!\!\!\!\!\!\!\!\!\!\!\!\quad + 2
\epsilon_0 \delta \left({a^2}_{\!\!\!1} h_v({\underline x}_1) +
\dots + {a^2}_{\!\!\!m} h_v({\underline x}_m)\right) + 2 \delta
\left(h_v({\underline y}_1) + \dots + h_v({\underline y}_{m -
1})\right) + o(\delta) , \end{split}
\end{equation*}
o\`u la derni\`ere inégalité provient du lemme \ref{c.22}. \\
Il ne reste qu'\`a majorer $\mid \underline{\ell} \mid$ par
$\frac{7}{2} \epsilon_1 \delta {a^2}_{\!\!\!1}$ (comme dans les
cas précédents) pour obtenir:
\begin{equation*}
\begin{split}
\log {\mid g_{\underline{\ell}} \mid}_v & \leq~ \left[(7 \epsilon_1 + 9) m_v + (42 \epsilon_1 + 28)\right] \delta {a^2}_{\!\!\!1} + h_v(P) \\
& \!\!\!\!\!\!\!\!\!\!\!\!\!\!\!\! + 2 \epsilon_0 \delta
\!\left({a^2}_{\!\!\!1} h_v({\underline x}_1) + \dots +
{a^2}_{\!\!\!m} h_v({\underline x}_m)\right) + 2 \delta
\!\left(h_v({\underline y}_1) + \dots + h_v({\underline y}_{m -
1})\right) + o(\delta) .
\end{split}
\end{equation*}
Et a fortiori:
\begin{equation*}
\begin{split}
\mathbf{(d)}~~
\log {\mid g_{\underline{\ell}} \mid}_v & \leq~ \left[(7 \epsilon_1 + 9) m_v + (48 \epsilon_1 + 28)\right] \delta {a^2}_{\!\!\!1} + h_v(P) \\
& \!\!\!\!\!\!\!\!\!\!\!\!\!\!\!\!\!\!\!\!\!\!\!\!\!\!\!\!\!\! + 2
\epsilon_0 \delta \!\left({a^2}_{\!\!\!1} h_v({\underline x}_1) +
\dots + {a^2}_{\!\!\!m} h_v({\underline x}_m)\right) + 2 \delta
\!\left(h_v({\underline y}_1) + \dots + h_v({\underline y}_{m -
1})\right) + o(\delta) .
\end{split}
\end{equation*}
En majorant maintenant la somme $\sum_{v \in M_K} \frac{[K_v : {\mathbb Q}_v]}{[K : \mathbb Q]} \log {\mid g_{\underline{\ell}} \mid}_v$ à l'aide des inégalités $\mathbf{(a)} , \mathbf{(b)} , \mathbf{(c)}$ et $\mathbf{(d)}$ -selon le cas correspondant \`a chaque place $v$- et en remarquant les égalités: \\
$\displaystyle \sum_{v \in S} \frac{[K_v : {\mathbb Q}_v]}{[K : \mathbb Q]}\lambda_v ~=~ 1 ~~~~~~~~~~ \text{(qui est une hypoth\`ese)} ,$ \\
$\displaystyle \sum_{v \in S} \frac{[K_v : {\mathbb Q}_v]}{[K : \mathbb Q]} \eta_v ~=:~ \eta , ~~~~ \sum_{v \in M_K} \frac{[K_v : {\mathbb Q}_v]}{[K : \mathbb Q]} h_v(P) ~=:~ \widetilde{h}(P) ,$ \\
$\displaystyle \sum_{v \in M_K} \frac{[K_v : {\mathbb Q}_v]}{[K : \mathbb Q]} h_v({\underline x}_i) ~=:~ h({\mathbf x}_i) ~~~~~~~~~~ (\forall i \in \{1 , \dots , m\}) ,$ \\
$\displaystyle \sum_{v \in M_K} \frac{[K_v : {\mathbb Q}_v]}{[K : \mathbb Q]} h_v({\underline y}_i) ~=:~ h({\mathbf y}_i) ~~~~~~~~~~ (\forall i \in \{1 , \dots , m - 1\}) $ \\
$\displaystyle \text{et} ~~ \sum_{v \in M_{K}^{\infty}} \frac{[K_v : {\mathbb Q}_v]}{[K : \mathbb Q]} ~=~ 1 ,$ \\
on obtient:
\begin{equation*}
\begin{split}
\sum_{v \in M_K} \frac{[K_v : {\mathbb Q}_v]}{[K : \mathbb Q]}\log {\mid g_{\underline{\ell}} \mid}_v & \leq~ \left[- \frac{7}{4} \epsilon \epsilon_1 h({\mathbf x}_1) + (7 \epsilon_1 + 9) \eta + (48 \epsilon_1 + 28)\right] \!\!\delta {a^2}_{\!\!\!1} \\
&
\!\!\!\!\!\!\!\!\!\!\!\!\!\!\!\!\!\!\!\!\!\!\!\!\!\!\!\!\!\!\!\!\!\!\!\!\!\!\!\!\!\!\!\!\!\!\!\!\!\!\!\!\!\!\!\!\!\!\!\!
+ \widetilde{h}(P) + 2 \epsilon_0 \delta \!\left({a^2}_{\!\!\!1}
h({\mathbf x}_1) + \dots + {a^2}_{\!\!\!m} h({\mathbf x}_m)\right)
+ 2 \delta \!\left(h({\mathbf y}_1) + \dots + h({\mathbf y}_{m -
1})\right) + o(\delta) .
\end{split}
\end{equation*}
En majorant finalement dans le deuxi\`eme membre de cette
derni\`ere inégalité $\widetilde{h}(P)$ par son estimation donnée
\`a la proposition \ref{c.12} qui est:
$$\widetilde{h}(P) ~\leq~ [14(\eta + 6) \epsilon_1 + 4 \eta + 12] \delta {a^2}_{\!\!\!1} + o(\delta) ,$$
${a^2}_{\!\!\!1} h({\mathbf x}_1) + \dots + {a^2}_{\!\!\!m}
h({\mathbf x}_m)$ par $2 m {a^2}_{\!\!\!1} h({\mathbf x}_1)$
(d'apr\`es (\ref{3.23})) et $h({\mathbf y}_1) + \dots + h({\mathbf
y}_{m - 1})$ par $\alpha [4 m h({\mathbf x}_1) + \frac{5}{2} \eta
+ \frac{27}{2}] {a^2}_{\!\!\!1}$ (d'apr\`es (\ref{3.24})), on a:
\begin{equation*}
\begin{split}
\sum_{v \in M_K} \frac{[K_v : {\mathbb Q}_v]}{[K : \mathbb Q]}\log {\mid g_{\underline{\ell}} \mid}_v & \leq~ [(4 m (\epsilon_0 + 2 \alpha)- 7\!/\!4 . \epsilon \epsilon_1) h({\mathbf x}_1) + (21 \eta + 132) \epsilon_1 \\
& \quad + (5 \eta + 27) \alpha + (13 \eta + 40)] \delta
{a^2}_{\!\!\!1} + o(\delta) .
\end{split}
\end{equation*}
Or, d'apr\`es (\ref{3.34}) le premier membre de cette inégalité
est nul; donc on doit avoir:
$$
\left[\!\!\left(\!4 m (\epsilon_0 + 2 \alpha)- \frac{7}{4}
\epsilon \epsilon_1 \!\!\right)\!h({\mathbf x}_1) + (21 \eta +
132) \epsilon_1 + (5 \eta + 27) \alpha + (13 \eta +
40)\!\right]\!\!\delta {a^2}_{\!\!\!1} + o(\delta) \geq 0 .$$
\\ En divisant les deux membres de cette derni\`ere inégalité par $\delta {a^2}_{\!\!\!1}$ et en faisant tendre $\delta$ vers l'infini on aura:
$$
\left(\!4 m (\epsilon_0 + 2 \alpha)- \frac{7}{4} \epsilon
\epsilon_1 \!\!\right)\!h({\mathbf x}_1) + (21 \eta + 132)
\epsilon_1 + (5 \eta + 27) \alpha + (13 \eta + 40) ~\!\geq~\! 0
;$$ c'est-\`a-dire:
$$
\left(\!\frac{7}{4} \epsilon \epsilon_1 - 4 m (\epsilon_0 + 2
\alpha)\!\!\right)\!h({\mathbf x}_1) ~\!\leq~\! (21 \eta + 132)
\epsilon_1 + (5 \eta + 27) \alpha + (13 \eta + 40) .$$ Comme
maintenant -d'apr\`es l'hypoth\`ese ${\mbox{H}}_3$- on a:
$\frac{7}{4} \epsilon \epsilon_1 - 4 m (\epsilon_0 + 2 \alpha)
\geq \frac{3}{4} \epsilon \epsilon_1$ on déduit de cette
derni\`ere inégalité:
\begin{eqnarray*}
h({\mathbf x}_1) & \leq & \frac{1}{\epsilon \epsilon_1} \frac{4}{3} \left\{(21 \eta + 132) \epsilon_1 + (5 \eta + 27) \alpha + (13 \eta + 40)\right\} \\
& < & \frac{1}{\epsilon \epsilon_1} \left\{(28 \eta + 176)
\epsilon_1 + (7 \eta + 36) \alpha + (18 \eta + 54)\right\} .
\end{eqnarray*}
En utilisant finalement le théor\`eme \ref{c.14} du formulaire il
vient:
\begin{eqnarray*}
\widehat{h}({\mathbf x}_1) & \leq & h({\mathbf x}_1) + \frac{3}{4} \eta + 5 \\
& < & \frac{1}{\epsilon \epsilon_1} \left\{(28 \eta + 176)
\epsilon_1 + (7 \eta + 36) \alpha + (18 \eta + 54)\right\} +
\frac{3}{4} \eta + 5 .
\end{eqnarray*}
Ce qui contredit l'hypoth\`ese ${\mbox{H}}_4$ et ach\`eve cette
démonstration.  $~~~~\blacksquare$\vspace{2mm}

Notons respectivement par $- {\mathbf x}_1 , \dots , {- \mathbf x}_m$ les
opposés des points ${\mathbf x}_1 , \dots , {\mathbf x}_m$ de $E$
et par $- {\underline x}_1 := (- x_1 : 1 : - z_1) , \dots , -
{\underline x}_m := (- x_m : 1 : - z_m)$ leurs représentants dans
${\mathbb P}_2$.
\begin{corollary}\label{c.27}
Notre forme $Q$ de $K[{\underline X}_1 , \dots , {\underline
X}_m]$ définie au §$6$ par:
$$
Q({\underline X}_1 , \dots , {\underline X}_m) :=
P\!\left({\underline X}_1 , \dots , {\underline X}_m , \underline
D \!\left({\underline F}^{(a_1)}({\underline X}_1) , {\underline
X}_m\right) , \dots , \underline D \!\left({\underline F}^{(a_{m -
1})}({\underline X}_{m - 1}) , {\underline X}_m\right)\right)
$$
est non identiquement nulle modulo $I(E^m)$, de multidegré $((2 +
\epsilon_0) \delta {a^2}_{\!\!\!1} , \dots , (2 + \epsilon_0)
\delta {a^2}_{\!\!\!m - 1} ,$ \\ $(2 m - 2 + \epsilon_0) \delta)$,
satisfait $\tau^m(Q) (- {\underline x}_1 , \dots , - {\underline
x}_m) \in {\mathfrak q}_{\delta/2}$ et est de hauteur de
Gauss-Weil majorée par:
$$\widetilde{h}(Q) \leq \left[14 (\eta + 6) \epsilon_1 + 8 \eta + 22\right] \delta {a^2}_{\!\!\!1} + o(\delta) .$$
\end{corollary}
{\bf Démonstration.---} Le fait que $Q$ est une forme de
$K[{\underline X}_1 , \dots , {\underline X}_m]$ non identiquement
nulle modulo l'idéal $I(E^m)$, de multidegré $((2 + \epsilon_0)
\delta {a^2}_{\!\!\!1} , \dots , (2 + \epsilon_0) \delta
{a^2}_{\!\!\!m - 1} , (2 m - 2 + \epsilon_0) \delta)$ et de
hauteur de Gauss-Weil majorée par $[14 (\eta + 6) \epsilon_1 + 8
\eta + 22] \delta {a^2}_{\!\!\!1} + \circ(\delta)$ suit du
corollaire \ref{c.13}. Par ailleurs la proposition \ref{c.26} est
équivalente à dire que $\Omega_{\underline
a}({\tau^*}_{\!\!\!\!\!\!- \mathbf y}P)(\underline 0 , \dots ,
\underline 0) \in {\mathfrak q}_{\delta/2}$ ce qui est équivalent
à dire aussi que $\tau^m(Q)(- {\underline x}_1 , \dots , -
{\underline x}_m)$ \\ $\in {\mathfrak q}_{\delta/2}$. Ceci ach\`eve la
démonstration du corollaire \ref{c.27}.
$~~~~\blacksquare$\vspace{1mm}
\section{Inégalité de la hauteur à la Vojta}
Le théor\`eme qui va suivre s'énonce à l'aide des param\`etres
$\epsilon_0 , \epsilon_1 , \alpha$ et $m$ et sous les contraintes
énoncées dans les §§$6$ et $8$. Nous choisirons ensuite les
param\`etres $\epsilon_0 , \epsilon_1$ et $\alpha$ en fonction de
$\epsilon$ et $m$ pour en déduire un théor\`eme ne dépendant que
des param\`etres $\epsilon , m$ et des param\`etres liés à la
courbe elliptique $E$.
\subsection{Premier théor\`eme:}
\begin{theorem}\label{c.38}
Soit $E$ une courbe elliptique définie sur un corps de nombres $K$
de degré $D$, plongée dans ${\mathbb P}_2$ à la Weierstrass,
d'équation projective $Y^2 Z = 4 X^3 - g_2 X Z^2 - g_3 Z^3 ~ (g_2
, g_3 \in K)$ et d'élément neutre (en tant que groupe) le point à
l'infini $\mathbf 0$ représenté dans ${\mathbb P}_2$ par les
coordonnées projectives $(0 : 1 : 0)$. On munit $E$ de la hauteur
de Néron-Tate, notée $\widehat{h} := {\mid.\mid}^2$, et on se
permet de plonger $E(\overline{\mathbb Q})$ dans l'espace
euclidien $E(\overline{\mathbb Q}) \otimes \mathbb R$. Soient
aussi $S$ un ensemble fini de places de $K$, $m_v ~ (v \in M_K)$
et $\eta$ les réels positifs définis au §$2$ et ${(\lambda_v)}_{v
\in S}$ une famille de réels positifs satisfaisant:
$$\sum_{v \in S} \frac{[K_v : {\mathbb Q}_v]}{[K : \mathbb Q]} \lambda_v ~=~ 1 .$$
Soient enfin $\epsilon$ un réel strictement positif, $m$ un entier
$\geq 2$ et $\epsilon_0 , \epsilon_1 , \alpha$ des réels positifs
satisfaisant les contraintes des §§$6$ et $8$, c'est-\`a-dire:
$$\epsilon_0 \leq 1/2 , ~~\frac{m - 1}{m!} \left(\frac{7}{3}\right)^{\!\!\!m} \!\!\!\frac{\epsilon_{1}^{m}}{\epsilon_0 (m + \epsilon_0) (1 + \epsilon)^{m - 2}} \leq \frac{1}{2} ~~\mbox{et}~~ 4 m
(\epsilon_0 + 2 \alpha) \leq \epsilon \epsilon_1 .$$ Pour tout
$m$-uplet $({\mathbf x}_1 , \dots , {\mathbf x}_m) \in E^m$
satisfaisant:
$$\cos({\mathbf x}_i , {\mathbf x}_j) \geq 1 -
\frac{\alpha}{4} ~~ \forall i , j = 1 , \dots , m ,$$
$$\widehat{h}({\mathbf x}_m) \geq \widehat{h}({\mathbf x}_{m - 1}) \geq \dots \geq \widehat{h}({\mathbf x}_1)
\geq \!\left(\!\frac{6 m^2}{\epsilon_1}\right)^{\!\!m} [(71
\epsilon_1 + 55) \eta + 800 \epsilon_1 + 272]$$
$$\mbox{et}~~~~~~~~~ {\mbox{dist}}_v({\mathbf x}_i , \mathbf 0) <
e^{- \lambda_v \epsilon h({\mathbf x}_i) - 2 m_v - c_v}
~~~~~~~~~~(\forall i \in \{1 , \dots , m\} , \forall v \in S)$$ on
a l'une au moins des inégalités suivantes:
\begin{eqnarray*}
\widehat{h}({\mathbf x}_2) & < & \sqrt{2} \left(\frac{6
m^2}{\epsilon_1}\right)^m \widehat{h}({\mathbf x}_1) \\
\vdots \\
\widehat{h}({\mathbf x}_{m - 1}) & < & \sqrt{2} \left(\frac{6
m^2}{\epsilon_1}\right)^m \widehat{h}({\mathbf x}_{m - 2}) \\
\widehat{h}({\mathbf x}_m) & < & \sqrt{2} (m - 1)\left(\frac{6
m^2}{\epsilon_1}\right)^m \widehat{h}({\mathbf x}_{m - 1})~.
\end{eqnarray*}
\end{theorem}
La démonstration de ce théorème utilise un lemme de Roth sur une
puissace d'une courbe elliptique. On utilise ici le théorème suivant qui est une conséquence du lemme de Roth
de \cite{300}:
\begin{theorem}\label{c70}
Soit $E$ une courbe elliptique définie sur un corps de nombres
$K$, plongée à la Weierstrass dans le plan projectif ${\mathbb
P}_2$ et d'équation projective: $Y^2 Z = 4 X^3 - g_2 X Z^2 - g_3
Z^3 ~(g_2 , g_3 \in K)$. On prend la variable $X$ comme
param\`etre pour $E$ au voisinage du point à l'infini $\mathbf 0$
de $E$, représenté dans ${\mathbb P}_2$ par $(0 : 1 : 0)$, qu'on
consid\`ere aussi comme élément neutre de $E$ et on pose $\eta :=
h(1 : g_2 : g_3)$. Soient aussi $m \geq 2$ un entier positif,
$\epsilon'$ un réel strictement positif, $\underline{\delta} =
(\delta_1 , \dots , \delta_m) \in {{\mathbb N}^*}^m$ et ${\mathbf
x}_1 , \dots , {\mathbf x}_m$ des points de $E(K)$ représentés
respectivement dans ${\mathbb P}_2$ par les syst\`emes de
coordonnées projectives: ${\underline x}_1 = (x_1 : 1 : z_1) ,
\dots , {\underline x}_m = (x_m : 1 : z_m)$. Posons $\mathbf x$ le
point de $E^m(K) : \mathbf x = ({\mathbf x}_1 , \dots , {\mathbf
x}_m) , W(\underline{\delta} , m \epsilon')$ le dessous d'escalier
de ${\mathbb N}^m$:
$$W(\underline{\delta} , m \epsilon') := \left\{(\tau_1 , \dots , \tau_m) \in {\mathbb N}^m / \frac{\tau_1}{\delta_1}
+ \dots + \frac{\tau_m}{\delta_m} < m \epsilon'\right\}$$ et
$\Sigma(\underline{\delta} , \epsilon')$ l'ensemble pondéré:
$$\Sigma(\underline{\delta} , \epsilon') := W(\underline{\delta} ,
m \epsilon') \times \{\mathbf x\} .$$ Posons par ailleurs
${\mathfrak q}_{\underline{\delta} , m \epsilon'}$ l'idéal de
l'anneau $K[[u_1 , \dots , u_m]]$ défini par: $${\mathfrak
q}_{\underline{\delta} , m \epsilon'} := \left\{f \in K[[u_1 ,
\dots , u_m]] / f_{i_1 , \dots , i_m} = 0 ,~ \forall (i_1 , \dots
, i_m) \in W(\underline{\delta} , m \epsilon')\right\}$$ o\`u
$f_{i_1 , \dots , i_m}$ désigne le coefficient de $u_{1}^{i_1}
\dots u_{m}^{i_m}$ dans la série $f$ de $\!K[[u_1 , \dots ,
u_m]]$. Soit enfin $Q$ une forme multihomog\`ene de $K[{\underline
X}_1 , \dots , {\underline X}_m]$ (o\`u ${\underline X}_i = (X_{i
0} , X_{i 1} , X_{i 2})$ pour $i = 1 , \dots , m$), non
identiquement nulle modulo $I(E^m)$, de multidegré
$\underline{\delta}$. On suppose qu'on a:
$$~~~~~~~~\frac{\delta_i}{\delta_{i + 1}} > \left(\frac{6 m}{\epsilon'}\right)^m~~~~ \mbox{pour}~ i = 1 , \dots , m - 1 .
$$ Alors, si $Q$ s'annule sur l'ensemble pondéré $\Sigma(\underline{\delta} , m
\epsilon')$, c'est-\`a-dire si:
$$\tau^m(Q)({\underline x}_1 , \dots , {\underline x}_m) \in {\mathfrak q}_{\underline{\delta} , m \epsilon'} ,$$
il existe une sous-variété irréductible, propre et
produit $V = V_1 \times \dots \times V_m$ de $E^m$, définie sur
$K$, contenant le point $\mathbf x$ de $E^m$ et satisfaisant:
\begin{equation*}
\begin{split}
d(V_1) \dots d(V_m) \!\sum_{\l = 1}^m \!\delta_{\l} \frac{h(V_{\l})}{d(V_{\l})} \leq 3^m \!\!\left(\!\frac{2(m + 1)}{\epsilon'}\!\right)
^{\!m - \dim V} \!\!\!\left[\widetilde{h}(Q) + m (3 \log m + 7) + 12 \right. \\
\left.+ [(4 m \epsilon' + 4) \eta + 70 m \epsilon' + 26] (\delta_1 +
\dots + \delta_m)\right] .
\end{split}
\end{equation*}
\end{theorem}
{\bf Démonstration.---} On applique le théor\`eme 5.3 de
\cite{300} pour $p = m$, les groupes algébriques $G_1 , \dots ,
G_p$ sont tous égaux \`a $E$ et plongés dans ${\mathbb P}_2$ à la
Weierstrass, $P = Q , \epsilon = \epsilon'$ et ${\mathbf x}_1 ,
\dots , {\mathbf x}_p , \delta_1 , \dots , \delta_p , K$ restent
tels qu'ils sont dans le corollaire. On a ainsi $g_i = 1$, $n_i =
2$ (pour tout $i = 1 , \dots , m$), $g = m$ et $d(G_1) = \dots =
d(G_m) = d(E) = 3$. D'apr\`es le formulaire (§$13$), pour tout
point fixé de $E^2$, il existe une famille de formes $\underline
A$ de $K[X_1 , Y_1 , Z_1 , X_2 , Y_2 , Z_2]$ de bidegré $(2 , 2)$
et de hauteur de Gauss-Weil majorée par $3 \eta + 5$, représentant
l'addition sur $E$ au voisinage de ce point. On peut prendre alors
dans cette application du théor\`eme 5.3 de \cite{300} $c = c' =
2$ et $\widetilde{h}(\Asoul_i) \leq 3 \eta + 5$ pour tout $i = 1 ,
\dots m$. Il est maintenant clair que toutes les hypoth\`eses du
théor\`eme 5.3 de \cite{300} sont satisfaites, d'o\`u l'existence
d'une sous-variété irréductible propre $V$ de $E^m$, définie sur
$K$, contenant le point $\mathbf x$ de $E^m(K)$ dont toutes les
composantes irréductibles sur $\Kbar$ sont des sous-variétés
produit de ${{\mathbb P}_2}^m$ et en désignant par $\widetilde V =
V_1 \times \dots \times V_m$ une des composantes irréductibles sur
$\Kbar$ de $V$ contenant le point $\mathbf x$, on peut écrire:
$$V = \bigcup_{\sigma \in \mbox{Gal}(\overline K/K)}
\!\!\!\!\!\!\!\!\!\sigma(\widetilde V) .$$
De plus, en posant $D := \frac{d(V)}{d(\widetilde V)}$ on a:
\begin{align}
D d(V_1) \dots d(V_m) & \leq~ \left(\frac{2 m}{\epsilon'}\right)^{m -
\dim V} 3^m ~~~~~~~~~~\mbox{et} \label{3.36} \\
D d(V_1) \dots d(V_m) \sum_{i = 1}^{m} \delta_i \frac{h(V_i)}{d(V_i)} & \leq~
\left(\frac{2 (m + 1)}{\epsilon'}\right)^{m - \dim V} \!\!\!3^m
\!\!\left[\widetilde{h}(Q) + \sum_{i = 1}^{m} R_i \delta_i + S\right] \label{3.37}
\end{align}
o\`u
$$R_i := \frac{h(E)}{6} + (m - 1) \epsilon' \left\{4 h(E) + 14 \log{3} + 26\right\} + \widetilde{h}(\Asoul_i) + 16 \log{3} + 2 \log{2}$$
pour tout $i = 1 , \dots , m$ et:
$$S := m \left(3 \log{m} + 3 \log{3} + 5 \log{2}\right) + 12 .$$
 Nous montrons d'abord que dans notre situation on a $V =
\widetilde V = V_1 \times \dots \times V_m$. Comme $\mathbf x =
({\mathbf x}_1 , \dots , {\mathbf x}_m) \in \widetilde V = V_1
\times \dots \times V_m$, c'est-à-dire ${\mathbf x}_i \in V_i$
pour tout $i \in \{1 , \dots , m\}$ et comme pour tout $i = 1 ,
\dots , m$, $V_i$ est une $\Kbar$- sous-variété irréductible de
$E$ (car $\widetilde V$ est une $\Kbar$- sous-variété irréductible
de $E^m$) et que $E$ est de dimension $1$ alors, pour tout $i \in
\{1 , \dots , m\}$, soit $V_i = \{{\mathbf x}_i\}$ ou bien $V_i =
E$. On remarque que dans ces deux cas $V_i$$(i = 1 , \dots , m)$
est définie sur $K$, donc $\widetilde V = V_1 \times \dots \times
V_m$ est aussi définie sur $K$. Il s'ensuit que: $\forall \sigma
\in \mbox{Gal}(\Kbar/K) ; \sigma(\widetilde V) = \widetilde V$,
d'o\`u $V = \widetilde V = V_1 \times \dots \times V_m$,
c'est-à-dire que $V$ est elle m\^eme une sous-variété produit de
$E^m$. L'entier $D$ du théor\`eme 5.3 de \cite{300} vaut alors
dans cette application $D := d(V)/d(\widetilde V) = 1$ (car $V =
\widetilde V$). Nous remarquons ainsi que l'inégalité (\ref{3.36})
est triviale. En effet, puisque $D = 1$ et $d(V_i) \in \{1 , 3\}
~\forall i = 1 , \dots , m$ (car, comme on vient de le voir
ci-dessus, $\forall i \in \{1 , \dots , m\}$: soit $V_i =
\{{\mathbf x}_i\}$ ou bien $V_i = E$) on a: $D d(V_1) \dots d(V_m)
\leq 3^m \leq (\frac{2 m}{\epsilon'})^{m - \dim V} 3^m$. C'est
l'inégalité (\ref{3.37}) par contre qui est intéressante dans
notre application. Afin d'en déduire l'inégalité du corollaire,
montrons que la hauteur unitaire de $E$, $h(E)$, peut \^etre
majorée par $\eta + 7$. Remarquons pour cela que $E$ est une
hypersurface de ${\mathbb P}_2$ définie par la forme $U$ de $K[X ,
Y , Z]$ (o\`u $U(X , Y , Z) = Y^2 Z + g_2 X Z^2 + g_3 Z^3 - 4
X^3$) qui est de degré $3$, donc d'apr\`es [Ph4]3 (page 347) on a:
$h(E) = h(U) + 3 \sum_{i = 1}^{1} \sum_{j = 1}^{i} \frac{1}{2 j}$,
c'est-à-dire $h(E) = h(U) + 3/2$. Par ailleurs, d'après le théor\`eme $4$
de [Le] fournissant la comparaison entre hauteur unitaire et
hauteur de Gauss-Weil d'une forme donnée, on a: $h(U) \leq \widetilde{h}(U) +
\frac{1}{2} \log {{3 + 2}\choose{2}} + 3 \sum_{j = 1}^{2}
\frac{1}{2 j}$, c'est-à-dire $h(U) \leq \widetilde{h}(U) +
\frac{1}{2} \log 10 + \frac{9}{4}$. Enfin on a: $\widetilde{h}(U)
= h(1 : -4 : g_2 : g_3) \leq h(1 : g_2 : g_3) + \log 4$,
c'est-à-dire $\widetilde{h}(U) \leq \eta + \log 4$. Il s'ensuit
finalement de toutes ces comparaisons qu'on a: $h(E) \leq \eta +
\log 4 + \frac{1}{2} \log 10 + \frac{15}{4} < \eta + 7$.
L'inégalité qui conclut le corollaire découle maintenant de
(\ref{3.37}) en majorant simplement $h(E)$ par $\eta + 7$,
$\widetilde{h}(\underline A)$ par $3 \eta + 5$ et en utilisant des
majorations numériques triviales. La démonstration est achevée.
$~~~~\blacksquare$\vspace{2mm}\\
{\bf Démonstration du théorème \ref{c.38}.---} Nous procédons par
l'absurde, c'est-\`a-dire que nous supposons en plus de toutes les
hypoth\`eses du théor\`eme \ref{c.38} qu'on a:
\begin{equation}
\begin{split}
\widehat{h}({\mathbf x}_2) & ~\geq~ \sqrt{2} \left(\frac{6
m^2}{\epsilon_1}\right)^m \widehat{h}({\mathbf x}_1) , \\
\widehat{h}({\mathbf x}_3) & ~\geq~ \sqrt{2} \left(\frac{6
m^2}{\epsilon_1}\right)^m \widehat{h}({\mathbf x}_2) , \\
\vdots \\
\widehat{h}({\mathbf x}_{m - 1}) & ~\geq~ \sqrt{2} \left(\frac{6
m^2}{\epsilon_1}\right)^m \widehat{h}({\mathbf x}_{m - 2}) \\
\!\!\!\!\!\!\!\!\!\mbox{et}~~ \widehat{h}({\mathbf x}_m) & ~\geq~
\sqrt{2} (m - 1)\left(\frac{6 m^2}{\epsilon_1}\right)^m
\widehat{h}({\mathbf x}_{m - 1}) . \label{3.38}
\end{split}
\end{equation}
On applique maintenant le théorème \ref{c70} pour la
courbe elliptique $E$ et les $m$ points $- {\mathbf x}_1 , \dots ,
- {\mathbf x}_m$ de $E$. On pose pour tout $i \in \{1 , \dots ,
m\}$ : \\ $a_i := [\mid{\mathbf x}_m \mid \!/\! \mid{\mathbf x}_i
\mid]$ et on prend $\epsilon' := \frac{\epsilon_1}{m}$, $\delta_i
:= (2 + \epsilon_0) \delta {a^2}_{\!\!\!i}$ pour $i = 1 , \dots ,
m - 1$ et $\delta_m := (2 m - 2 + \epsilon_0) \delta$. Le
corollaire \ref{c.27} du §$8$ précédent nous fournit bien une
forme $Q$ de $K[{\underline X}_1 , \dots , {\underline X}_m]$
(avec ${\underline X}_i := (X_{i 0} , X_{i 1} , X_{i 2})$), non
identiquement nulle modulo l'idéal $I(E^m)$, de multidegré
$(\delta_1 , \dots , \delta_m)$ et s'annulant sur l'ensemble
pondéré $T_{\delta/2} \times \{- \mathbf x\}$. Or, un simple
calcul montre que le dessous d'escalier:
$$T_{\delta/2} := \left\{(\tau_1 , \dots , \tau_m) \in {\mathbb N}^m / \frac{\tau_1}{{a^2}_{\!\!\!1}} + \dots +
\frac{\tau_{m - 1}}{{a^2}_{\!\!\!m - 1}} + \frac{\tau_m}{m - 1}
\leq \frac{7}{2} \epsilon_1 \delta \right\}$$ contient le dessous
d'escalier:
$$W(\underline{\delta} , m \epsilon') = W(\underline{\delta} , \epsilon_1) :=
\left\{(\tau_1 , \dots , \tau_m) \in {\mathbb N}^m /
\frac{\tau_1}{\delta_1} + \dots + \frac{\tau_m}{\delta_m} <
\epsilon_1\right\} ;$$ donc la forme $Q$ s'annule aussi -a
fortiori- sur l'ensemble pondéré $\Sigma(\underline{\delta} , m
\epsilon') := W(\underline{\delta} , m \epsilon') \times \{-
\mathbf x\}$. Il nous reste alors -pour vérifier toutes les
hypoth\`eses du théorème \cite{c70}- qu'\`a montrer pour
tout $i = 1 , \dots , m - 1 : \delta_i / \delta_{i + 1} > (6 m /
\epsilon')^m$.\\ Pour $i \in \{1 , \dots , m - 2\}$ on a:
\begin{eqnarray*}
\frac{\delta_i}{\delta_{i + 1}} =
\frac{{a^2}_{\!\!\!i}}{{a^2}_{\!\!\!i + 1}} & \geq &
\frac{1}{\sqrt{2}} \frac{\widehat{h}({\mathbf x}_{i +
1})}{\widehat{h}({\mathbf x}_i)} ~~~~~~~~~~~~ (\text{d'apr\`es (\ref{3.25}) du lemme \ref{c.30}}) \\
 & > & \left(\frac{6 m^2}{\epsilon_1}\right)^{\!\!m} = \left(\frac{6
 m}{\epsilon'}\right)^{\!\!m} ~~~(\text{d'apr\`es nos hypoth\`eses (\ref{3.38})}) ;
 \end{eqnarray*}
 donc $\delta_i / \delta_{i + 1} > (6 m /
\epsilon')^m ~\forall i = 1 , \dots , m - 2$, et pour $i = m - 1$
on a: $$\frac{\delta_i}{\delta_{i + 1}} = \frac{\delta_{m -
1}}{\delta_m} = \frac{2 + \epsilon_0}{2 m - 2 + \epsilon_0}
{a^2}_{\!\!\!m - 1} > \frac{{a^2}_{\!\!\!m - 1}}{m -
1}~~~~~~~~~~~~~~~~~~~~~~$$
\begin{eqnarray*}
\mbox{et}~~~~ {a^2}_{\!\!\!m - 1} = \frac{{a^2}_{\!\!\!m -
1}}{{a^2}_{\!\!\!m}} & \geq & \frac{1}{\sqrt{2}}
\frac{\widehat{h}({\mathbf x}_m)}{\widehat{h}({\mathbf x}_{m -
1})} ~~~~~~ (\text{d'apr\`es (\ref{3.25}) du lemme \ref{c.30}})~~~~~~ \\
 & > & (m - 1) \left(\frac{6 m^2}{\epsilon_1}\right)^{\!\!m}
 ~
 (\text{d'apr\`es nos hypoth\`eses (\ref{3.38})}) ;
\end{eqnarray*}
donc $\delta_i / \delta_{i + 1} > (6 m^2 / \epsilon_1)^m = (6 m /
\epsilon')^m$ aussi pour $i = m - 1$, d'o\`u: $\forall i \in \{1 ,
\dots , m - 1\} : \delta_i / \delta_{i + 1} > (6 m /
\epsilon')^m$. Ainsi toutes les hypoth\`eses du théorème
\ref{c70} sont satisfaites, d'o\`u l'existence d'une sous-variété
irréductible, propre et produit $V = V_1 \times \dots \times V_m$
de $E^m$, définie sur $K$, contenant le point $(- {\mathbf x}_1 ,
\dots , - {\mathbf x}_m)$ de $E^m(K)$ et vérifiant:
\begin{equation}
\begin{split}
d(V_1) \dots d(V_m) \sum_{i =1}^{m} \delta_i \frac{h(V_i)}{d(V_i)}
& \leq~ 3^m \left(\frac{2 (m + 1)}{\epsilon'}\right)
^{\!\!\!m - \dim V} \left[\widetilde{h}(Q) + m (3 \log m + 7) + 12\right. \\
& \quad \left. + [(4 m \epsilon' + 4) \eta + 70 m \epsilon' + 26]
(\delta_1 + \dots +
\delta_m)\phantom{\widetilde{h}(Q)}\!\!\!\!\!\!\!\!\!\!\!\!\!\right]
. \label{3.39}
\end{split}
\end{equation}
En majorant $\delta_1 + \dots + \delta_m$ par:
\begin{eqnarray*}
\delta_1 + \dots + \delta_m & \!\!= & \!\!(2 + \epsilon_0) \delta
({a^2}_{\!\!\!1} + \dots + {a^2}_{\!\!\!m - 1}) + (2 m - 2 +
\epsilon_0) \delta \\
 & \!\!\leq & \!\!(2 + \epsilon_0) \delta ({a^2}_{\!\!\!1} + \dots + {a^2}_{\!\!\!m - 1} + m - 1) \\
& \!\!\leq & \!\!(2 + \epsilon_0)(1 + 1/24) \delta {a^2}_{\!\!\!1}
~~
(\mbox{d'apr\`es $4)$ et $5)$ du lemme \ref{c.31}}) \\
& \!\!< & \!\!3 \delta {a^2}_{\!\!\!1} ;
\end{eqnarray*}
$\widetilde{h}(Q)$ par $[14(\eta + 6) \epsilon_1 + 8 \eta + 22]
\delta {a^2}_{\!\!\!1} + o(\delta)$ (d'apr\`es le corollaire
\ref{c.27}), l'exposant $m - \dim V$ par $m$ et $m(3 \log m + 7) +
12$ par $o(\delta)$, un simple calcul montre que le deuxi\`eme
membre de (\ref{3.39}) est majoré par:
$$\left(\frac{6 m (m + 1)}{\epsilon_1}\right)^{\!\!\!m} \left[(26 \epsilon_1 + 20) \eta + 294 \epsilon_1 + 100\right]
\delta {a^2}_{\!\!\!1} + o(\delta) .$$ Par ailleurs comme la
sous-variété $V = V_1 \times \dots \times V_m$ de $E^m$ est
irréductible, chacune des sous-variétés $V_i ~(i = 1 , \dots , m)$
sera de m\^eme; par conséquent pour tout $i \in \{1 , \dots , m\}$
on a soit $V_i$ est réduite à un point (donc dans ce cas $V_i =
\{- {\mathbf x}_i\}$, car $- {\mathbf x}_i \in V_i$ puisque $(-
{\mathbf x}_1 , \dots , - {\mathbf x}_m) \in V$) ou bien $V_i$ est
égale à $E$ tout entier. Comme maintenant on ne peut pas avoir
$V_i = E ~\forall i \in \{1 , \dots , m\}$ car $V$ est une
sous-variété propre de $E^m$, alors il existe au moins un $i_0 \in
\{1 , \dots , m\}$ tel que $V_{i_0} = \{- {\mathbf x}_{i_0}\}$. le
premier membre de (\ref{3.39}) peut \^etre alors minoré par:
\begin{eqnarray*}
d(V_1) \dots d(V_m) \sum_{i = 1}^m \delta_i \frac{h(V_i)}{d(V_i)}
& \geq & \delta_{i_0} h(\{- {\mathbf x}_{i_0}\}) = \delta_{i_0}
h({\mathbf x}_{i_0}) \\
 & \geq & (2 + \epsilon_0) \delta {a^2}_{\!\!\!i_0} h({\mathbf
 x}_{i_0}) \\
  & > & 2 \delta {a^2}_{\!\!\!i_0} h({\mathbf x}_{i_0}) \\
 & \geq & \delta {a^2}_{\!\!\!1} h({\mathbf x}_1) ~~~~ \text{(d'apr\`es (\ref{3.27}) du lemme \ref{c.30})} .
\end{eqnarray*}
Il suit de cette minoration du premier membre de l'inégalité
(\ref{3.39}) et de la majoration précédente pour son deuxi\`eme
membre qu'on a:
$$\delta {a^2}_{\!\!\!1} h({\mathbf x}_1) < \left(\frac{6 m (m +
1)}{\epsilon_1}\right)^{\!\!\!m} \!\!\!\left[(26 \epsilon_1 + 20)
\eta + 294 \epsilon_1 + 100\right] \delta {a^2}_{\!\!\!1} +
o(\delta) .$$ Cette derni\`ere inégalité entra{\sf\^\i}ne, en
divisant ces deux membres par $\delta {a^2}_{\!\!\!1}$ et en
faisant tendre $\delta$ vers l'infini:
\begin{eqnarray*}
h({\mathbf x}_1) & < & \left(\frac{6 m (m +
1)}{\epsilon_1}\right)^{\!\!\!m} \!\!\!\left[(26 \epsilon_1 + 20)
\eta + 294 \epsilon_1 + 100\right] \\
 & < & \left(\frac{6 m^2}{\epsilon_1}\right)^{\!\!\!m} \!\!\!\left[(26 e \epsilon_1 + 20 e)
\eta + 294 e \epsilon_1 + 100 e\right] ~~ \text{(car $(m + 1)^m <
e
m^m$)} \\
 & < & \left(\frac{6 m^2}{\epsilon_1}\right)^{\!\!\!m} \!\!\!\left[(71 \epsilon_1 + 54,5)
\eta + 800 \epsilon_1 + 271,9\right] .
\end{eqnarray*}
Par le théor\`eme \ref{c.14} du formulaire on obtient enfin:
$$\widehat{h}({\mathbf x}_1) < \left(\frac{6 m^2}{\epsilon_1}\right)^{\!\!\!m} \!\!\!\left[(71 \epsilon_1 + 55)
\eta + 800 \epsilon_1 + 272\right] ,$$ ce qui est en contradiction
avec notre hypoth\`ese:
$$\widehat{h}({\mathbf x}_m) \geq \widehat{h}({\mathbf x}_{m - 1}) \geq \dots \geq \widehat{h}({\mathbf x}_1)
\geq \!\left(\!\frac{6 m^2}{\epsilon_1}\right)^{\!\!m}
\!\!\!\!\left[(71 \epsilon_1 + 55) \eta + 800 \epsilon_1 +
272\right] .$$ La démonstration est achevée.
$~~~~\blacksquare$\vspace{1mm}
\subsection{Choix des param\`etres $\epsilon_0 , \epsilon_1$ et $\alpha$:}
Soit $E$ une courbe elliptique comme dans le théor\`eme
\ref{c.38}, $\epsilon$ un réel strictement positif et $m$ un
entier $\geq 2$. Afin d'avoir un énoncé dépendant de moins de
param\`etres, nous choisissons les réels strictement positifs
$\epsilon_0 , \epsilon_1$ et $\alpha$ en fonction de $\epsilon$ et
$m$ de façon à satisfaire les contraintes du théor\`eme
\ref{c.38}, qui sont: $\epsilon_0 \leq 1/2,$
\begin{eqnarray}
\frac{m - 1}{m!} \left(\frac{7}{3}\right)^{\!\!m}
\!\!\!\!\!\frac{{\epsilon_1}^{\!m}}{\epsilon_0 (m + \epsilon_0) (1
+
\epsilon_0)^{m - 2}} & \leq & \frac{1}{2} ~~\text{et} \label{3.40} \\
4 m (\epsilon_0 + 2 \alpha) & \leq & \epsilon \epsilon_1 .
\label{3.41}
\end{eqnarray}
La procédure pour ce choix est la suivante:
\\ Nous remarquons que la relation (\ref{3.40}) est impliquée par la
relation suivante:
\begin{equation}
\epsilon_0 ~\geq~ \frac{2
\left(\frac{7}{3}\right)^{\!m}}{m!}~\!{{\epsilon_1}^{\!m}}
\label{3.42}
\end{equation}
et que la relation (\ref{3.41}) est impliquée par les deux
relations suivantes:
\begin{eqnarray}
\epsilon_0 & \leq & \frac{1}{8 m} ~\!\epsilon
\epsilon_1 \label{3.43} \\
\alpha & \leq & \frac{1}{16 m} ~\!\epsilon \epsilon_1 .
\label{3.44}
\end{eqnarray}
Les deux relations (\ref{3.42}) et (\ref{3.43}) s'assemblent en la
relation suivante donnant un encadrement pour le param\`etre
$\epsilon_0$ en fonction des autres param\`etres:
\begin{equation}
\frac{2 \left(\frac{7}{3}\right)^{\!m}}{m!}~\!{{\epsilon_1}^{\!m}}
~\leq~ \epsilon_0 ~\leq~ \frac{1}{8 m} ~\!\epsilon \epsilon_1 .
\label{3.45}
\end{equation}
Ainsi nous venons de montrer que le théor\`eme \ref{c.38} reste -a
fortiori- valable quand on remplace les deux contraintes
(\ref{3.40}) et (\ref{3.41}) par les deux relations (\ref{3.44})
et (\ref{3.45}). Il nous suffit donc de choisir $\epsilon_0 ,
\epsilon_1$ et $\alpha$ (en fonction de $m$ et $\epsilon$) de
façon \`a satisfaire (\ref{3.44}) et (\ref{3.45}). Nous remarquons
que la relation (\ref{3.45}) exige d'avoir $\frac{2
\left(\frac{7}{3}\right)^{\!m}}{m!}~\!{{\epsilon_1}^{\!m}} \leq
\frac{1}{8 m} ~\!\epsilon \epsilon_1 $, c'est-à-dire d'avoir:
\begin{equation}
\epsilon_1 ~\leq~
\frac{3}{7}.\left(\frac{1}{38}\right)^{\!\!\!\frac{1}{m -
1}}\!\!.\left((m - 1)!\right)^{\frac{1}{m -
1}}.\epsilon^{\frac{1}{m - 1}} . \label{3.46}
\end{equation}
Maintenant comme $\left(\frac{1}{38}\right)^{\frac{1}{m - 1}} \geq
\frac{1}{38}$ (car $m \geq 2$) et $\left((m -
1)!\right)^{\frac{1}{m - 1}} \geq \frac{m - 1}{e}$ (car $\forall n
\in {\mathbb N}^* : \frac{n}{e} \leq (n!)^{\frac{1}{n}} \leq n$),
alors (\ref{3.46}) est impliquée par la relation:
$$\epsilon_1 ~\leq~ \frac{3}{7}.\frac{1}{38}.\frac{1}{e} (m - 1).\epsilon^{\frac{1}{m - 1}}$$
qui est de nouveau impliquée par la suivante:
\begin{equation}
\epsilon_1 ~\leq~ \frac{1}{484} m \epsilon^{\frac{1}{m - 1}}
\label{3.47}
\end{equation}
puisque $\frac{3}{7} \frac{1}{38} \frac{1}{e} > \frac{1}{242}$ et
$m - 1 \geq \frac{m}{2}$ (car $m \geq 2$). \\ En regardant
maintenant l'énoncé du théor\`eme \ref{c.38} (inégalité de Vojta),
on voit qu'on a plut\^ot intér\^et à prendre $\epsilon_1$ le plus
grand possible, donc le meilleur choix pour $\epsilon_1$ sous
(\ref{3.47}) est:
\[ \boxed{\epsilon_1 ~=~ \frac{1}{484} m \epsilon^{\frac{1}{m - 1}}}~.\] En reportant maintenant cette valeur de $\epsilon_1$ dans (\ref{3.45}), on
a, pour $\epsilon_0$, l'encadrement suivant:
$$2 \left(\frac{7}{1452}\right)^{\!\!m} \frac{m^m}{m!} \epsilon^{\frac{m}{m - 1}} ~\leq~ \epsilon_0 ~\leq~ \frac{1}{3872} \epsilon^{\frac{m}{m - 1}} .$$ Or, on peut vérifier qu'on a: $\forall m \geq 2 : 2
\left(\frac{7}{1452}\right)^{\!\!m} \frac{m^m}{m!} \leq
\frac{1}{4000}$ (pour montrer cela, il suffit de faire le calcul
numérique pour $m = 2$ et d'utiliser la majoration $\frac{m^m}{m!}
\leq \sum_{n = 0}^{\infty} \frac{m^n}{n!} = e^m$ quand $m \geq
3$). Le dernier encadrement pour $\epsilon_0$ est satisfait pour
la valeur:
\[ \boxed{\epsilon_0 ~=~ \frac{1}{4000} \epsilon^{\frac{m}{m - 1}}}~.\] En reportant aussi la valeur choisie pour $\epsilon_1$ dans (\ref{3.44}), on obtient:
\begin{equation}
\alpha ~\leq~ \frac{1}{7744} \epsilon^{\frac{m}{m - 1}} .
\label{3.48}
\end{equation}
Le param\`etre $\alpha$ choisi, on doit recouvrir l'espace
euclidien $E(\overline{\mathbb Q}) \otimes \mathbb R$ par des
c\^ones d'angle $\simeq \arccos(1 - \alpha/4)$. Or, notre souhait
est d'avoir un nombre minimal de tels c\^ones pour ce
recouvrement, donc on a plut\^ot intér\^et à prendre ces c\^ones
les plus larges possible, ce qui revient à choisir $\alpha$ le
plus grand possible. D'o\`u, le meilleur choix pour $\alpha$ sous
(\ref{3.48}): \[ \boxed{\alpha ~=~ \frac{1}{7744}
\epsilon^{\frac{m}{m - 1}}}~.\] Ce choix des param\`etres
$\epsilon_0 , \epsilon_1$ et $\alpha$ en fonction de $\epsilon$ et
$m$ satisfait bien les relations (\ref{3.44}) et (\ref{3.45}), par
conséquent il satisfait -a fortiori- les contraintes (\ref{3.40})
et (\ref{3.41}) du théor\`eme \ref{c.38}. Ce dernier reste valable
lorsqu'on remplace $\epsilon_0 , \epsilon_1$ et $\alpha$ par les
choix précédents et on obtient finalement le théor\`eme suivant:
\begin{theorem}\label{c.39}
Soit $E$ une courbe elliptique définie sur un corps de nombres $K$
de degré $D$, plongée dans ${\mathbb P}_2$ à la Weierstrass,
d'équation projective $Y^2 Z = 4 X^3 - g_2 X Z^2 - g_3 Z^3 ~ (g_2
, g_3 \in K)$ et d'élément neutre (en tant que groupe) le point à
l'infini $\mathbf 0$ représenté dans ${\mathbb P}_2$ par les
coordonnées projectives $(0 : 1 : 0)$. On munit $E$ de la hauteur
de Néron-Tate, notée $\widehat{h} := {\nb{.}}^2$, et on se permet
de plonger $E(\overline{\mathbb Q})$ dans l'espace euclidien
$E(\overline{\mathbb Q}) \otimes \mathbb R$. Soient aussi $S$ un
ensemble fini de places de $K$, $m_v ~ (v \in M_K)$ et $\eta$ les
réels positifs définis au §$2$ et ${(\lambda_v)}_{v \in S}$ une
famille de réels positifs satisfaisant:
$$\sum_{v \in S} \frac{[K_v : {\mathbb Q}_v]}{[K : \mathbb Q]} \lambda_v ~=~ 1 .$$
Soient enfin $\epsilon$ un réel strictement positif et $m$ un
entier $\geq 2$. Pour tout $m$-uplet $({\mathbf x}_1 , \dots ,$ \\
${\mathbf x}_m) \in E^m$ satisfaisant:
\begin{equation*}
\cos({\mathbf x}_i , {\mathbf x}_j) ~\geq~ 1 - \frac{1}{30976}
\epsilon^{\frac{m}{m - 1}} ~~~~~~~~~~\text{pour $i , j = 1 , \dots
, m$,}~~~~~~~~~~
\end{equation*}
\begin{equation*}
\begin{split}
\widehat{h}({\mathbf x}_m) ~\geq~ \widehat{h}({\mathbf x}_{m - 1}) ~\geq~ \dots ~\geq~ \widehat{h}({\mathbf x}_1) ~\geq & \\
&\!\!\!\!\!\!\!\!\!\!\!\!\!\!\!\!\!\!\!\!\!\!\!\!\!\!\!\!\!\!\!\!\!\!\!\!\!\!\!\!\!\!\!\!\!\!\!\!\!\!\!\!\!\!\!\!\!\!\!\!\!\!\!\!\!\!\!\!\!\!\!\!\!\!\!\left(2904
m\right)^m \left[\left\{55 \epsilon^{- \frac{m}{m - 1}} + m
\epsilon^{- 1}\right\} \eta + \left\{272 \epsilon^{- \frac{m}{m -
1}} + 2 m \epsilon^{- 1}\right\}\right]
\end{split}
\end{equation*}
\begin{equation*}
\text{et}~~~~ {\mbox{dist}}_v({\mathbf x}_i , \mathbf 0) ~<~ e^{-
\lambda_v \epsilon h({\mathbf x}_i) - 2 m_v - c_v}
~~\text{$\forall i \in \{1 , \dots , m\}$ et $\forall v \in
S$;}~~~~~~~~~
\end{equation*}
on a l'une au moins des inégalités suivantes:
\begin{eqnarray*}
\widehat{h}({\mathbf x}_2) & < & \sqrt{2} \left(2904 m\right)^m
\epsilon^{- \frac{m}{m - 1}} \widehat{h}({\mathbf x}_1) \\
\vdots \\
\widehat{h}({\mathbf x}_{m - 1}) & < & \sqrt{2} \left(2904
m\right)^m \epsilon^{- \frac{m}{m - 1}} \widehat{h}({\mathbf x}_{m - 2}) \\
\widehat{h}({\mathbf x}_m) & < & \sqrt{2} (m - 1) \left(2904
m\right)^m \epsilon^{- \frac{m}{m - 1}} \widehat{h}({\mathbf x}_{m
- 1})~.
\end{eqnarray*}
\end{theorem}
\section{Inégalité de la hauteur à la Mumford}
Notre but dans cette section est de démontrer le théor\`eme
suivant:
\begin{theorem}\label{c.40}
Soit $E$ une courbe elliptique définie sur un corps de nombres $K$
de degré $D$, plongée dans ${\mathbb P}_2$ à la Weierstrass,
d'équation projective $Y^2 Z = 4 X^3 - g_2 X Z^2 - g_3 Z^3 ~ (g_2
, g_3 \in K)$ et d'élément neutre (en tant que groupe) le point à
l'infini $\mathbf 0$ représenté dans ${\mathbb P}_2$ par les
coordonnées projectives $(0 : 1 : 0)$. On munit $E$ de la hauteur
de Néron-Tate, notée $\widehat{h} := {\mid.\mid}^2$, et on se
permet de plonger $E(\overline{\mathbb Q})$ dans l'espace
euclidien $E(\overline{\mathbb Q}) \otimes \mathbb R$. Soient
aussi $S$ un ensemble fini de places de $K$, $m_v ~ (v \in M_K)$
et $\eta$ les réels positifs définis au §$2$ et ${(\lambda_v)}_{v
\in S}$ une famille de réels positifs satisfaisant:
$$\sum_{v \in S} \frac{[K_v : {\mathbb Q}_v]}{[K : \mathbb Q]} \lambda_v ~=~ 1 .$$
Soient enfin $\epsilon$ un réel strictement positif et $\beta$,
$\theta$ les deux réels positifs: $\beta := \frac{\epsilon}{2}$ et
$\theta := \frac{1}{3} \sqrt{\epsilon}$. Alors pour tout couple
$({\mathbf x}_1 , {\mathbf x}_2) \in E^2(K) , {\mathbf x}_1 \neq
{\mathbf x}_2$, satisfaisant:
\begin{equation}
\cos({\mathbf x}_1 , {\mathbf x}_2) ~\geq~ 1 - \frac{\beta}{4}~,
\label{3.49}
\end{equation}
\begin{equation}
\widehat{h}({\mathbf x}_2) ~\geq~ \widehat{h}({\mathbf x}_1)
~\geq~ \frac{1}{\epsilon} (54 \eta + 204) ~~\text{et} \label{3.50}
\end{equation}
\begin{equation}
{\mbox{dist}}_v({\mathbf x}_j , \mathbf 0) ~<~ e^{- \lambda_v
\epsilon h({\mathbf x}_j) - 2 m_v - c_v} ~~\text{pour $j = 1 , 2$
et $v \in S$;} \label{3.51}
\end{equation}
on a: $$\widehat{h}({\mathbf x}_2) \geq (1 + \theta)
\widehat{h}({\mathbf x}_1) .$$
\end{theorem}
Avant de rentrer dans les étapes de la démonstration, expliquons
grosso-modo le schéma qu'on va suivre. Nous allons procéder de la
m\^eme mani\`ere que pour le théor\`eme \ref{c.38} (l'inégalité de
la hauteur à la Vojta) sauf qu'ici la situation est beaucoup plus
simple du fait que d'une part la construction des fonctions
auxiliaires est explicite -donc ne nécessite pas l'utilisation du
lemme de Siegel- et d'autre part qu'on n'aura pas besoin
d'utiliser un lemme de zéros pour la conclusion finale! Plus
précisément, voilà ce qu'on va faire: on proc\`ede par l'absurde,
c'est-à-dire on suppose qu'il existe un couple $({\mathbf x}_1 ,
{\mathbf x}_2)$ de points de $E^2(K)$, vérifiant toutes les
hypoth\`eses du théor\`eme \ref{c.40} et tel que:
$\widehat{h}({\mathbf x}_2) < (1 + \theta) \widehat{h}({\mathbf
x}_1)$. On construit explicitement deux formes $Q_1$ et $Q_2$ de
$A(E^2)$ qui ont $({\mathbf x}_1 , {\mathbf x}_2)$ comme point
d'intersection isolé. Ensuite on fait l'extrapolation qui, du fait
que nos deux points ${\mathbf x}_1$ et ${\mathbf x}_2$ sont
supposés assez proches de l'origine (au sens de la distance
projective), oblige nos fonctions auxiliaires $Q_1$ et $Q_2$ à
s'annuler aussi au point $(\mathbf 0 , \mathbf 0)$ de $E^2(K)$.
Finalement l'annulation de $Q_1$ et $Q_2$ au point $(\mathbf 0 ,
\mathbf 0)$ entra{\sf\^\i}ne facilement (sans utiliser un lemme de
zéros) qu'on a ${\mathbf x}_1 = {\mathbf x}_2$, ce qui est en
contradiction avec notre hypoth\`ese et achevera la démonstration
du théor\`eme \ref{c.40}.\\ Soient ${\mathbf x}_1$ et ${\mathbf
x}_2$ deux points de $E(K)$ satisfaisant toutes les hypoth\`eses
du théor\`eme \ref{c.40}. Comme les hauteurs normalisées de
${\mathbf x}_1$ et ${\mathbf x}_2$ sont non nulles (puisqu'on a
fait l'hypoth\`ese $\widehat{h}({\mathbf x}_2) \geq
\widehat{h}({\mathbf x}_1) \geq \mbox{fct}(D , \epsilon , \eta ,
m_w) > 0$) alors ${\mathbf x}_1$ et ${\mathbf x}_2$ sont tous les
deux différents des $4$ points de $2$-torsion de $E$, et donc on
peut représenter respectivement ${\mathbf x}_1$ et ${\mathbf x}_2$
dans ${\mathbb P}_2(K)$ par des syst\`emes de coordonnées
projectives de la forme: ${\underline x}_1 = (x_1 : 1 : z_1)$ et
${\underline x}_2 = (x_2 : 1 : z_2)$. Soient aussi $\mathbf y$ le
point de $E(K) : \mathbf y := {\mathbf x}_1 - {\mathbf x}_2$ et
$\underline y$ un de ses représentants dans ${\mathbb P}_2(K)$.
\subsection{Construction des fonctions auxiliaires:}
Les hypoth\`eses (\ref{3.51}) du théor\`eme \ref{c.40}
entra{\sf\^\i}nent d'apr\`es le lemme \ref{c.29} du formulaire
qu'il existe une famille de formes totalement explicites
$\underline D := (D_0 , D_1 , D_2)$ de $K[X_1 , Y_1 , Z_1 , X_2 ,
Y_2 , Z_2]$, de bidegrés $(2 , 2)$, représentant la différence
dans $E$ au voisinage de $\{\mathbf 0\} \times \{\mathbf 0\}$
ainsi qu'au voisinage de $\{{\mathbf x}_1\} \times \{{\mathbf
x}_2\}$. De plus, on peut vérifier (puisque c'est explicite) que
cette famille $\underline D$ est de hauteur logarithmique locale
$v$-adique majorée par $3 m_v$ lorsque $v$ est une place finie sur
$K$ et de longueur logarithmique locale $v$-adique majorée par $3
m_v + \log 425$ lorsque $v$ est une place infinie sur $K$ et que
la famille constituée seulement des deux formes $D_0$ et $D_2$ est
de hauteur logarithmique locale $v$-adique majorée par $2 m_v$
lorsque $v$ est une place finie sur $K$ et de longueur
logarithmique locale $v$-adique majorée par $2 m_v + \log 152$
lorsque $v$ est une place infinie sur $K$.
\\ Considérons maintenant les deux formes suivantes $P_1$ et $P_2$
de $K[\underline Y]$ (avec $\underline Y := (Y_0 , Y_1 ,$ \\
$Y_2)$) définies par:
\begin{eqnarray*}
P_1(\underline Y) & := & D_0(\underline Y , \underline y) \\
P_2(\underline Y) & := & D_2(\underline Y , \underline y)
\end{eqnarray*}
et les deux formes suivantes $Q_1$ et $Q_2$ de $K[{\underline X}_1
, {\underline X}_2]$ (avec ${\underline X}_1 := (X_{1 0} , X_{1 1}
, $ \\ $X_{1 2})$ et ${\underline X}_2 := (X_{2 0} , X_{2 1} ,
X_{2 2})$) définies par:
\begin{eqnarray*}
Q_1({\underline X}_1 , {\underline X}_2) & := &
P_1\!\left(\underline D({\underline X}_1 , {\underline
X}_2)\right)
 ~=~  D_0\!\!\left(\underline D({\underline X}_1 , {\underline X}_2) ,
\underline y\right) \\
Q_2({\underline X}_1 , {\underline X}_2) & := &
P_2\!\left(\underline D({\underline X}_1 , {\underline
X}_2)\right)
 ~=~  D_2\!\!\left(\underline D({\underline X}_1 , {\underline X}_2) ,
\underline y\right) .
\end{eqnarray*}
On voit clairement que nos deux formes $P_1$ et $P_2$ de
$K[\underline Y]$ sont non identiquement nulles modulo $I(E)$ et
qu'elles s'annulent, toutes les deux, au point $\mathbf y$ de
$E(K)$, par conséquent les deux formes déduites $Q_1$ et $Q_2$
sont non identiquement nulles modulo $I(E^2)$ et s'annulent toutes
les deux au point $({\mathbf x}_1 , {\mathbf x}_2)$ de $E^2(K)$.
\\ Pour faire le lien avec le §$5$, ce qui va nous permettre
d'appliquer les résultats qui y sont obtenus, il suffit de prendre
$m = 2$ et $a_1 = a_2 = 1$, c'est-à-dire $\underline a = (1 , 1)$.
En posant $\psi_{\underline a} = \psi , \varphi_{\underline a} =
\varphi , \Omega_{\underline a} = \Omega$ et $i$ le plongement de
Weierstrass de $E$ dans ${\mathbb P}_2$, on a:
$$\begin{array}{lcccc}
\varphi : ~~~E^2 & \stackrel{\psi}{\rightarrow} & E^2 \times E &
\stackrel{i^3}{\hookrightarrow} & {{\mathbb P}_2}^3 \\
~~~~({\mathbf p}_1 , {\mathbf p}_2) & \mapsto & ({\mathbf p}_1 ,
{\mathbf p}_2 , {\mathbf p}_1 - {\mathbf p}_2) & \mapsto &
\left(i({\mathbf p}_1) , i({\mathbf p}_2) , i({\mathbf p}_1 -
{\mathbf p}_2)\right) .
\end{array}$$
Ainsi, pour toute forme $P$ de $K[{\underline X}_1 , {\underline
X}_2 , \underline Y]$ on a: $\Omega(P) = \tau^2(Q)$, o\`u $Q$
désigne la forme de $K[{\underline X}_1 , {\underline X}_2]$
définie en fonction de $P$ par:
$$Q({\underline X}_1 , {\underline
X}_2) := P\left({\underline X}_1 , {\underline X}_2 , \underline
D({\underline X}_1 , {\underline X}_2)\right) .$$ En considérant
les formes $P_1$ et $P_2$ de $K[\underline Y]$ comme des formes de
$K[{\underline X}_1 , {\underline X}_2 , \underline Y]$ on a:
$$\Omega(P_1) = \tau^2(Q_1) ~~ \mbox{et}~~ \Omega(P_2) =
\tau^2(Q_2) .$$ Par ailleurs, comme la famille $\underline D$ est
constituée des formes $D_0 , D_1$ et $D_2$, toutes de bidegré $(2
, 2)$, on voit clairement que nos deux formes $P_1$ et $P_2$ de
$K[\underline Y]$ sont, toutes les deux, de degré $2$ et que nos
deux formes $Q_1$ et $Q_2$ de $K[{\underline X}_1 , {\underline
X}_2]$ sont, toutes les deux, de bidegré $(4 , 4)$. Les
estimations des hauteurs (ou longueurs) logarithmiques locales de
$P_1 , P_2 , Q_1$ et $Q_2$ sont données par le lemme suivant:
\begin{lemma}\label{c.41}
Pour toute place $v$ de $K$ on a:
\\ $\bullet$ lorsque $v$ est finie:
\begin{eqnarray*}
h_v\!\!\left(\{P_1 , P_2\}\right) & \leq & 2 m_v + 2
h_v(\underline y)
\\
h_v\!\!\left(\{Q_1 , Q_2\}\right) & \leq & 8 m_v + 2
h_v(\underline y)
\end{eqnarray*}
$\bullet$ et lorsque $v$ est infinie:
\begin{eqnarray*}
\ell_v\!\!\left(\{P_1 , P_2\}\right) & \leq & 2 m_v + 2
h_v(\underline y) + \log 152
\\
\ell_v\!\!\left(\{Q_1 , Q_2\}\right) & \leq & 8 m_v + 2
h_v(\underline y) + 18 .
\end{eqnarray*}
\end{lemma}
{\bf Démonstration.---} D'apr\`es les propriétés des hauteurs et
longueurs logarithmiques locales, on a pour toute place $v$ de
$K$:
\begin{eqnarray*}
h_v\!\!\left(\{P_1 , P_2\}\right) & \leq & h_v\!\!\left(\{D_0 ,
D_2\}\right) + 2 h_v(\underline y) , \\
\ell_v\!\!\left(\{P_1 , P_2\}\right) & \leq &
\ell_v\!\!\left(\{D_0 ,
D_2\}\right) + 2 h_v(\underline y) , \\
h_v\!\!\left(\{Q_1 , Q_2\}\right) & \leq & 2 h_v(\underline D) +
h_v\!\!\left(\{P_1 , P_2\}\right) \\
\mbox{et}~ \ell_v\!\!\left(\{Q_1 , Q_2\}\right) & \leq & 2
\ell_v(\underline D) + \ell_v\!\!\left(\{P_1 , P_2\}\right) .
\end{eqnarray*}
Le reste est un simple calcul.  $~~~~\blacksquare$\vspace{1mm}
\begin{proposition}\label{c.42}
Pour $j = 1 , 2$ on a:
$$\tau^2(Q_j) = \frac{1}{X_{1 1}^{4} X_{2 1}^{4}} \sum_{i_1 \geq 0 , i_2 \geq 0} \frac{X_{1 1}^{i_1} X_{2 1}^{i_2}
\partial^{(i_1 , i_2)}\!P_j}{{\Delta({\underline X}_1)}^{f(i_1)} {\Delta({\underline X}_2)}^{f(i_2)}} u_{1}^{i_1}
u_{2}^{i_2}$$ o\`u les $\partial^{(i_1 , i_2)} P_j ~ ((i_1 , i_2)
\in {\mathbb N}^2)$ sont des formes de $K[{\underline X}_1 ,
{\underline X}_2]$ de degrés:
\begin{eqnarray*}
{d°}_{\!\!\!\!{\underline X}_1} \partial^{(i_1 , i_2)} P_j & \leq
& g(i_1)
+ 4 \\
{d°}_{\!\!\!\!{\underline X}_2} \partial^{(i_1 , i_2)} P_j & \leq
& g(i_2) + 4
\end{eqnarray*}
et, pour $T \in \mathbb N$ et $v$ une place sur $K$, la famille
des formes $\partial^{(i_1 , i_2)} P_j , i_1 + i_2 \leq T$, est de
hauteur logarithmique locale $h_v$ majorée par:
$$(2 T + 14) m_v + 2 h_v(\underline y)$$ lorsque $v$ est finie et
elle est de longueur logarithmique locale $\ell_v$ majorée par:
$$(2 T + 14) m_v + 2 h_v(\underline y) + 12 T + 46$$ lorsque $v$
est infinie.
\end{proposition}
{\bf Démonstration.---} Il suffit de s'adapter à la situation du
§$5$ comme expliqué auparavant et d'appliquer le corollaire
\ref{c.11} de cette derni\`ere et d'utiliser les estimations des
hauteurs (et longueurs) locales des formes $P_1$ et $P_2$ données
par le lemme \ref{c.41}.  $~~~~\blacksquare$\vspace{1mm}
\begin{proposition}\label{c.43}
Pour $j = 1 , 2$ on a:
$$\tau^2(Q_j)({\underline x}_1 , {\underline x}_2) = \!\!\!\sum_{i_1 \geq 0 , i_2 \geq 0} \!\!\!\!f_{i_1 , i_2}^{(j)} u_{1}^{i_1} u_{2}^{i_2}$$
o\`u $f_{i_1 , i_2}^{(j)} := \frac{\left(\partial^{(i_1 ,
i_2)}\!P_j\right)({\underline x}_1 , {\underline x}_2)}
{{\Delta({\underline X}_1)}^{f(i_1)} {\Delta({\underline
X}_2)}^{f(i_2)}} , ~ (i_1 , i_2) \in {\mathbb N}^2$, sont des
nombres de $K$, nuls pour $i_1 = i_2 = 0$ et satisfaisant pour
tout $(i_1 , i_2) \in {\mathbb N}^2$ et pour toute place $v \in
S$:
$${\mid f_{i_1 , i_2}^{(j)}~\mid}_v \leq \begin{cases}
\exp\{14 m_v + 2 h_v(\underline y)\}.(e^{2 m_v})^{i_1 + i_2} & \text{si $v$ est finie} \\
\exp\{14 m_v + 2 h_v(\underline y) + 46\}.(e^{2 m_v + 13})^{i_1 +
i_2} & \text{si $v$ est infinie}
\end{cases} .$$
\end{proposition}
{\bf Démonstration.---} Soit $j \in \{1 , 2\}$ fixé. Le fait que
$f_{i_1 , i_2}^{(j)} = 0$ pour $i_1 = i_2 = 0$ est simplement d\^u
au fait que notre forme $Q_j$ s'annule au point $({\underline x}_1
, {\underline x}_2)$. Soit maintenant $(i_1 , i_2)$ un couple fixé
de ${\mathbb N}^2$ et montrons les majorations de la proposition
\ref{c.43} pour les ${\mid f_{i_1 , i_2}^{(j)}~\mid}_v ~ (v \in
S)$. En désignant par $(d_1 , d_2)$ le bidegré de la forme
$\partial^{(i_1 , i_2)}P_j$ de $K[{\underline X}_1 , {\underline
X}_2]$, on a clairement pour toute place $v$ de $K$:
$${\mid\left(\partial^{(i_1 , i_2)}P_j\right)({\underline x}_1 ,
{\underline x}_2)\mid}_v \leq
\begin{cases}
H_v\!\!\left(\partial^{(i_1 , i_2)}P_j\right).{H_v({\underline
x}_1)}^{d_1} {H_v({\underline
x}_2)}^{d_2} & \!\!\!\text{si $v$ est finie} \\
L_v\!\!\left(\partial^{(i_1 , i_2)}P_j\right).{H_v({\underline
x}_1)}^{d_1} {H_v({\underline x}_2)}^{d_2} & \!\!\!\text{si $v$
est infinie}
\end{cases} .$$
Or, lorsque $v \in S$, l'hypoth\`ese principale (\ref{3.51}) du
théor\`eme \ref{c.40} entra{\sf\^\i}ne, d'apr\`es la propriété
$ii)$ du §$2$ pour les distances ${\mbox{dist}}_v ~ (v \in M_K)$
qu'on a: $H_v({\underline x}_1) = H_v({\underline x}_2) = 1$.
D'o\`u, pour toute place $v \in S$:
$${\mid\left(\partial^{(i_1 , i_2)}P_j\right)({\underline x}_1 , {\underline x}_2)\mid}_v \leq
\begin{cases}
H_v\!\!\left(\partial^{(i_1 ,
i_2)}P_j\right) & \text{si $v$ est finie} \\
L_v\!\!\left(\partial^{(i_1 , i_2)}P_j\right) & \text{si $v$ est
infinie}
\end{cases} .$$
En appliquant par ailleurs la proposition \ref{c.42} pour $T = i_1
+ i_2$ on a, pour toute place $v$ de $K$:
\\ $\bullet$ lorsque $v$ est finie:
\begin{eqnarray*}
H_v\!\!\left(\partial^{(i_1 , i_2)}P_j\right) & = &
\exp\left\{h_v\!\!\left(\partial^{(i_1 , i_2)}P_j\right)\right\} \\
 & \leq & \exp\left\{2 m_v (i_1 + i_2) + 14 m_v + 2 h_v(\underline
 y)\right\} ,
\end{eqnarray*}
\\ $\bullet$ et lorsque $v$ est infinie:
\begin{eqnarray*}
L_v\!\!\left(\partial^{(i_1 , i_2)}P_j\right) & = &
\exp\left\{\ell_v\!\!\left(\partial^{(i_1 , i_2)}P_j\right)\right\} \\
 & \leq & \exp\left\{(2 m_v + 12)(i_1 + i_2) + 14 m_v + 2 h_v(\underline y) + 46\right\} .
\end{eqnarray*}
D'o\`u, pour toute place $v \in S$:
\begin{equation}
\begin{split}
{\mid\left(\partial^{(i_1 , i_2)}P_j\right)({\underline x}_1 , {\underline x}_2)\mid}_v & ~~ \\
&
\!\!\!\!\!\!\!\!\!\!\!\!\!\!\!\!\!\!\!\!\!\!\!\!\!\!\!\!\!\!\!\!\!\!\!\!\!\!\!\!\!\!\!\!\!\!\!\!\!\!\!
\leq
\begin{cases}
\exp\!\left\{2 m_v (i_1 + i_2) + 14
m_v + 2 h_v(\underline y)\right\} & \text{si $v$ est finie} \\
\exp\!\left\{(2 m_v + 12)(i_1 + i_2) + 14 m_v + 2 h_v(\underline
y) + 46\right\} & \text{si $v$ est infinie}
\end{cases} . \label{3.52}
\end{split}
\end{equation}
D'autre part, l'hypoth\`ese (\ref{3.51}) du théor\`eme \ref{c.40}
entra{\sf\^\i}ne, d'apr\`es le lemme \ref{c.35} de l'appendice,
qu'on a pour toute place $v$ dans $S$:
\begin{equation}
\frac{1}{{\mid{\Delta({\underline x}_1)}^{f(i_1)}
{\Delta({\underline x}_2)}^{f(i_2)}\mid}_v} \leq
\begin{cases} 1 & \text{si $v$ est finie} \\
{\sqrt{e}}^{f(i_1) + f(i_2)} \leq e^{i_1 + i_2} & \text{si $v$ est
infinie}
\end{cases} . \label{3.53}
\end{equation}
Il résulte finalement de (\ref{3.52}) et (\ref{3.53}) qu'on a pour
toute place $v$ de $S$:
\begin{eqnarray*}
{\mid{f_{i_1 , i_2}}^{\!\!\!\!\!\!\!\!(j)}~\mid}_v ~:=~ {\mid
\left(\partial^{(i_1 , i_2)}P_j\right)({\underline x}_1 ,
{\underline x}_2) \mid}_v . \frac{1}{{\mid {\Delta({\underline
x}_1)}^{f(i_1)}
{\Delta({\underline x}_2)}^{f(i_2)}\mid}_v}~~~~~~~~~~~~~~~~~~~~ \\
\leq~ \begin{cases}
\exp\left\{2 m_v (i_1 + i_2) + 14 m_v + 2 h_v(\underline y)\right\} & \text{si $v$ est finie} \\
\exp\left\{(2 m_v + 13) (i_1 + i_2) + 14 m_v + 2 h_v(\underline y)
+ 46\right\} & \text{si $v$ est infinie}
\end{cases} ,
\end{eqnarray*}
ce qui n'est rien d'autre que la majoration de la proposition
\ref{c.43} pour les ${\mid f_{i_1 , i_2}^{(j)} \mid}_v ~ (v \in
S)$. La démonstration est achevée. $~~~~\blacksquare$\vspace{1mm}
\begin{corollary}[Immédiat]\label{c.44}
Pour toute place $v \in S$, les deux séries de la proposition
\ref{c.43} précédente sont absolument convergentes en valeur
absolue $v$-adique d\`es que:
$${\mid u_k \mid}_v < \begin{cases}
e^{- 2 m_v} & \text{si $v$ est finie} \\
e^{- 2 m_v - 13} & \text{si $v$ est infinie}
\end{cases} ~~~~ k = 1 , 2 .$$
En particulier, pour toute place $v \in S$, les deux séries
sus-citées convergent absolument en valeur absolue $v$-adique au
voisinage du point $(- x_1 , - x_2)$.
\end{corollary}
{\bf Démonstration.---} Etant donné une place $v \in S$, la
condition suffisante pour la convergence absolue $v$-adique des
deux séries de la proposition \ref{c.43} est claire et
entra{\sf\^\i}ne -d'apr\`es l'hypoth\`ese principale (\ref{3.51})
du théor\`eme \ref{c.40} et la propriété $ii)$ du §$2$ pour les
distances ${\mbox{dist}}_w ~ (w \in M_K)$- la convergence absolue
$v$-adique de ces m\^eme séries au voisinage du point $(- x_1 , -
x_2)$.  $~~~~\blacksquare$\vspace{1mm}
\subsection{Extrapolation:}
\begin{proposition}\label{c.45}
Si l'un des deux nombres $Q_1(\underline 0 , \underline 0)$ et
$Q_2(\underline 0 , \underline 0)$ est non nul, alors on a:
$$h(\mathbf y) ~\geq~ \frac{\epsilon}{2} \min\!\left\{h({\mathbf x}_1) , h({\mathbf x}_2)\right\} - 7 \eta - \frac{47}{2} .$$
\end{proposition}
{\bf Démonstration.---} Supposons que pour un $j \in \{1 , 2\}$ on
a: $Q_j(\underline 0 , \underline 0) \neq 0$. Dans ce cas, comme
le nombre $Q_j(\underline 0 , \underline 0)$ appartient \`a $K$
(car $Q_j$ est une forme \`a coefficients dans $K$ et $\underline
0 \in K^3$) alors ce dernier doit satisfaire la formule du
produit:
\begin{equation}
\sum_{v \in K} \frac{[K_v : {\mathbb Q}_v]}{[K : \mathbb Q]} \log
{\mid Q_j(\underline 0 , \underline 0) \mid}_v ~=~ 0 .
\label{3.54}
\end{equation}
Nous majorons maintenant astucieusement le membre de gauche de (\ref{3.54}) en fonction de $\epsilon , \eta , h({\mathbf x}_1) , h({\mathbf x}_2)$ et $h(\mathbf y)$. Pour ce faire nous distinguons les deux cas suivants: \\
$\underline{{\mbox{1}}^{\mbox{er}} cas:}$ (si $v \in S$) \\
Dans ce cas le corollaire \ref{c.44} nous dit que la série
numérique:
$$\tau^2(Q_j)({\underline x}_1 , {\underline x}_2)(- x_1 , - x_2) = \sum_{i_1 \geq 0 , i_2 \geq 0} f_{i_1 , i_2}^{(j)} (- x_1)^{i_1} (- x_2)^{i_2}$$
est absolument convergente en valeur absolue $v$-adique. Or,
d'apr\`es la définition m\^eme du monomorphisme $\tau$, cette
derni\`ere ne peut converger que vers $Q_j(\underline 0 ,
\underline 0)$, par conséquent on doit avoir:
\begin{eqnarray*}
Q_j(\underline 0 , \underline 0) & = & \sum_{i_1 \geq 0 , i_2 \geq
0} \!\!f_{i_1 , i_2}^{(j)} (- x_1)^{i_1} (- x_2)^{i_2} \\
 & = & \!\!\!\!\!\!\!\!\sum_{(i_1 , i_2) \in {\mathbb N}^2 \setminus \{(0 , 0)\}} \!\!\!\!\!\!\!\!\!\!\!f_{i_1 , i_2}^{(j)} (- x_1)^{i_1} (- x_2)^{i_2} ~~~~ (\mbox{car}~ f_{i_1 , i_2}^{(j)} = 0 ~\mbox{pour}~ (i_1 , i_2) = (0 , 0)) .
\end{eqnarray*}
D'o\`u, la majoration:
$${\mid Q_j(\underline 0 , \underline 0) \mid}_v \leq \left\{\!\!\!\begin{array}{ll}
\max_{(i_1 , i_2) \in {\mathbb N}^2 \setminus \{(0 , 0)\}}
\left\{{\mid f_{i_1 , i_2}^{(j)} \mid}_v.{{\mid x_1
\mid}_v}^{\!\!\!\!i_1}.{{\mid x_2 \mid}_v}^{\!\!\!\!i_2}\right\}
& \mbox{si}~ v ~ \mbox{est finie} \\
\sum_{(i_1 , i_2) \in {\mathbb N}^2 \setminus \{(0 , 0)\}}
\left\{{\mid f_{i_1 , i_2}^{(j)} \mid}_v.{{\mid x_1
\mid}_v}^{\!\!\!\!i_1}.{{\mid x_2 \mid}_v}^{\!\!\!\!i_2}\right\} &
\mbox{si}~ v ~ \mbox{est infinie}
\end{array}\right.\!\!.$$
En utilisant les estimations de la proposition \ref{c.43} pour les
${\mid f_{i_1 , i_2}^{(j)} \mid}_v ~ (i_1 , i_2 \in \mathbb N)$ et
en majorant les ${\mid x_k \mid}_v ~ (k = 1 , 2)$ -gr\^ace \`a
l'hypoth\`ese (\ref{3.51}) du théor\`eme \ref{c.40} et \`a la
propriété $ii)$ du §$2$ pour les distances ${\mbox{dist}}_w ~ (w
\in M_K)$- par:
$${\mid x_k \mid}_v < \left\{\!\!\!\begin{array}{ll}
e^{- \lambda_v \epsilon h({\mathbf x}_k) - 2 m_v} & \mbox{si}~ v
~\mbox{est
finie} \\
e^{- \lambda_v \epsilon h({\mathbf x}_k) - 2 m_v - 16} &
\mbox{si}~ v ~\mbox{est infinie}
\end{array}\right.~~~~ k = 1 , 2~;$$
la derni\`ere majoration de ${\mid Q_j(\underline 0 , \underline
0) \mid}_v$ entra{\sf\^\i}ne:
$$
{\mid Q_j(\underline 0 , \underline 0) \mid}_v ~\leq~ \exp\{14 m_v
+ 2 h_v(\underline y)\} \!\!\!\max_{(i_1 , i_2) \in {\mathbb N}^2
\setminus \{(0 , 0)\}} \!\!\! e^{- \lambda_v \epsilon \left\{i_1
h({\mathbf x}_1) + i_2 h({\mathbf x}_2)\right\}}~~~~~~~~~~
$$
si $v$ est finie et:
$$
{\mid Q_j(\underline 0 , \underline 0) \mid}_v ~\leq~ \exp\{14 m_v
+ 2 h_v(\underline y) + 46\} \!\!\!\!\!\!\sum_{(i_1 , i_2) \in
{\mathbb N}^2 \setminus \{(0 , 0)\}} \!\!\!\!\!\!\!\!\!e^{-
\lambda_v \epsilon \left\{i_1 h({\mathbf x}_1) + i_2 h({\mathbf
x}_2)\right\}} (e^{-3})^{i_1 + i_2}
$$
si $v$ est infinie. \\
Comme maintenant pour tout $(i_1 , i_2) \in {\mathbb N}^2
\setminus \{(0 , 0)\}$ on a: $i_1 h({\mathbf x}_1) + i_2
h({\mathbf x}_2) \geq (i_1 + i_2) \min\left\{h({\mathbf x}_1) ,
h({\mathbf x}_2)\right\} \geq \min\left\{h({\mathbf x}_1) ,
h({\mathbf x}_2)\right\}$ et que:
$$\sum_{(i_1 , i_2) \in {\mathbb N}^2 \setminus \{(0 , 0)\}} (e^{-3})^{i_1 + i_2} = \frac{1}{(1 - e^{-3})^2} - 1 < e ,$$ alors:
$${\mid Q_j(\underline 0 , \underline 0) \mid}_v
\leq \left\{\!\!\!\begin{array}{ll} \exp\{14 m_v + 2
h_v(\underline y)\} e^{- \lambda_v \epsilon \min\left\{h({\mathbf
x}_1) , h({\mathbf
x}_2)\right\}} & \!\!\mbox{si}~ v ~\mbox{est finie} \\
\exp\{14 m_v + 2 h_v(\underline y) + 47\} e^{- \lambda_v \epsilon
\min\left\{h({\mathbf x}_1) , h({\mathbf x}_2)\right\}} &
\!\!\mbox{si}~ v ~\mbox{est infinie}
\end{array}\right.\!\!\!.$$
En prenant finalement les logarithmes des deux membres de cette
derni\`ere inégalité, on obtient:
$$\mathbf{(e)} : \log {\mid Q_j(\underline 0 , \underline 0) \mid}_v \!\leq \!\begin{cases} \!
- \lambda_v \epsilon \min\{h({\mathbf x}_1) , h({\mathbf x}_2)\} + 14 m_v + 2 h_v(\underline y) & \!\!\!\!\text{si $v \!\!\nmid \!\!\infty$} \\
\! - \lambda_v \epsilon \min\{h({\mathbf x}_1) , h({\mathbf
x}_2)\} + 14 m_v + 2 h_v(\underline y) + 47 & \!\!\!\!\text{si $v
\!\!\mid \!\!\infty$}
\end{cases} \!.$$
$\underline{{\mbox{2}}^{\mbox{i\`eme}} \mbox{cas:}}$ (si $v \not\in S$) \\
Dans ce deuxi\`eme cas, on majore na{\sf\"\i}vement le nombre
$\log {\mid Q_j(\underline 0 , \underline 0) \mid}_v$. On remarque
que $Q_j(\underline 0 , \underline 0)$ est un coefficient de la
forme $Q_j$, donc on doit avoir:
\begin{equation*}
\begin{split}
\log {\mid Q_j(\underline 0 , \underline 0) \mid}_v &\leq~ h_v(Q_j) \\
&\leq
\begin{cases}
8 m_v + 2 h_v(\underline y) & \text{si $v$ est finie} \\
8 m_v + 2 h_v(\underline y) + 18 & \text{si $v$ est infinie}
\end{cases}
\end{split}
\end{equation*}
d'apr\`es le lemme \ref{c.41}. D'o\`u à fortiori:
$$\mathbf{(f)} : \log {\mid Q_j(\underline 0 , \underline 0) \mid}_v \leq
\begin{cases}
14 m_v + 2 h_v(\underline y) & \text{si $v$ est finie} \\
14 m_v + 2 h_v(\underline y) + 47 & \text{si $v$ est infinie}
\end{cases} .~~~~~~~~~~~~~~~~~~~~$$ \\
En majorant maintenant, pour toute place $v$ de $K$, le nombre
$\log {\mid Q_j(\underline 0 , \underline 0) \mid}_v$ en utilisant
$\mathbf{(e)}$ si $v \in S$ et $\mathbf{(f)}$ si $v \in S$ et en
tenant compte des identités:
\begin{align}
&~\!\sum_{v \in S} \frac{[K_v : {\mathbb Q}_v]}{[K : \mathbb Q]} \lambda_v = 1 ~~~~\text{(par hypoth\`ese),} ~~~~\!\!\sum_{v \in M_K} \frac{[K_v : {\mathbb Q}_v]}{[K : \mathbb Q]} m_v = \eta , \notag \\
&\sum_{v \in M_K} \frac{[K_v : {\mathbb Q}_v]}{[K : \mathbb Q]}
h_v(\underline y) = h(\mathbf y) ~~~~\text{et}
~~~~~~~~~~~~\!\sum_{v \in M_{K}^{\infty}} \frac{[K_v : {\mathbb
Q}_v]}{[K : \mathbb Q]} = 1 ; \notag
\end{align}
on a la majoration:$$\sum_{v \in M_K} \frac{[K_v : {\mathbb
Q}_v]}{[K : {\mathbb Q}]} \log {\mid Q_j(\underline 0 , \underline
0) \mid}_v \leq - \epsilon \min\{h({\mathbf x}_1) , h({\mathbf
x}_2)\} + 14 \eta + 2 h(\mathbf y) + 47 .$$ Or, d'apr\`es
(\ref{3.54}) le nombre $\sum_{v \in M_K} \frac{[K_v : {\mathbb
Q}_v]}{[K : {\mathbb Q}]} \log {\mid Q_j(\underline 0 , \underline
0) \mid}_v$ est nul. La derni\`ere inégalité est donc équivalente
\`a:
$$- \epsilon \min\{h({\mathbf x}_1) , h({\mathbf x}_2)\} + 14 \eta + 2 h(\mathbf y) + 47 \geq 0 ;$$ ce qui donne finalement:
$$h(\mathbf y) \geq \frac{\epsilon}{2} \min\{h({\mathbf x}_1) , h({\mathbf x}_2)\} - 7 \eta - \frac{47}{2} .$$
La démonstration est achevée.  $~~~~\blacksquare$\vspace{1mm}
\subsection{Démonstration du théor\`eme \ref{c.40}:}
Nous procédons par l'absurde. Supposons qu'il existe un couple
$({\mathbf x}_1 , {\mathbf x}_2)$ de points de $E(K)$,
satisfaisant toutes les hypoth\`eses du théor\`eme \ref{c.40} mais
ne satisfaisant pas sa conclusion, c'est-à-dire qu'on a:
$\widehat{h}({\mathbf x}_1) \leq \widehat{h}({\mathbf x}_2) < (1 +
\theta) \widehat{h}({\mathbf x}_1)$. On a donc, en fonction de la
norme de Néron-Tate:
$$\mid{\mathbf x}_1\mid \leq \mid{\mathbf x}_2\mid < (1
+ \theta)^{\frac{1}{2}} \mid{\mathbf x}_1\mid < (1 +
\frac{\theta}{2}) \mid{\mathbf x}_1\mid ,$$ d'o\`u:
\begin{equation}
\mid{\mathbf x}_1\mid \leq \mid{\mathbf x}_2 \mid < \left(1 +
\frac{\theta}{2}\right) \mid{\mathbf x}_1\mid . \label{3.55}
\end{equation}
Gr\^ace \`a (\ref{3.55}) et à l'hypoth\`ese (\ref{3.49}) du
théor\`eme \ref{c.40}, on a la série d'inégalités suivante:
\begin{eqnarray*}
{\mid \mathbf y \mid}^2 = {\mid{\mathbf x}_1 - {\mathbf x}_2
\mid}^2 & = & {\mid {\mathbf x}_2 \mid}^2 + {\mid{\mathbf x}_1
\mid}^2 - 2 <{\mathbf x}_1 ,
{\mathbf x}_2> \\
 & = & \left(\mid{\mathbf x}_2 \mid - \mid{\mathbf x}_1 \mid\right)^2 + 2
 \left(\mid{\mathbf x}_1 \mid.\mid{\mathbf x}_2 \mid - <{\mathbf x}_1 , {\mathbf
 x}_2>\right) \\
 & = & \left(\mid{\mathbf x}_2 \mid - \mid{\mathbf x}_1 \mid\right)^2 + 2
 \left(1 - \cos({\mathbf x}_1 , {\mathbf x}_2)\right) \mid{\mathbf
 x}_1 \mid.\mid{\mathbf x}_2 \mid \\
 & < & \left(\frac{\theta}{2} \mid{\mathbf x}_1 \mid\right)^2 + 2
 \left(\frac{\beta}{4}\right) \mid{\mathbf x}_1 \mid.\left(1 +
 \frac{\theta}{2}\right)\mid{\mathbf x}_1 \mid ,
\end{eqnarray*}
c'est-\`a-dire:
$$\widehat{h}(\mathbf y) < \left(\frac{\beta}{2} + \frac{\theta^2}{4} + \frac{\beta \theta}{4}\right) \widehat{h}
({\mathbf x}_1) .$$ Or $\frac{\beta}{2} + \frac{\theta^2}{4} +
\frac{\beta \theta}{4} = \frac{\epsilon}{4} + \frac{\epsilon}{36}
+ \frac{1}{24} \epsilon \sqrt{\epsilon} \leq (\frac{1}{4} +
\frac{1}{36} + \frac{1}{24}) \epsilon \leq \frac{\epsilon}{3}$,
d'o\`u:
\begin{equation}
\widehat{h}(\mathbf y) < \frac{\epsilon}{3} \widehat{h}({\mathbf
x}_1) . \label{3.56}
\end{equation}
Maintenant il vient:
\begin{equation}
\begin{split}
h(\mathbf y) &\leq~ \widehat{h}(\mathbf y) + \frac{3}{2} \eta +
8 \\
&<~ \frac{\epsilon}{3} \widehat{h}({\mathbf x}_1) +
\frac{3}{2} \eta + 8 \\
&\leq~ \frac{\epsilon}{2} \left(\widehat{h}({\mathbf x}_1) -
\frac{3}{4} \eta - 5\right) - 7 \eta - \frac{47}{2} \label{3.57}
\end{split}
\end{equation}
o\`u, dans cette série d'inégalités, la premi\`ere inégalité suit
du théor\`eme \ref{c.14} du formulaire, la deuxi\`eme de
(\ref{3.56}) et la troisi\`eme de l'hypoth\`ese (\ref{3.50}) du
théor\`eme \ref{c.40}. Finalement, comme on a d'apr\`es le
théor\`eme \ref{c.14} du formulaire: $\widehat{h}({\mathbf x}_1) -
\frac{3}{4} \eta - 5 = \min\{\widehat{h}({\mathbf x}_1) ,
\widehat{h}({\mathbf x}_2)\} - \frac{3}{4} \eta - 5 \leq
\min\{h({\mathbf x}_1) , h({\mathbf x}_2)\}$, on déduit de
(\ref{3.57}):
$$h(\mathbf y) < \frac{\epsilon}{2} \min\{h({\mathbf x}_1) , h({\mathbf x}_2)\} - 7 \eta - \frac{47}{2} .$$
Cette derni\`ere inégalité entra{\sf\^\i}ne, d'apr\`es la
proposition \ref{c.45}, que nos deux formes $Q_1$ et $Q_2$ de
$K[{\underline X}_1 , {\underline X}_2]$ s'annulent, toutes les
deux, en $(\underline 0 , \underline 0)$, ce qui implique que nos
deux formes $P_1$ et $P_2$ de $K[\underline Y]$ s'annulent, toutes
les deux, en $\underline 0$. Ce dernier fait montre qu'on a: $-
\mathbf y = \mathbf 0$. En effet, le point $- \mathbf y = \mathbf
0 - \mathbf y$ de $E$ peut \^etre représenté dans ${\mathbb P}_2$
par le syst\`eme de coordonnées projectives $\left(D_0(\underline
0 , \underline y) : D_1(\underline 0 , \underline y) ,
D_2(\underline 0 , \underline y)\right) = \left(P_1(\underline 0)
: D_1(\underline 0 , \underline y) : P_2(\underline 0)\right)$, or
lorsque $P_1(\underline 0) = P_2(\underline 0) = 0$, cette
représentation est identique à $(0 : 1 : 0) = \underline 0$, par
conséquent le point $- \mathbf y$ sera identique à l'origine
$\mathbf 0$ de $E$. Enfin $- \mathbf y = \mathbf 0$ entraine
$\mathbf y = \mathbf 0$, puis ${\mathbf x}_1 = {\mathbf x}_2$, ce
qui donne la contradiction cherchée avec les hypoth\`eses et
ach\`eve cette démonstration.
\section{Démonstration des théor\`emes \ref{c.1}, \ref{c.2} et \ref{c.3}}
Notons dans toute la suite, $R$ l'expression dépendant de
$\epsilon , m$ et $\eta$ suivante:
$$
R ~\!:=~\! \left(2904 m\right)^m \left[\left\{55 \epsilon^{-
\frac{m}{m - 1}} + m \epsilon^{- 1}\right\} \eta + \left\{272
\epsilon^{- \frac{m}{m - 1}} + 2 m \epsilon^{- 1}\right\}\right] ,
$$
et désignons respectivement par $(I)$ et $(II)$ les deux
syst\`emes d'inégalités simultanées (en $\mathbf x \in E(K)$):
\begin{align}
{\mbox{dist}}_v(\mathbf x , \mathbf 0) &<~  e^{- \lambda_v \epsilon h(\mathbf x) - 2 m_v - c_v} ~~~~~~~~~~ (v \in S) \tag{$I$} \\
{\mbox{dist}}_v(\mathbf x , \mathbf 0) &<~ e^{- \lambda_v \epsilon
h(\mathbf x) - \lambda_v R - 2 m_v - c_v} ~~~~~ (v \in S) .
\tag{$II$}
\end{align}
Afin de démontrer nos théor\`emes principaux sur le décompte des
points de $E(K)$ satisfaisant le syst\`eme d'inégalités $(I)$,
nous allons estimer le nombre de tels points qui sont de hauteurs
$> R$, lesquels sont appelés ci-dessous ``points de hauteurs assez
grandes'', puis le nombre de ces points qui sont de hauteurs $\leq
R$, lesquels sont nommés ci-dessous ``points de hauteurs assez
petites''. \\
Notons que l'estimation du nombre de points de $E(K)$,
satisfaisant le syst\`eme d'inégalités $(I)$ et qui sont de
hauteurs assez grandes s'obtient à partir des théor\`emes
\ref{c.39} et \ref{c.40}, alors que pour estimer le nombre de
points de $E(K)$ qui sont de hauteurs assez petites, nous
disposons de deux méthodes différentes qui sont:
\subsubsection*{La premi\`ere} consiste à remplacer le syst\`eme d'inégalités $(I)$ par le syst\`eme d'inégalités $(II)$ de sorte que seul le point $\mathbf 0$ puisse satisfaire $(II)$ tout en étant de hauteur $\leq R$.
\subsubsection*{La deuxi\`eme} consiste à estimer le nombre de tous
les points de $E(K)$ de hauteurs $\leq R$ sans tenir compte du
syst\`eme d'inégalités $(I)$. Cette deuxi\`eme méthode est moins
élémentaire que la premi\`ere dans le sens o\`u elle nécessite
l'utilisation d'un résultat fournissant une borne inférieure
(strictement positive) pour l'ensemble des hauteurs de Néron-Tate
des points non de torsion de $E(K)$ et d'un résultat fournissant
une majoration pour le nombre de points de torsion de $E(K)$. De
tels résultats se trouvent respectivement dans [H-S1] et [Me].
\subsection{Décompte des points de hauteurs assez grandes:}
Le décompte des points de $E(K)$ satisfaisant le syst\`eme
d'inégalités $(I)$ et qui sont de hauteurs $> R$ est donné par le
théor\`eme suivant:
\begin{theorem}\label{c.46}
Soit $E$ une courbe elliptique définie sur un corps de nombres $K$
de degré $D$, plongée dans ${\mathbb P}_2$ à la Weierstrass,
d'équation projective $Y^2 Z = 4 X^3 - g_2 X Z^2 - g_3 Z^3 ~ (g_2
, g_3 \in K)$ et d'élément neutre (en tant que groupe) le point à
l'infini $\mathbf 0$ représenté dans ${\mathbb P}_2$ par les
coordonnées projectives $(0 : 1 : 0)$. On munit $E$ de la hauteur
de Néron-Tate définie au §$13$, notée $\widehat{h}$. Soient aussi
$S$ un ensemble fini de places sur $K$, $m_v ~ (v \in M_K)$,
$\eta$ les réel positifs définis au §$2$ et ${(\lambda_v)}_{v \in
S}$ une famille de réels positifs satisfaisant:
$$\sum_{v \in S} \frac{[K_v : {\mathbb Q}_v]}{[K : \mathbb Q]} \lambda_v = 1 .$$ Soient enfin $\epsilon$ un réel strictement positif et $m$ un entier $\geq 2$. On a:
\begin{equation*}
\begin{split}
\sharp \left\{\mathbf x \in E(K) , ~\mathbf x ~\text{satisfait $(I)$ et}~ \widehat{h}(\mathbf x) > R\right\} & \\
&\!\!\!\!\!\!\!\!\!\!\!\!\!\!\!\!\!\!\!\!\!\!\!\!\!\!\!\!\!\!\!\!\!\!\!\!\!\!\!\!\!\!\!\!\!\!\!\!\!\!\!\!\!\!\!\!\!\!\!\!\!\!\leq~
4 \epsilon^{- \frac{1}{2}} \left[m (m - 1) (\log m + 8.8) + m
\nb{\log{\epsilon}} \right].\left(499 \epsilon^{- \frac{m}{2 (m -
1)}}\right)^r ,
\end{split}
\end{equation*}
o\`u $r$ désigne le rang du groupe de Mordell-Weil de $E(K)$.
\end{theorem}
{\bf Démonstration.---} Remarquons d'abord que la conclusion du
théor\`eme \ref{c.40} reste toujours valable lorsqu'on remplace
les hypoth\`eses de ce dernier par celles du théor\`eme \ref{c.39}
pour $m = 2$. Notons toujours par $\alpha$ le param\`etre choisi
au sous-paragraphe $10.2$, par $\theta$ le param\`etre du
théor\`eme \ref{c.40} et désignons par $r$ le rang du groupe de
Mordell-Weil de $E(K)$. Notre procédé consiste à recouvrir
l'espace euclidien $E(K) \otimes \mathbb R \simeq {\mathbb R}^r$
par un nombre fini de c\^ones de centres $\mathbf 0$ et d'angles
$\leq \arccos(1 - \alpha/4)$, ensuite en utilisant les deux
théor\`emes \ref{c.39} et \ref{c.40} nous aurons un décompte des
points de $E(K)$ satisfaisant le syst\`eme d'inégalités $(I)$ et
qui sont de hauteurs de Néron-Tate $> R$, dans chacun de ces
petits c\^ones. Finalement, en multipliant le décompte obtenu dans
un tel petit c\^one par le nombre de c\^ones que comporte notre
recouvrement, on aura l'estimation du théor\`eme \ref{c.46}.
D'apr\`es le lemme \ref{b.21} de l'appendice, le
nombre minimal de c\^ones d'angles $\leq \arccos(1 - \alpha/4)$
suffisant pour recouvrir l'espace euclidien $E(K) \otimes \mathbb
R \simeq {\mathbb R}^r$ est majoré par $\left(1 + \frac{8}{\sqrt{2
\alpha}}\right)^{\!\!r} < \left(499 \epsilon^{- \frac{m}{2 (m
-1)}}\right)^{\!\!r}$ (en vertu du choix du param\^etre $\alpha$
effectué au sous-paragraphe $9.2$). Dans un tel recouvrement,
supposons qu'on a $\ell$ points de $E(K) $$(\ell \geq 1)$
satisfaisant le syst\`eme d'inégalités $(I)$, appartenant à un
m\^eme c\^one et qui sont de hauteurs de Néron-Tate $> R$. Notons
par ${\mathbf x}_1 , {\mathbf x}_2 , \dots , {\mathbf x}_{\ell}$
ces $\ell$ points, ordonnés selon l'ordre croissant de leurs
hauteurs de Néron-Tate:
$$R ~<~ \widehat{h}({\mathbf x}_1) ~\leq~ \widehat{h}({\mathbf x}_2) ~\leq~ \dots ~\leq~ \widehat{h}({\mathbf x}_{\ell}) .$$
En effectuant la division euclidienne de l'entier positif $\ell -
1$ sur l'entier positif non nul $m - 1$:
$$\ell - 1 ~=~ k (m - 1) + r ~,~ 0 \leq r \leq m - 2  ~\!;$$
considérons parmis les $\ell$ points précédents, seulement les
points: ${\mathbf x}_1 , {\mathbf x}_2 , \dots , {\mathbf x}_{k(m
- 1) + 1}$ et soient parmis ces derniers ${\mathbf y}_1 , {\mathbf
y}_2 , \dots , {\mathbf y}_m$ les $m$ points:
$${\mathbf y}_i ~:=~ {\mathbf x}_{k (i - 1) + 1} ~~~~~~~~~~ i = 1 , \dots , m .$$
On a évidemment aussi:
$$R ~<~ \widehat{h}({\mathbf y}_1) ~\leq~ \widehat{h}({\mathbf y}_2) ~\leq~ \dots ~\leq~ \widehat{h}({\mathbf y}_m) .$$
Maintenant, d'une part le $m$-uplet $({\mathbf y}_1 , {\mathbf
y}_2 , \dots , {\mathbf y}_m)$ de $E^m(K)$ ainsi considéré
satisfait clairement toutes les hypoth\`eses du théor\`eme
\ref{c.39}, donc -d'apr\`es ce dernier- il doit exister un entier
positif $j \in \{2 , \dots , m\}$ tel que l'on ait:
\begin{equation}
\widehat{h}({\mathbf y}_j) ~<~ \sqrt{2} (m - 1) (2904 m)^m
\left(\frac{1}{\epsilon}\right)^{\!\!\frac{m}{m - 1}}
\widehat{h}({\mathbf y}_{j - 1}) \label{3.58}
\end{equation}
et d'autre part -d'apr\`es la remarque faite au début de cette
démonstration- les hypoth\`ese du théor\`eme \ref{c.40} sont
clairement satisfaites pour tout couple $({\mathbf x}_n , {\mathbf
x}_{n + 1}) ,~ n = 1 , \dots , \ell - 1$, donc d'apr\`es le
théor\`eme \ref{c.40} on a pour tout $n \in \{1 , \dots , \ell -
1\} : \widehat{h}({\mathbf x}_{n + 1}) \geq (1 + \theta)
\widehat{h}({\mathbf x}_n)$ et plus généralement pour tout $n_1 ,
n_2 ~(n_1 \leq n_2)$ dans $\{1 , \dots , \ell - 1\} :
\widehat{h}({\mathbf x}_{n_2}) \geq (1 + \theta)^{n_2 - n_1}
\widehat{h}({\mathbf x}_{n_1})$. En particulier pour un entier $j
\in \{2 , \dots , m\}$ satisfaisant (\ref{3.58}) on a:
$$\widehat{h}({\mathbf y}_j) ~=~ \widehat{h}({\mathbf x}_{k(j - 1) + 1}) ~\geq~ (1 + \theta)^k \widehat{h}({\mathbf x}
_{k(j - 2) + 1}) ~=~ (1 + \theta)^k \widehat{h}({\mathbf y}_{j -
1}) ,$$ c'est-à-dire:
\begin{equation}
\widehat{h}({\mathbf y}_j) ~\geq~ (1 + \theta)^k
\widehat{h}({\mathbf y}_{j - 1}) . \label{3.59}
\end{equation}
Comme les points ${\mathbf y}_1 , \dots , {\mathbf y}_m$ sont de
hauteurs de Néron-Tate suppérieures à $R$ (donc non nulles), les
deux relations (\ref{3.58}) et (\ref{3.59}) entra{\sf\^\i}nent:
$$(1 + \theta)^k ~<~ \sqrt{2} (m - 1) (2904 m)^m \epsilon^{- \frac{m}{m - 1}} ,$$
$$\mbox{ce qui donne:}~~~~~~~~~~~~~~~~~~~ k ~<~ \frac{\log \left[\sqrt{2} (m - 1) (2904 m)^m \epsilon^{- \frac{m}{m - 1}}\right]}{\log(1 + \theta)} .~~~~~~~~~~~~~~~~~~~~~~$$
Puis, comme $\ell - 1 = k(m - 1) + r$ et $r \leq m - 2$, on a:
$\ell - 1 \leq k(m - 1) + m - 2$, c'est-à-dire $\ell \leq (m -
1)(k + 1)$, d'o\`u, d'apr\`es la majoration précédente pour $k$:
$$\ell ~\leq~ (m - 1) \left(\frac{\log \left[\sqrt{2} (m - 1) (2904 m)^m
\epsilon^{- \frac{m}{m - 1}}\right]}{\log(1 + \theta)} + 1\right)
.$$ En majorant maintenant $\sqrt{2} (m - 1)$ par $2^m$ et en
minorant $\log(1 + \theta)$ par $\frac{1}{4} \sqrt{\epsilon}$ (qui
vient simplement du fait que la fonction réelle $x \mapsto
\frac{\log(1 + x)}{x}$ est décroissante sur $]0 , + \infty[$ et
que $\theta := \frac{1}{3} \sqrt{\epsilon} \leq \frac{1}{3}$) on
a:
\begin{equation*}
\begin{split}
\frac{\log \left[\sqrt{2} (m - 1) (2904 m)^m \epsilon^{- \frac{m}{m - 1}}\right]}{\log(1 + \theta)} + 1 &\leq 4 \epsilon^{- \frac{1}{2}} \log \!\left[(5808 m)^m \epsilon^{- \frac{m}{m - 1}}\right] + 1 \\
&\!\!\!\!\!\!\leq 4 \epsilon^{- \frac{1}{2}} \left\{m \left(\log m + \log 5808\right) + \frac{m}{m - 1} \!\mid \!\log \epsilon \!\mid \right\} + 1 \\
&\!\!\!\!\!\!\leq 4 \epsilon^{- \frac{1}{2}} \left[m \left(\log m
+ 8.8\right) + \frac{m}{m - 1} \!\mid \!\log \epsilon \!\mid
\right] .
\end{split}
\end{equation*}
D'o\`u la majoration finale de $\ell$:
$$\ell ~\leq~ 4 \epsilon^{- \frac{1}{2}} \left[m(m - 1)(\log m + 8.8) + m \!\mid \!\log \epsilon \!\mid \right] .$$
On vient ainsi de montrer que le nombre de points de $E(K)$
satisfaisant le syst\`eme d'inégalité $(I)$, qui sont de hauteur
de Néron-Tate $> R$ et contenus dans un m\^eme c\^one de notre
recouvrement, ne peut dépasser la quantité $4 \epsilon^{-
\frac{1}{2}} [m(m - 1)(\log m + 8.8) +$ \\ $m \!\mid \!\log
\epsilon \!\mid]$. Il ne reste qu'à multiplier cette derni\`ere
quantité par l'estimation $\!(\!499 \epsilon^{\!- \frac{m}{2(m -
1)}}\!)^r\!\!$ du nombre de c\^ones constituants notre
recouvrement de l'espace euclidien $E(K) \otimes \mathbb R$ pour
aboutir finalement à l'estimation du théor\`eme \ref{c.46} pour le
nombre de points de $E(K)$ satisfaisant le syst\`eme d'inégalité$
(I)$ et qui sont de hauteurs de Néron-Tate $> R$. La démonstration
est achevée.  $~~~~\blacksquare$\vspace{1mm}
\begin{corollary}[Immédiat]\label{c.47}
Sous les hypoth\`eses du théor\`eme \ref{c.46} précédent, on a: \\ \\
$\displaystyle \sharp \left\{\mathbf x \in E(K) , ~\mathbf x ~\text{satisfait $(II)$ et}~ \widehat{h}(\mathbf x) > R\right\}$ \\ \\
$\displaystyle ~~~~~~~~~~~~~~~~~~~~\leq~ 4 \epsilon^{-
\frac{1}{2}} \left[m (m - 1) (\log m + 8.8) + m \!\mid \!\log
\epsilon \!\mid \right].\left(499
\epsilon^{- \frac{m}{2 (m - 1)}}\right)^r$, \\ \\
o\`u $r$ désigne le rang du groupe de Mordell-Weil de $E(K)$.
\end{corollary}
{\bf Démonstration.---} Il suffit de remarquer que le syst\`eme
d'inégalités $(II)$ entra{\sf\^\i}ne trivialement le syst\`eme
d'inégalité $(I)$.  $~~~~\blacksquare$\vspace{1mm}
\subsection{Décompte des points de hauteurs assez petites:}
\subsubsection{Premi\`ere méthode:}
Comme déjà dit, elle consiste à montrer que l'ensemble des points
de $E(K)$ satisfaisant l'inégalité $(II)$ et qui sont de hauteurs
de Néron-Tate $\leq R$, est réduit au singleton $\{\mathbf 0\}$.
C'est le théor\`eme suivant:
\begin{theorem}\label{c.48}
Sous les hypoth\`eses du théor\`eme \ref{c.46}, on a:
$$\displaystyle \sharp \left\{\mathbf x \in E(K) , ~\mathbf x ~\text{satisfait $(II)$ et}~ \widehat{h}(\mathbf x) \leq
R\right\} ~=~ \left\{\mathbf 0\right\} .$$
\end{theorem}
{\bf Démonstration.---} Nous procédons par l'absurde, c'est-à-dire
que nous supposons qu'il existe un point $\mathbf x$ de $E(K)$
différent de l'origine $\mathbf 0$ de $E$, de hauteur de
Néron-Tate $\widehat{h}(\mathbf x) \leq R$ et satisfaisant le
syst\`eme d'inégalités $(II)$, c'est-à-dire satisfaisant:
$${\mbox{dist}}_v(\mathbf x , \mathbf 0) ~<~ e^{- \lambda_v \epsilon h(\mathbf x) - \lambda_v R - 2 m_v - c_v} ~<~ 1 ~~~~~~~~~~ \text{pour $v \in S$} .$$
Ces inégalités entra{\sf\^\i}nent -d'apr\`es la propriété $ii)$ du
§$2$ pour les distances ${\mbox{dist}}_w ~ (w \in M_K)$- que
$\mathbf x \in E(K) \setminus \{Y = 0\}$ et en désignant par
$\underline x = (x , 1 , z) \in K^3$ son représentant d'ordonnée
$Y = 1$ dans ${\mathbb P}_2$, on a:
$${\mbox{dist}}_v(\mathbf x , \mathbf 0) ~=~ \max\left({\mid x \mid}_v , {\mid z \mid}_v\right) ~~~~~~~~~~ \forall v \in S .$$
Ceci nous permet de déduire du syst\`eme d'inégalité $(II)$, le
syst\`eme d'inégalités:
$$\max\left({\mid x \mid}_v , {\mid z \mid}_v\right) ~<~ e^{- \lambda_v \epsilon h(\mathbf x) - \lambda_v R - 2 m_v - c_v} ~~~~~~~~~~ \text{pour $v \in S$} .$$ En reportant toutes ces inégalités dans la somme:
$$\sum_{v \in S} \frac{[K_v : {\mathbb Q}_v]}{[K : \mathbb Q]} \log \max\left({\mid x \mid}_v , {\mid z \mid}_v\right)$$ et en utilisant l'hypoth\`ese $\sum_{v \in S} \frac{[K_v : {\mathbb Q}_v]}{[K : \mathbb Q]} \lambda_v = 1$, on obtient:
\begin{equation}
\begin{split}
\sum_{v \in S} \frac{[K_v : {\mathbb Q}_v]}{[K : \mathbb Q]} \log \max\left({\mid x \mid}_v , {\mid z \mid}_v \right) &< - \epsilon h(\mathbf x) - R - \!\sum_{v \in S} \frac{[K_v : {\mathbb Q}_v]}{[K : \mathbb Q]} (2 m_v + c_v) \\
&\leq - \epsilon h(\mathbf x) - R . \label{3.60}
\end{split}
\end{equation}
Par ailleurs, on a par définition:
$$\sum_{v \in M_K} \frac{[K_v : {\mathbb Q}_v]}{[K : \mathbb Q]} \log \max\left({\mid x \mid}_v , {\mid z \mid}_v\right) ~\geq~ 0 ,$$
ce qui donne:
\begin{equation*}
\begin{split}
\sum_{v \in S} \frac{[K_v : {\mathbb Q}_v]}{[K : \mathbb Q]} \log \max\left({\mid x \mid}_v , {\mid z \mid}_v \right) & \geq - \sum_{v \not\in S} \frac{[K_v : {\mathbb Q}_v]}{[K : \mathbb Q]} \log\max\left({\mid x \mid}_v , {\mid z \mid}_v \right) \\
& \!\!\!\!\!\!\!\!\!\!\!\!\!\!\!\!\!\!\!\!\!\!\!\!\!\!\!\!\!\!\!\!\!\!\!\geq - \sum_{v \not\in S} \frac{[K_v : {\mathbb Q}_v]}{[K : \mathbb Q]} \log\max\left({\mid x \mid}_v , 1 , {\mid z \mid}_v \right) \\
&
\!\!\!\!\!\!\!\!\!\!\!\!\!\!\!\!\!\!\!\!\!\!\!\!\!\!\!\!\!\!\!\!\!\!\!\geq
- \sum_{v \in M_K} \frac{[K_v : {\mathbb Q}_v]}{[K : \mathbb Q]}
\log\max\left({\mid x \mid}_v , 1 , {\mid z \mid}_v \right) = -
h(\mathbf x) .
\end{split}
\end{equation*}
D'o\`u:
\begin{equation}
\sum_{v \in S} \frac{[K_v : {\mathbb Q}_v]}{[K : \mathbb Q]} \log
\max\left({\mid x \mid}_v , {\mid z \mid}_v\right) \geq -
h(\mathbf x) . \label{3.61}
\end{equation}
En comparant (\ref{3.60}) et (\ref{3.61}) on obtient:
$$- \epsilon h(\mathbf x) - R > - h(\mathbf x) ,$$
c'est-à-dire: $\displaystyle ~~~~~~~~~~ h(\mathbf x) > \frac{R}{1 - \epsilon} > (1 + \epsilon) R = R + \epsilon R$. \\
On conclut, gr\^ace au théor\`eme \ref{c.14} du formulaire:
\begin{equation*}
\begin{split}
\widehat{h}(\mathbf x) & \geq h(\mathbf x) - \frac{3}{2} \eta - 8 \\
& > R + \epsilon R - \frac{3}{2} \eta - 8 \\
& > R ~~~~~~~~~~ \left(\text{car $R > \frac{1}{\epsilon}
(3\!/\!2.\eta + 8)$}\right) .
\end{split}
\end{equation*}
Ce qui donne la contradiction cherchée et ach\`eve cette
démonstration.  $~~~~\blacksquare$\vspace{1mm}
\subsubsection{Deuxi\`eme méthode:}
Comme déjà dit, elle consiste simplement à estimer le nombre de
points de $E(K)$ qui sont de hauteurs de Néron-Tate $\leq R$, sans
tenir compte du syst\`eme d'inégalités $(I)$. On dispose plus
généralement d'un lemme estimant le nombre de points (à
coordonnées dans un corps de nombre donné $K$) d'une variété
abelienne $A$ (définie sur $K$) qui sont de hauteurs de Néron-Tate
majorées par une constante positive donnée $R$. Notons que cette
estimation dépend de deux quantités liées à $A$: la premi\`ere est
la plus petite valeur non nulle des hauteurs des points de $A(K)$,
qu'on désigne par ${\widehat{h}}_{\min}$ et la deuxi\`eme est le
nombre de points de torsion de $A(K)$, qu'on désigne par $\sharp
{A(K)}_{\rm{tor}}$. Dans notre cas ($A = E$ est une courbe
elliptique), afin d'avoir un résultat compl\`etement explicite, on
se ref\`ere à [Me] pour majorer $\sharp {A(K)}_{\rm{tor}}$ et à
[H-S1] pour minorer ${\widehat{h}}_{\min}$, en fonction des
invariants habituels de la courbe elliptique $E$. Notre lemme est
le suivant:
\begin{lemma}\label{c.49}
Soit $A$ une variété abélienne définie sur un corps de nombres $K$
et munie d'une hauteur de Néron-Tate $\widehat{h}$ et soit $r$ le
rang du groupe de Mordell-Weil $A(K)$ qu'on suppose non nul.
Alors, pour tout réel positif $R$ on a:
$$\sharp\left\{\mathbf x \in A(K) ~/~ \widehat{h}(\mathbf x) \leq R\right\}
~\leq~ \sharp {A(K)}_{\rm{tor}}~.\!\left(1 + \sqrt{\frac{4
R}{{\widehat{h}}_{\min}}}\right)^{\!\!r} ,$$ o\`u
${A(K)}_{\rm{tor}}$ désigne le groupe fini des points de torsion
de $A(K)$ et ${\widehat{h}}_{\min}$ désigne la plus petite valeur
des hauteurs des points de non torsion de $A(K)$.
\end{lemma}
{\bf Démonstration.---} Soient $\Gamma$ le groupe de Mordell-Weil
de $A(K)$ et ${\mathbf p}_1 , \dots , {\mathbf p}_r$ des
générateurs de $\Gamma$. Le groupe de Mordell-Weil $A(K)$ s'écrit
alors:
$$A(K) ~=~ {A(K)}_{\rm{tor}} \oplus {\mathbf p}_1 \mathbb Z \oplus \dots \oplus {\mathbf p}_r \mathbb Z .$$
Posons, pour un réel positif donné $R$:
$$\Gamma_R ~:=~ \left\{\mathbf x \in {\mathbf p}_1 \mathbb Z \oplus \dots \oplus {\mathbf p}_r \mathbb Z ~/~ \mathbf x \neq \mathbf 0 ~\mbox{et}~ \widehat{h}(\mathbf x) \leq R\right\} .$$
Il est clair qu'on a:
$$\sharp\left\{\mathbf x \in A(K) ~/~ \widehat{h}(\mathbf x) \leq R\right\} ~=~ \sharp {A(K)}_{\rm{tor}} . \sharp \Gamma_R .$$
Pour aboutir à la conclusion du lemme \ref{c.49}, on doit alors
montrer qu'on a:
$$\sharp \Gamma_R ~\leq~ \left(1 + \sqrt{\frac{4 R}{{\widehat{h}}_{\min}}}\right)^{\!\!r} .$$
Nous introduisons pour cela l'entier $N \geq 1$ défini par:
$$N ~:=~ \left[1 + \sqrt{\frac{4 R}{{\widehat{h}}_{\min}}}\right] ~~~~\text{(o\`u $[.]$ désigne la partie enti\`ere)}$$ et nous considérons l'application $\cl\rest{\Gamma_R}$ de $\Gamma_R$ dans $\Gamma/N \Gamma$, restriction de l'homomorphisme surjectif $\mbox{cl}$ de $\Gamma$ dans $\Gamma/N \Gamma$ associant à chaque point du groupe $\Gamma$ sa classe modulo le groupe $N \Gamma$:
\[
\begin{array}{rcl}
\Gamma_R \subset \Gamma & \stackrel{\cl\rest{\Gamma_R}}{\longrightarrow} & \Gamma/N \Gamma \\
\mathbf x & \longmapsto & \mbox{cl}(\mathbf x) ~\mbox{mod}~ N
\Gamma
\end{array}.
\]
Le choix de l'entier $N$ pousse l'application $\cl\rest{\Gamma_R}$
à \^etre injective. En effet, si ceci n'était pas le cas, il
devrait exister au moins un couple $(\mathbf x , \mathbf y)$ de
${\Gamma_R}^{2}$ tel que l'on ait $\mathbf x \neq \mathbf y$ et
$\mbox{cl}(\mathbf x) = \mbox{cl}(\mathbf y)$. Pour un tel couple,
le point $\mathbf x - \mathbf y$ vérifie: $\mbox{cl}(\mathbf x -
\mathbf y) = \mbox{cl}(\mathbf 0)$ et $\mathbf x - \mathbf y \neq
\mathbf 0$, ce qui revient à dire que $\mathbf x - \mathbf y \in
(N \Gamma)^*$ ou en d'autres termes que le point $\mathbf x -
\mathbf y$ s'écrit sous la forme:
$$\mathbf x - \mathbf y ~=~ N \mathbf z ~~~~\text{avec $\mathbf z \in \Gamma \!\setminus \!\{\mathbf 0\}$} .$$
Comme $\mathbf z \in \Gamma \!\setminus \!\{\mathbf 0\}$, $\mathbf
z$ est un point non de torsion, donc sa hauteur de Néron-Tate est
minorée par ${\widehat{h}}_{\min}$. Ce qui permet de minorer la
hauteur de Néron-Tate du point $\mathbf x - \mathbf y$ par:
$$\widehat{h}(\mathbf x - \mathbf y) = \widehat{h}(N \mathbf z) = N^2 \widehat{h}(\mathbf z) \geq N^2 {\widehat{h}}_{\min} .$$
Par ailleurs, la hauteur de Néron-Tate du point $\mathbf x -
\mathbf y$ est majorée par:
$$\widehat{h}(\mathbf x - \mathbf y) = {\mid \mathbf x - \mathbf y \mid}^2 \leq (\mid \mathbf x \mid + \mid \mathbf y \mid)^2 \leq (\sqrt{R} + \sqrt{R})^2 ~~~~~~ \text{(car $\mathbf x , \mathbf y \in \Gamma_R$)},$$
$$\mbox{c'est-\`a-dire:}~~~~~~~~~~ \widehat{h}(\mathbf x - \mathbf y) ~\leq~ 4 R .~~~~~~~~~~~~~~~~~~~~~~~~~~~~~~~~~~~~~~~~~~~~~~~~~~~~~~~~~~~~~~$$
Il résulte de la majoration et de la minoration de
$\widehat{h}(\mathbf x - \mathbf y)$ qu'on a:
$$N^2 {\widehat{h}}_{\min} ~\leq~ 4 R .$$
Ce qui conduit à la contradiction $N \leq \sqrt{\frac{4
R}{{\widehat{h}}_{\min}}}$ et prouve que l'application
$\cl\rest{\Gamma_R}$ est bien une injection. Finalement
l'injectivité de $\cl\rest{\Gamma_R}$ entra{\sf\^\i}ne:
$$\sharp \Gamma_R ~\leq~ \sharp (\Gamma/N \Gamma) ~=~ N^r ~\leq~ \left(1 + \sqrt{\frac{4 R}{{\widehat{h}}_{\min}}}\right)^{\!\!r} ,$$
$$\mbox{c'est-à-dire:}~~~~~~~~~~ \sharp \Gamma_R ~\leq~ \left(1 + \sqrt{\frac{4 R}{{\widehat{h}}_{\min}}}\right)^{\!\!r} .~~~~~~~~~~~~~~~~~~~~~~~~~~~~~~~~~~~~~~~~~~~~~~~~~~~~~~$$
La démonstration du lemme \ref{c.49} est achevée.
$~~~~\blacksquare$\vspace{1mm}
\subsection{Démonstration des deux premiers théor\`emes principaux:}
Nous déduisons d'abord immédiatement du corollaire \ref{c.47} et
du théor\`eme \ref{c.48} le théor\`eme suivant:
\begin{theorem}\label{c.50}
Sous les hypoth\`eses du théor\`eme \ref{c.46}, on a:
\begin{equation*}
\begin{split}
\sharp\left\{\mathbf x \in E(K) ~/~ \mathbf x ~\mbox{satisfait}~
(II)\right\}
~ \leq & \\
&\!\!\!\!\!\!\!\!\!\!\!\!\!\!\!\!\!\!\!\!\!\!\!\!\!4 \epsilon^{-
\frac{1}{2}} \left[m (m - 1) (\log m + 9) + m
\nb{\log{\epsilon}}\right] \left(499 \epsilon^{- \frac{m}{2(m -
1)}}\right)^{\!\!r} .
\end{split}
\end{equation*}
\end{theorem}
De ce théor\`eme \ref{c.50} nous apparait deux façons de procéder
pour optimiser notre résultat. La premi\`ere consiste à choisir
l'entier $m \geq 2$ (en fonction de $\epsilon$ et $r$) de façon à
rendre la quantité $R = R(m)$ minimale afin d'affaiblir le mieu
possible le syst\`eme d'inégalité $(II)$. La deuxi\`eme par contre
consiste à choisir $m$ pour rendre plut\^ot la quantité du
théor\`eme \ref{c.50} estimant le nombre de points de $E(K)$
satisfaisant $(II)$, minimale. Le corollaire \ref{c.51} de
l'appendice montre que pour $\epsilon$ suffisamment petit (plus
précisément $\epsilon \leq \frac{1}{15788}$), l'entier:
$$m_0 := \left[\sqrt{\frac{2 \nb{\log{\epsilon}}}{\log{\nb{\log{\epsilon}}} - \log{\log{\nb{\log{\epsilon}}}} + 16}} + 2\right]$$
(avec $[.]$ désigne la partie enti\`ere), minimise presque $R$ et
on a:
$$R(m_0) ~\leq~ 56(\eta + 5) \epsilon^{- 1 - \frac{183}{\log
\mid \log \epsilon \!\mid}} .$$ D'autre part, un simple calcul
montre que la quantité:
$$ 4 \epsilon^{- \frac{1}{2}} \left[m_0 (m_0 - 1) (\log{m_0} + 9) + m_0 \!\mid \!\log \epsilon \!\mid\right] \left(499 \epsilon^{- \frac{m_0}{2(m_0 - 1)}}\right)^{\!\!r}$$ est majorée par:
$$34 \epsilon^{- 1/2}~\!{\!\mid \!\log \epsilon \!\mid}^{3/2}\!\!\left(\log \!\mid \!\log \epsilon \!\mid\right)^{\!\!- 1/2} \times \left[499
\epsilon^{- 1/2} \exp\!\left(\!\sqrt{\!\mid \!\log \epsilon \!\mid
\log \!\mid \!\log \epsilon \!\mid}\right)\right]^{r} .$$
Le premier théor\`eme principal résulte ainsi simplement du théor\`eme \ref{c.50} pour $m = m_0$.\\
Par ailleurs, le corollaire \ref{c.53} de l'appendice montre que
l'entier:
$$m_1 ~:=~ \left[\frac{r}{4} \nb{\log{\epsilon}} + 2\right] $$ (avec $[.]$ désigne la partie enti\`ere) minimise presque la quantité:
$$ 4 \epsilon^{- \frac{1}{2}} \left[m (m - 1) (\log m + 9) + m \nb{\log{\epsilon}}\right] \!\!\left(499 \epsilon^{- \frac{m}{2(m - 1)}}\right)^{\!\!r}$$ et lui donne une valeur
$$\leq~ 2 r^2 \epsilon^{- \frac{1}{2}} \left(\nb{\log{\epsilon}}\right)^{\!2} \!\left(\log r + \log{\nb{\log{\epsilon}}} + 82\right)
\!\!\left(499 \epsilon^{- \frac{1}{2}}\right)^{\!\!r} .$$ D'autre
part, un simple calcul montre que la valeur de $R$ lorsque $m$
vaut $m_1$ est majorée par: $$R(m_1) ~\leq~
\exp\!\left\{\!\left(\frac{r}{4} \nb{\log{\epsilon}} +
2\right)\!\! \left(\log{\nb{\log{\epsilon}}} + \log r + 16\right)
+ \log(\eta + 5)\!\right\}\! .$$ Apr\`es toutes ces majorations,
le deuxi\`eme théor\`eme principal résulte du théor\`eme
\ref{c.50}, en prenant $m = m_1$ dans ce dernier.
\subsection{Démonstration du troisi\`eme théor\`eme principal:}
Du théor\`eme \ref{c.46} et du lemme \ref{c.49} pour $A = E$
résulte immédiatement l'estimation:
\begin{equation*}
\begin{split}
\sharp\left\{\mathbf x \in E(K) ~/~ \text{$\mathbf x$ satisfait $(I)$}\right\}~& \leq~ \sharp {E(K)}_{\rm{tor}}~\!\!\left(\!1 + \sqrt{\frac{4 R(m)}{{\widehat{h}}_{\min}}}\right)^{\!\!r} \\
&
\!\!\!\!\!\!\!\!\!\!\!\!\!\!\!\!\!\!\!\!\!\!\!\!\!\!\!\!\!\!\!\!\!\!\!\!\!\!\!\!\!\!\!\!\!\!\!\!\!\!+
4 \epsilon^{- \frac{1}{2}} \left[m (m - 1) (\log m + 9) + m
\nb{\log{\epsilon}}\right] \!\!\left(\!499 \epsilon^{-
\frac{m}{2(m - 1)}}\right)^{\!\!r} .
\end{split}
\end{equation*}
Il reste à choisir l'entier $m$ de façon à optimiser cette
estimation. Pour ce faire, on remarque que lorsque $m$ est assez
grand, la quantité:
$$4 \epsilon^{- \frac{1}{2}} \left[m (m - 1) (\log m + 9) + m \nb{\log{\epsilon}}\right]\!\!\left(\!499 \epsilon^{- \frac{m}{2(m - 1)}}\right)^{\!\!r}$$
devient négligeable devant la quantité:
$$\sharp {E(K)}_{\rm{tor}}~\!\!\left(\!1 + \sqrt{\frac{4 R(m)}{{\widehat{h}}_{\min}}}\right)^{\!\!r} \gg \ll R(m)^{\frac{r}{2}} \geq \left((2904 m)^m \epsilon^{- \frac{m}{m - 1}}\right)^{\!\frac{r}{2}} .$$
Donc, optimiser l'estimation précédente revient presque à
optimiser $R(m)$, or ceci a été déja fait dans la démonstration du
premier théor\`eme principal. En prenant comme dans ce dernier:
$$m = m_0 := \left[\sqrt{\frac{2 \nb{\log{\epsilon}}}{\log{\nb{\log{\epsilon}}} - \log{\log{\nb{\log{\epsilon}}}} + 16}} + 2\right]$$
(avec $[.]$ désigne la partie enti\`ere) et en tenant compte des
estimations faites au cour de sa démonstration, le théor\`eme
\ref{c.3} suit.
\section{Démonstration des corollaires \ref{c.18}, \ref{c.19} et \ref{c.20}}
L'argument qu'on utilise pour déduire le corollaire \ref{c.18}
(resp \ref{c.19} et \ref{c.20}) du théor\`eme \ref{c.1} (resp
\ref{c.2} et \ref{c.3}) est l'objet du lemme suivant:
\begin{lemma}\label{c.21}
Soit $E$ une courbe elliptique définie sur un corps de nombres $K$
de degré $D$, plongée dans ${\mathbb P}_2$ à la Weierstrass,
d'équation projective $Y^2 Z = 4 X^3 - g_2 X Z^2 - g_3 Z^3 ~ (g_2
, g_3 \in K)$ et d'élément neutre (en tant que groupe) le point à
l'infini $\mathbf 0$ représenté dans ${\mathbb P}_2$ par les
coordonnées projectives $(0 : 1 : 0)$. \\ Supposons que pour tout
réel strictement positif $\epsilon$, pour tout ensemble fini $S$
de places sur $K$ et pour toute famille ${(\lambda_v)}_{v \in S}$
de réels positifs satisfaisant:
$$\sum_{v \in S} \frac{[K_v : {\mathbb Q}_v]}{[K : \mathbb Q]} \lambda_v = 1 ,$$ il existe une application:
$$
\begin{array}{rcl}
F : M_k & \rightarrow & {\mathbb R}_+ \\
v & \mapsto & F_v
\end{array}
$$
tel que le syst\`eme d'inégalités simultanées:
$${\mbox{dist}}_v(\mathbf x , \mathbf 0) \leq e^{- \lambda_v \epsilon h(\mathbf x) - F_v} ~~~~~~~~~~ (v \in S)$$
-dont les inconnues sont les points $\mathbf x$ de $E(K)$- n'admet
qu'un nombre fini de solutions constituant un ensemble de cardinal
majoré par une fonction positive:
$${fct}_{E , D}(\epsilon , S)$$
en $\epsilon$ et $S$. \\
Alors pour tout réel $\epsilon > 0$ et tout sous-ensemble fini $S$
de places sur $K$; pour tout choix de réels $0 < \epsilon'(T) <
\epsilon$, associés aux sous-ensembles $T$ de $S$, l'inégalité:
$$\prod_{v \in S} {{\mbox{dist}}_v(\mathbf x , \mathbf 0)}^{\frac{[K_v : {\mathbb Q}_v]}{[K : \mathbb Q]}} \leq e^{- \epsilon h(\mathbf x) - \sum_{v \in S} \frac{[K_v : {\mathbb Q}_v]}{[K : \mathbb Q]} F_v}$$
n'est satisfaite que par un nombre fini de points $\mathbf x$ de
$E(K)$ constituant un ensemble de cardinal majoré par:
$$\sum_{T \in \mathcal P (S)} \!\!\binom{A(T) + \card{~\!(T)} - 1}{\card{~\!(T)} - 1} {fct}_{E , D}\left(\epsilon'(T) , T\right)$$
o\`u $\mathcal P (S)$ désigne l'ensemble de toutes les parties de
$S$ et $A(T)$ le plus petit entier $\geq \frac{\epsilon'(T)
\card(T)}{\epsilon - \epsilon'(T)}$.
\end{lemma}
{\bf Démonstration.---} Soient $\epsilon$ un réel strictement
positif et $S$ un ensemble fini de places de $K$. Désignons par
$X$ l'ensemble des points $\mathbf x$ de $E(K)$ satisfaisant
l'inégalité:
\begin{gather}
\prod_{v \in S} {{\mbox{dist}}_v(\mathbf x , \mathbf
0)}^{\frac{[K_v : {\mathbb Q}_v]}{[K : \mathbb Q]}} ~\leq~ e^{-
\epsilon h(\mathbf x) - \sum_{v \in S} \frac{[K_v : {\mathbb
Q}_v]}{[K : \mathbb Q]} F_v} \label{3.62}
\end{gather}
et par $\pi$ l'application de $X$ dans $\mathcal P (S)$ associant
\`a chaque point $\mathbf x \in X$ l'ensemble $\pi(\mathbf x)$ des
places $v$ de $S$ pour lesquelles $\mathbf x$ satisfait
l'inégalité:
$${{\mbox{dist}}_v(\mathbf x , \mathbf 0)} ~\leq~ e^{- F_v} .$$
Il est clair que $\pi$ est bien définie, de plus on peut remarquer
d'apr\`es (\ref{3.62}) que pour tout $\mathbf x \in X$, le
sous-ensemble $\pi(\mathbf x)$ de $S$ n'est jamais vide, ce qui
revient à dire qu'on a: $\pi^{-1}(\{\emptyset\}) = \emptyset$. \\
Soit maintenant $T$ un sous-ensemble quelconque de $S$ tel que
$\pi^{-1}(\{T\}) \neq \emptyset$, $T$ est donc non vide puisque on
vient de remarquer que $\pi^{-1}(\{\emptyset\}) = \emptyset$. Par
définition m\^eme de l'application $\pi$, on a:
\begin{align}
\forall \mathbf x \in \pi^{-1}(\{T\}) , \forall v \in T : ~~~~
{{\mbox{dist}}_v(\mathbf x , \mathbf 0)} & \leq~ e^{- F_v} \label{3.63} \\
\forall \mathbf x \in \pi^{-1}(\{T\}) , \forall v \in S \setminus
T :
{{\mbox{dist}}_v(\mathbf x , \mathbf 0)} & >~ e^{- F_v} . \label{3.64} \\
\intertext{Des deux inégalités (\ref{3.62}) et (\ref{3.64}), on
déduit l'inégalité importante:} \forall \mathbf x \in
\pi^{-1}(\{T\}) : \prod_{v \in T} {{\mbox{dist}}_v(\mathbf x ,
\mathbf 0)}^{\frac{[K_v : {\mathbb Q}_v]}{[K : \mathbb Q]}} &
\leq~ e^{- \epsilon h(\mathbf x) - \sum_{v \in T} \frac{[K_v :
{\mathbb Q}_v]}{[K : \mathbb Q]} F_v}. \label{3.65}
\end{align}
Définissons maintenant pour tout $\mathbf x \in \pi^{-1}(T)$ et
pour toute place $v$ dans $T$ le réel positif $\xi_v(\mathbf x)$
par:
\begin{equation}
\begin{split}
{{\mbox{dist}}_v(\mathbf x , \mathbf 0)} & ~=~ e^{- \xi_v(\mathbf x) \epsilon h(\mathbf x) - F_v} ~~~~~~~~~~\text{si $h(\mathbf x) \neq 0$} \\
\xi_v(\mathbf x) & :=~ \frac{1}{\card{~\!\!T}} \frac{[K : \mathbb
Q]}{[K_v : {\mathbb Q}_v]} ~~~~~~\text{si $h(\mathbf x) = 0$}
\end{split} \!\!\!\!\!\!\!\!\!\!\!\!\!\!\!\!\!\!\!\!\!\!\!\!\!\!\!\!\!\!\!\!\!\!\!\!\!\!\!\!\!\!\!\!\!\!\!\!\!\!\!\!\!\!\!\!\!\!\!\!. \label{3.66}
\end{equation}
(lorsque $h(\mathbf x) \neq 0$, l'existence de $\xi_v(\mathbf x)$ est justifiée par (\ref{3.63})). \\
D'apr\`es (\ref{3.65}) on doit avoir:
\begin{equation}
\forall \mathbf x \in \pi^{-1}(\{T\}) : ~~\sum_{v \in T}
\frac{[K_v : {\mathbb Q}_v]}{[K : \mathbb Q]} \xi_v(\mathbf x) ~
\geq~ 1 . \label{3.67}
\end{equation}
Maintenant on a pour tout $\mathbf x \in \pi^{-1}(\{T\})$:
\begin{eqnarray*}
A(T) + \mbox{card}~\!(T) & \leq & A(T) \frac{\epsilon}{\epsilon'(T)} \\
& \leq & \sum_{v \in T} A(T) \frac{\epsilon}{\epsilon'(T)} \frac{[K_v : {\mathbb Q}_v]}{[K : \mathbb Q]} \xi_w(\mathbf x) ~~~~~~~~~~\text{(d'apr\`es (\ref{3.67}))} \\
& \leq & \sum_{v \in T} \left[ A(T) \frac{\epsilon}{\epsilon'(T)}
\frac{[K_v : {\mathbb Q}_v]}{[K : \mathbb Q]} \xi_v(\mathbf
x)\right] + \mbox{card}~\!(T) ;
\end{eqnarray*}
d'o\`u:
\begin{equation}
\forall \mathbf x \in \pi^{-1}(\{T\}) : ~~A(T) ~\leq~ \sum_{v \in
T} \left[ A(T) \frac{\epsilon}{\epsilon'(T)} \frac{[K_v : {\mathbb
Q}_v]}{[K : \mathbb Q]} \xi_v(\mathbf x)\right] . \label{3.68}
\end{equation}
Cette derni\`ere inégalité (\ref{3.68}) entra{\sf\^\i}ne -pour
tout $\mathbf x \in \pi^{-1}(\{T\})$- l'existence d'une famille
d'entiers positifs ${(a_v(\mathbf x))}_{v \in T}$ satisfaisant
pour toute place $v \in T$:
\begin{equation}
\begin{split}
a_v(\mathbf x) &\leq~ \left[ A(T) \frac{\epsilon}{\epsilon'(T)} \frac{[K_v : {\mathbb Q}_v]}{[K : \mathbb Q]} \xi_v(\mathbf x)\right] \\
&\leq~ A(T) \frac{\epsilon}{\epsilon'(T)} \frac{[K_v : {\mathbb
Q}_v]}{[K : \mathbb Q]} \xi_v(\mathbf x) \label{3.69}
\end{split}
\end{equation}
et
\begin{equation}
\sum_{v \in T} a_v(\mathbf x) ~=~ A(T) . \label{3.70}
\end{equation}
Or, l'équation (\ref{3.70}) admet exactement $\binom{A(T) +
\card{~\!\!(T)} - 1}{\card{~\!\!(T)} - 1}$ solutions en entiers
positifs ${(a_v)}_{v \in T}$, donc il doit exister un
sous-ensemble $E_T$ de $\pi^{-1}(\{T\})$ de cardinal:
\begin{equation}
\card{~\!\!(E_T)} ~\geq~
\frac{\card{\left(\pi^{-1}(\{T\})\right)}}{\binom{A(T) +
\card{~\!\!(T)} - 1}{\card{~\!\!(T)} - 1}} \label{3.71}
\end{equation}
tel que pour toute place $v \in T$, les entiers positifs
$a_v(\mathbf x) ,~ \mathbf x \in E_T$ soient tous égaux. Notons
alors simplement l'entier positif $a_v(\mathbf x) ~ (\mathbf x \in
E_T)$ par $a_v$ (pour toute place $v\in T$). Ces ${(a_v)}_{v \in
T}$ satisfont d'apr\`es (\ref{3.69}) et (\ref{3.70}):
\begin{align}
\forall \mathbf x \in E_T : ~~~~ a_v & \leq~ A(T) \frac{\epsilon}{\epsilon'(T)} \frac{[K_v : {\mathbb Q}_v]}{[K : \mathbb Q]} \xi_v(\mathbf x) \label{3.72} \\
\intertext{et}
\sum_{v \in T} a_v & =~ A(T) . \label{3.73} \\
\intertext{Posons aussi pour toute place $v \in T$:}
\lambda_v & :=~ \frac{a_v}{A(T)} \frac{[K : \mathbb Q]}{[K_v : {\mathbb Q}_v]}. \label{3.74} \\
\intertext{D'apr\`es (\ref{3.72}) et (\ref{3.73}) ces réels
positifs ${(\lambda_v)}_{v \in T}$ satisfont pour toute place $v
\in T$ et tout $\mathbf x \in E_T$:}
\lambda_v \epsilon'(T) & \leq~ \epsilon \xi_v(\mathbf x) \label{3.75} \\
\intertext{et}
\sum_{v \in T} \frac{[K_v : {\mathbb Q}_v]}{[K : \mathbb Q]} \lambda_v & =~ 1 . \label{3.76} \\
\intertext{On vérifie aisément gr\^ace \`a l'inégalité
(\ref{3.75}) et à l'égalité (de définition) (\ref{3.66}) (lorsque
$h(\mathbf x) \neq 0$) et gr\^ace à l'inégalité (\ref{3.63})
(lorsque $h(\mathbf x) = 0$) qu'on a:}
\forall \mathbf x \in E_T : ~~~~{\mbox{dist}}_v(\mathbf x , \mathbf 0) & \leq~ e^{- \lambda_v \epsilon'(T) h(\mathbf x) - F_v} ~~~~~~~~~~ \forall v \in T . \label{3.77} \\
\intertext{Or d'apr\`es l'hypoth\`ese du lemme \ref{c.21}, le
nombre de solutions $\mathbf x \in E(K)$ du syst\`eme (\ref{3.77})
est majoré par ${fct}_{E , D}(\epsilon'(T) , T)$. D'o\`u:}
\card{~\!\!(E_T)} & \leq~ {fct}_{E , D}(\epsilon'(T) , T) \label{3.78} \\
\intertext{et puis, par (\ref{3.71}), on en déduit qu'on a:}
\card{\left(\pi^{-1}(\{T\})\right)} & \leq~ \binom{A(T) +
\card{~\!\!(T)} - 1}{\card{~\!\!(T)} - 1} {fct}_{E ,
D}(\epsilon'(T) , T). \label{3.79}
\end{align}
Remarquons enfin que cette derni\`ere inégalité (\ref{3.79}) est
trivialement satisfaite lorsque $\pi^{-1}(\{T\}) = \emptyset$,
donc (\ref{3.79}) est valable pour tout $T$ dans $\mathcal P (S)$.
\\ En additionant membre \`a membre les inégalités (\ref{3.79})
correspondant \`a chaque $T \in \mathcal P (S)$, on obtient:
\begin{align}
\sum_{T \in \mathcal P (S)} \card{\left(\pi^{-1}(\{T\})\right)} &
\leq~ \sum_{T \in \mathcal P (S)} \binom{A(T) + \card{~\!\!(T)} -
1}{\card{~\!\!(T)} - 1} {fct}_{E , D}(\epsilon'(T) , T) . \notag
\end{align}
Pour conclure, il ne reste qu'\`a remarquer que $\sum_{T \in
\mathcal P (S)} \card{\left(\pi^{-1}(\{T\})\right)} =
\card{~\!\!(X)}$. La démonstration est achevée.
$~~~~\blacksquare$\vspace{2mm}

Dans ce qui suit, nous allons en déduire le corollaire \ref{c.18} du théor\`eme \ref{c.1} en utilisant le lemme \ref{c.21} précédent; les deux autres corollaires \ref{c.19} et \ref{c.20} se déduisent respectivement des deux théor\`emes \ref{c.2} et \ref{c.3} exactement de la m\^eme mani\`ere. \\
Nous appliquons le lemme \ref{c.21} avec:
\begin{eqnarray*}
F_v &:=& 56(\eta + 5) \epsilon^{- 1 - \frac{183}{\log \mid \log \epsilon \!\mid}} \lambda_v + 2 m_v + c_v ~~~~ \forall v \in M_K \\
{fct}_E(\epsilon , S) &=& 34 \epsilon^{- 1/2} \!\!\left(\!\mid
\!\log \epsilon \!\mid \right)^{\!3/2}
\!\!\left(\log \!\mid \!\log \epsilon \!\mid \right)^{\!\!- 1/2} \left[499 \epsilon^{- 1/2} \!\!\exp\!\left(\sqrt{\!\mid \!\log \epsilon \!\mid \log \!\mid \!\log \epsilon \!\mid}\right)\right]^{r} \\
& = & {fct}_E(\epsilon) ~~~~\text{(indépendante de $S$)}
\end{eqnarray*}
$$\forall T \in \mathcal P (S) : ~~~~ \epsilon'(T) = \frac{\epsilon}{2} ~~\text{et}~~ A(T) = \card{~\!\!(T)} .~~~~~~~~~~~~~~~~~~~~~~~~~~~~~~~~~~~~~~~~~~~~~~~~$$
Ainsi, d'apr\`es le théor\`eme \ref{c.1}, l'hypoth\`ese principale
du lemme \ref{c.21} est bien vérifiée, donc ce dernier entraine
que l'ensemble des points $\mathbf x$ de $E(K)$ satisfaisant
l'inégalité:
$$\prod_{v \in S} {{\mbox{dist}}_v(\mathbf x , \mathbf 0)}^{\frac{[K_v : {\mathbb Q}_v]}{[K : \mathbb Q]}} \leq e^{- \epsilon h(\mathbf x) - \sum_{v \in S} \frac{[K_v : {\mathbb Q}_v]}{[K : \mathbb Q]} F_v}$$
est de cardinal:
$$\leq \sum_{T \in \mathcal P (S)} \binom{A(T) + \card{~\!\!(T)} - 1}{\card{~\!\!(T)} - 1} {fct}_E(\epsilon'(T) , T) ,$$
c'est-à-dire:
\begin{eqnarray*}
& \leq & {fct}_E\left(\frac{\epsilon}{2}\right) . \sum_{T \in \mathcal P (S)} \binom{2 \card{~\!\!(T)} - 1}{\card{~\!\!(T)} - 1} \\
& \leq & {fct}_E\left(\frac{\epsilon}{2}\right) . \sum_{n = 0}^{\card{~\!\!(S)}} \binom{\card{~\!\!(S)}}{n} \binom{2 n - 1}{n - 1} ~~~~~~~~~~\text{(en posant $n = \card{~\!\!(T)}$)} \\
& \leq & {fct}_E\left(\frac{\epsilon}{2}\right) . \sum_{n = 0}^{\card{~\!\!(S)}} \binom{\card{~\!\!(S)}}{n} 4^n \\
& \leq & 5^{\card{~\!\!(S)}} .
{fct}_E\left(\frac{\epsilon}{2}\right) .
\end{eqnarray*}
Par ailleurs, en utilisant l'hypoth\`ese:
$$\sum_{v \in S} \frac{[K_v : {\mathbb Q}_v]}{[K : \mathbb Q]} \lambda_v = 1$$
et en majorant $\sum_{v \in S} \frac{[K_v : {\mathbb Q}_v]}{[K :
\mathbb Q]} m_v$ par:
$$\sum_{v \in S} \frac{[K_v : {\mathbb Q}_v]}{[K : \mathbb Q]} m_v ~\leq~ \sum_{v \in M_K} \frac{[K_v : {\mathbb Q}_v]}{[K : \mathbb Q]} m_v ~=:~ \eta$$
et $\sum_{v \in S}\frac{[K_v : {\mathbb Q}_v]}{[K : \mathbb Q]}
c_v$ par:
$$\sum_{v \in S}\frac{[K_v : {\mathbb Q}_v]}{[K : \mathbb Q]} c_v ~\leq~ \sum_{v \in M_K} \frac{[K_v : {\mathbb Q}_v]}{[K : \mathbb Q]} c_v = \sum_{v \in {M_K}^{\!\!\infty}} \frac{[K_v : {\mathbb Q}_v]}{[K : \mathbb Q]} c_v = 16 ,$$
on en déduit pour la quantité: $\sum_{v \in S} \frac{[K_v :
{\mathbb Q}_v]}{[K : \mathbb Q]} F_v$ la majoration:
\begin{eqnarray*}
\sum_{v \in S} \frac{[K_v : {\mathbb Q}_v]}{[K : \mathbb Q]} F_v & \leq & 56(\eta + 5) \epsilon^{- 1 - \frac{183}{\log \mid \log \epsilon \!\mid}} + 2 \eta + 16 \\
& \leq & 57(\eta + 5) \epsilon^{- 1 - \frac{183}{\log \mid \log
\epsilon \!\mid}} ;
\end{eqnarray*}
ce qui entra{\sf\^\i}ne que l'inégalité:
$$\prod_{v \in S} {{\mbox{dist}}_v(\mathbf x , \mathbf 0)}^{\frac{[K_v : {\mathbb Q}_v]}{[K : \mathbb Q]}} \leq e^{- \epsilon h(\mathbf x) - \sum_{v \in S} \frac{[K_v : {\mathbb Q}_v]}{[K : \mathbb Q]} F_v}$$
est impliquée par l'inégalité:
$$\prod_{v \in S} {{\mbox{dist}}_v(\mathbf x , \mathbf 0)}^{\frac{[K_v : {\mathbb Q}_v]}{[K : \mathbb Q]}} \leq e^{- \epsilon h(\mathbf x) - 57(\eta + 5) \epsilon^{- 1 - \frac{183}{\log \mid \log \epsilon \!\mid}}} .$$
Le corollaire \ref{c.18} s'ensuit. La démonstration est achevée. $~~~~\blacksquare$

\section{Formulaire}
Dans tous ce paragraphe, soit $E$ une courbe elliptique définie
sur un corps de nombres $K$ et plongée à la Weierstrass dans
l'espace projectif ${\mathbb P}_2$ et soit:
$$Y^2 Z = 4 X^3 - g_2 X Z^2 - g_3 Z^3~~~~~~ g_2 , g_3 \in K$$
son équation projective dans ce plongement. \\Nous donnons dans le
théor\`eme \ref{c.8} qui suit trois familles de formules
d'addition explicites sur $E$ ainsi que les cartes o\`u chacune de
ces familles de formules est valable, de mani\`ere à avoir un
atlas constitué de trois cartes de $E^2$. \\Dans le théor\`eme
\ref{c.9} nous montrons, par un procédé de récurence, l'existence
d'une famille de formes représentant globalement la multiplication
d'un point de $E$ par un entier positif $n$ donné. Le degré et la
hauteur de ces formes est bien controlé en fonction de $n$ et de
la hauteur de $E$. Notre référence principale pour ce théor\`eme
\ref{c.9} est le chapitre $2$ de [La3]. \\Enfin dans le théor\`eme
\ref{c.14} nous donnons une valeur explicite pour la constante de
Néron-Tate de $E \hookrightarrow {\mathbb P}_2$ qui est
essentiellement celle de [Zi-Sch] et [Da1].
\subsection{Formules d'addition sur $E$:}
Nous appelons syst\`eme complet de familles de formes représentant
l'addition sur $E \hookrightarrow {\mathbb P}_2 ,$ la donnée d'un
nombre fini de familles de formes ${\underline A}_i = (A_{i 0} ,
A_{i 1} , A_{i 2}) , i \in I$ ($~I$ fini ) de $K[(X_1 , Y_1 , Z_1)
, (X_2 , Y_2 , Z_2)]$ et d'un nombre fini d'ouverts non vides
${\Omega}_i , i \in I$ de $E^2$ formant un recouvrement pour $E^2$
tel que pour tout $i \in I$ et pour tout couple de points
$({\mathbf p}_1 , {\mathbf p}_2) \in {\Omega}_i \subset E^2 ,$ le
point ${\mathbf p}_1 + {\mathbf p}_2$ de $E$ peut \^etre
représenté dans ${\mathbb P}_2$ par le syst\`eme de coordonnées
projectives ${\underline A}_i ({\mathbf p}_1 , {\mathbf p}_2) =
(A_{i 0}({\mathbf p}_1 , {\mathbf p}_2) : A_{i 1}({\mathbf p}_1 ,
{\mathbf p}_2) : A_{i 2}({\mathbf p}_1 , {\mathbf p}_2))$.\\
Notons que dans un syst\`eme complet de familles de formes
représentant l'addition sur $E,$ il n'est pas vraiment nécessaire
de préciser les ouverts ${\Omega}_i , i \in I$ de $E^2$
correspondant à chacune des familles de formes ${\underline A}_i ,
i \in I$. En effet si $\underline A = (A_0 , A_1 , A_2)$ est une
famille de formes de $K[(X_1 , Y_1, Z_1) , (X_2 , Y_2 , Z_2)]$
représentant l'addition sur $E$ sur un ouvert non vide $\Omega$ de
$E^2$ alors $\underline A$ doit certainement représenter
l'addition sur $E$ sur tout le sous-ensemble de $E^2$ là o\`u les
trois formes $A_0 , A_1$ et $A_2$ ne s'annulent pas simultanément.
Ce dernier étant un ouvert de $E^2$, ne dépend que de la famille
de formes $\underline A$ et contient $\Omega$, qu'il peut
avantageusement remplacer. On a le théor\`eme suivant:
\begin{theorem}\label{c.8}
Un syst\`eme complet de familles de formes représentant l'addition
sur $E$ est donné par:
\begin{eqnarray*}
\rm{1}) ~~ A_0 & := & {Y_1}^2 X_2 Z_2 - 2 X_1 Y_1 Y_2 Z_2 + 2 Y_1
Z_1 X_2 Y_2
- 3 g_3 X_1 Z_1 {Z_2}^2 - g_2 {X_1}^2 {Z_2}^2 \\
 & & + 3 g_3 {Z_1}^2 X_2 Z_2 + g_2 {Z_1}^2 {X_2}^2 - X_1 Z_1
 {Y_2}^2 \\
 A_1 & := & {Y_1}^2 Y_2 Z_2 + 3 g_3 Y_1 Z_1 {Z_2}^2 + g_2 X_1 Y_1
 {Z_2}^2 + 2 g_2 Y_1 Z_1 X_2 Z_2 - Y_1 Z_1 {Y_2}^2 \\
  & & -12 X_1 Y_1 {X_2}^2 - 3 g_3 {Z_1}^2 Y_2 Z_2 - 2 g_2 X_1 Z_1
  Y_2 Z_2 - g_2 {Z_1}^2 X_2 Y_2 \\
   & & + 12 {X_1}^2 X_2 Y_2
\end{eqnarray*}
\begin{eqnarray*}
 A_2 & := & {Y_1}^2 {Z_2}^2 + g_2 X_1 Z_1 {Z_2}^2 - g_2 {Z_1}^2 X_2 Z_2 - 12 {X_1}^2 X_2 Z_2 - {Z_1}^2 {Y_2}^2 \\
  & & + 12 X_1 Z_1 {X_2}^2 \\
\end{eqnarray*}
\begin{eqnarray*}
\rm{2}) ~~ A_0 & := & 4 {Y_1}^2 {X_2}^2 + {g_2}^2 X_1 Z_1 {Z_2}^2
+ 12 g_3
{X_1}^2 {Z_2}^2 - {g_2}^2 {Z_1}^2 X_2 Z_2 \\
 & & + 4 g_2 {X_1}^2 X_2 Z_2 - 12 g_3 {Z_1}^2 {X_2}^2 - 4 g_2 X_1 Z_1
{X_2}^2 - 4 {X_1}^2
 {Y_2}^2 \\
 A_1 & := & 4 {Y_1}^2 X_2 Y_2 - {g_2}^2 Y_1 Z_1 {Z_2}^2 - 12 g_3
 X_1 Y_1 {Z_2}^2 - 24 g_3 Y_1 Z_1 X_2 Z_2 \\
  & & - 8 g_2 X_1 Y_1 X_2 Z_2 - 4 g_2 Y_1 Z_1 {X_2}^2 - 4 X_1 Y_1
  {Y_2}^2 + {g_2}^2 {Z_1}^2 Y_2 Z_2 \\
  & & + 24 g_3 X_1 Z_1 Y_2 Z_2 + 4 g_2 {X_1}^2 Y_2 Z_2 + 12 g_3
  {Z_1}^2 X_2 Y_2 + 8 g_2 X_1 Z_1 X_2 Y_2 \\
 A_2 & := & 4 {Y_1}^2 X_2 Z_2 + 8 X_1 Y_1 Y_2 Z_2 - 8 Y_1 Z_1 X_2
 Y_2 - 12 g_3 X_1 Z_1 {Z_2}^2 \\
 & & - 4 g_2 {X_1}^2 {Z_2}^2 + 12 g_3 {Z_1}^2 X_2 Z_2 + 4 g_2
 {Z_1}^2 {X_2}^2 - 4 X_1 Z_1 {Y_2}^2 \\
\end{eqnarray*}
\begin{eqnarray*}
\rm{3}) ~~ A_0 & := & 4 {Y_1}^2 X_2 Y_2 + {g_2}^2 Y_1 Z_1 {Z_2}^2
+ 12 g_3
X_1 Y_1 {Z_2}^2 + 24 g_3 Y_1 Z_1 X_2 Z_2 \\
 & & + 8 g_2 X_1 Y_1 X_2 Z_2 + 4 g_2 Y_1 Z_1 {X_2}^2 + 4 X_1 Y_1
 {Y_2}^2 + {g_2}^2 {Z_1}^2 Y_2 Z_2 \\
 & & + 24 g_3 X_1 Z_1 Y_2 Z_2 + 4 g_2 {X_1}^2 Y_2 Z_2 + 12 g_3
 {Z_1}^2 X_2 Y_2 + 8 g_2 X_1 Z_1 X_2 Y_2 \\
A_1 & := & 4 {Y_1}^2 {Y_2}^2 + \left({g_2}^3 - 36 {g_3}^2\right)
{Z_1}^2 {Z_2}^2 - 12 g_2 g_3 X_1 Z_1 {Z_2}^2 - 4 {g_2}^2 {X_1}^2
{Z_2}^2 \\
 & & - 12 g_2 g_3 {Z_1}^2 X_2 Z_2 - 16 {g_2}^2 X_1 Z_1 X_2 Z_2 -
 144 g_3 {X_1}^2 X_2 Z_2 \\
 & & - 4 {g_2}^2 {Z_1}^2 {X_2}^2 - 144 g_3 X_1 Z_1 {X_2}^2 - 48
 g_2 {X_1}^2 {X_2}^2 \\
 A_2 & := & 4 {Y_1}^2 Y_2 Z_2 - 12 g_3 Y_1 Z_1 {Z_2}^2 - 4 g_2 X_1
 Y_1 {Z_2}^2 - 8 g_2 Y_1 Z_1 X_2 Z_2 \\
  & & + 4 Y_1 Z_1 {Y_2}^2 + 48 X_1 Y_1 {X_2}^2 - 12 g_3 {Z_1}^2
  Y_2 Z_2 - 8 g_2 X_1 Z_1 Y_2 Z_2 \\
  & & - 4 g_2 {Z_1}^2 X_2 Y_2 + 48 {X_1}^2 X_2 Y_2 \\
\end{eqnarray*}
Ainsi, ces trois familles sont constituées de formes de bidegré
$(2 , 2)$ et sont -pour toute place $v$ de $K$- de hauteur logarithmique $v$-adique $h_v(\underline A)$ majorée par: \\ \\
pour $\rm{1})$ ~~~~ $\displaystyle h_v(\underline A) \leq
\begin{cases}
m_v + \log 12 & \text{si $v$ est infinie} \\
m_v & \text{si $v$ est finie}
\end{cases}$ \\ \\
pour $\rm{2})$ ~~~~ $\displaystyle h_v(\underline A) \leq
\begin{cases}
2 m_v + \log 24 & \text{si $v$ est infinie} \\
2 m_v & \text{si $v$ est finie}
\end{cases}$ \\ \\
pour $\rm{3})$ ~~~~ $\displaystyle h_v(\underline A) \leq
\begin{cases}
3 m_v + \log 144 & \text{si $v$ est infinie} \\
3 m_v & \text{si $v$ est finie}
\end{cases}$. \\ \\
Par conséquent, ces trois familles de $\rm{1}) , \rm{2})$ et $\rm{3})$ sont de hauteurs de Gauss-Weil majorées respectivement par $\eta + \log 12 , 2 \eta + \log 24$ et $3 \eta + \log 144$. \\
Par ailleurs, quand $v$ est une place infinie sur $K$, ces trois
familles de $\rm{1}) , \rm{2})$ et $\rm{3})$ sont de longueurs
logarithmiques $v$-adique $\ell_v(\underline A)$ majorées
respectivement par: $m_v + \log 38 , 2 m_v + \log 106$ et $3 m_v +
\log 425$.
\end{theorem}
{\bf Démonstration.---} Ce syst\`eme complet de formules
d'additions est celui donné dans [La-Ru], on n'a fait que
développer les calculs de cette référence.
$~~~~\blacksquare$\vspace{1mm}
\begin{remarque}\label{c.54}
Les trois familles de formules d'additions données par le
théor\`eme \ref{c.8} sont aussi de hauteurs de Gauss-Weil majorées
par $h(1 : {g_2}^3 : {g_3}^2) + \log 144$.
\end{remarque}
\subsection{Formules de multiplication d'un point de $E$ par un entier positif donné:}
Nous consacrons ce sous-paragraphe à la démonstration du
théor\`eme suivant:
\begin{theorem}\label{c.9}
Pour tout entier $n \geq 1,$ il existe une famille de formes
${\underline F}^{(n)} := ({F_0}^{(n)} , {F_1}^{(n)} ,
{F_2}^{(n)})$ de $K[X , Y , Z]$ de degré $n^2$ chacune,
représentant globalement la multiplication par $n$ sur $E
\hookrightarrow {\mathbb P}_2$ tel que pour toute place finie
(resp infinie) $v$ de $K$, la hauteur logarithmique locale
$v$-adique $h_v$ (resp la longueur logarithmique locale $v$-adique
$\ell_v$) de la famille ${\underline F}^{(n)}$ est majorée par
$\frac{3}{2} m_v . n^2$ (resp par $\frac{3}{2} (m_v + 3) . n^2$).
Par conséquent cette famille de formes ${\underline F}^{(n)}$ est
de hauteur de Gauss-Weil majorée par:
$$\widetilde{h}({\underline F}^{(n)}) \leq \frac{3}{2} (\eta + 3).n^2 .$$
\end{theorem}
Afin de démontrer ce théor\`eme, on se ref\`ere dans toute la
suite de ce sous-paragraphe au chapitre $2$ de [La3]. Soit $\wp$
la fonction de Weierstrass associée à $E$ et $\Lambda$ le réseau
de périodes de $E$. La courbe elliptique $E$ est alors isomorphe
au groupe abélien $\mathbb C / \Lambda$ et l'isomorphisme en
question est donné par:
\begin{eqnarray*}
\rho : \mathbb C / \Lambda & \longrightarrow & E \\
 z & \longmapsto & (\wp(z) : {\wp}'(z) : 1) ~~\mbox{si}~~ z \neq 0 \\
 0 & \longmapsto & (0: 1 : 0)~.
\end{eqnarray*}
Pour tout $n \in \mathbb N$, notons ${(\mathbb C / \Lambda)}_n$ le
sous-groupe des points de $n$-torsion de $\mathbb C / \Lambda$.
D'apr\`es le chapitre $2$ de [La3], il existe pour tout entier $n
\geq 1$ une fonction elliptique $f_n$ vérifiant:
\begin{equation}
{f_n(z)}^2 ~=~ n^2 \!\!\!\!\!\!\!\!\!\!\prod_{u \in {(\mathbb C /
\Lambda)}_n , u \neq 0} \!\!\!\!\!\!\!\!\!\!\left(\wp(z) -
\wp(u)\right) . \label{3.80}
\end{equation}
Pour tout $n \geq 1$, les zéros de la fonction elliptique $f_n$
sont les points de $n$-torsion de $\mathbb C / \Lambda$ autre que
$0$ et se sont tous des zéros simples, par ailleurs ses p\^oles
sont les points du réseau $\Lambda$ et sont évidemment d'ordre
$(n^2 - 1)$ chacun. On sait aussi d'apr\`es [La3] que les
fonctions $f_n , (n \geq 1)$ s'écrivent comme des polyn\^omes à
coefficients dans $K$ en $\wp(z)$ et ${\wp}'(z)$, c'est-à-dire
qu'il existe une suite de polyn\^omes ${(Q_n)}_{n \geq 1}$ de $K[X
, Y]$ tels que:
$$f_n(z) = Q_n\left(\wp(z) ~,~ {\wp}'(z)\right)~~ \forall n \geq 1 .$$
On peut m\^eme préciser qu'on a:
\\ $\bullet$ si $n$ est impair: $Q_n(X , Y) = P_n(X)$ o\`u $P_n$
est un polyn\^ome de $K[X]$ de degré $\frac{n^2 - 1}{2}$ et de
coefficient dominant $n$.
\\ $\bullet$ si $n$ est pair: $Q_n(X , Y) = \frac{1}{2} Y P_n(X)$
o\`u $P_n$ est un polyn\^ome de $K[X]$ de degré $\frac{n^2 -
4}{2}$ et de coefficient dominant $n$.
\\ Il est vérifié dans [La3] qu'on a:
$$Q_1(X , Y) = 1~,~ Q_2(X , Y) = Y~,~ Q_3(X , Y) = 3 X^4 - \frac{3}{2} g_2 X^2 -3 g_3 X - \frac{1}{16} {g_2}^2$$
$$\mbox{et}~~ Q_4(X , Y) = \frac{1}{2} Y \left(4 X^6 \!-\! 5 g_2 X^4 \!-\! 20 g_3 X^3 \!-\! \frac{5}{4} {g_2}^2 X^2
\!-\! g_2 g_3 X \!-\! 2 {g_3}^2 \!+\! \frac{{g_2}^3}{16}\right).$$
Posons aussi, par convention, $Q_{-1} \equiv -1$ et $Q_0 \equiv
0$. Le théor\`eme $1.3$ du chapitre $2$ de [La3] donne des
formules de récurence permettant de calculer de proche en proche
ces polyn\^omes $Q_n, n \geq 1.$ Ces formules sont:
\begin{equation}
\begin{split}
\forall n \geq 1 : ~~~~ \left\{
\begin{array}{lcl}
Q_{2 n + 1} & = & Q_{n + 2} Q_{n}^{3} - Q_{n - 1} Q_{n + 1}^{3} \\
 Y Q_{2n} & = & Q_n \left( Q_{n + 2} Q_{n - 1}^{2} - Q_{n - 2} Q_{n + 1}^{2}\right)
\end{array}
\right..
\end{split} \label{3.81}
\end{equation}
 Gr\^ace à ces formules de récurence, on peut donner des estimations pour les hauteurs (ou longueurs) locales des polyn\^omes $Q_n , n \geq 1$. On obtient le lemme suivant:
\begin{lemma}\label{c.55}
Soit $v$ une place sur $K$ et $\xi := \sqrt[3]{\frac{1}{2}}$, on
a:
\\ $\bullet$ quand $v$ est infinie et $n \geq 2 :$
$$L_v(Q_n) \leq \left\{
\begin{array}{cc}
\xi {(4 M_v)}^{\frac{n^2 - 1}{4}} & \mbox{si} ~n~ \mbox{est
impair}\\
\xi {(4 M_v)}^{\frac{n^2 - 4}{4}} & \mbox{si} ~n~ \mbox{est pair}
\end{array}
\right.~,$$
\\ $\bullet$ quand $v$ est finie, $v \mid 2$ et $n \geq 1 :$
$$H_v(Q_n) \leq \left\{
\begin{array}{cc}
{(4 M_v)}^{\frac{n^2 - 1}{4}} & \mbox{si} ~n~ \mbox{est
impair}\\
{(4 M_v)}^{\frac{n^2 - 4}{4}} & \mbox{si} ~n~ \mbox{est pair}
\end{array}
\right.~,$$
\\ $\bullet$ et quand $v$ est finie, $v \nmid 2$ et $n \geq 1 :$
$$H_v(Q_n) \leq \left\{
\begin{array}{cc}
{M_v}^{\frac{n^2 - 1}{4}} & \mbox{si} ~n~ \mbox{est
impair}\\
{M_v}^{\frac{n^2 - 4}{4}} & \mbox{si} ~n~ \mbox{est pair}
\end{array}
\right.~.$$
\end{lemma}
{\bf Démonstration.---}
On proc\`ede par récurrence. \\
$~~~$-Dans le cas $v$ infinie, on vérifie l'estimation du lemme
\ref{c.55} pour $n = 2 , 3 , 4 , 5 , 6$ et puis on utilise les
formules (\ref{3.81}) en distinguant les cas: $n = 4 k , 4 k + 1 ,
4 k + 2 , 4 k + 3
~~(k \in {\mathbb N}^{*})$ pour la récurrence. \\
$~~~~$-Dans le cas $v$ finie, on vérifie l'estimation du lemme
\ref{c.55} pour $n = 1 ,$ \\ $2 , 3 , 4$ et puis on utilise aussi
les formules (\ref{3.81}) en distinguant les m\^emes cas que
précédemment pour établir la récurrence.
$~~~~\blacksquare$\vspace{2mm}

Du lemme \ref{c.55} découle
immédiatement le corollaire suivant:
\begin{corollary}\label{c.56}
Soit $v$ une place sur $K$ et $\xi := \sqrt[3]{\frac{1}{2}}$. On
a:
\\ $\bullet$ quand $v$ est infinie et $n \geq 2 :$
$$L_v(P_n) \leq \left\{
\begin{array}{cc}
\xi {(4 M_v)}^{\frac{n^2 - 1}{4}} & \mbox{si} ~n~ \mbox{est
impair}\\
2 \xi {(4 M_v)}^{\frac{n^2 - 4}{4}} & \mbox{si} ~n~ \mbox{est
pair}
\end{array}
\right.~,$$
\\ $\bullet$ quand $v$ est finie, $v \mid 2$ et $n \geq 1 :$
$$H_v(P_n) \leq \left\{
\begin{array}{cc}
{(4 M_v)}^{\frac{n^2 - 1}{4}} & \mbox{si} ~n~ \mbox{est
impair}\\
\frac{1}{2}{(4 M_v)}^{\frac{n^2 - 4}{4}} & \mbox{si} ~n~ \mbox{est
pair}
\end{array}
\right.~,$$
\\ $\bullet$ et quand $v$ est finie, $v \nmid 2$ et $n \geq 1 :$
$$H_v(P_n) \leq \left\{
\begin{array}{cc}
{M_v}^{\frac{n^2 - 1}{4}} & \mbox{si} ~n~ \mbox{est
impair}\\
{M_v}^{\frac{n^2 - 4}{4}} & \mbox{si} ~n~ \mbox{est pair}
\end{array}
\right.~.$$
\end{corollary}

Enfin, on trouve aussi dans [La3] des formules exprimant les
fonctions elliptiques $\wp(n z)$ et ${\wp}'(n z)$ (pour un entier
$n \in {\mathbb N}^{*}$) en fonction de $\wp(z) , {\wp}'(z)$ et
des $f_m , m \in {\mathbb N}^{*}.$ Ces formules sont:
\begin{equation}
\begin{split}
\left\{
\begin{array}{lcl}
\wp(n z) & = & \wp(z) - \frac{f_{n + 1}(z).f_{n - 1}(z)}{{f_n (z)}^2} \\
{\wp}'(n z) & = & \frac{f_{2 n}(z)}{{f_n (z)}^4} =
\frac{1}{{\wp}'(z)}.\frac{f_{n + 2}(z){f_{n - 1}(z)}^2 - f_{n -
2}(z) {f_{n + 1}(z)}^2}{{f_n(z)}^3}
\end{array}
\right..
\end{split}\label{3.82}
\end{equation}
Ces formules (\ref{3.82}) donnent -a priori- des formes
représentant la multiplication d'un point de $E$ par un entier $n
\geq 1$, mais écrites en fonction des polyn\^omes $Q_m , m \geq
1$. On pourra donc estimer leurs degrés et hauteurs gr\^ace au
lemme \ref{c.55}. En effet, soient pour $n \geq 1$ les polyn\^omes
suivants de $K[X , Y]$:
\begin{equation}
\begin{split}
\begin{array}{ccc}
{\tilde F}_{0}^{(n)} (X , Y) & := & X Q_{n}^{3} (X , Y) -
Q_n Q_{n + 1} Q_{n - 1} (X , Y) \\
{\tilde F}_{1}^{(n)} (X , Y) & := & \frac{1}{Y} \left(Q_{n + 2}
Q_{n - 1}^{2} - Q_{n - 2} Q_{n + 1}^{2} \right) (X , Y) \\
{\tilde F}_{2}^{(n)} (X , Y) & := & Q_{n}^{3} (X , Y) .
\end{array}
\end{split}\label{3.83}
\end{equation}
On peut énoncer:
\begin{lemma}\label{c.57}
Pour tout point $\mathbf p := (x : y : 1)$ de $E \setminus \{0\}$
et tout entier $n \geq 1 ,$ le point $n.\mathbf p$ de $E$ est
représenté dans ${\mathbb P}_2$ par les coordonnées projectives:
$$n.\mathbf p = \left({{\tilde F}_{\!0}}^{(n)} (x, y) : {{\tilde F}_{\!1}}^{(n)} (x, y) : {{\tilde F}_{\!2}}^{(n)} (x, y)
\right).$$
\end{lemma}
{\bf Démonstration.---} Soit $\mathbf p := (x : y : 1) := (\wp(z)
: {\wp}'(z) : 1)$ avec $z \in {(\mathbb C / \Lambda)}^{*}$, un
point de $E \setminus \{0\}$ et soit $n$ un entier $\geq 1$. On
peut représenter le point $n.\mathbf p$ de $E$ dans l'espace
projectif ${\mathbb P}_2$ par:
\begin{equation*}
\begin{split}
n.\mathbf p  &=  \left(f_{n}^{3} (z) \wp(nz) : f_{n}^{3}
(z) {\wp}'(nz) : f_{n}^{3}(z)\right) \\
&\!\!\!\!\!\!\!\!\!= \left(f_{n}^{3} (z) \wp (z) - f_n (z) f_{n +
1}(z) f_{n - 1} (z) : \frac{f_{n + 2}(z) {f_{n - 1}(z)}^2 -
  f_{n - 2}(z) {f_{n + 1}(z)}^2} {{\wp}'(z)} : f_{n}^{3}(z)\!\!\right) \\
&\!\!\!\!\!\!\!\!\!=  \left(\!x Q_{n}^{3}(x , y) - Q_n Q_{n + 1} Q_{n - 1}(x , y) : \frac{1}{y} (Q_{n + 2} Q_{n - 1}^{2} - Q_{n - 2} Q_{n + 1}^{2})(x , y) : Q_{n}^{3} (x , y)\!\!\right) \\
&\!\!\!\!\!\!\!\!\!= \left({\tilde F}_{0}^{(n)}(x , y) : {\tilde
F}_{1}^{(n)}(x , y) : {\tilde F}_{2}^{(n)}(x , y)\right)
\end{split}
\end{equation*}
o\`u la deuxi\`eme égalité vient des formules (\ref{3.82}). Pour
compléter la preuve de notre lemme, il ne reste qu'à vérifier que
cette représentation du point $n.\mathbf p$ est bien définie dans
${\mathbb P}_2 ,$ c'est-à-dire que les trois expressions ${{\tilde
F}_{\!0}}^{(n)} (x , y) , {{\tilde F}_{\!1}}^{(n)} (x , y)$ et
${{\tilde F}_{\!2}}^{(n)} (x , y)$ ne peuvent s'annuler
simultanément pour un point $\mathbf p = (x : y : 1)$ de $E
\setminus \{0\}.$ Procédons par l'absurde. Supposons que pour un
certain point $\mathbf p = (x : y : 1) = (\wp(z) : {\wp}'(z) : 1)$
de $E \setminus \{0\}$ et pour un certain $n \geq 1$ on a :
${\tilde F}_{0}^{(n)} (x , y) = {\tilde F}_{1}^{(n)} (x , y) =
{\tilde F}_{2}^{(n)} (x , y) = 0,$ ceci revient à dire qu'on a:
$$Q_n(x , y) = \frac{1}{y} (Q_{n + 2} Q_{n - 1}^{2} -
    Q_{n - 2} Q_{n + 1}^{2}) (x , y) = 0,$$
    qui s'écrit en fonction de $z$:
    $$f_n(z) = \frac{1}{{\wp}'(z)} \left(f_{n + 2} f_{n - 1}^{2} - f_{n - 2}
    f_{n + 1}^{2}\right)(z) = 0$$
    et qu'on peut écrire aussi d'apr\`es les formules (\ref{3.82}):
    $$f_n(z) = \left(\frac{f_{2 n}}{f_n}\right)(z) = 0 .$$
Maintenant, on a d'une part:
$$f_n(z) = 0 ~\Leftrightarrow~ z \in {(\mathbb C / \Lambda)}_n, ~z \neq 0$$
et d'autre part, en reprenant l'expression (\ref{3.80}) pour
${f_n(z)}^2$ on a:
$$\left(\frac{f_{2 n}}{f_n}\right)^{\!\!2} \!\!(z) = 4\!\!\!\!\!\!\!\!\!\! \prod_{
\begin{array}{c}
\scriptstyle u \in {(\mathbb C / \Lambda)}_{2 n} \\
\scriptstyle u \not\in {(\mathbb C / \Lambda)}_n
\end{array}
}\!\!\!\!\!\!\!\!\! \left(\wp(z) - \wp(u)\right)$$ donc:
$$\frac{f_{2 n}}{f_n} (z) = 0 ~\Leftrightarrow~ z \in {(\mathbb C / \Lambda)}_{2 n} ,~ z \not\in {(\mathbb C / \Lambda)}_n
.$$ On voit ainsi que les deux fonctions elliptiques $f_n$ et
$\frac{f_{2 n}}{f_n}$ ne peuvent pas s'annuler simultanément, par
conséquent ${{\tilde F}_{\!0}}^{(n)} (X , Y) , {{\tilde
F}_{\!1}}^{(n)} (X , Y)$ et ${{\tilde F}_{\!2}}^{(n)} (X , Y)$ ne
peuvent pas s'annuler simultanément aussi. La représentation $n
\mathbf p = ({{\tilde F}_{\!0}}^{(n)} (X , Y) : {{\tilde
F}_{\!1}}^{(n)} (X , Y) : {{\tilde F}_{\!2}}^{(n)} (X , Y))$ est
effectivement bien définie sur ${\mathbb P}_2$ pour tout point
$\mathbf p = (x : y : 1) \in E \setminus \{0\}$ et tout entier $n
\geq 1$ ce qui ach\`eve cette démonstration.
$~~~~\blacksquare$\vspace{2mm}\\
$\underline{\text{\bf{Probl\`eme}}}$: $~~$ On peut voir facilement
que les polyn\^omes ${{\tilde F}_{\!0}}^{(n)} , {{\tilde
F}_{\!1}}^{(n)}$ et ${{\tilde F}_{\!2}}^{(n)} (n \geq 1)$ de $K[X
, Y]$, représentant la multiplication d'un point $\mathbf p$ de $E
\setminus \{0\}$ par l'entier positif $n$, sont de degré $\leq
\frac{3}{2} n^2$ et de hauteur de Gauss-Weil majorée par
$\frac{3}{4} (\eta + 4).n^2$ (qu'on estime gr\^ace au lemme
\ref{c.55}). Le probl\`eme est que lorsqu'on les homogénéise en
des formes $F_{0}^{(n)} , F_{1}^{(n)}$ et $F_{2}^{(n)}$ de $K[X ,
Y , Z]$, les formes homog\`enes obtenues sont non définies à
l'origine, ce qui est indésirable. Pour régler ce probl\`eme, nous
allons réduire (en un certain sens) les polyn\^omes ${{\tilde
F}_{\!0}}^{(n)} , {{\tilde F}_{\!1}}^{(n)}$ et ${{\tilde
F}_{\!2}}^{(n)} ~(n \geq 1)$ modulo l'équation affine de la courbe
elliptique $E$, qui est $Y^2 = 4 X^3 - g_2 X - g_3$ et nous
homogénéisons ensuite les polyn\^omes réduits ainsi obtenus qui
donnerons cette fois-ci des formes représentant globalement la
multiplication d'un point de $E$ par $n$. Il est important de
signaler que dans cette réduction on gagne un peu sur les degrés
(on obtient des formes de degré $n^2$, ce qui est optimal) mais on
perd sur l'estimation de la hauteur (on obtient une famille de
formes de hauteurs majorées par $\frac{3}{2} (\eta + 3).n^2$).
\subsubsection*{Réduction des polyn\^omes ${{\tilde
F}_{\!0}}^{(n)} , {{\tilde F}_{\!1}}^{(n)}$ et ${{\tilde
F}_{\!2}}^{(n)}$ modulo l'équation affine de $E$:} Nous
définissons $T$ l'application de réduction modulo l'équation
affine de la courbe elliptique $E$ par:
$$T : K[X , Y] \longrightarrow K[X , Y]$$
associant à chaque polyn\^ome $P$ de $K[X , Y]$ le polyn\^ome
$T(P)$ de $K[X , Y]$ dont le degré en $X$ est $\leq 2$ et qui
équivaut à $P$ modulo l'équation affine de $E$. Il est clair que
$T$ est unique et bien définie. Etant donné $P$ un polyn\^ome de
$K[X , Y]$, pour calculer concr\`etement $T(P)$ on proc\`ede comme
suit: \\ $\bullet$ si ${d°}_{\!\!\!\!X} P \leq 2$, on prend $T(P)
= P$,
\\ $\bullet$ sinon, $X^3$ divise forcément l'un au moins des
mon\^omes de $P,$ on remplace un tel $X^3$ par $\frac{1}{4} (Y^2 +
g_2 X + g_3)$ et on réit\`ere cette opération jusqu'à l'obtention
d'un polyn\^ome de degré en $X$ inférieur ou égal à $2$. C'est ce
dernier qu'on prend pour $T(P)$.
\\ Posons pour tout $n \in \mathbb N: T_n := T(X^n)$, les $T_n, n \in \mathbb
N$, sont donc des polyn\^omes de $K[X , Y]$ dont le degré en $X$
est $\leq 2$. Les premiers $T_n$ sont: $$T_0(X , Y) = 1 ,~ T_1(X ,
Y) = X ,~ T_2(X , Y) = X^2 ,$$$$ T_3(X , Y) = \frac{1}{4} (Y^2 +
g_2 X + g_3) = \frac{1}{4} g_2 X + \frac{1}{4} (Y^2 + g_3) ,
\ldots ~\mbox{etc}$$ Pour $n \in \mathbb N$, écrivons $T_n$ comme
polyn\^ome en $X$ (de degré $\leq 2$) à coefficients polyn\^omes
de $K[Y]$:
\begin{equation}
T_n(X , Y) = A_n(Y) X^2 + B_n(Y) X + C_n(Y) .\label{3.84}
\end{equation}
Les premiers éléments des suites ${(A_n)}_n , {(B_n)}_n$ et
${(C_n)}_n$ sont alors:$$
\begin{array}{ccl}
A_0 \equiv 0 , & A_1 \equiv 0 , & A_2 \equiv 1 ; \\
B_0 \equiv 0 , & B_1 \equiv 1 , & B_2 \equiv 0 ; \\
C_0 \equiv 1 , & C_1 \equiv 0 , & C_2 \equiv 0 .
\end{array}
$$
Nous établissons maintenant des relations de récurrence permettant
de calculer de proche en proche les $A_n , B_n , C_n ~(n \in
\mathbb N)$. On a pour tout $n \in \mathbb N$:
\begin{eqnarray*}
T_{n + 1} & := & T(X^{n + 1}) \\
          & = & T(X T_n) \\
          & = & T\left(A_n(Y) X^3 + B_n(Y) X^2 + C_n(Y) X \right)
          \\
          & = & A_n(Y) \left[\frac{1}{4} \left(Y^2 + g_2 X +
          g_3\right)\right] + B_n(Y) X^2 + C_n(Y) X \\
          & = & B_n(Y) X^2 + \left(\frac{1}{4} g_2 A_n(Y) +
          C_n(Y)\right) X + \frac{1}{4} \left(Y^2 + g_3\right) A_n(Y) ,\\
\end{eqnarray*}
d'o\`u:
\begin{equation}
\forall n \in \mathbb N : \left\{
\begin{array}{lcl}
A_{n + 1} & = & B_n \\
B_{n + 1} & = & \frac{1}{4} g_2 A_n + C_n \\
C_{n + 1} & = & \frac{1}{4} (Y^2 + g_3) A_n
\end{array}
\right.. \label{3.85}
\end{equation}
En utilisant ces formules, on déduira d'autres formules de
récurrence liant les termes de chacune des trois suites ${(A_n)}_n
, {(B_n)}_n$ et ${(C_n)}_n$ indépendamment des termes des deux
autres. On montre qu'on a pour tout $n \in \mathbb N$:
\begin{equation}
\left\{
\begin{array}{lcl}
A_{n + 3} & = & \frac{1}{4} g_2 A_{n + 1} + \frac{1}{4} (Y^2 +
g_3) A_n \\
B_{n + 3} & = & \frac{1}{4} g_2 B_{n + 1} + \frac{1}{4} (Y^2 +
g_3) B_n \\
C_{n + 3} & = & \frac{1}{4} g_2 C_{n + 1} + \frac{1}{4} (Y^2 +
g_3) C_n
\end{array}
\right.. \label{3.86}
\end{equation}
Maintenant, gr\^ace à ces formules de récurrence (\ref{3.85}) et
(\ref{3.86}), nous allons estimer les degrés et les hauteurs
locales des polyn\^omes $A_n , B_n$ et $C_n ~(n \in \mathbb N)$ de
$K[Y]$ et nous en déduisons plus généralement des estimations pour
les degrés et les hauteurs locales d'une réduction $T(P)$ d'un
polyn\^ome $P \in K[X]$ modulo l'équation affine de notre courbe
elliptique $E,$ en fonction du degré et des hauteurs locales de
$P;$ il ne reste apr\`es ça qu'à appliquer ces estimations aux
polyn\^omes ${{\tilde F}_{\!0}}^{(n)} , {{\tilde F}_{\!1}}^{(n)}$
et ${{\tilde F}_{\!2}}^{(n)} ~(n \in {\mathbb N}^*)$ pour en
déduire des estimations pour le degré, les hauteurs locales ainsi
que la hauteur de Gauss-Weil de leurs réductions $T({{\tilde
F}_{\!0}}^{(n)}) , T({{\tilde F}_{\!1}}^{(n)})$ et $T({{\tilde
F}_{\!2}}^{(n)}) ~(n \in {\mathbb N}^*)$ modulo l'équation affine
de $E$.
\subsubsection*{a)-Estimations sur les degrés:}
Le lemme suivant estime les degrés des polyn\^omes $A_n , B_n$ et
$C_n ~(n \in \mathbb N)$ et donne m\^eme les degrés exacts de ces
derniers ainsi que la valeur du coefficient dominant de chacun
lorsque $g_2 \neq 0$.
\begin{lemma}\label{c.58}
Pour tout $n \in {\mathbb N}^*$;
\begin{description}
\item[$\bullet$] $A_n$ est de degré $\leq 2 (n - 2 [\frac{n}{3}] -
2)$ et le coefficient de $Y^{2 (n - 2 [\frac{n}{3}] - 2)}$ dans
son écriture canonique vaut $ \left(\!
\begin{array}{c}
[\frac{n}{3}] \\
2 + 3 [\frac{n}{3}] - n
\end{array}
\!\right). {(\frac{1}{4})}^{[\frac{n}{3}]}\!.~ {g_2}^{2 + 3
[\frac{n}{3}] - n}$ \item[$\bullet$] $B_n$ est de degré $\leq 2 (n
- 2 [\frac{n + 1}{3}] - 1)$ et le coefficient de $Y^{2 (n - 2
[\frac{n + 1}{3}] - 1)}$ dans son écriture canonique vaut $
\left(\!
\begin{array}{c}
[\frac{n + 1}{3}] \\
1 + 3 [\frac{n + 1}{3}] - n
\end{array}
\!\right). {(\frac{1}{4})}^{[\frac{n + 1}{3}]}\!.~ {g_2}^{1 + 3
[\frac{n + 1}{3}] - n}$ \item[$\bullet$] $C_n$ est de degré $\leq
2 (n - 2 [\frac{n - 1}{3}] - 2)$ et le coefficient de $Y^{2 (n - 2
[\frac{n - 1}{3}] - 2)}$ dans son écriture canonique vaut $
\left(\!
\begin{array}{c}
[\frac{n - 1}{3}] \\
3 + 3 [\frac{n - 1}{3}] - n
\end{array}
\!\right). {(\frac{1}{4})}^{[\frac{n - 1}{3}] + 1}\!.~ {g_2}^{3 +
3 [\frac{n - 1}{3}] - n}$.
\end{description}
\end{lemma}
{\bf Démonstration.---} Pour obtenir les estimations concernant
les $A_n ~(n \in \mathbb N),$ on proc\`ede par récurrence en
utilisant la relation: $A_{n + 3} = \frac{1}{4} g_2 A_{n + 1} +
\frac{1}{4} (Y^2 + g_3) A_n ~(n \in \mathbb N)$ de (\ref{3.86}) et
cette récurrence se fait en distinguant les trois cas: $n = 3 k ,~
n = 3 k + 1$ et $n = 3 k + 2 ~(k \in \mathbb N).$ Les estimations
concernant les $B_n$ et les $C_n ~(n \in \mathbb N)$ peuvent se
déduire ensuite directement de celles des $A_n$ -sans refaire la
récurrence- gr\^ace aux deux relations: $B_n = A_{n + 1}$ et $C_n
= \frac{1}{4} (Y^2 + g_3) A_{n - 1} ~(n \in {\mathbb N}^*)$ de
(\ref{3.85}).  $~~~~\blacksquare$\vspace{2mm}

Le lemme qui suit est une conséquence du lemme \ref{c.58} précédent, il donne le
degré total et le degré en $Y$ d'une réduction $T(P)$ d'un
polyn\^ome $P \in K[X]$ modulo l'équation affine de $E$ et donne
aussi un mon\^ome significatif de l'écriture canonique de $T(P)$.
\begin{lemma}\label{c.59}
Pour tout polyn\^ome $P$ de $K[X]$ de degré $n \in \mathbb N$ et
de coefficient dominant $a_n \in K^*$, sa réduction $T(P)$ modulo
l'équation affine de $E$ vérifie:
$${d°}_{\!\!\! \rm{tot}} T(P) = n - [\frac{n}{3}] ~,~ {d°}_{\!\!\! \rm{Y}} T(P) = 2 [\frac{n}{3}]$$
et $T(P)$ contient dans son écriture canonique le mon\^ome $$ a_n
\left(\frac{1}{4}\right)^{\!\![\frac{n}{3}]} X^{n - 3
[\frac{n}{3}]} Y^{2 [\frac{n}{3}]} .$$
\end{lemma}
{\bf Démonstration.---} En utilisant le lemme \ref{c.58} et en
distinguant les trois cas: $n = 3 k ,~ n = 3 k + 1$ et $n = 3 k +
2$ on montre qu'on a pour tout $n \in \mathbb N:$
$${d°}_{\!\!\! \rm{tot}} T_n \leq n - [\frac{n}{3}] ~,~ {d°}_{\!\!\! \rm{Y}} T_n \leq 2 [\frac{n}{3}]$$
et que $T_n$ contient dans son écriture canonique le mon\^ome:
$$\left(\frac{1}{4}\right)^{\!\![\frac{n}{3}]} X^{n - 3 [\frac{n}{3}]} Y^{2 [\frac{n}{3}]} .$$
Une remarque, utile pour la suite de cette démonstration, qu'on
peut tirer immédiatement de ce qui préc\`ede est que pour tout
couple $(k , \ell) \in {\mathbb N}^2 , ~(k < \ell)$, le
coefficient de $X^{\ell - 3 [\frac{\ell}{3}]} Y^{2
[\frac{\ell}{3}]}$ dans l'écriture canonique de $T_k$ est nul. En
effet si ce coefficient était non nul, on aurait: $${d°}_{\!\!\!
\rm{tot}} T_k \geq \ell -
[\frac{\ell}{3}]~~\mbox{et}~~{d°}_{\!\!\! \rm{Y}} T_k \geq 2
[\frac{\ell}{3}] ,$$ c'est-à-dire: $$k - [\frac{k}{3}] \geq \ell -
[\frac{\ell}{3}]~~\mbox{et}~~2 [\frac{k}{3}] \geq 2
[\frac{\ell}{3}]$$ ce qui est clairement impossible puisque $k <
\ell$. Donc ce coefficient est effectivement nul. En écrivant
maintenant $\displaystyle P(X) =: \sum_{i = 0}^{n} a_i X^i ~(a_n
\neq 0)$, on a:
$$
T(P) ~=~ \sum_{i = 0}^{n} a_i T(X^i) ~=~ \sum_{i = 0}^{n} a_i T_i
$$
$$
\begin{array}{rclcl}
\mbox{d'o\`u:} ~~~~~~~~~\!\!~ {d°}_{\!\!\! \rm{tot}} T(P) & \leq &
\max\left({d°}_{\!\!\!
\rm{tot}} T_i ~,~ i = 0 , \ldots , n\right) & \leq & n - [\frac{n}{3}] , \\
{d°}_{\!\!\! \rm{Y}} T(P)~ & \leq & \max\left({d°}_{\!\!\!Y} T_i
~,~ i = 0 , \ldots , n\right) & \leq & 2[\frac{n}{3}]
\end{array}~~~~~~~~
$$
et, de plus, d'apr\`es la remarque précédente $T(P)$ contient dans
son écriture canonique le mon\^ome $a_n
{(\frac{1}{4})}^{[\frac{n}{3}]} X^{n - 3 [\frac{n}{3}]} Y^{2
[\frac{n}{3}]}$, donc on a bien ${d°}_{\!\!\! \rm{tot}} T(P) = n -
[\frac{n}{3}]$ et ${d°}_{\!\!\! \rm{Y}} T(P) = 2 [\frac{n}{3}]$ ce
qui ach\`eve cette démonstration.$~~~~\blacksquare$\vspace{2mm}

Nous passons enfin, dans le corollaire qui suit, à l'estimation du
degré total et du degré partiel en $Y$ des réductions respectives
$T({{\tilde F}_{\!0}}^{(n)}) , T({{\tilde F}_{\!1}}^{(n)})$ et
$T({{\tilde F}_{\!2}}^{(n)})$ des polyn\^omes ${{\tilde
F}_{\!0}}^{(n)} , {{\tilde F}_{\!1}}^{(n)}$ et ${{\tilde
F}_{\!2}}^{(n)} ~(n \in \mathbb N)$ modulo l'équation affine de
$E$, en donnant en m\^eme temps un mon\^ome significatif
l'écriture canonique de $T({{\tilde F}_{\!1}}^{(n)})$.
\begin{corollary}\label{c.60}
Pour tout $n \in \mathbb N$, les polyn\^omes réduits $T({{\tilde
F}_{\!0}}^{(n)}) , T({{\tilde F}_{\!1}}^{(n)})$ et $T({{\tilde
F}_{\!2}}^{(n)})$ des polyn\^omes ${{\tilde F}_{\!0}}^{(n)} ,
{{\tilde F}_{\!1}}^{(n)}$ et ${{\tilde F}_{\!2}}^{(n)}$
respectivement, modulo l'équation affine de $E$, sont de degré
total $\leq n^2$ et de degré partiel en $Y$ strictement inférieur
à $n^2$ pour $T({{\tilde F}_{\!0}}^{(n)})$ et $T({{\tilde
F}_{\!2}}^{(n)})$ et égal à $n^2$ pour $T({{\tilde
F}_{\!1}}^{(n)})$. De plus, $T({{\tilde F}_{\!1}}^{(n)})$ contient
dans son écriture canonique le mon\^ome ${(\frac{1}{2})}^{n^2 - 1}
Y^{n^2}$.
\end{corollary}
{\bf Démonstration.---} On remplace dans les formules (\ref{3.83})
de définition des polyn\^omes ${{\tilde F}_{\!0}}^{(n)} , {{\tilde
F}_{\!1}}^{(n)}$ et ${{\tilde F}_{\!2}}^{(n)}$ -selon les cas $n$
impair et $n$ pair- les polyn\^omes $Q_k ~(k \in \mathbb N)$ de
$K[X , Y]$ par leurs expressions en fonction des polyn\^omes $P_k
~(k \in \mathbb N)$ de $K[X]$. Les degrés et les coefficients
dominants de ces derniers étant connus, il suffit d'utiliser le
lemme \ref{c.59} pour avoir toutes les assertions du corollaire
\ref{c.60}.  $~~~~\blacksquare$\vspace{1mm}
\subsubsection*{b)-Estimations des hauteurs locales et de la hauteur de Gauss-Weil:}
Les estimations arithmétiques des polyn\^omes de $K[Y]$ ou de $K[X
, Y]$, utilisent dans le cas fini la hauteur $v$-adique $H_v$
alors que dans le cas infini nous prenons la longueur $v$-adique
$L_v$, qui est plus facile à gérer pour ce cas et qui
entra{\sf\^\i}ne une estimation pour $H_v$ puisque on a: $H_v \leq
L_v$. Le lemme suivant donne des estimations pour les hauteurs
$v$-adique (resp longueurs $v$-adique) des polyn\^omes $A_n , B_n$
et $C_n ~(n \in {\mathbb N}^*)$ en fonction de $M_v$ et de $n$
quand $v$ est une place finie (resp infinie).
\begin{lemma}\label{c.61}
Soit $v$ une place de $K$. Pour tout $n \in {\mathbb N}^*$ on a:
\\ $\bullet$ quand $v$ est infinie:
\begin{eqnarray*}
L_v(A_n) & \leq & {M_v}^{[\frac{n - 2}{2}]} \\
\mbox{et}~~ \max\left(L_v(B_n) , L_v(C_n)\right) & \leq &
{M_v}^{[\frac{n - 1}{2}]}, \\
\end{eqnarray*}
$\bullet$ quand $v$ est finie et $v \mid 2:$
\begin{eqnarray*}
H_v(A_n) & \leq & {4 M_v}^{[\frac{n - 2}{2}]} \\
\mbox{et}~~ \max\left(H_v(B_n) , H_v(C_n)\right) & \leq &
{4 M_v}^{[\frac{n - 1}{2}]} \\
\end{eqnarray*}
$\bullet$ et quand $v$ est finie et $v \nmid 2:$
\begin{eqnarray*}
H_v(A_n) & \leq & {M_v}^{[\frac{n - 2}{2}]} \\
\mbox{et}~~ \max\left(H_v(B_n) , H_v(C_n)\right) & \leq &
{M_v}^{[\frac{n - 1}{2}]} . \\
\end{eqnarray*}
\end{lemma}
{\bf Démonstration.---} Les estimations des hauteurs $v$-adiques
(resp des longueurs $v$-adiques) des polyn\^omes $A_n ~(n \geq 1)$
de $K[Y]$ quand $v$ est finie (resp quand $v$ est infinie)
s'obtiennent par un procédé de récurrence sur $n$ en utilisant la
relation: $A_{n + 3} = \frac{1}{4} g_2 A_{n + 1} + \frac{1}{4}
(Y^2 + g_3) A_n ~(n \in \mathbb N)$ donnée par (\ref{3.86}).
Ensuite, pour en déduire les estimations analogues concernant les
polyn\^omes $B_n$ et $C_n ~(n \geq 1)$, il suffit d'utiliser les
deux relations: $B_n = A_{n + 1}$ et $C_n = \frac{1}{4} (Y^2 +
g_3) A_{n - 1} ~(n \geq 1)$ données par (\ref{3.85}) et lemme
s'ensuit.$~~~~\blacksquare$\vspace{2mm}

Comme conséquence du lemme \ref{c.61}, le lemme qui suit donne -lorsque $P$ est un
polyn\^ome de $K[Y]$ et $v$ une place finie (resp infinie) de $K$-
une estimation pour la hauteur $v$-adique (resp la longueur
$v$-adique) de sa réduction $T(P)$ modulo l'équation affine de
$E$, en fonction de sa hauteur $v$-adique (resp de sa longueur
$v$-adique), de son degré et de $M_v$.
\begin{lemma}\label{c.62}
Soit $v$ une place de $K$, pour tout polyn\^ome $P$ de $K[X]$ de
degré $n \geq 1,$ sa réduction $T(P)$ modulo l'équation affine de
$E$ satisfait:
\\ $\bullet$ quand $v$ est infinie:
$$L_v(T(P)) \leq 3 L_v(P) {M_v}^{[\frac{n - 1}{2}]},$$
$\bullet$ quand $v$ est finie et $v \mid 2:$
$$H_v(T(P)) \leq H_v(P) {(4 M_v)}^{[\frac{n - 1}{2}]}$$
$\bullet$ et quand $v$ est finie et $v \nmid 2:$
$$H_v(T(P)) \leq H_v(P) {M_v}^{[\frac{n - 1}{2}]}.$$
\end{lemma}
{\bf Démonstration.---} Ecrivons:
$$P(X) =: \sum_{i = 0}^{n} a_i X^i ~~(a_n \neq 0),$$
d'o\`u
\begin{eqnarray*}
T(P) & = & \sum_{i = 0}^{n} a_i T(X^i) \\
     & = & \sum_{i = 0}^{n} a_i T_i \\
     & = & \sum_{i = 0}^{n} a_i \left[A_i(Y) X^2 + B_i(Y) X + C_i(Y)\right]
\end{eqnarray*}
et d'apr\`es les propriétés des longueurs et des hauteurs locales:
\\ $\bullet$ quand $v$ est infinie:
$$L_v(T(P)) \leq L_v(P) . \max_{0 \leq i \leq n} \left(L_v(A_i) + L_v(B_i) + L_v(C_i)\right)$$
$\bullet$ et quand $v$ est finie:
$$H_v(T(P)) \leq H_v(P) . \max_{0 \leq i \leq n} \left(\max\{H_v(A_i) , H_v(B_i) , H_v(C_i)\}\right).$$
Il ne reste qu'à appliquer le lemme \ref{c.61} pour conclure.$~~~~\blacksquare$\vspace{2mm}

Comme application du lemme \ref{c.62}, le lemme qui suit estime les hauteurs $v$-adiques
(resp les longueurs $v$-adiques) des réductions $T({{\tilde
F}_{\!0}}^{(n)}) , T({{\tilde F}_{\!1}}^{(n)})$ et $T({{\tilde
F}_{\!2}}^{(n)})$ des polyn\^omes ${{\tilde F}_{\!0}}^{(n)} ,
{{\tilde F}_{\!1}}^{(n)}$ et ${{\tilde F}_{\!2}}^{(n)} ~(n \in
\mathbb N)$, modulo l'équation affine de $E$, quand $v$ est une
place finie (resp infinie).
\begin{lemma}\label{c.63}
Soit $v$ une place de $K$, pour tout $n \in {\mathbb N}^*$ on a:
\\ $\bullet$ quand $v$ est infinie:
\begin{equation*}
\begin{split}
\max\!\left\{L_v\left(T({{\tilde F}_{\!0}}^{(n)})\right) \right.&
\left., L_v\left(T({{\tilde F}_{\!1}}^{(n)})\right) ,
 L_v\left(T({{\tilde F}_{\!2}}^{(n)})\right)\right\}~~~~~~~~~~~~~~~~~~~~ \\
 &\leq \left\{\!\!
 \begin{array}{ll}
 3 . {(2 M_v)}^{\frac{3 n^2 - 3}{2}} & \mbox{si} ~n~ \mbox{est
 impair} \\
 6 . {(2 M_v)}^{\frac{3 n^2 - 2}{2}} & \mbox{si} ~n~ \mbox{est pair}
 \end{array}
 \right.,
\end{split}
\end{equation*}
 $\bullet$ quand $v$ est finie et $v \mid 2$:
\begin{equation*}
\begin{split}
\max\!\left\{H_v\left(T({{\tilde F}_{\!0}}^{(n)})\right) \right.&
\left., H_v\left(T({{\tilde F}_{\!1}}^{(n)})\right) ,
 H_v\left(T({{\tilde F}_{\!2}}^{(n)})\right)\right\}~~~~~~~~~~~~~~~~~~~~ \\
&\leq \left\{\!\!
 \begin{array}{ll}
 {(4 M_v)}^{\frac{3 n^2 - 3}{2}} & \mbox{si} ~n~ \mbox{est
 impair} \\
 {(4 M_v)}^{\frac{3 n^2 - 2}{2}} & \mbox{si} ~n~ \mbox{est pair}
 \end{array}
 \right.
\end{split}
\end{equation*}
 $\bullet$ et quand $v$ est finie et $v \nmid 2$:
\begin{equation*}
\begin{split}
\max\!\left\{H_v\left(T({{\tilde F}_{\!0}}^{(n)})\right) \right.&
\left., H_v\left(T({{\tilde F}_{\!1}}^{(n)})\right) ,
 H_v\left(T({{\tilde F}_{\!2}}^{(n)})\right)\right\}~~~~~~~~~~~~~~~~~~~~ \\
&\leq \left\{\!\!
 \begin{array}{ll}
 {M_v}^{\frac{3 n^2 - 3}{2}} & \mbox{si} ~n~ \mbox{est
 impair} \\
 {M_v}^{\frac{3 n^2 - 2}{2}} & \mbox{si} ~n~ \mbox{est pair}
 \end{array}
 \right..
\end{split}
\end{equation*}
\end{lemma}
{\bf Démonstration.---} Pour $n = 1, 2 , 3$ on vérifie les
estimations du lemme directement sur les formules de $T({{\tilde
F}_{\!0}}^{(n)}) , T({{\tilde F}_{\!1}}^{(n)})$ et $T({{\tilde
F}_{\!2}}^{(n)})$ (données à la fin de ce sous-chapitre) et pour
$n \geq 4$, on remplace dans les formules de définition
(\ref{3.83}) des ${{\tilde F}_{\!0}}^{(n)} , {{\tilde
F}_{\!1}}^{(n)}$ et ${{\tilde F}_{\!2}}^{(n)}$ -selon le cas $n$
impair ou $n$ pair- les polyn\^omes $Q_m ~(m \in \mathbb N)$ de
$K[X , Y]$ par leurs expressions en fonction des polyn\^omes $P_m
~(m \in \mathbb N)$. Les hauteurs $v$-adiques (resp longueurs
$v$-adiques) des $P_m ~(m \in \mathbb N)$ étant estimées dans le
corollaire \ref{c.56} quand $v$ est une place finie (resp
infinie), les degrés des $P_m ~(m \in \mathbb N)$ étant aussi
connus, l'application du lemme \ref{c.62} permet alors de
conclure.$~~~~\blacksquare$\vspace{2mm}

Enfin étant donné un entier $n \geq 1,$ nous tirons du lemme \ref{c.63} une estimation
pour la hauteur de Gauss-Weil de la famille des polyn\^omes
$T({{\tilde F}_{\!0}}^{(n)}) , T({{\tilde F}_{\!1}}^{(n)})$ et
$T({{\tilde F}_{\!2}}^{(n)})$ réduits des polyn\^omes ${{\tilde
F}_{\!0}}^{(n)} , {{\tilde F}_{\!1}}^{(n)}$ et ${{\tilde
F}_{\!2}}^{(n)}$ respectivement, modulo l'équation affine de $E$.
Cette estimation dépend de $n$ et de $\eta$.
\begin{corollary}\label{c.64}
Pour tout entier $n \geq 1$ on a:
$$\widetilde h \left(T({{\tilde F}_{\!0}}^{(n)}) , T({{\tilde F}_{\!1}}^{(n)}) , T({{\tilde F}_{\!2}}^{(n)})\right)
 \leq \frac{3}{2} (\eta + 3) n^2 .$$
\end{corollary}
{\bf Démonstration.---} Il suffit d'utiliser la définition de la
hauteur de Gauss-Weil d'une famille finie de polyn\^omes:
\begin{equation*}
\begin{split}
\widetilde h \left(T({{\tilde F}_{\!0}}^{(n)}) , T({{\tilde
F}_{\!1}}^{(n)}) ,
T({{\tilde F}_{\!2}}^{(n)})\right) & \\
&\!\!\!\!\!\!\!\!\!\!\!\!\!\!\!\!\!\!\!\!\!\!\!\!\!\!\!\!\!\!\!\!\!\!\!\!\!\!\!\!\!\!\!\!\!\!\!\!\!\!\!\!\!\!:=
\sum_{v \in M_K} \frac{[K_v : {\mathbb Q}_v]}{[K : \mathbb Q]} .
\log \max\left\{ H_v\left(T({{\tilde F}_{\!0}}^{(n)})\right) ,
H_v\left(T({{\tilde F}_{\!1}}^{(n)})\right) , H_v\left(T({{\tilde
F}_{\!2}}^{(n)})\right)\right\} ,
 \end{split}
\end{equation*}
de majorer dans le cas $v$ infinie $H_v$ par $L_v$ et d'utiliser
les estimation du lemme \ref{c.63} en remarquant qu'on a:
$\displaystyle \sum_{v \in M_{K}^{\infty}} [K_v : {\mathbb Q}_v] =
\sum_{v \mid 2} [K_v : {\mathbb Q}_v] = [K : \mathbb Q]$ et
$\displaystyle \sum_{v \in M_K} \frac{[K_v : {\mathbb Q}_v]}{[K :
\mathbb Q]} . \log M_v$ \\ $=\eta$. Nous avons ainsi:
\begin{equation*}
\begin{split}
\widetilde h \left(T({{\tilde F}_{\!0}}^{(n)}) , T({{\tilde
F}_{\!1}}^{(n)}) , T({{\tilde F}_{\!2}}^{(n)})\right) &\leq
\left\{\!\!
\begin{array}{ll}
\frac{3 n^2 - 3}{2} (\eta + \log 8) + \log 3 & \mbox{si} ~n~ \mbox{est impair} \\
\frac{3 n^2 - 2}{2} (\eta + \log 8) + \log 6 & \mbox{si} ~n~
\mbox{est pair}
\end{array}
\right. \\
&\leq \frac{3}{2} (\eta + 3) n^2 ~~~~ (\forall n \geq 1) .
\end{split}
\end{equation*}
Ceci démontre l'estimation du corollaire.$~~~~\blacksquare$\vspace{1mm}
\subsubsection*{Démonstration du théor\`eme \ref{c.9}:}
Soit $n$ un entier $\geq 1$, les polyn\^omes de $K[X , Y] :
T({{\tilde F}_{\!0}}^{(n)}) , T({{\tilde F}_{\!1}}^{(n)})$ et
$T({{\tilde F}_{\!2}}^{(n)})$ réduits des polyn\^omes ${{\tilde
F}_{\!0}}^{(n)} , {{\tilde F}_{\!1}}^{(n)}$ et ${{\tilde
F}_{\!2}}^{(n)}$ respectivement, modulo l'équation affine de $E$,
sont -d'apr\`es le corollaire \ref{c.60}- de degré total $\leq
n^2$. Prenons respectivement pour $F_{0}^{(n)} , F_{1}^{(n)}$ et
$F_{2}^{(n)}$ leurs homogénéisés en des formes de $K[X , Y , Z]$
de degré $n^2$, c'est-à-dire posons:
\begin{eqnarray*}
F_{0}^{(n)} (X , Y , Z) & := & Z^{n^2}. T({{\tilde
F}_{\!0}}^{(n)})\left(\frac{X}{Z} ~,~ \frac{Y}{Z}\right), \\
F_{1}^{(n)} (X , Y , Z) & := & Z^{n^2}. T({{\tilde
F}_{\!1}}^{(n)})\left(\frac{X}{Z} ~,~ \frac{Y}{Z}\right) \\
\mbox{et}~~ F_{2}^{(n)} (X , Y , Z) & := & Z^{n^2}. T({{\tilde
F}_{\!2}}^{(n)})\left(\frac{X}{Z} ~,~ \frac{Y}{Z}\right) .
\end{eqnarray*}
Posons aussi ${\underline F}^{(n)} := (F_{0}^{(n)} , F_{1}^{(n)} ,
F_{2}^{(n)})$ la famille constituée de ces formes, il est clair
que ${\underline F}^{(n)}$ est une famille de formes de $K[X , Y ,
Z]$ de degré $n^2$ chacune et pour tout point $\mathbf p := (x : y
: z)$ de $E \hookrightarrow {\mathbb P}_2$, le point
$\left(F_{0}^{(n)} (x , y , z) : F_{1}^{(n)} (x , y , z) :
F_{2}^{(n)} (x , y , z)\right)$ représente dans ${\mathbb P}_2$ le
point $n.\mathbf p$ de $E$ sauf peut \^etre quand $\mathbf p$ est
le point à l'infini de $E$. Ceci résulte du lemme \ref{c.57} et du
fait que les polyn\^omes $T({{\tilde F}_{\!0}}^{(n)}) , T({{\tilde
F}_{\!1}}^{(n)})$ et $T({{\tilde F}_{\!2}}^{(n)})$ sont équivalent
modulo l'équation affine de $E$, aux polyn\^omes ${{\tilde
F}_{\!0}}^{(n)} , {{\tilde F}_{\!1}}^{(n)}$ et ${{\tilde
F}_{\!2}}^{(n)}$ respectivement. Montrons maintenant que m\^eme
lorsque $\mathbf p := (x : y : z)$ est le point à l'infini
$\underline 0 := (0 : 1 : 0)$ de $E$, le point projectif
$\left(F_{0}^{(n)} (x , y , z) : F_{1}^{(n)} (x , y , z) :\right.$ \\ $\left.
F_{2}^{(n)} (x , y , z)\right)$ est bien défini et représente
toujours dans ${\mathbb P}_2$ le point $n.\mathbf p$ de $E$. Comme
$n.\mathbf{0} = \mathbf{0}$, cela revient simplement à montrer
qu'on a: $F_{0}^{(n)} (0 , 1 , 0) = F_{2}^{(n)} (0 , 1 , 0) = 0$
et $F_{2}^{(n)} (0 , 1 , 0) \neq 0$, or ceci résulte immédiatement
du corollaire \ref{c.60}. En effet d'apr\`es ce dernier les deux
polyn\^omes $T({{\tilde F}_{\!0}}^{(n)})$ et $T({{\tilde
F}_{\!2}}^{(n)})$ sont de degré total $\leq n^2$ et de degré
partiel en $Y$ strictement inférieur à $n^2$. Donc leurs
homogénéisés $F_{0}^{(n)}$ et $F_{2}^{(n)}$ ne contiennent -dans
leurs écritures canoniques- que des mon\^omes de $K[X , Y , Z]$
dont l'un au moins des degrés en $X$ ou en $Z$ est non nul, ce qui
entra{\sf\^\i}ne qu'on a bien: $F_{0}^{(n)} (0 , 1 , 0) =
F_{2}^{(n)} (0 , 1 , 0) = 0$. Quant à $T({{\tilde
F}_{\!1}}^{(n)}),$ il est de degré total égal à son degré partiel
en $Y$ et égal à $n^2$, de plus il contient -dans son écriture
canonique- le mon\^ome ${(\frac{1}{2})}^{n^2 - 1} Y^{n^2}$, par
conséquent son homogénéisé $F_{1}^{(n)}$ ne contient -dans son
écriture canonique- que le mon\^ome ${(\frac{1}{2})}^{n^2 - 1}
Y^{n^2}$ et d'autre mon\^omes de $K[X , Y , Z]$ dont l'un au moins
des degrés en $X$ ou en $Z$ est non nul, ce qui montre bien qu'on
a: $F_{1}^{(n)} (0 , 1 , 0) = {(\frac{1}{2})}^{n^2 - 1} \neq 0$,
d'o\`u: $(F_{0}^{(n)} (\underline 0) : F_{1}^{(n)} (\underline 0)
: F_{2}^{(n)} (\underline 0)) = (0 : 1 : 0) = \underline 0$. On
vient de montrer que la famille de formes ${\underline F}^{(n)} =
(F_{0}^{(n)} , F_{1}^{(n)} , F_{2}^{(n)})$ représente globalement
la multiplication par $n$ d'un point de $E \hookrightarrow
{\mathbb P}_2$. Par ailleurs, la hauteur de Gauss-Weil de la
famille de formes ${\underline F}^{(n)}$ est évidemment la m\^eme
que celle des déshomogénéisés $T({{\tilde F}_{\!0}}^{(n)}) ,
T({{\tilde F}_{\!1}}^{(n)})$ et $T({{\tilde F}_{\!2}}^{(n)})$
qu'on a majoré avant dans
le corollaire \ref{c.64} par $\frac{3}{2} (\eta + 3) n^2$. Ceci conclut la démonstration du théor\`eme. \\

Avant de conclure ce sous-chapitre, donnons explicitement les
formes constituant les deux familles ${\underline F}^{(2)}$ et
${\underline F}^{(3)}$; cela nous donne ainsi des formules de
duplication et de triplication sur $E$. Pour aboutir à ces
formules on a besoin seulement des polyn\^omes $Q_i ~(i = 0 , 1 ,
\dots , 4)$ qui sont donnés tout au début de ce sous-chapitre et
du polyn\^ome $Q_5$ qui est égal à $Q_4 Q_{2}^{3} - Q_1 Q_{3}^{3}$
d'apr\`es les relations (\ref{3.81}). Gr\^ace au formules de
définition (\ref{3.83}), on calcule les polyn\^omes de $K[X , Y] :
{{\tilde F}_{\!0}}^{(n)} , {{\tilde F}_{\!1}}^{(n)}$ et ${{\tilde
F}_{\!2}}^{(n)}$ pour $n = 2$ et $n = 3$, on les réduit ensuite
modulo l'équation affine de $E$ par l'application de la réduction
$T$ définie à la page $83$ et on homogénéise enfin les polyn\^omes
réduits obtenus $T({{\tilde F}_{\!0}}^{(n)}) , T({{\tilde
F}_{\!1}}^{(n)})$ et $T({{\tilde F}_{\!2}}^{(n)})$ pour $n = 2$ et
$n = 3$. Voilà ce qu'on obtient en faisant tout ces calculs:
\subsubsection*{$\bullet$ Formules de duplication sur $E$:}
On peut représenter globalement la duplication sur $E$ par les
formes suivantes:
\begin{eqnarray*}
F_{0}^{(2)} (X , Y , Z) & := & \frac{1}{4} X Y^3 + \frac{3}{4} g_2
X^2 Y Z + \frac{9}{4} g_3 X Y Z^2 +
\frac{1}{16} g_{2}^{2} Y Z^3 \\
F_{1}^{(2)} (X , Y , Z) & := & \frac{1}{8} Y^4 - \frac{3}{8} g_2 X
Y^2 Z - \frac{9}{4} g_3 Y^2 Z^2 - \frac{9}{8} g_{2}^{2} X^2
Z^2 \\
 & & - \frac{27}{8} g_2 g_3 X Z^3 - \frac{27}{8} g_{3}^{2} Z^4 +
\frac{1}{32} g_{2}^{3} Z^4 \\
F_{2}^{(2)} (X , Y , Z) & := & Y^3 Z .
\end{eqnarray*}
\subsubsection*{$\bullet$ Formules de triplication sur $E$:}
On peut représenter globalement la triplication sur $E$ par les
formes suivantes:
\begin{eqnarray*}
F_{0}^{(3)} (X , Y , Z) & := & \frac{3}{256} X Y^8 +
\frac{21}{128} g_2 X^2 Y^6 Z + \frac{9}{8} g_3 X Y^6 Z^2 +
\frac{125}{1024} g_{2}^{2} Y^6 Z^3 \\
 & & + \frac{81}{128} g_2 g_3 X^2 Y^4 Z^3 - \left(\frac{3}{256} g_{2}^{3} + \frac{81}{128} g_{3}^{2}\right) X
 Y^4 Z^4 \\
  & & - \frac{837}{1024} g_{2}^{2} g_3 Y^4 Z^5 - \left(\frac{147}{1024} g_{2}^{4} + \frac{1053}{128} g_2 g_{3}^{2}\right)
  X^2 Y^2 Z^5 \\
   & & - \left(\frac{45}{32} g_{2}^{3} g_3 + \frac{243}{32}
   g_{3}^{3}\right) X Y^2 Z^6 - \frac{27}{64} g_{2}^{3} X^7 Z^2 \\
   & &  - \left(\frac{405}{256} g_{2}^{3} {g_3}^2 + \frac{37}{4096} g_{2}^{6}
    + \frac{729}{256} g_{3}^{4}\right) X Z^8 - \left(\frac{729}{128} g_2 g_{3}^{3} \right.\\
    & & + \left.\frac{459}{1024} g_3 g_{2}^{4}\right) X^2 Z^7 - \left(\frac{7}{1024} g_{2}^{5} + \frac{2241}{1024} g_{2}^{2}
    g_{3}^{2}\right) Y^2 Z^7  \\
    & & - \left(\frac{1215}{1024} g_{2}^{2} g_{3}^{3} +
    \frac{9}{1024} g_{2}^{5} g_3\right) Z^9
\end{eqnarray*}
\begin{eqnarray*}
F_{1}^{(3)} (X , Y , Z) & := & \frac{1}{256} Y^9 - \frac{9}{128}
g_2 X Y^7 Z - \frac{27}{32} g_3 Y^7 Z^2 -
\frac{225}{256} g_{2}^{2} X^2 Y^5 Z^2 \\
 & & - \frac{621}{128} g_2 g_3 X Y^5 Z^3 - \left(\frac{75}{256} g_{2}^{3} + \frac{1215}{128} g_{3}^{2}\right)Y^5 Z^4 \\
  & & - \frac{1215}{128} g_{2}^{2} g_3 X^2 Y^3 Z^4 - \left(\frac{579}{1024} g_{2}^{4} + \frac{2511}{128} g_2 g_{3}^{2}\right)
   X Y^3 Z^5 \\
   & & - \left(\frac{423}{256} g_{2}^{3} g_3 + \frac{243}{16}
   g_{3}^{3}\right) Y^3 Z^6 - \left(\frac{3645}{256} g_{2}^{2} g_{3}^{2}
   \right. \\
   & & - \left.\frac{189}{1024} g_{2}^{5}\right) X^2 Y Z^6 -\left(\frac{1539}{1024} g_{2}^{4} g_3 + \frac{2187}{128}
    g_2 g_{3}^{3}\right) X Y Z^7 \\
    & & - \left(\frac{3}{4096} g_{2}^{6} + \frac{2187}{256} g_{3}^{4} + \frac{81}{64} g_{2}^{3} g_{3}^{2}\right) Y Z^8
\end{eqnarray*}
\begin{eqnarray*}
F_{2}^{(3)} (X , Y , Z) & := & \frac{27}{256} Y^8 Z -
\frac{27}{128} g_2 X Y^6 Z^2 - \frac{27}{32} g_3 Y^6 Z^3 -
\frac{27}{256} g_{2}^{2} X^2 Y^4 Z^3 \\
 & & + \frac{81}{128} g_2 g_3 X Y^4 Z^4 + \left(\frac{243}{128} g_{3}^{2} + \frac{27}{256} g_{2}^{3}\right) Y^4 Z^5 \\
  & & + \frac{243}{128} g_{2}^{2} g_3 X^2 Y^2 Z^5 + \left(\frac{63}{1024} g_{2}^{4} + \frac{243}{128} g_2
  g_{3}^{2}\right) X Y^2 Z^6 \\
  & & - \frac{81}{256} g_{2}^{3} g_3 Y^2 Z^7 - \left(\frac{1215}{256} g_{2}^{2} g_{3}^{2} + \frac{63}{1024} g_{2}^{5}\right)
   X^2 Z^7 \\
   & & - \left(\frac{729}{128} g_2 g_{3}^{3} + \frac{513}{1024} g_{2}^{4}
   g_3\right) X Z^8 - \left(\frac{1}{4096} g_{2}^{6}\right. \\
   & & + \left.\frac{27}{64} g_{2}^{3} g_{3}^{2} + \frac{729}{256}
   g_{3}^{4}\right)Z^9 .\\
\end{eqnarray*}
\subsection{Une valeur admissible pour la constante de comparaison entre hauteur projective et hauteur de Néron-Tate sur $E$:}
Pour tout point $\mathbf p \in E(K),$ désignons par $h(\mathbf p)$
la hauteur logarithmique absolue de $\mathbf p$ (voir §$2$) et par
$\widehat{h} (\mathbf p)$ la hauteur normalisée (ou de Néron-Tate)
de $\mathbf p$ qui est définie par:
$$\widehat{h} (\mathbf p) := \lim_{n \rightarrow + \infty} \frac{h(n.\mathbf p)}{n^2} .$$
On a le théor\`eme suivant:
\begin{theorem}\label{c.14}
Pour tout point $\mathbf p \in E(K)$ on a:
$$- \frac{3}{4} \eta - 5 ~\leq~ h(\mathbf p) - \widehat{h} (\mathbf p) ~\leq~ \frac{3}{2} \eta + 8 .$$
\end{theorem}
{\bf Démonstration.---} C'est une conséquence immédiate du lemme
$3.4$ de [Da1]. En effet, on majore simplement la constante $h :=
\max\{1 , h(1 : g_2 : g_3)\}$ introduite dans cette référence par
$h(1 : g_2 : g_3) + 1$, puis on majore les constantes numériques
$5 \log 2 + \frac{3}{4}$ et $8 \log 2 + \frac{3}{2}$ par $5$ et
$8$ respectivement et le théor\`eme \ref{c.14} en découle.
$~~~~\blacksquare$\vspace{1mm}
\section{Appendice}
Dans ce paragraphe nous démontrons les quelques lemmes
élémentaires utilisés depuis le §$7$.
\begin{lemma}[G. Rémond]\label{b.21}
Soient $c_1$ un réel $>1$ et $r$ un entier $\geq 1$. L'espace euclidien ${\mathbb R}^r$ peut \^etre recouvert par un nombre
$\leq [(1 + \sqrt{8 c_1})^r]$ de c\^ones de sommet $0$ et d'angle $\arccos{(1 - \frac{1}{c_1})}$ chacun.
\end{lemma}
{\bf Démonstration.---}
Durant toute cette démonstration, quand $x$ et $y$ sont deux points de ${\mathbb R}^r$, on note $(x , y)$ le réel appartenant à l'intervalle $]-\pi , \pi]$
représentant l'angle $(\overrightarrow{0x} , \overrightarrow{0y})$ en radians.
Comme tout point de ${\mathbb R}^r$ s'obtient par une homothétie de centre $0$ d'un point de la sph\`ere unité $S := S(0 , 1)$ alors: recouvrir l'espace euclidien ${\mathbb R}^r$ par des c\^ones de sommet $0$ et d'angle $\arccos{(1 - \frac{1}{c_1})}$ revient à recouvrir $S$ par de tels (m\^emes) c\^ones. Pour faire cela en ayant une estimation du nombre de c\^ones du recouvrement, on se base sur le fait suivant: \og $\forall x \in S$ et $\forall \varphi \in ]0 , \frac{\pi}{4}[$ le sous-ensemble $B(x , 2 \sin{\varphi}) \cap S$ de ${\mathbb R}^r$ est contenu dans un certain c\^one de centre $0$ et d'angle $4 \varphi$ \fg. Montrons d'abord ce fait: soient $y_1$ et $y_2$ deux points quelconques de $B(x , 2 \sin{\varphi}) \cap S$, on a d'une part:
\begin{align*}
{\nb{y_i - x}}^2 &= {\nb{y_i}}^2 + {\nb{x}}^2 - 2 \nb{y_i} . \nb{x} . \cos{(x , y_i)} \\
&= 2 - 2 \cos{(x , y_i)} ~~~~~~~~~ i = 1 , 2
\end{align*}
et d'autre part:
\begin{align*}
\nb{y_i - x} \leq 2 \sin{\varphi} &\Rightarrow {\nb{y_i - x}}^2 \leq 4 {\sin}^2{\varphi} ~~~~~~ i = 1 , 2 \\
\intertext{d'o\`u on déduit:}
2 - 2 \cos{(x , y_i)} &\leq 4 {\sin}^2 \varphi ~~~~~~~~~~~~~~~~~~~~~ i = 1 , 2 \\
\mbox{i.e}~~~~~~~~~~~~~~~~~~~~~~~~~~~~~~~~~ \cos{(x , y_i)} &\geq 1 - 2 {\sin}^2 \varphi = \cos{2 \varphi} ~~~~ i = 1 , 2 ~~~~~~~~~~~~~~~~~~~~~\\
\mbox{d'o\`u}~~~~~~~~~~~~~~~~~~~~~~~~~~~~~~~~~~ \nb{(x , y_i)} &\leq 2 \varphi ~~~~~~~~~~~~~~~~~~~~~~~~~~ i = 1 , 2 . ~~~~~~~~~~~~~~~~~~~~~~~~~~~~~~~~~~~~~~~~~~~~~~~~~~~~~~~~~~~~~~~~
\end{align*}
Ainsi
$$\nb{(y_1 , y_2)} \leq \nb{(x , y_1)} + \nb{(x , y_2)} \leq 2 \varphi + 2 \varphi = 4 \varphi ,$$
or ceci est vrai pour tout couple de points $(y_1 , y_2) \in {\left(B(x , 2 \sin{\varphi}) \cap S \right)}^2$ ce qui montre qu'effectivement le sous-ensemble $B(x , 2 \sin{\varphi}) \cap S$ de ${\mathbb R}^r$ est contenu dans un certain c\^one de centre $0$ et d'angle $4 \varphi$.
\\ Posons maintenant $\varphi := \frac{1}{4} \arccos{(1 - \frac{1}{c_1})}$, d'apr\`es le fait ci-dessus le probl\`eme revient à recouvrir la sph\`ere unité $S$ par un nombre fini de boules ouvertes de rayon $2 \sin{\varphi}$ chacune et dont les centres appartiennent à $S$, en ayant une estimation du nombre de boules de ce recouvrement (qui est aussi le nombre de c\^ones du recouvrement déduit par le fait précédent). Mais, il suffit d'appliquer le lemme \ref{b.20} pour $S$ avec $\rho$ remplacé par $1$ et $\gamma$ par $\frac{1}{2 \sin{\varphi}}$; ce qui montre que $S$ peut \^etre recouverte par un nombre $\leq [(\frac{1}{\sin{\varphi}} + 1)^r]$ de petites boules de rayon $2 \sin{\varphi}$ chacune et dont les centres ont dans $S$. Par conséquent, en remontant le raisonnement, ${\mathbb R}^r$ est recouvert par un nombre $\leq [(\frac{1}{\sin{\varphi}} + 1)^r]$ de c\^ones de sommet $0$ et d'angle $4 \varphi = \arccos{(1 - \frac{1}{c_1})}$ chacun. Il nous reste à vérifier que $\frac{1}{\sin{\varphi}} + 1 \leq 1 + \sqrt{8 c_1}$. On a par définition de $\varphi$:
$$1 - \frac{1}{c_1} = \cos{4 \varphi} = 2 {\cos}^2 2 \varphi - 1 = 2 {\left(1 - 2 {\sin}^2 \varphi \right)}^2 - 1$$
d'o\`u on tire:
\begin{equation*}
\begin{split}
\sin{\varphi} &= \sqrt{\frac{1}{2} \left(1 - \sqrt{1 - \frac{1}{2 c_1}}\right)} \\
&= \frac{1}{2 \sqrt{c_1} \sqrt{1 + \sqrt{1 - \frac{1}{2 c_1}}}} \\
&\geq \frac{1}{\sqrt{8 c_1}}
\end{split}
\end{equation*}
d'o\`u $\displaystyle ~~\frac{1}{\sin{\varphi}} + 1 \leq 1 + \sqrt{8 c_1}~~$ ce qui ach\`eve la preuve. $~~~~\blacksquare$\vspace{1mm}
\begin{lemma}\label{c.29}
Soient $X_1 , X_2 , Y_1$ et $Y_2$ des réels strictement positifs,
$(a , b)$ un couple de réels positifs différent de $(0 , 0)$ et
$k$ un réel $> 1$. Supposons qu'on ait:
$$- a ~\leq~ Y_i - X_i ~\leq~ b ~~~~\text{pour $i = 1 , 2$.}$$
Alors, la double inégalité:
$$\frac{1}{k} \frac{Y_1}{Y_2} ~\leq~ \frac{X_1}{X_2} ~\leq~ k \frac{Y_1}{Y_2}$$
est satisfaite d\`es que:
\begin{align}
X_i &\geq~ \frac{a k + b}{k - 1} ~~~~\text{pour $i = 1 , 2$} \notag \\
\intertext{ou d\`es que:} Y_i &\geq~ \frac{b k + a}{k - 1}
~~~~\text{pour $i = 1 , 2$.} \notag
\end{align}
\end{lemma}
{\bf Démonstration.---} $\bullet$ Supposons qu'on a: $X_i \geq
\frac{a k + b}{k - 1} > a$ pour $i = 1 , 2$. Des inégalités:
$$- a \leq Y_i - X_i \leq b ~~~~(i = 1 , 2)$$
$\displaystyle \mbox{vient:}~~~~~~~~~~~~~ \frac{X_1 - a}{X_2 + b} \leq \frac{Y_1}{Y_2} \leq \frac{X_1 + b}{X_2 - a}$. \\
Comme on a:
\begin{align}
\frac{X_1 + b}{X_2 - a} \leq k \frac{X_1}{X_2} & \Longleftrightarrow \left(X_1 - \frac{b}{k - 1}\right) \left(X_2 - \frac{k a}{k - 1}\right) \geq \frac{a b k}{(k - 1)^2} , \notag \\
\frac{X_1 - a}{X_2 + b} \geq \frac{1}{k} \frac{X_1}{X_2} &
\Longleftrightarrow \left(X_1 - \frac{a k}{k - 1}\right) \left(X_2
- \frac{b}{k - 1}\right) \geq \frac{a b k}{(k - 1)^2} \notag
\end{align}
et que les deux inégalités de droite sont satisfaites, puisqu'on a
supposé que $X_i \geq \frac{a k + b}{k - 1}$ pour $i = 1 , 2$, on
a bien:
$$\frac{X_1 + b}{X_2 - a} \leq k \frac{X_1}{X_2} ~~\mbox{et}~~ \frac{X_1 - a}{X_2 + b} \geq \frac{1}{k} \frac{X_1}{X_2}$$
\begin{align}
&\text{d'o\`u:} ~~~~~~~~~~ &\frac{1}{k} \frac{X_1}{X_2} \leq \frac{X_1 - a}{X_2 + b} \leq \frac{Y_1}{Y_2} \leq \frac{X_1 + b}{X_2 - a} \leq k \frac{X_1}{X_2}~~~~~~~~~~~~~~ \notag \\
&\text{puis:} ~~~~~~~~~~&\frac{1}{k} \frac{Y_1}{Y_2} \leq
\frac{X_1}{X_2} \leq k \frac{Y_1}{Y_2}
.~~~~~~~~~~~~~~~~~~~~~~~~~~~~~~~~~~~ \notag
\end{align}
ce qui démontre le premier cas du lemme \ref{c.29}. \\
$\bullet$ Le deuxi\`eme cas du lemme \ref{c.29} suit simplement du
premier cas en permutant $X_1$ avec $Y_1$, $X_2$ avec $Y_2$ et $a$
avec $b$. La démonstration est achevée.
$~~~~\blacksquare$\vspace{1mm}
\begin{lemma}\label{c.30}
Soit $E$ une courbe elliptique plongée dans ${\mathbb P}_2$ à la Weierstrass, définie sur un corps de nombres $K$, d'équation projective: \\
$Y^2 Z = 4 X^3 - g_2 X Z^2 - g_3 Z^3 ~ (g_2 , g_3 \in K)$ et soit
$m \geq 2$ un entier positif, ${\mathbf x}_1 , \dots , {\mathbf
x}_m$ des points de $E$ satisfaisant les hypoth\`eses
${\mbox{H}}_1$ du §$7$. Posons pour $i = 1 , \dots , m$:
$$a_i := \left[\frac{\mid {\mathbf x}_m \mid}{\mid {\mathbf x}_i \mid}\right] .$$
Alors, pour tout $i , j$ dans $\{1 , \dots , m\}$, on a:
\begin{eqnarray}
\frac{1}{\sqrt{2}} a_{j}^{2} \widehat{h}({\mathbf x}_j) & \leq a_{i}^{2} \widehat{h}({\mathbf x}_i) & \leq \sqrt{2} a_{j}^{2} \widehat{h}({\mathbf x}_j) , \label{3.25} \\
\frac{1}{\sqrt{2}} \frac{h({\mathbf x}_i)}{h({\mathbf x}_j)} & \leq \frac{\widehat{h}({\mathbf x}_i)}{\widehat{h}({\mathbf x}_j)} & \leq \sqrt{2} \frac{h({\mathbf x}_i)}{h({\mathbf x}_j)} ,\label{3.26} \\
\frac{1}{2} a_{j}^{2} h({\mathbf x}_j) & \leq a_{i}^{2} h({\mathbf
x}_i) & \leq 2 a_{j}^{2} h({\mathbf x}_j) . \label{3.27}
\end{eqnarray}
\end{lemma}
{\bf Démonstration.---}
Remarquons d'abord que les trois doubles inégalités (\ref{3.25}), (\ref{3.26}) et (\ref{3.27}) sont symétriques par rapport à $i$ et $j$ et qu'elles sont triviales lorsque $i = j$. Donc, on peut supposer -sans restreindre la généralité- qu'on a $1 \leq i < j \leq m$. En tenant compte de cette hypoth\`ese, commençons par démontrer la double inégalité (\ref{3.25}). Pour cela on distingue les deux cas suivants: \\
$\bullet \underline{{\mbox{1}}^{\mbox{er}} \mbox{cas}}$: (si $j =
m$) Dans ce premier cas, on applique le lemme \ref{c.29} avec:
$X_1 = a_i$, $X_2 = Y_1 = Y_2 = \frac{\mid {\mathbf x}_m
\mid}{\mid {\mathbf x}_i \mid}$, $a = 0$, $b = 1$ et $k =
\sqrt[4]{2}$. Les hypoth\`eses du lemme \ref{c.29} sont clairement
vérifiées, par conséquent -d'apr\`es le deuxi\`eme cas de ce
dernier- on a:
\begin{equation}
\frac{1}{\sqrt[4]{2}} \leq a_i \frac{\mid {\mathbf x}_i \mid}{\mid
{\mathbf x}_m \mid} \leq \sqrt[4]{2} \label{3.28}
\end{equation}
d\`es que:
$$\frac{\mid {\mathbf x}_m \mid}{\mid {\mathbf x}_i \mid} \geq \frac{\sqrt[4]{2}}{\sqrt[4]{2} - 1} = 6,29\dots .$$
Or, cette derni\`ere inégalité est bien vérifiée puisqu'on a,
d'apr\`es les hypoth\`eses ${\mbox{H}}_1$:
$$\frac{\mid {\mathbf x}_m \mid}{\mid {\mathbf x}_i \mid} \geq \frac{\mid {\mathbf x}_m \mid}{\mid {\mathbf x}_{m - 1} \mid} \geq 7 ~~~~ \text{(car $i \leq j - 1 \leq m - 1$)} .$$
Donc, la double inégalité (\ref{3.28}) est établie (pour ce premier cas). Il suffit d'en prendre les carrés des trois termes et de remarquer que $a_j = a_m = 1$ pour en déduire la double inégalité (\ref{3.25}) dans ce premier cas. \\
$\bullet \underline{{\mbox{2}}^{\mbox{\`eme}} \mbox{cas}}$: (si $j
\leq m - 1$) Dans ce deuxi\`eme cas, on applique le lemme
\ref{c.29} avec: $X_1 = a_i$, $X_2 = a_j$, $Y_1 = \frac{\mid
{\mathbf x}_m \mid}{\mid {\mathbf x}_i \mid}$, $Y_2 = \frac{\mid
{\mathbf x}_m \mid}{\mid {\mathbf x}_j \mid}$, $a = 0$, $b = 1$ et
$k = \sqrt[4]{2}$. Les hypoth\`eses du lemme \ref{c.29} se
vérifient tr\`es facilement, par conséquent -d'apr\`es le
deuxi\`eme cas de ce dernier- on a:
\begin{equation}
\frac{1}{\sqrt[4]{2}} \frac{\mid {\mathbf x}_j \mid}{\mid {\mathbf
x}_i \mid} \leq \frac{a_i}{a_j} \leq \sqrt[4]{2} \frac{\mid
{\mathbf x}_j \mid}{\mid {\mathbf x}_i \mid} \label{3.29}
\end{equation}
d\`es que:
$$\min\!\left(\frac{\mid {\mathbf x}_m \mid}{\mid {\mathbf x}_i \mid} , \frac{\mid {\mathbf x}_m \mid}{\mid {\mathbf x}_j \mid}\right) \geq \frac{\sqrt[4]{2}}{\sqrt[4]{2} - 1} ,$$
c'est-à-dire d\`es que:
$$\frac{\mid {\mathbf x}_m \mid}{\mid {\mathbf x}_j \mid} \geq \frac{\sqrt[4]{2}}{\sqrt[4]{2} - 1} = 6,29\dots .$$
Or, cette derni\`ere est bien vérifiée puisqu'on a, d'apr\`es les
hypoth\`eses ${\mbox{H}}_1$:
$$\frac{\mid {\mathbf x}_m \mid}{\mid {\mathbf x}_j \mid} \geq \frac{\mid {\mathbf x}_m \mid}{\mid {\mathbf x}_{m - 1} \mid} \geq 7 ~~~~\text{(car $j \leq m - 1$)} .$$
Donc, la double inégalité (\ref{3.29}) est établie (pour ce deuxi\`eme cas). Il ne reste qu'à prendre les carrés des trois termes de cette derni\`ere pour aboutir à la double inégalité (\ref{3.25}) dans ce deuxi\`eme cas. En conclusion la double inégalité (\ref{3.25}) est vraie pour tout $i , j$ dans $\{1 , \dots , m\}$. \\
Montrons maintenant la double inégalité (\ref{3.26}). Pour cela on
applique le lemme \ref{c.29} avec: $X_1 = \widehat{h}({\mathbf
x}_i)$, $X_2 = \widehat{h}({\mathbf x}_j)$, $Y_1 = h({\mathbf
x}_i)$, $Y_2 = h({\mathbf x}_j)$, $a = \frac{3}{4} \eta + 5$, $b =
\frac{3}{2} \eta + 8$ et $k = \sqrt{2}$. Le théor\`eme \ref{c.14}
montre ainsi que les hypoth\`eses du lemme \ref{c.29} sont bien
vérifiées, donc d'apr\`es le premier cas de ce dernier, la double
inégalité (\ref{3.26}) est satisfaite d\`es que:
$$\min\left(\widehat{h}({\mathbf x}_i) , \widehat{h}({\mathbf x}_j)\right) \geq \frac{\left(\frac{3}{4} \eta + 5\right) \sqrt{2} + \frac{3}{2} \eta + 8}{\sqrt{2} - 1} ,$$
c'est-à-dire d\`es que:
$$\widehat{h}({\mathbf x}_i) \geq \left(3 + \frac{9}{4} \sqrt{2}\right) \eta + 18 + 13 \sqrt{2} = 6,18\dots \eta + 36,38\dots .$$
Or, cette derni\`ere est bien satisfaite puisque, d'apr\`es les hypoth\`eses ${\mbox{H}}_1$, on a m\^eme: $\widehat{h}({\mathbf x}_i) \geq \widehat{h}({\mathbf x}_1) \geq 7 \eta + 37$. Ce qui démontre la double inégalité (\ref{3.26}) pour tout $i , j$ dans $\{1 , \dots , m\}$. \\
Finalement, nous remarquons que la double inégalité (\ref{3.26})
(qu'on vient de démontrer) est équivalente à la double inégalité:
$$\frac{1}{\sqrt{2}} \frac{h({\mathbf x}_j)}{\widehat{h}({\mathbf x}_j)} \leq \frac{h({\mathbf x}_i)}{\widehat{h}({\mathbf x}_i)} \leq \sqrt{2} \frac{h({\mathbf x}_j)}{\widehat{h}({\mathbf x}_j)} .$$
Le produit membre à membre de cette derni\`ere avec (\ref{3.25})
(déja démontrée ci-dessus) nous am\`ene à l'inégalité (\ref{3.27})
et ach\`eve cette démonstration. $~~~~\blacksquare$\vspace{1mm}
\begin{lemma}\label{c.31}
Sous les hypoth\`ese du lemme \ref{c.30} on a:
\begin{description}
\item[1)] $\forall i \in \{1 , \dots , m - 1\} : ~~~~ a_i \geq
\frac{1}{\sqrt{\alpha}}$, \item[2)] $\forall i \in \{1 , \dots ,
m\} : ~~~~ a_i \geq 7^{m - i} \geq 7 (m - i)$, \item[3)] $\forall
i \in \{2 , \dots , m\} : ~~~~ a_{i - 1} \geq 7 a_i$, \item[4)]
$a_{1}^{2} + \dots + a_{m}^{2} \leq (1 + 1/48) a_{1}^{2}$ et
\item[5)] $m - 1 \leq (1 + 1/48) \alpha a_{1}^{2}$.
\end{description}
\end{lemma}
{\bf Démonstration.---} $\bullet$ Démontrons $1)$:
soit $i \in \{1 , \dots , m - 1\}$. On a: $a_i := [\frac{\mid x_m \mid}{\mid x_i \mid}] > \frac{\mid x_m \mid}{\mid x_i \mid} - 1$. Comme d'apr\`es les hypoth\`eses ${\mbox{H}}_1$ on a: $\frac{\mid x_m \mid}{\mid x_i \mid} \geq \frac{1}{\sqrt{\alpha}} + 1$, c'est-à-dire: $\frac{\mid x_m \mid}{\mid x_i \mid} - 1 \geq \frac{1}{\sqrt{\alpha}}$ alors $a_i > \frac{1}{\sqrt{\alpha}}$. Ce qui démontre $1)$ du lemme \ref{c.31}. \\
$~~~~~\!\bullet$ Démontrons $2)$:
soit $i \in \{1 , \dots , m\}$. D'apr\`es l'hypoth\`ese ${\mbox{H}}_1$, affirmant: \og $\forall j \in \{2 , \dots , m\}$, $\mid x_j \mid \geq 7 \mid x_{j - 1} \mid$ \fg, on a $\mid x_m \mid \geq 7^{m - i} \mid x_i \mid$, c'est-à-dire: $\frac{\mid x_m \mid}{\mid x_i \mid} \geq 7^{m - i}$. En prenant les parties enti\`eres pour cette derni\`ere inégalité on a finalement $a_i \geq 7^{m - i}$. Par ailleurs la minoration: $7^{m - i} \geq 7 (m - i)$ est claire. D'o\`u s'ensuit $2)$ du lemme \ref{c.31}. \\
$~~~~\!\bullet$ Démontrons $3)$: Soit $i \in \{2 , \dots , m\}$. D'apr\`es l'hypoth\`ese ${\mbox{H}}_1$ on a: $\mid x_i \mid \geq 7 \mid x_{i - 1} \mid$; d'o\`u $\frac{\mid x_m \mid}{\mid x_{i - 1} \mid} \geq 7 \frac{\mid x_m \mid}{\mid x_i \mid} \geq 7 [\frac{\mid x_m \mid}{\mid x_i \mid}] = 7 a_i$. En prenant les parties enti\`eres des deux membres de l'inégalité $\frac{\mid x_m \mid}{\mid x_{i - 1} \mid} \geq 7 a_i$ on a finalement: $a_{i - 1} \geq 7 a_i$ ce qui démontre $3)$ du lemme \ref{c.31}. \\
$~~~~~\!\bullet$ Démontrons $4)$: d'apr\`es $3)$ -qu'on vient de
démontrer- on a: $$a_1 \geq 7 a_2 \geq 7^2 a_3 \geq \dots \geq
7^{i - 1} a_i \geq \dots \geq 7^{m - 1} a_m$$ donc:
$$a_{1}^{2} + \dots + a_{m}^{2} \leq \left(1 + \frac{1}{49} + \frac{1}{49^2} + \dots + \frac{1}{49^{m - 1}}\right) a_{1}^{2} \leq \left(1 + \frac{1}{48}\right) a_{1}^{2}$$
ce qui démontre $4)$ du lemme \ref{c.31}. \\
$~~~~\!\bullet$ Démontrons $5)$: d'apr\`es $1)$ -déja démontré- on
a pour tout $i$ dans $\{1 , \dots , m - 1\}$ l'inégalité $1 \leq
\alpha {a^2}_{\!\!\!i}$ d'o\`u:
$$\sum_{i = 1}^{m - 1} 1 \leq \sum_{i = 1}^{m - 1} \alpha {a^2}_{\!\!\!i}$$
$\displaystyle \text{c'est-à-dire:} ~~~~~~~~~~ m - 1 \leq \alpha \left({a^2}_{\!\!\!1} + \dots + {a^2}_{\!\!\!m - 1}\right)$. \\
comme $a_{1}^{2} + \dots + a_{m - 1}^{2} \leq a_{1}^{2} + \dots +
a_{m}^{2} \leq (1 + 1/48) a_{1}^{2}$ d'apr\`es $4)$ -qu'on vient
de démontrer- on en déduit finalement qu'on a $m - 1 \leq (1 +
1/48) \alpha a_{1}^{2}$ ce qui démontre $5)$ du lemme \ref{c.31}
et ach\`eve cette démonstration. $~~~~\blacksquare$\vspace{1mm}
\begin{lemma}\label{c.32}
Sous les hypoth\`eses du lemme \ref{c.30} on a, pour tout $i \in
\{1 , \dots , m - 1\}$:
$$\widehat{h}(a_i {\mathbf x}_i - {\mathbf x}_m) \leq \alpha \left(a_{1}^{2} \widehat{h}({\mathbf x}_i) + \widehat{h}({\mathbf x}_m)\right) .$$
\end{lemma}
{\bf Démonstration.---} Pour tout $i \in \{1 , \dots , m - 1\}$ on
écrit:
\begin{eqnarray*}
{\mid \!a_i {\mathbf x}_i - {\mathbf x}_m \!\mid}^2 & = & a_{i}^{2} {\mid \!{\mathbf x}_i \!\mid}^2 + {\mid \!{\mathbf x}_m \!\mid}^2 - 2 a_i <{\mathbf x}_i , {\mathbf x}_m> \\
& = & \left(a_i \mid \!{\mathbf x}_i \!\mid - \mid \!{\mathbf x}_m \!\mid\right)^2 + 2 a_i \left(\mid \!{\mathbf x}_i \!\mid . \mid \!{\mathbf x}_m \!\mid - <{\mathbf x}_i , {\mathbf x}_m>\right) \\
& \leq & {\mid \!{\mathbf x}_i \!\mid}^2 + 2 a_i \frac{\alpha}{4} \mid {\mathbf x}_i \!\mid . \mid \!{\mathbf x}_m \!\mid \\
& \leq & {\mid \!{\mathbf x}_i \!\mid}^2 + \frac{\alpha}{2} a_i \mid \!{\mathbf x}_i \!\mid . \mid \!{\mathbf x}_m \!\mid \\
& \leq & \frac{{\mid \!{\mathbf x}_m \!\mid}^2}{a_{i}^{2}} + \alpha \frac{a_i (a_i + 1)}{2} {\mid \!{\mathbf x}_i \!\mid}^2 \\
& \leq & \alpha \left(a_{i}^{2} {\mid \!{\mathbf x}_i \!\mid}^2 +
{\mid \!{\mathbf x}_m \!\mid}^2\right)
\end{eqnarray*}
o\`u, pour l'obtention de cette série d'inégalités, on a utilisé:
la premi\`ere inégalité de l'hypoth\`ese ${\mbox{H}}_1$ qui
entra{\sf\^\i}ne $\mid \!{\mathbf x}_i \!\mid \!. \!\mid
\!{\mathbf x}_m \!\mid - <{\mathbf x}_i , {\mathbf x}_m> \leq
\frac{\alpha}{4} \!\mid \!{\mathbf x}_i \!\mid \!. \!\mid
\!{\mathbf x}_m \!\mid$; $a_i := [\frac{\mid {\mathbf x}_m
\mid}{\mid {\mathbf x}_i \mid}]$ qui entra{\sf\^\i}ne $\frac{\mid
{\mathbf x}_m \mid}{\mid {\mathbf x}_i \mid} - 1 < a_i \leq
\frac{\mid {\mathbf x}_m \mid}{\mid {\mathbf x}_i \mid}$;
l'assertion $1)$ du lemme \ref{c.31} qui entra{\sf\^\i}ne:
$\frac{1}{a_{i}^{2}} < \alpha$ et enfin la majoration triviale
$\frac{a_i (a_i + 1)}{2} \leq a_{i}^{2}$. Ceci ach\`eve cette
démonstration.$~~~~\blacksquare$\vspace{2mm}

Avant de passer aux autres lemmes, il est intéressant de remarquer que le lemme \ref{c.30} est vrai m\^eme si on se
restreint seulement aux troisi\`eme et quatri\`eme inégalités de l'hypoth\`ese ${\mbox{H}}_1$.
Le lemme \ref{c.31} est vrai m\^eme si on se restreint seulement aux deuxi\`eme et troisi\`eme inégalités
de l'hypoth\`ese ${\mbox{H}}_1$ et le lemme \ref{c.32} reste valable m\^eme si on se restreint seulement
aux premi\`ere et deuxi\`eme inégalités de l'hypoth\`ese ${\mbox{H}}_1$.

Le lemme qui suit est utilisé dans le §$6$ pour affirmer:
$\card{~\!\!(T_{\delta})} = (1 + o(\delta))
\mbox{vol}~\!(T_{\delta})$.
\begin{lemma}\label{c.33}
Soient $\delta$ un entier positif, $r_1 , \dots , r_n ~ (n \in
{\mathbb N}^*)$ des entiers strictement positifs et $T$ un entier
positif satisfaisant:
$$T \leq \frac{1}{r_1} + \dots + \frac{1}{r_n} .$$
Alors, l'ensemble:
$$C := \left\{(\alpha_1 , \dots , \alpha_n) \in {\mathbb N}^n / \frac{\alpha_1}{r_1} + \dots + \frac{\alpha_n}{r_n} \leq \delta \right\}$$
est de cardinal encadré par:
$$r_1 \dots r_n \binom{\delta + T}{n} ~\leq~ \card{~\!\!(C)} ~\leq~ r_1 \dots r_n \binom{\delta + n}{n} .$$
\end{lemma}
{\bf Démonstration.---} L'idée consiste à encadrer $C$ (au sens de
l'inclusion) par deux ensembles de cardinaux faciles à calculer.
Les deux ensembles dont il s'agit ici sont du type de la famille
d'ensembles $C_R ~ (R \in \mathbb N)$ définis par:
$$C_R := \left\{(\alpha_1 , \dots , \alpha_n) \in {\mathbb N}^n / \left[\frac{\alpha_1}{r_1}\right] + \dots + \left[\frac{\alpha_n}{r_n}\right] \leq R\right\} ,$$
o\`u $[.]$ désigne la partie enti\`ere. On a ainsi la double
inclusion:
$$
C_{\delta - n + T} ~\subset~ C ~\subset~ C_{\delta} .
$$
En effet, l'inclusion $C \subset C_{\delta}$ est immédiate et
l'inclusion $C_{\delta - n + T} \subset C$ s'obtient en remarquant
que pour tout rationnel positif $p/q ~ (p \in \mathbb N , q \in
{\mathbb N}^*)$ on a: $p/q \leq [p/q] + 1 - 1/q$. La double
inclusion précédente entra{\sf\^\i}ne pour $\card{~\!\!(C)}$
l'encadrement:
\begin{equation}
\card{~\!\!(C_{\delta - n + T})} ~\leq~ \card{~\!\!(C)} ~\leq~
\card{~\!\!(C_{\delta})} . \label{3.87}
\end{equation}
Par ailleurs, nous affirmons que pour tout entier positif $R$, le
cardinal de l'ensemble $C_R$ vaut exactement:
$$\card{~\!\!(C_R)} = r_1 \dots r_n \binom{R + n}{n} .$$
En effet, un élément $(\alpha_1 , \dots , \alpha_n)$ de $C_R$ est
forcément une solution d'un syst\`eme:
\begin{gather}
\left\{\!\!\!
\begin{array}{cl}
\left[\frac{\alpha_1}{r_1}\right]  & \!\!\!=~ x_1 \\
\vdots \\
\left[\frac{\alpha_n}{r_n}\right] & \!\!\!=~ x_n
\end{array}
\right. \tag{$\Xi$}
\end{gather}
pour un certain $(x_1 , \dots , x_n) \in {\mathbb N}^n$
satisfaisant: $x_1 + \dots + x_n \leq R$. Comme l'ensemble des
$(x_1 , \dots , x_n)$ de ${\mathbb N}^n$ satisfaisant $x_1 + \dots
+ x_n \leq R$ est de cardinal exactement $\binom{R + n}{n}$ et que
pour un tel point $(x_1 , \dots , x_n)$, l'ensemble des solutions
du syst\`eme $(\Xi)$ est exactement $r_1 \dots r_n$ (car ces
solutions sont: $(r_1 x_1 + i_1 , \dots , r_n x_n + i_n) /~ 0 \leq
i_1 \leq r_1 - 1 , \dots , 0 \leq i_n \leq r_n - 1$, et sont donc
au nombre de $r_1 \dots r_n$) alors le cardinal de $C_R$ est
effectivement $r_1 \dots r_n \binom{R + n}{n}$. En particulier:
$$\card{~\!\!(C_{\delta - n + T})} = r_1 \dots r_n \binom{\delta +
T}{n} ~~\text{et}~~ \card{~\!\!(C_{\delta})} = r_1 \dots r_n
\binom{\delta + n}{n} .$$ Ce qui conclut, d'apr\`es (\ref{3.87}),
cette démonstration.  $~~~~\blacksquare$\vspace{1mm}
\begin{remarque}\label{c.34}
Lorsque $r_1 = \dots = r_n = 1$, on peut prendre $T = n$ et dans
ce cas le lemme \ref{c.33} est précis et on retrouve le fait que
$\card~\!\!\{(\alpha_1 , \dots , \alpha_n) \in {\mathbb N}^n /~
\alpha_1 + \dots + \alpha_n \leq \delta\} = \binom{\delta + n}{n}
~~ \forall \delta \in \mathbb N$.
\end{remarque}
\begin{lemma}\label{c.35}
Soient $E$ la courbe elliptique du §$2$, $\mathbf x$ un point de
$E(K)$ représenté dans ${\mathbb P}_2$ par le syst\`eme de
coordonnées projectives $\underline x = (x , 1 , z) \in K^3$, $v$
une place de $K$ et $\epsilon$ un réel strictement positif.
Supposons qu'on a:
\begin{gather}
{\mbox{dist}}_v(\mathbf x , \mathbf 0) ~<~ e^{-2 m_v - c_v} ,
\tag{$*$}
\end{gather}
alors on a: \\
$\bullet$ si $v$ est finie:
$$\max\left({\mid x \mid}_v , {\mid z \mid}_v \right) ~<~ \min\left\{1 , \frac{1}{{\mid g_2 \mid}_v} , \frac{1}{{\mid g_3 \mid}_v}\right\} ~~\mbox{et}~~ {\mid \Delta(\underline x) \mid}_v = 1 $$
$\bullet$ et si $v$ est infinie:
$$\max\left({\mid x \mid}_v , {\mid z \mid}_v \right) ~<~ e^{- 16} \min\left\{1 , \frac{1}{{\mid g_2 \mid}_v} , \frac{1}{{\mid g_3 \mid}_v}\right\} ~,~ {\mid \Delta(\underline x) \mid}_v > \frac{1}{\sqrt e}$$
et pour tout entier positif $d$:
$$\sum_{\begin{array}{c}
\scriptstyle \underline{\alpha} \in {\mathbb N}^3 \\
\scriptstyle \mid \underline{\alpha} \mid = d
\end{array}}{\mid {{\underline x}^{\underline{\alpha}}} \mid}_v ~\leq~ e .$$
\end{lemma}
{\bf Démonstration.---} D'apr\`es la propriété $ii)$ du §$2$ pour
la distance ${\mbox{dist}}_v$, l'hypoth\`ese $(*)$
entra{\sf\^\i}ne qu'on a:
$${\mbox{dist}}_v(\mathbf x , \mathbf 0) ~=~ \max\left({\mid x \mid}_v , {\mid z \mid}_v\right)$$ et par conséquent on a:
$$\max\left({\mid x \mid}_v , {\mid z \mid}_v\right) ~<~ e^{-2 m_v - c_v} .$$
Cette inégalité entra{\sf\^\i}ne bien les estimations du lemme
\ref{c.35} concernant la quantité $\max({\mid x \mid}_v , {\mid z
\mid}_v)$ puisqu'on a: $e^{-2 m_v} = {M_v}^{-2} \leq {M_v}^{-1} =
\min\{1 , 1/{\mid g_2 \mid}_v , 1/{\mid g_3 \mid}_v\}$. Démontrons
maintenant -en distinguant les deux cas ``$v$ finie'' et ``$v$
infinie''- l'assertion du lemme \ref{c.35} pour ${\mid
\Delta(\underline x) \mid}_v$. On a $\Delta(\underline x) = 3 g_3
z^2 + 2 g_2 x z + 1$.
\\

$\underline{{\mbox{1}}^{\mbox{er}} \mbox{cas}}$ (si $v$ est finie): \\
Dans ce cas on a:
$$ {\mid \Delta(\underline x) \mid}_v ~\leq~ \max\left(1 , {\mid 3 g_3 z^2 + 2 g_2 x z \mid}_v\right) ,$$
o\`u cette derni\`ere inégalité devient une égalité lorsque ${\mid
3 g_3 z^2 + 2 g_2 x z \mid}_v < 1$. Or, c'est bien le cas dans
notre situation. En effet, on a:
\begin{eqnarray*}
{\mid 3 g_3 z^2 + 2 g_2 x z \mid}_v & = & {\mid z \mid}_v.{\mid 3 g_3 z + 2 g_2 x \mid}_v \\
& \leq & {\mid 3 g_3 z + 2 g_2 x \mid}_v ~~~~~~~~~~\text{(car
${\mid z \mid}_v < 1$)}
\end{eqnarray*}
puis:
\begin{eqnarray*}
{\mid 3 g_3 z^2 + 2 g_2 x z \mid}_v & \leq & \max\left({\mid 3 g_3 z \mid}_v , {\mid 2 g_2 x \mid}_v\right) \\
& \leq & \max\left({\mid g_3 \mid}_v.{\mid z \mid}_v , {\mid g_2 \mid}_v.{\mid x \mid}_v\right) \\
& < & 1 ~~~~~~~~~~\left(\text{car}~ {\mid z \mid}_v <
\frac{1}{{\mid g_3 \mid}_v} ~\text{et}~ {\mid x \mid}_v <
\frac{1}{{\mid g_2 \mid}_v}\right) .
\end{eqnarray*}
D'o\`u: $\displaystyle {\mid \Delta(\underline x) \mid}_v ~=~ 1 .$
\\

$\underline{{\mbox{2}}^{\mbox{i\`eme}} \mbox{cas}}$ (si $v$ est infinie): \\
Pour ce deuxi\`eme cas on a:
$${\mid \Delta(\underline x) \mid}_v ~\geq~ 1 - {\mid 3 g_3 z^2 + 2 g_2 x z \mid}_v .$$
Or:
\begin{eqnarray*}
{\mid 3 g_3 z^2 + 2 g_2 x z \mid}_v & = & {\mid z \mid}_v.{\mid 3 g_3 z + 2 g_2 x \mid}_v \\
& \leq & {\mid z \mid}_v.\left(3 {\mid g_3 \mid}_v.{\mid z \mid}_v + 2 {\mid g_2 \mid}_v.{\mid x \mid}_v\right) \\
& \leq & e^{- 16} \left(3 e^{- 16} + 2 e^{- 16}\right) \\
& \leq & 5 e^{- 32} \\
& < & 1 - \frac{1}{\sqrt e} .
\end{eqnarray*}
D'o\`u alors: $\displaystyle {\mid \Delta(\underline x) \mid}_v >
\frac{1}{\sqrt e}$. \\ Finalement, pour achever la démonstration
du lemme \ref{c.35}, il ne reste qu'\`a vérifier que dans ce
deuxi\`eme cas ($v$ infinie), pour tout entier positif $d$ on a
$\sum_{\mid \underline{\alpha} \mid = d} {\mid {\underline
x}^{\underline{\alpha}} \mid}_v \leq e$. En effet, on a pour tout
$d \in \mathbb N$:
\begin{eqnarray*}
\sum_{\!\!\!\!\!\! \underline{\alpha} \in {\mathbb N}^3 , \mid
\underline{\alpha} \mid = d}
{\mid {\underline x}^{\underline{\alpha}} \mid}_v & = & \sum_{\!\!\!\!\!\! \alpha_0 + \alpha_1 + \alpha_2 = d} {{\mid x \mid}^{\alpha_0}}_{\!\!\!\!\!\!v} ~ {{\mid z \mid}^{\alpha_2}}_{\!\!\!\!\!\!v} ~~~~~~~~~~\text{(avec $(\alpha_0 , \alpha_1 , \alpha_2) := \underline{\alpha}$)} \\
& = & \sum_{\!\!\!\!\!\! \alpha_0 + \alpha_2 \leq d} {{\mid x \mid}^{\alpha_0}}_{\!\!\!\!\!\!v} ~ {{\mid z \mid}^{\alpha_2}}_{\!\!\!\!\!\!v} \\
& \leq & \sum_{\alpha_0 = 0}^{\infty} \sum_{\alpha_2 = 0}^{\infty} {{\mid x \mid}^{\alpha_0}}_{\!\!\!\!\!\!v} ~ {{\mid z \mid}^{\alpha_2}}_{\!\!\!\!\!\!v} \\
& \leq & \left(\sum_{\alpha_0 = 0}^{\infty} {{\mid x \mid}^{\alpha_0}}_{\!\!\!\!\!\!v}\right) \!\!\! \left(\sum_{\alpha_2 = 0}^{\infty} {{\mid z \mid}^{\alpha_2}}_{\!\!\!\!\!\!v}\right) \\
& \leq & \frac{1}{1 - {\mid x \mid}_v}.\frac{1}{1 - {\mid z \mid}_v} \\
& \leq & \left(\frac{1}{1 - e^{- 16}}\right)^2 ~~~~~~\left(\text{car ${\mid x \mid}_v < e^{- 16}$ et ${\mid z \mid}_v < e^{- 16}$}\right) \\
& < & e .
\end{eqnarray*}
La démonstration du lemme \ref{c.35} est ainsi compl\`ete.
$~~~~\blacksquare$\vspace{1mm}
\begin{lemma}\label{c.54}
Soit $a$ un réel positif $> a_0 := 15788$ et $f$ une fonction
réelle sur $[2 , + \infty[$ définie par:
$$f(x) := (2904 x)^x a^{\frac{x}{x - 1}} .$$
Alors $f$ atteint son minimum en un point unique $\xi = \xi(a)$
vérifiant l'encadrement:
$$\sqrt{\frac{2 \log a}{\log\!\left(\frac{\log a}{\log\log a}\right)+ 21}} + 1 ~<~ \xi ~<~
\sqrt{\frac{2 \log a}{\log\!\left(\frac{\log a}{\log\log
a}\right)+ 16}} + 1 .$$ De plus en désignant par $m_0$ l'entier
positif:
$$m_0 ~:=~ \left[\sqrt{\frac{2 \log a}{\log\log a - \log\log\log a + 16}} + 2\right] ,$$
o\`u $[.]$ désigne la partie enti\`ere ($m_0$ minimise donc \og
presque \fg la fonction $f$), on a:
$$f(m_0) \leq a^{1 + \frac{183}{\log\log a}} .$$
\end{lemma}
{\bf Démonstration.---} Minimiser $f$ revient à minimiser son
logarithme népérien:
$$g(x) := \log f(x) = x(\log x + \log 2904) + \frac{x}{x - 1} \log a .$$
Etudions les variations de $g$: On a $\forall x \in [2 , +
\infty[$: $~\! g'(x) = \log x + 1 + \log 2904 - \frac{\log a}{(x -
1)^2}$. Il est clair que $g'$ est strictement croissante sur $[2 ,
+ \infty[$ et, puisque on a: $g'(2) = \log 2 + 1 + \log 2904 -
\log a = \log \left(\frac{5808 e}{a}\right) < 0$ et $\lim_{x
\rightarrow + \infty} g'(x) = + \infty$, alors, d'apr\`es le
théor\`eme de la bijection, $g'$ réalise une bijection de $[2 , +
\infty[$ sur $[g'(2) , + \infty[$. D'o\`u l'existence d'un unique
$\xi \in [2 , + \infty[$ tel que $g'(\xi) = 0$. De plus, $g'$ est
strictement négative sur $[2 , \xi[$ et elle est strictement
positive sur $]\xi , + \infty[$. Ce qui entra{\sf\^\i}ne que la
fonction $g$ est strictement décroissante sur $[2 , \xi]$ et elle
est strictement croissante sur $[\xi , + \infty[$, par conséquent
$g$ atteint son minimum en $\xi$ et il en est de m\^eme pour $f$
puisque $f =: \exp(g)$. Démontrer l'encadrement du lemme
\ref{c.54} pour $\xi$ revient, d'apr\`es cette étude, à démontrer
les deux inégalités:
\begin{align}
g'\!\!\left(\sqrt{\frac{2 \log a}{\log\log a - \log\log\log a + 21}} + 1\right)&<~ 0 \notag \\
\intertext{et} g'\!\!\left(\sqrt{\frac{2 \log a}{\log\log a -
\log\log\log a + 16}} + 1\right)&>~ 0 . \notag
\end{align}
Commençons par démontrer la premi\`ere, on a:
\newpage
\begin{equation*}
\begin{split}
g'\!\!\left(\sqrt{\frac{2 \log a}{\log\log a - \log\log\log a + 21}} + 1\right) \\
&\!\!\!\!\!\!\!\!\!\!\!\!\!\!\!\!\!\!\!\!\!\!\!\!\!\!\!\!\!\!\!\!\!\!\!\!\!\!\!\!\!\!\!\!\!\!\!\!\!\!\!\!\!\!\!\!\!\!\!\!\!\!\!\!\!\!\!\!\!\!= \log\!\!\left(\sqrt{\frac{2 \log a}{\log\log a - \log\log\log a + 21}} + 1\right) + 1 +\log 2904 \\
&~~~~\!\!\quad- \frac{\log\log a - \log\log\log a + 21}{2} \\
&\!\!\!\!\!\!\!\!\!\!\!\!\!\!\!\!\!\!\!\!\!\!\!\!\!\!\!\!\!\!\!\!\!\!\!\!\!\!\!\!\!\!\!\!\!\!\!\!\!\!\!\!\!\!\!\!\!\!\!\!\!\!\!\!\!\!\!\!\!\!\leq~ \log\!\!\left(\sqrt{\frac{2 \log a}{\log\log a - \log\log\log a + 21}}\right) \\
&\!\!\!\!\!\!\!\!\!\!\!\!\!\!\!\!\!\!\!\!\!\!\!\!\!\!\!\!\!\!\!\!\!\!\!\!\!\quad+ \sqrt{\frac{\log\log a - \log\log\log a + 21}{2 \log a}} + 1 + \log 2904 \\
&~~~~~~~~~\!\quad- \frac{1}{2} \log{\!\left(\frac{\log{a}}{\log{\log{a}}}\right)} - \frac{21}{2} \\
&\!\!\!\!\!\!\!\!\!\!\!\!\!\!\!\!\!\!\!\!\!\!\!\!\!\!\!\!\!\!\!\!\!\!\!\!\!\!\!\!\!\!\!\!\!\!\!\!\!\!\!\!\!\!\!\!\!\!\!\!\!\!\!\!\!\!\!\!\!\!\leq~ \frac{1}{2} \log 2 + \frac{1}{2} \log\log a - \frac{1}{2} \log(\log\log a - \log\log\log a + 21) \\
&\!\!\!\!\!\!\!\!\!\!\!\!\!\!\!\!\!\!\!\!\!\!\!\!\!\!\!\!\!\!\!\!\!\!\!\!\!\!\!\!\!\!\!\!\!\!\!\!\!\!\!\quad+ 1 + \log 2904 - \frac{21}{2} + \sqrt{\frac{\log\log a - \log\log\log a + 21}{2 \log a}} \\
&~~~~~~~~~~~~~~~\!\quad- \frac{1}{2} \log{\left(\frac{\log{a}}{\log{\log{a}}}\right)} \\
&\!\!\!\!\!\!\!\!\!\!\!\!\!\!\!\!\!\!\!\!\!\!\!\!\!\!\!\!\!\!\!\!\!\!\!\!\!\!\!\!\!\!\!\!\!\!\!\!\!\!\!\!\!\!\!\!\!\!\!\!\!\!\!\!\!\!\!\!\!\!\leq~ \frac{1}{2} \log\!\!\left(\frac{\log\log a}{\log\log a - \log\log\log a + 21}\right) \\
&\!\!\!\!\!\!\!\!\!\!\!\!\!\!\!\!\!\!\!\!\!\!\!\!\!\!\!\!\!\!\!\!\!\!\!\!\!\!\!\!\!\!\!\!\!\!\!\!\!\!\!\!\!\!\!\!\!\!\!\!\!\!\!\!\!\quad+ \sqrt{\frac{\log\log a - \log\log\log a + 21}{2 \log a}} + \frac{1}{2} \log 2 + \log 2904 - \frac{19}{2} \\
&\!\!\!\!\!\!\!\!\!\!\!\!\!\!\!\!\!\!\!\!\!\!\!\!\!\!\!\!\!\!\!\!\!\!\!\!\!\!\!\!\!\!\!\!\!\!\!\!\!\!\!\!\!\!\!\!\!\!\!\!\!\!\!\!\!\!\!\!\!\!\leq~
\frac{1}{2} \log\!\!\left(\frac{e^{22}}{e^{22} - 1}\right) +
\sqrt{
\frac{\log\log a_0 - \log\log\log a_0 + 21}{2 \log a_0}} \\
&~~~~~~~~~\!\!\quad+ \frac{1}{2} \log 2 + \log 2904 - \frac{19}{2} \\
&\!\!\!\!\!\!\!\!\!\!\!\!\!\!\!\!\!\!\!\!\!\!\!\!\!\!\!\!\!\!\!\!\!\!\!\!\!\!\!\!\!\!\!\!\!\!\!\!\!\!\!\!\!\!\!\!\!\!\!\!\!\!\!\!\!\!\!\!\!\!<~
0 .
\end{split}
\end{equation*}
L'avant derni\`ere inégalité provient de:
\begin{eqnarray*}
\frac{e^{22} - 1}{e^{22}} \log\log a & \leq & \log\log a -
\log\log\log a + 21 \\   & \leq & \frac{\log\log a_0 -
\log\log\log a_0 + 21}{\log a_0} \log a ,
\end{eqnarray*}
qui se montre par une simple étude de fonctions (en $a$) en tenant
compte du fait \\ $a > a_0 := 15788$. La premi\`ere inégalité est
alors démontrée. \\ Montrons maintenant la deuxi\`eme inégalité,
on a:
\begin{equation*}
\begin{split}
g'\!\!\left(\sqrt{\frac{2 \log a}{\log\log a - \log\log\log a + 16}} + 1\right) \\
&\!\!\!\!\!\!\!\!\!\!\!\!\!\!\!\!\!\!\!\!\!\!\!\!\!\!\!\!\!\!\!\!\!\!\!\!\!\!\!\!\!\!\!\!\!\!\!\!\!\!\!\!\!\!\!\!\!\!\!\!\!\!\!\!\!\!\!\!\!\!\!\!\!\!\!\!\!\!\!\!=~ \log\!\!\left(\sqrt{\frac{2 \log a}{\log\log a - \log\log\log a + 16}} + 1\right) + 1 +\log 2904 \\
&~~~~~\!\!\!\!\quad- \frac{\log\log a - \log\log\log a + 16}{2} \\
&\!\!\!\!\!\!\!\!\!\!\!\!\!\!\!\!\!\!\!\!\!\!\!\!\!\!\!\!\!\!\!\!\!\!\!\!\!\!\!\!\!\!\!\!\!\!\!\!\!\!\!\!\!\!\!\!\!\!\!\!\!\!\!\!\!\!\!\!\!\!\!\!\!\!\!\!\!\!\!\!>~ \log\!\!\left(\sqrt{\frac{2 \log a}{\log\log a - \log\log\log a + 16}}\right) + 1 + \log 2904 \\
&~\!\quad- \frac{1}{2} \log\log a + \frac{1}{2} \log\log\log a - 8 \\
&\!\!\!\!\!\!\!\!\!\!\!\!\!\!\!\!\!\!\!\!\!\!\!\!\!\!\!\!\!\!\!\!\!\!\!\!\!\!\!\!\!\!\!\!\!\!\!\!\!\!\!\!\!\!\!\!\!\!\!\!\!\!\!\!\!\!\!\!\!\!\!\!\!\!\!\!\!\!\!\!\geq~ \frac{1}{2} \log 2 + \frac{1}{2} \log\log a - \frac{1}{2} \log(\log\log a - \log\log\log a + 16) \\
&\!\!\!\!\!\!\!\!\!\!\!\!\!\!\!\!\!\!\!\!\!\!\!\!\!\!\!\quad- \frac{1}{2} \log\log a + \frac{1}{2} \log\log\log a + \log 2904 - 7 \\
&\!\!\!\!\!\!\!\!\!\!\!\!\!\!\!\!\!\!\!\!\!\!\!\!\!\!\!\!\!\!\!\!\!\!\!\!\!\!\!\!\!\!\!\!\!\!\!\!\!\!\!\!\!\!\!\!\!\!\!\!\!\!\!\!\!\!\!\!\!\!\!\!\!\!\!\!\!\!\!\!\geq~ \frac{1}{2} \log\!\!\left(\frac{\log\log a}{\log\log a - \log\log\log a + 16}\right) + \frac{1}{2} \log 2 + \log 2904 - 7 \\
&\!\!\!\!\!\!\!\!\!\!\!\!\!\!\!\!\!\!\!\!\!\!\!\!\!\!\!\!\!\!\!\!\!\!\!\!\!\!\!\!\!\!\!\!\!\!\!\!\!\!\!\!\!\!\!\!\!\!\!\!\!\!\!\!\!\!\!\!\!\!\!\!\!\!\!\!\!\!\!\!\geq~ \frac{1}{2} \log\!\!\left(\frac{\log\log a_0}{\log\log a_0 - \log\log\log a_0 + 16}\right) + \frac{1}{2} \log 2 + \log 2904 - 7 \\
&\!\!\!\!\!\!\!\!\!\!\!\!\!\!\!\!\!\!\!\!\!\!\!\!\!\!\!\!\!\!\!\!\!\!\!\!\!\!\!\!\!\!\!\!\!\!\!\!\!\!\!\!\!\!\!\!\!\!\!\!\!\!\!\!\!\!\!\!\!\!\!\!\!\!\!\!\!\!\!\!>~
0 .
\end{split}
\end{equation*}
L'avant derni\`ere inégalité se montre en étudiant la fonction:
$$a \mapsto \frac{\log\log a}{\log\log a - \log\log\log a + 16} .$$
Ceci ach\`eve la démonstration de la deuxi\`eme inégalité et ach\`eve ainsi la démonstration de la premi\`ere partie du lemme. \\
Montrons maintenant la deuxi\`eme partie du lemme \ref{c.54}. En
posant:
$$t ~:=~ \sqrt{\frac{2 \log a}{\log\log a - \log\log\log a + 16}} ~,$$
on a $m_0 = [t + 2]$ (o\`u $[.]$ désigne la partie enti\`ere),
d'o\`u l'encadrement:
$$t + 1 ~<~ m_0 ~\leq~ t + 2$$
pour l'entier $m_0$. \\ De $a > a_0$ on déduit, gr\^ace à une
simple étude de fonction, qu'on a:
$$t ~>~ \sqrt{\frac{2 \log a_0}{\log\log a_0 - \log\log\log a_0 + 16}} ~>~ 1 .$$
Par ailleurs, la premi\`ere partie du présent lemme (déja
démontrée) affirme qu'on a $t + 1 > \xi$, ce qui entra{\sf\^\i}ne
qu'on a aussi $m_0 > \xi$, puis, d'apr\`es la croissance stricte
de la fonction $g'$, qu'on a $g'(m_0) > g'(\xi) = 0$,
c'est-à-dire:
$$g'(m_0) ~>~ 0 .$$
Pour majorer sans peine le réel positif $g(m_0)$, nous remarquons
que la fonction $g$ satisfait l'équation différentielle:
$$x (x - 1) y' + y ~=~ x^2 (\log x + 1 + \log 2904) - x .$$
On a ainsi:
\begin{eqnarray*}
g(m_0) & = & {m_0}^{\!\!2} \log{m_0} + (1 + \log{2904}) {m_0}^{\!\!2} - m_0 - ({m_0}^{\!\!2} - m_0) g'(m_0) \\
       & < & {m_0}^{\!\!2} \log{m_0} + 9 {m_0}^{\!\!2} - m_0 ~~~~~~~~~~~~~~~~~~~~~~~~~\text{(car $g'(m_0) > 0$)} \\
       & < & (t + 2)^2 \log{(t + 2)} + 9 (t + 2)^2 - (t + 2) ~~~~~~~~\text{(car $m_0 \leq t + 2$)} \\
       & < & (t^2 + 4 t + 4) (\log t + \frac{2}{t}) + 9 t^2 + 35 t + 34 \\
       & < & t^2 \log t + 9 t^2 + 4 t \log t + 37 t + 4 \log t + 42 + \frac{8}{t} \\
       & < & t^2 \log t + 99 t^2 ~~~~\text{(en utilisant l'inégalité: $\log x \leq \frac{x}{e}$, $~ \forall x > 0$).}
\end{eqnarray*}
En fonction de $a$, les calculs donnent:
$$t^2 \log t + 99 t^2 ~=~ \log{a} ~\!. \!\!\left[1 + s(
\log{s} + \log 2 + 182) . \frac{1}{\log\log a} \right] $$
$$\text{avec}~~~~~~~~~~s ~:=~ \frac{\log\log a}{\log\log a - \log\log\log a + 16} ~=~ t^2 . \frac{\log\log a}{2 \log a} .~~~~~~~~~~~~~~~~~~~~~~~~~~~$$
Or, une simple étude de fonction montre qu'on a:
$$s ~:=~ \frac{\log\log a}{\log\log a - \log\log\log a + 16} ~\leq~ \frac{e^{17}}{e^{17} - 1} .$$
D'o\`u:
\begin{eqnarray*}
t^2 \log{t} + 99 t^2 & \!\!\!\leq & \!\!\!\log{a}~\!.\!\left[1 +
\frac{e^{17}}{e^{17} - 1} \!\left(
\log\!\left(\frac{e^{17}}{e^{17} - 1}\right) + \log 2 + 182\right)\!\frac{1}{\log\log a}\right] \\
 &\!\!\! < & \!\!\!\log{a}~\!.\!\left[1 + \frac{183}{\log\log a}\right].
\end{eqnarray*}
Par conséquent:
$$g(m_0) ~\leq~ t^2 \log t + 99 t^2 ~<~ \log{a}~\!.\!\left[1 + \frac{183}{\log\log a}\right] .$$
D'o\`u finalement: $\displaystyle f(m_0) = \exp{g(m_0)} < a^{1 +
\frac{183}{\log\log a}}$, ce qui ach\`eve la démonstration du
lemme \ref{c.54}.  $~~~~\blacksquare$\vspace{1mm}
\begin{corollary}\label{c.51}
Soit $R$ l'expression définie au §$11$, c'est-\`a-dire:
$$R(m) ~:=~ c_1(2904 m)^m a^{\frac{m}{m - 1}} + c_2 m(2904 m)^m a ~~~~(m \in {\mathbb N}_{\geq 2}) ,$$
avec $c_1 := 55 \eta + 272 ~,~ c_2 := \eta + 2 < \frac{c_1}{55}~$
et $~ a := \frac{1}{\epsilon}$. Sous l'hypoth\`ese $\epsilon <
\frac{1}{15788}$, l'entier $m_0$ défini dans le lemme \ref{c.54}
précédent minimise presque l'expression $R$ et lui donne une
valeur majorée par:
$$R(m_0) ~\leq~ 56(\eta + 5) a^{1 + \frac{183}{\log\log a}} .$$
\end{corollary}
{\bf Démonstration.---}
Le fait que $m_0$ minimise presque $R$ est d\^u simplement au fait que $m_0$ minimise presque l'expression $f(m) := (2904 m)^m a^{\frac{m}{m - 1}}$ (d'apr\`es le lemme \ref{c.54} précédent) et au fait qu'on a $m(2904 m)^m a \ll f(m)$ lorsque l'entier $m$ est proche de $m_0$ (voir les calculs ci-dessous). \\ Montrons maintenant la majoration du corollaire \ref{c.51} pour $R(m_0)$. \\

Comme $c_2 < \frac{c_1}{55}$, on a pour tout $m \geq 2$: $R(m)
\leq c_1 f(m)(1 + \frac{1}{55} m a^{- \frac{1}{m - 1}})$, en
particulier, pour l'entier $m_0$ on a:
$$R(m_0) ~\leq~ c_1 f(m_0)\left(1 + \frac{1}{55} m_0 a^{- \frac{1}{m_0 - 1}}\right) .$$
D'une part, d'apr\`es le lemme \ref{c.54} précédent on a:
$$f(m_0) ~<~ a^{1 + \frac{183}{\log\log a}}$$
et d'autre part, en désignant par $t$ le réel introduit dans la
démonstration du lemme \ref{c.54}, on a:
$$m_0 a^{- \frac{1}{m_0 - 1}} ~\leq~ (t + 2) a^{- \frac{1}{t + 1}} ~<~ e^{-3}$$
o\`u dans cette double inégalité, l'inégalité de gauche sort de
l'encadrement $t + 1 < m_0 \leq t + 2$ pour l'entier $m_0$ et
l'inégalité de droite est équivalente à l'inégalité $\frac{\log
a}{t + 1} > \log{(t + 2)} + 3$ laquelle se démontre en utilisant
des estimations du type de celles déja utilisées dans la
démonstration du lemme \ref{c.54} en tenant compte seulement du
fait $a > a_0$.
\begin{eqnarray*}
\text{D'o\`u:}~~~~~~~~~~~~~~ R(m_0) & \leq & c_1 a^{1 + \frac{183}{\log\log a}}(1 + e^{-3}/55)~~~~~~~~~~~~~~~~~~~~~~~~~~~~~~~~~~~~~~~~~~~~~~~~~~~~~~ \\
 & < & (56 \eta + 273) a^{1 + \frac{183}{\log\log a}} ,
\end{eqnarray*}
et a fortiori: $\displaystyle R(m_0) ~\!<~\! 56(\eta + 5) a^{1 +
\frac{183}{\log\log a}}$, ce qui ach\`eve cette démonstration.
 $~~~~\blacksquare$\vspace{1mm}
\begin{lemma}\label{c.52}
Soit $a$ un réel $\geq e^2$, $r$ un entier $\geq 1$ et $q$ la
fonction réelle sur $[2 , + \infty[$ définie par:
$$q(x) := \left(x(x - 1)(\log x + 9) + \frac{2}{r} \log{a}.x\right) a^{\frac{x}{x - 1}} .$$
Alors, la valeur minimale de $q$ sur $[2 , + \infty[$ vérifie:
$$16 (\log a)^2 a ~<~ \min_{x \in [2 , + \infty[} q(x) ~<~ 2 (\log a)^2 (\log\log a + 82) a .$$
De plus, la valeur de $q$ en l'entier $m_1 := [\frac{\log a}{2} +
2]$ (o\`u $[.]$ désigne la partie enti\`ere) est majorée par:
$$q(m_1) ~\leq~ 2 (\log a)^2 (\log\log a + 82) a .$$
\end{lemma}
{\bf Démonstration.---} Il est clair qu'on a:
$$\forall x \geq 2 ~~~~~~~~q(x) ~\geq~ 9 (x - 1)^2 a^{\frac{x}{x - 1}} ,$$
$$\text{donc}~~~~~~~~~~~~~~ \min_{x \in [2 , + \infty[} q(x) ~\geq \min_{x \in [2 , + \infty[} 9 (x - 1)^2 a^{\frac{x}{x - 1}} .~~~~~~~~~~~~~~~~~~~~~$$
Or, une simple étude de fonction montre que la fonction $9 (x -
1)^2 a^{\frac{x}{x - 1}}$ (sur $[2 , + \infty[$) atteint sa valeur
minimale au point $x = \frac{\log a}{2} + 1$ et cette derni\`ere
vaut $\frac{9}{4} e^2 (\log a)^2 a > 16 (\log a)^2 a$ ce qui
entra{\sf\^\i}ne qu'on a:
$$\min_{x \in [2 , + \infty[} q(x) ~>~ 16 (\log a)^2 a .$$
Par ailleurs, l'entier $m_1 := [\frac{\log a}{2} + 2]$ vérifie
l'encadrement $\frac{\log a}{2} + 1 < m_1 \leq \frac{\log a}{2} +
2$ gr\^ace auquel on obtient les deux majorations:
$$m_1 (m_1 - 1) (\log m_1 + 9) + \frac{2}{r} \log{a}.m_1 ~\leq~ \left(\frac{\log a}{2}\right)^2 (\log\log a +82)$$
$$\text{et}~~~~~~~~~~~~~~~~~~~~~~~~~~~~~~~ a^{\frac{m_1}{m_1 - 1}} ~<~ e^2 a .~~~~~~~~~~~~~~~~~~~~~~~~~~~~~~~~$$
\begin{eqnarray*}
\text{Ce qui donne}~~~~~~ q(m_1) & \!\!\!\!:= & \!\!\!\!\left(\!m_1(m_1 - 1)(\log m_1 + 9) + \frac{2}{r} \log{a}.m_1\!\right)\!a^{\frac{m_1}{m_1 - 1}} ~~~~~ \\
&\!\!\!\!< & \!\!\!\!\frac{e^2}{4} (\log a)^2 (\log\log a + 82) a \\
& \!\!\!\!< & \!\!\!\!2 (\log a)^2 (\log\log a + 82) a .
\end{eqnarray*}
La démonstration du lemme \ref{c.52} se termine par la remarque
banale:
$$\min_{x \in [2 , + \infty[} q(x) ~\leq~ q(m_1) ~<~ 2 (\log a)^2 (\log\log a + 82) a .~~~~\blacksquare$$
\begin{corollary}\label{c.53}
Soit $m$ un entier $\geq 2$ (que nous prenons pour variable), $r$
un entier $\geq 1$ (que nous prenons pour param\`etre), $\epsilon$
un réel positif $\leq e^{- 4/r}$ (que nous prenons pour
param\`etre aussi) et $S$ l'expression définie par:
$$S_{r , \epsilon}(m) ~:=~
 4 \epsilon^{- \frac{1}{2}} \left[m(m - 1)(\log m + 9) + m \!\mid \!\log \epsilon \!\mid\right]\!\!\left(\!499 \epsilon^{- \frac{m}{2(m - 1)}}\!\right)^{\!\!r} .$$
Alors, $S$ (en tant que fonction en $m$) atteint presque sa valeur
minimale en l'entier $m_1 := [\frac{r}{4} \!\mid \!\log \epsilon
\!\mid + 2]$ (o\`u $[.]$ désigne la partie enti\`ere) et la valeur
de $S$ en $m_1$ est majorée par:
$$S(m_1) ~\leq~ 2 r^2 \epsilon^{- \frac{1}{2}} {\!\mid \!\log \epsilon \!\mid}^2 \!\!\left(\log r + \log \!\mid \!\log \epsilon \!\mid + 82\right) \!\!\left(\!\!499 \epsilon^{- \frac{1}{2}}\right)^{\!\!r} .$$
\end{corollary}
{\bf Démonstration.---} La fonction $q$ définie au lemme
\ref{c.52} s'écrit pour $a := \epsilon^{- r/2}$ ($a$ vérifie bien
$a \geq e^2$ puisque $\epsilon \leq e^{- 4/r}$):
$$q(x) = \left(x(x - 1)(\log x + 9) + x \!\mid \!\log \epsilon \!\mid \!\right)\epsilon^{- \frac{r x}{2 (x - 1)}} .$$
Ainsi $S(m)$ est un produit de $q(m)$ avec une expression
indépendante de $m$. En effet, on a:
$$S(m) = 4 \epsilon^{- 1/2}(499)^r q(m) .$$
Il s'ensuit que $S(m)$ et $q(m)$ atteignent leurs valeurs
minimales au m\^eme point. Or, on sait, d'apr\`es le lemme
\ref{c.52}, que $q(m)$ atteint presque sa valeur minimale en $m_1
:= [\frac{\log a}{2} + 2] = [\frac{r}{4} \!\mid \!\log \epsilon
\!\mid + 2]$ et que:
\begin{eqnarray*}
q(m_1) & \leq & 2 (\log a)^2 (\log\log a + 82) a \\
& \leq & 2 \!\left(\frac{r}{2} \!\mid \!\log \epsilon \!\mid \right)^{\!2} \!\!\left(\!\log{\frac{r}{2}} + \log \!\mid \!\log \epsilon \!\mid + 82\right) \epsilon^{- \frac{r}{2}} \\
& \leq & \frac{1}{2} r^2 \!\left(\mid \!\log \epsilon \!\mid
\right)^{\!2} \!\left(\log{r} + \log \!\mid \!\log \epsilon
\!\mid\! + 82\right) \epsilon^{- \frac{r}{2}} .
\end{eqnarray*}
Ainsi $S(m)$ atteint aussi presque sa valeur minimale au m\^eme
point $m_1 := [\frac{r}{4} \!\mid \!\log \epsilon \!\mid + 2]$ et:
\begin{eqnarray*}
S(m_1) & = & 4 \epsilon^{- \frac{1}{2}} (499)^r q(m_1) \\
& \leq & 2 r^2 \epsilon^{- \frac{1}{2}} \left(\mid \!\log \epsilon
\!\mid \right)^{\!2} \!\left(\log{r} + \log \!\mid \!\log \epsilon
\!\mid + 82\right) \!\!\left(\!499 \epsilon^{-
\frac{1}{2}}\right)^{\!\!r} .
\end{eqnarray*}
Ce qui termine cette démonstration. $~~~~\blacksquare$

\end{document}